\documentclass[a4paper,12pt,amsart,frenchb]{article}

\usepackage{amsmath,amsbsy,amsfonts,amssymb}
\usepackage[french]{babel}
\newcommand{\ass}[2]{\vskip0.3cm\noindent
{\bf {#1}}. { \sl {#2}}\vskip0.3cm\noindent
}
  
\begin{document}
 
   \title{ Stabilisation de la formule des traces tordue V:   int\'egrales orbitales et endoscopie sur le corps r\'eel}
\author{J.-L. Waldspurger}
\date{2 avril 2014}
\maketitle

La stabilisation de la formule des traces tordue n\'ecessite quelques travaux pr\'eliminaires. Voici l'un d'eux. Dans deux  articles pr\'ec\'edents ([II] et [III]), on a trait\'e des int\'egrales orbitales pond\'er\'ees et de leur stabilisation, le corps de base \'etant local non-archim\'edien. On consid\`ere ici la m\^eme question, le corps de base \'etant cette fois r\'eel. Quelles sont les diff\'erences? De temps en temps, on doit faire un peu plus de topologie que sur un corps local non-archim\'edien. En effet, dans le cas r\'eel, les distributions sont des formes lin\'eaires continues sur des espaces de fonctions. Cette continuit\'e est essentielle et intervient dans la plupart  des d\'emonstrations. Mais ceci n'induit gu\`ere de changement dans la structure des preuves. Il y a une diff\'erence beaucoup plus  perturbante. Utilisons les notations des articles pr\'ecedents: $G$ est un groupe r\'eductif connexe d\'efini sur un corps local $F$, $\tilde{G}$ est un espace tordu sous $G$ et ${\bf a}$ est un \'el\'ement de $H^1(W_{F};Z(\hat{G}))$ dont se d\'eduit un caract\`ere $\omega $ de $G(F)$. On d\'efinit l'espace $D_{g\acute{e}om}(\tilde{G}(F),\omega)$ des distributions $\omega$-\'equivariantes \`a support contenu dans un nombre fini de classes de conjugaison. On note $D_{orb}(\tilde{G}(F),\omega)$ le sous-espace engendr\'e par les int\'egrales orbitales. Si $F$ est non-archim\'edien, ces deux espaces sont \'egaux. Dans le pr\'esent article, on a $F={\mathbb R}$ et ces deux espaces ne sont plus du tout \'egaux. Pour un espace de Levi $\tilde{M}$,  
 on d\'efinit \`a la suite d'Arthur les int\'egrales orbitales pond\'er\'ees tordues et leurs variantes $\omega$-\'equivariantes.   On peut  consid\'erer ces derni\`eres comme des formes bilin\'eaires $(\boldsymbol{\gamma},f)\mapsto I_{\tilde{M}}^{\tilde{G}}(\boldsymbol{\gamma},f)$, o\`u  $f$ appartient \`a $C_{c}^{\infty}(\tilde{G}({\mathbb R}))$ et $\boldsymbol{\gamma}$ appartient \`a $D_{orb}(\tilde{M}({\mathbb R}),\omega)$. La d\'efinition ne s'\'etend pas de fa\c{c}on simple \`a  $\boldsymbol{\gamma}\in D_{g\acute{e}om}(\tilde{M}({\mathbb R}),\omega)$. Quand on veut "stabiliser" ces int\'egrales, on est conduit \`a transf\'erer les distributions $\boldsymbol{\gamma}$. Consid\'erons une donn\'ee endoscopique ${\bf M}'=(M',{\cal M}',\tilde{\zeta})$ de $(M,\tilde{M},{\bf a})$. Supposons que l'on soit dans une situation simple o\`u le recours \`a des donn\'ees auxiliaires ne soit pas n\'ecessaire. On peut alors d\'efinir un transfert
 $$D_{g\acute{e}om}^{st}(\tilde{M}'({\mathbb R}))\to D_{g\acute{e}om}(\tilde{M}({\mathbb R}),\omega),$$
 o\`u $D_{g\acute{e}om}^{st}(\tilde{M}'({\mathbb R}))$ est le sous-espace des \'el\'ements de $D_{g\acute{e}om}(\tilde{M}'({\mathbb R}))$ qui sont stables. Notons de m\^eme $D_{orb}^{st}(\tilde{M}'({\mathbb R}))$ le sous-espace des \'el\'ements de $D_{orb}(\tilde{M}'({\mathbb R}))$ qui sont stables. Comme le montre un exemple d\^u \`a Magdy Assem (je remercie Kottwitz de me l'avoir indiqu\'e), le transfert n'envoie pas en g\'en\'eral l'espace $D_{orb}^{st}(\tilde{M}'({\mathbb R}))$ dans $D_{orb}(\tilde{M}({\mathbb R}),\omega)$. Les constructions que l'on a faites en [II] sur un corps local non-archim\'edien s'effondrent  si on se limite aux espaces $D_{orb}$. La m\'ethode utilis\'ee par Arthur pour r\'esoudre ce probl\`eme consiste \`a g\'en\'eraliser la d\'efinition des int\'egrales $I_{\tilde{M}}^{\tilde{G}}(\boldsymbol{\gamma},{\bf f})$  aux \'el\'ements $\boldsymbol{\gamma}$ de $D_{g\acute{e}om}(\tilde{M}({\mathbb R}),\omega)$ tout entier, cf. [A1]. Nous utilisons une autre m\'ethode. Dans un premier temps, on se limite aux \'el\'ements de $D_{g\acute{e}om}(\tilde{M}({\mathbb R}),\omega)$ dont le support est form\'e d'\'el\'ements $\gamma\in \tilde{M}({\mathbb R})$ qui sont $\tilde{G}$-\'equisinguliers, c'est-\`a-dire tels que $M_{\gamma}=G_{\gamma}$.   Pour de tels \'el\'ements, il est facile d'\'etendre la d\'efinition des termes $I_{\tilde{M}}^{\tilde{G}}(\boldsymbol{\gamma},{\bf f})$. La th\'eorie, restreinte \`a ces distributions, est alors tr\`es semblable \`a celle du cas non-archim\'edien. C'est l'objet de la section 1. Pour traiter les distributions \`a support quelconque, on d\'efinit par r\'ecurrence un espace $D_{tr-orb}(\tilde{M}({\mathbb R}),\omega)$, cf. section 2. C'est grosso-modo l'espace engendr\'e par $D_{orb}(\tilde{M}({\mathbb R}),\omega)$ et par les images par transfert endoscopique d'espaces $D_{tr-orb}^{st}(\tilde{M}'({\mathbb R}))$ quand ${\bf M}'$ parcourt toutes les donn\'ees endoscopiques de $(M,\tilde{M},{\bf a})$. Cette d\'efinition assure que ces espaces sont "stables par transfert", en un sens assez clair. Le point est de montrer que l'on peut \'etendre la d\'efinition des termes $I_{\tilde{M}}^{\tilde{G}}(\boldsymbol{\gamma},{\bf f})$ aux \'el\'ements $\boldsymbol{\gamma}\in D_{tr-orb}(\tilde{M}({\mathbb R}),\omega)$. On \'enonce pr\'ecis\'ement en 2.4 toutes les propri\'et\'es que ces termes doivent v\'erifier. Les preuves, dont le seul int\'er\^et est d'exister,  sont 
  l'objet des sections 3, 4 et 5. En gros, on montre que les termes $I_{\tilde{M}}^{\tilde{G}}(\boldsymbol{\gamma},{\bf f})$ peuvent \^etre d\'efinis par un proc\'ed\'e de limite similaire \`a celui qu'utilise Arthur pour les d\'efinir   dans le cas o\`u $\boldsymbol{\gamma}$ est une int\'egrale orbitale \`a support singulier. En fait, dans le cas g\'en\'eral, on ne m\`ene \`a bien ce programme qu'en admettant une hypoth\`ese. Celle-ci est l'assertion principale de la stabilisation, \`a savoir l'\'egalit\'e  $I_{\tilde{M}}^{\tilde{G},{\cal E}}(\boldsymbol{\gamma},{\bf f})=I_{\tilde{M}}^{\tilde{G}}(\boldsymbol{\gamma},{\bf f})$  quand $\boldsymbol{\gamma}$ est une int\'egrale orbitale \`a support fortement $\tilde{G}$-r\'egulier. A ce gros probl\`eme de d\'efinition pr\`es, les preuves sont tr\`es similaires \`a celles du cas non-archim\'edien. On ne les reprendra que rapidement. Le point qui nous int\'eressera surtout sera de v\'erifier que l'on ne sort jamais d'espaces pour lesquels tous les termes sont d\'efinis.  On obtient les m\^emes r\'esultats qu'en [III], c'est-\`a-dire que tous les r\'esultats esp\'er\'es r\'esultent de l'assertion principale \'evoqu\'ee ci-dessus. 
 Comme on vient de le dire, la diff\'erence est que non seulement les r\'esultats eux-m\^emes d\'ependent de l'assertion principale mais m\^eme la d\'efinition de certains termes en d\'epend. Toutefois, on obtient des d\'efinitions et r\'esultats non conditionnels dans le cas o\`u $G$ est quasi-d\'eploy\'e, $\tilde{G}$ est \`a torsion int\'erieure et ${\bf a}=1$. Dans ce cas, on montre par la m\^eme m\'ethode qu'en [III] que l'assertion principale est cons\'equence des r\'esultats d'Arthur.

Supposons $(G,\tilde{G},{\bf a})$ quasi-d\'eploy\'e et \`a torsion int\'erieure. L'int\'er\^et des constructions des  paragraphes  2 \`a  5 est que le germe en $a=1$ de la fonction $a\mapsto I_{\tilde{M}}^{\tilde{G}}(a\boldsymbol{\gamma},f)$ appartient \`a un espace de germes que l'on  contr\^ole bien. C'est ce qui nous permet d'obtenir les r\'esultats voulus. Toutefois, pour certaines applications ult\'erieures, la d\'efinition qu'elles nous fournissent de $I_{\tilde{M}}^{\tilde{G}}(\boldsymbol{\gamma},f)$ n'a pas la meilleure forme possible. Il s'av\`ere que l'on aura besoin d'une autre approximation de cette int\'egrale, que nous \'etablissons dans la section 6. Dans cette nouvelle approximation, on perd sur l'espace des germes de fonctions en $a$, qui est moins contr\^olable. Par contre, il n'intervient plus dans la formule que des int\'egrales $I_{\tilde{L}}^{\tilde{G}}(\boldsymbol{\gamma}_{\tilde{L}}(a),f)$, o\`u $\boldsymbol{\gamma}_{\tilde{L}}(a)$  est une distribution sur $\tilde{L}({\mathbb R})$ \`a support $\tilde{G}$-\'equisingulier. On obtient ainsi une approximation de $I_{\tilde{M}}^{\tilde{G}}(\boldsymbol{\gamma},f)$ par de telles int\'egrales. 

Dans la derni\`ere section, on dira quelques mots du cas o\`u le corps de base n'est plus ${\mathbb R}$ mais ${\mathbb C}$.

\bigskip

\section{Int\'egrales orbitales pond\'er\'ees}
\bigskip

\subsection{La situation}
Dans cet article, le corps de base est ${\mathbb R}$. On consid\'erera soit un triplet $(G,\tilde{G},{\bf a})$ comme en [IV] 1.1, soit un "$K$-triplet" $(KG,K\tilde{G},{\bf a})$ comme en [I] 1.11. Dans le premier cas, on fixe un espace de Levi minimal $\tilde{M}_{0}$ de $\tilde{G}$ et un sous-groupe compact maximal $K$ de $G({\mathbb R})$ en bonne position relativement \`a $M_{0}$. Dans le second, on \'ecrit $K\tilde{G}=(\tilde{G}_{p})_{p\in \Pi}$. Pour tout $p\in \Pi$, on fixe un espace de Levi minimal $\tilde{M}_{p,0}$ de $\tilde{G}_{p}$ et un sous-groupe compact maximal $K_{p}$ de $G_{p}({\mathbb R})$ qui soit en bonne position relativement \`a $M_{p,0}$. On d\'efinit le $K$-espace de Levi minimal $K\tilde{M}_{0}$ de $K\tilde{G}$ et 
 l'ensemble ${\cal L}(K\tilde{M}_{0})$ comme en  [I] 3.5. On notera symboliquement $K=(K_{p})_{p\in \Pi}$. D'une fa\c{c}on g\'en\'erale, on simplifie la notation en supprimant la lettre $K$ de $K\tilde{G}$ pour les objets ind\'ependants de $p\in \Pi$, par exemple  on pose $dim(G)=dim(G_{p})$ pour tout $p$.

On raisonne par r\'ecurrence sur $dim(G_{SC})$. Pour d\'emontrer une assertion concernant un triplet $(G,\tilde{G},{\bf a})$ quasi-d\'eploy\'e et \`a torsion int\'erieure, on suppose connues toutes les assertions concernant des triplets quasi-d\'eploy\'es et \`a torsion int\'erieure $(G',\tilde{G}',{\bf a}')$ tels que $dim(G'_{SC})< dim(G_{SC})$. Pour d\'emontrer une assertion concernant un triplet $(G,\tilde{G},{\bf a})$ qui n'est pas  quasi-d\'eploy\'e et \`a torsion int\'erieure, on suppose connues toutes les assertions concernant des triplets quasi-d\'eploy\'es et \`a torsion int\'erieure $(G',\tilde{G}',{\bf a}')$ tels que $dim(G'_{SC})\leq dim(G_{SC})$, ainsi que toutes les assertions concernant un triplet $(G',\tilde{G}',{\bf a}')$ quelconque tel que $dim(G'_{SC})< dim(G_{SC})$. Si une assertion concerne un espace de Levi $\tilde{M}$ de $\tilde{G}$, on suppose connues toutes les assertions concernant le m\^eme triplet $(G,\tilde{G},{\bf a})$ et un espace de Levi $\tilde{L}\in {\cal L}(\tilde{M})$ tel que $\tilde{L}\not=\tilde{M}$. Pour d\'emontrer une assertion concernant un $K$-triplet $(KG,K\tilde{G},{\bf a})$, on  suppose connues toutes les assertions concernant des triplets quasi-d\'eploy\'es et \`a torsion int\'erieure $(G',\tilde{G}',{\bf a}')$ tels que $dim(G'_{SC})\leq dim(G_{SC})$ et toutes les assertions concernant des $K$-triplets $(KG',K\tilde{G}',{\bf a})$ tels que $dim(G'_{SC})<dim(G_{SC})$ (c'est-\`a-dire que, pour nous, la notion de $K$-triplet quasi-d\'eploy\'e et \`a torsion int\'erieure n'existe pas). Si une assertion concerne un $K$-espace de Levi $K\tilde{M}\in {\cal L}(K\tilde{M}_{0})$, on suppose connues toutes les assertions concernant le m\^eme triplet $(KG,K\tilde{G},{\bf a})$ et un $K$-espace de Levi $K\tilde{L}\in {\cal L}(K\tilde{M})$ tel que $K\tilde{L}\not=K\tilde{M}$.

\subsection{L'application $\phi_{\tilde{M}}$}
Jusqu'en 1.7, on consid\`ere un triplet $(G,\tilde{G},{\bf a})$. Soit $\tilde{M}\in {\cal L}(\tilde{M}_{0})$. En [W1] 6.4, on a d\'efini, en suivant Arthur, une application $\phi_{\tilde{M}}$ qui, dans le cas qui nous occupe o\`u le corps de base est ${\mathbb R}$, n'\'etait d\'efinie que sur l'espace $C_{c}^{\infty}(\tilde{G}({\mathbb R}),K)$ des fonctions $K$-finies. Cela parce que l'on utilisait le th\'eor\`eme de Paley-Wiener de Delorme-Mezo, qui s'applique \`a ces fonctions. En utilisant le th\'eor\`eme de Renard, cf. [IV] 1.4, on va \'etendre cette application \`a tout $C_{c}^{\infty}(\tilde{G}({\mathbb R}))$. Pour simplifier, on fixe des mesures de Haar sur tous les groupes qui apparaissent.

Notons $C_{ac}^{\infty}(\tilde{G}({\mathbb R}))$ l'espace des fonctions $f:\tilde{G}({\mathbb R})\to {\mathbb C}$ telles que, pour toute fonction $b\in C_{c}^{\infty}({\cal A}_{\tilde{G}})$, la fonction $f(b\circ H_{\tilde{G}}):\gamma\mapsto b(H_{\tilde{G}}(\gamma))f(\gamma)$ appartienne \`a $C_{c}^{\infty}(\tilde{G}({\mathbb R}))$. Consid\'erons l'immense produit
$$\underline{P}=C_{c}^{\infty}(\tilde{G}({\mathbb R}))^{C_{c}^{\infty}({\cal A}_{\tilde{G}})}.$$
Ses \'el\'ements sont les familles $(\phi_{b})_{b\in C_{c}^{\infty}({\cal A}_{\tilde{G}})}$, o\`u $\phi_{b}\in C_{c}^{\infty}(\tilde{G}({\mathbb R}))$ pour tout $b$. L'espace $C_{c}^{\infty}(\tilde{G}({\mathbb R}))$ \'etant muni de la topologie usuelle, on munit $P$ de la topologie produit. Notons $P$ le sous-espace de $\underline{P}$ form\'e des familles $(\phi_{b})_{b\in C_{c}^{\infty}({\cal A}_{\tilde{G}})}$ telles que, pour tous $b,b'\in C_{c}^{\infty}({\cal A}_{\tilde{G}})$, on ait l'\'egalit\'e
$$\phi_{b}(b'\circ H_{\tilde{G}})=\phi_{b'}(b\circ H_{\tilde{G}}).$$
On v\'erifie que c'est un sous-espace ferm\'e de $P$ et on le munit de la topologie induite. On a une application lin\'eaire
$$\begin{array}{ccc}C_{ac}^{\infty}(\tilde{G}({\mathbb R}))&\to& P\\ f&\mapsto& (f(b\circ H_{\tilde{G}}))_{b\in C_{c}^{\infty}({\cal A}_{\tilde{G}})}\\ \end{array}$$
On v\'erifie qu'elle est injective et que son image est $P$. On munit $C_{ac}^{\infty}(\tilde{G}({\mathbb R}))$ de la topologie pour laquelle cette application devient un hom\'eomorphisme de cet espace sur $P$. Concr\`etement, un voisinage de $0$ dans $C_{ac}^{\infty}(\tilde{G}({\mathbb R}))$ est un sous-ensemble $U$ pour lequel il existe un nombre fini  d'\'el\'ements $b_{1},...,b_{n}\in C_{c}^{\infty}({\cal A}_{\tilde{G}})$ et, pour tout $i$, un voisinage $U_{i}$ de $0$ dans $C_{c}^{\infty}(\tilde{G}({\mathbb R}))$ de sorte que $U$ contienne les $f\in C_{ac}^{\infty}(\tilde{G}({\mathbb R}))$ tels que, pour tout $i$, $f(b_{i}\circ H_{\tilde{G}})$ appartienne \`a $U_{i}$. 

On v\'erifie que, pour tout $\gamma\in \tilde{G}({\mathbb R})$, l'int\'egrale orbitale $f\mapsto I^{\tilde{G}}(\gamma,\omega,f)$ se prolonge en une forme lin\'eaire continue sur $C_{ac}^{\infty}(\tilde{G}({\mathbb R}))$. De m\^eme, pour un espace de Levi $\tilde{M}\in {\cal L}(\tilde{M}_{0})$, pour une $\omega$-repr\'esentation temp\'er\'ee  $\tilde{\pi}$ de $\tilde{M}({\mathbb R})$ (de longueur finie) et pour $X\in {\cal A}_{\tilde{M}}$, l'application $f\mapsto J_{\tilde{M}}^{\tilde{G}}(\tilde{\pi},X,f)$ d\'efinie en [W1] 6.4 se prolonge en une forme lin\'eaire continue sur $C_{ac}^{\infty}(\tilde{G}({\mathbb R}))$. 

Modifiant la d\'efinition que l'on avait donn\'ee en [W1] 6.4, on note $I_{ac}(\tilde{G}({\mathbb R}),\omega)$ le quotient de $C_{ac}^{\infty}(\tilde{G}({\mathbb R}))$ par le sous-espace des $f\in C_{ac}^{\infty}(\tilde{G}({\mathbb R}))$ telles que $I^{\tilde{G}}(\gamma,\omega,f)=0$ pour tout $\gamma$ fortement r\'egulier. On voit que ce sous-espace est ferm\'e et on munit le quotient de la topologie quotient. On peut d\'efinir des espaces $\underline{IP}$ et $IP$ en rempla\c{c}ant dans les d\'efinitions de $\underline{P}$ et $P$ les espaces $C_{c}^{\infty}(\tilde{G}({\mathbb R}))$ par $I(\tilde{G}({\mathbb R}),\omega)$. On v\'erifie qu'une application analogue \`a celle construite ci-dessus identifie hom\'eomorphiquement $I_{ac}(\tilde{G}({\mathbb R}),\omega)$ \`a $IP$.   Remarquons que, pour toute $\omega$-repr\'esentation de longueur finie $\tilde{\pi}$ de $\tilde{G}({\mathbb R})$ et pour tout $X\in {\cal A}_{\tilde{G}}$, la forme lin\'eaire $f\mapsto I^{\tilde{G}}(\tilde{\pi},X,f)$ sur $C_{ac}^{\infty}(\tilde{G}({\mathbb R}))$ se factorise en une forme lin\'eaire continue sur $I_{ac}(\tilde{G}({\mathbb R}),\omega)$. Cela r\'esulte de la locale int\'egrabilit\'e des caract\`eres. 

\ass{Proposition}{Soit $\tilde{M}\in {\cal L}(\tilde{M}_{0})$. Il existe une unique application lin\'eaire continue
$$\begin{array}{ccc}C_{ac}^{\infty}(\tilde{G}({\mathbb R}))&\to&I_{ac}(\tilde{M}({\mathbb R}),\omega)\\  f&\mapsto& \phi_{\tilde{M}}(f)\\ \end{array}$$
telle que, pour tout $f\in C_{ac}^{\infty}(\tilde{G}({\mathbb R}))$,  pour toute $\omega$-repr\'esentation temp\'er\'ee et de longueur finie $\tilde{\pi}$ de $\tilde{M}({\mathbb R})$ et pour tout $X\in {\cal A}_{\tilde{M}}$, on ait l'\'egalit\'e
$$I^{\tilde{M}}(\tilde{\pi},X,\phi_{\tilde{M}}(f))=J_{\tilde{M}}^{\tilde{G}}(\tilde{\pi},X,f).$$}

{\bf Remarque.} Comme toujours, il est plus canonique de voir l'application $\phi_{\tilde{M}}$ comme une application lin\'eaire de $C_{ac}^{\infty}(\tilde{G}({\mathbb R}))\otimes Mes(G({\mathbb R}))$ dans $I_{ac}(\tilde{M}({\mathbb R})\otimes Mes(M({\mathbb R}))$.
\bigskip

Preuve. On reprend la d\'emonstration donn\'ee en [W1] 6.4, en supprimant des conjugaisons complexes inopportunes. Fixons un sous-ensemble compact $\tilde{C}$ de $\tilde{G}({\mathbb R})$ et notons $C^{\infty}(\tilde{C})$ le sous-espace des fonctions $f\in C_{c}^{\infty}(\tilde{G}({\mathbb R}))$ telles que $Supp(f)\subset \tilde{C}$. Appelons semi-norme pour $\tilde{C}$ toute fonction $\alpha:C^{\infty}(\tilde{C})\to {\mathbb R}_{\geq0}$ telle qu'il existe un nombre fini d'op\'erateurs diff\'erentiels  $D_{1}$,...,$D_{m}$ sur $\tilde{G}({\mathbb R})$ invariants par translations \`a gauche de sorte que, pour tout $f\in C^{\infty}(\tilde{C})$, on ait l'\'egalit\'e
$$\alpha(f)=\sum_{i=1,...,m}sup_{\gamma\in \tilde{C}}\vert (D_{i}f)(\gamma)\vert .$$

On va commencer par \'etablir l'existence d'une application lin\'eaire
$$\begin{array}{ccc}C^{\infty}(\tilde{C})&\to&I_{ac}(\tilde{M}({\mathbb R}),\omega)\\  f&\mapsto& \phi_{\tilde{M}}(f)\\ \end{array}$$
v\'erifiant les propri\'et\'es de l'\'enonc\'e. Soient $f\in C^{\infty}( \tilde{C})$ et $b\in C_{c}^{\infty}({\cal A}_{\tilde{M}})$. Pour une $\omega$-repr\'esentation temp\'er\'ee et de longueur finie $\tilde{\pi}$ de $\tilde{M}({\mathbb R})$, posons
$$\varphi_{f,b}(\tilde{\pi})=\int_{{\cal A}_{\tilde{M}}}J_{\tilde{M}}^{\tilde{G}}(\tilde{\pi},X,f)b(X)\,dX.$$
 Par inversion de Fourier, on a l'\'egalit\'e
 $$\varphi_{f,b}(\tilde{\pi})=\int_{i{\cal A}_{\tilde{M}}^*}J_{\tilde{M}}^{\tilde{G}}(\tilde{\pi}_{\lambda},f)\hat{b}(-\lambda)\,d\lambda,$$
 o\`u
 $$\hat{b}(\lambda)=\int_{{\cal A}_{\tilde{M}}}b(X)e^{<X,\lambda>}\,dX.$$
 Fixons $\tilde{Q}=\tilde{L}U_{Q}\in {\cal F}^{\tilde{M}}(\tilde{M}_{0})$. Pour tout $\boldsymbol{\tau}\in {\cal E}_{ell,0}(\tilde{L},\omega)$, on d\'efinit une fonction $\varphi_{f,b,\boldsymbol{\tau}}$ sur ${\cal A}_{\tilde{L},{\mathbb C}}^*$ par
 $$\varphi_{f,b,\boldsymbol{\tau}}(\nu)=\varphi_{f,b}(Ind_{\tilde{Q}}^{\tilde{M}}(\tilde{\pi}_{\boldsymbol{\tau}_{\nu}})).$$
 On va prouver les propri\'etes suivantes
 
 (1) il existe $r>0$ ind\'ependant de $f$ et, pour tout $N>0$ et tout  $\boldsymbol{\tau}\in {\cal E}_{ell,0}(\tilde{L},\omega)$, il existe $C>0$ tel que, pour tout $\nu\in {\cal A}_{\tilde{L},{\mathbb C}}^*$, on a la majoration
 $$\vert \varphi_{f,b,\boldsymbol{\tau}}(\nu)\vert \leq C(1+\vert \nu\vert )^{-N}e^{r\vert Re(\nu)\vert };$$
 
 (2) pour tout $N>0$, il existe une semi-norme $\alpha$ pour $\tilde{C}$ de sorte que, pour tout  $\boldsymbol{\tau}\in {\cal E}_{ell,0}(\tilde{L},\omega)$ et tout $\nu\in i{\cal A}_{\tilde{L}}^*$, on ait la majoration
  $$\vert \varphi_{f,b,\boldsymbol{\tau}}(\nu)\vert \leq \alpha(f)(1+\vert \mu(\boldsymbol{\tau})+\vert \nu\vert )^{-N},$$
  o\`u $\mu(\boldsymbol{\tau})$ est le param\`etre infinit\'esimal de $\boldsymbol{\tau}$.

  Notons que le $r$ de (1) et la semi-norme $\alpha$ de (2) peuvent d\'ependre de $b$.
  On a une formule de descente famili\`ere
  $$J_{\tilde{M}}^{\tilde{G}}(Ind_{\tilde{Q}}^{\tilde{M}}((\tilde{\pi}_{\boldsymbol{\tau}})_{\nu+\lambda}),f)=\sum_{\tilde{M}'\in {\cal L}(\tilde{L})}d_{\tilde{L}}^{\tilde{G}}(\tilde{M},\tilde{M}')J_{\tilde{L}}^{\tilde{M}'}((\tilde{\pi}_{\boldsymbol{\tau}})_{\nu+\lambda},f_{\tilde{P}',\omega}),$$
  o\`u les $\tilde{P}'$ sont des \'el\'ements de ${\cal P}(\tilde{M}')$ d\'etermin\'es par le choix d'un param\`etre auxiliaire. Pour tout $\tilde{M}'$, il existe un sous-ensemble compact $\tilde{C}'$ de $\tilde{M}'({\mathbb R})$ tel que l'application $f\mapsto f_{\tilde{P}',\omega}$ envoie continuement $C^{\infty}(\tilde{C})$ dans $C^{\infty}(\tilde{C}')$. On ne perd rien \`a fixer $\tilde{M}'$, \`a d\'efinir une fonction $\phi_{f',b,\boldsymbol{\tau}}$ sur ${\cal A}_{\tilde{L},{\mathbb C}}^*$ pour $f'\in C^{\infty}(\tilde{C}')$ par
  $$\phi_{f',b,\boldsymbol{\tau}}(\nu)=\int_{i{\cal A}^*_{\tilde{M}}} J_{\tilde{L}}^{\tilde{M}'}((\tilde{\pi}_{\boldsymbol{\tau}})_{\nu+\lambda},f')\hat{b}(-\lambda)\,d\lambda,$$
  et \`a d\'emontrer pour cette fonction des propri\'et\'es similaires \`a (1) et (2). On a prouv\'e en [W1] 6.4 que la fonction $\phi_{f',b,\boldsymbol{\tau}}$ s'\'ecrivait sous la forme
  $$\phi_{f',b,\boldsymbol{\tau}}(\nu)=\int_{{\cal A}_{\tilde{L}}}\Psi(Y)e^{<\nu,Y>}\,dY,$$
  o\`u $\Psi$ est une fonction $C^{\infty}$ sur ${\cal A}_{\tilde{L}}$ dont le support est inclus dans un compact qui ne d\'epend que de $b$ et de $\tilde{C}'$. L'assertion (1) en r\'esulte. On a prouv\'e en [W1]  5.2 (en copiant une fois de plus Arthur) que, pour tout $N$, il existait une semi-norme $\alpha$ pour $\tilde{C}'$  telle que, pour tout  $\boldsymbol{\tau}\in {\cal E}_{ell,0}(\tilde{L},\omega)$ et tout $\nu\in i{\cal A}_{\tilde{L}}^*$, on ait la majoration
$$\vert J_{\tilde{L}}^{\tilde{M}'}((\tilde{\pi}_{\boldsymbol{\tau}})_{\nu},f')\vert \leq \alpha(f')(1+\vert \mu(\boldsymbol{\tau})\vert)^{-N}(1+\vert \nu\vert )^{-N}.$$
 Pour $\nu\in i{\cal A}_{\tilde{L}}^*$ et $\lambda\in i{\cal A}_{\tilde{M}}^*$, on a la majoration
$$\vert \nu\vert \leq \vert \nu+\lambda\vert +\vert \lambda\vert ,$$
d'o\`u
$$1+\vert \nu\vert \leq (1+\vert \nu+\lambda\vert )(1+\vert \lambda\vert ),$$
puis
$$(1+\vert \nu+\lambda\vert )^{-N}\leq (1+\vert \nu\vert )^{-N}(1+\vert \lambda\vert )^{N}.$$
La fonction $\hat{b}$ est de Schwartz, donc v\'erifie une majoration
$$\vert \hat{b}(\lambda)\vert  \leq c(1+\vert \lambda\vert )^{-N-a_{\tilde{M}}-1}.$$
Il r\'esulte de ces majorations que l'on a l'in\'egalit\'e
 $$\vert \phi_{f',b,\boldsymbol{\tau}}(\nu) \vert \leq cc'\alpha(f')(1+\vert \mu(\boldsymbol{\tau})\vert)^{-N}(1+\vert \nu\vert )^{-N},$$
 o\`u
 $$c'=\int_{i{\cal A}_{\tilde{M}}^*}(1+\vert \lambda\vert )^{-a_{\tilde{M}}-1}\,d\lambda.$$
 Cela d\'emontre la majoration (2) cherch\'ee. 
 
Fixons un ensemble de repr\'esentants $\underline{{\cal E}}_{ell,0}(\tilde{L},\omega)$ comme en [IV] 1.4 (on utilise dans ce qui suit les notations de cette r\'ef\'erence).  Les propri\'et\'es (1) et (2) jointes au lemme [IV] 1.3 montrent que, pour tout $f\in C^{\infty}(\tilde{C})$, la famille $(\varphi_{f,b,\boldsymbol{\tau}})_{\boldsymbol{\tau}\in \underline{{\cal E}}_{ell,0}(\tilde{L},\omega)}$ appartient \`a $PW^{\infty}_{ell}(\tilde{L},\omega)$. De plus, l'application qui, \`a $f$, associe cette famille, est continue. En sommant ces applications sur tout $\tilde{L}\in {\cal L}^{\tilde{M}}(\tilde{M}_{0})$, on obtient une application continue de $C^{\infty}(\tilde{C})$ dans $PW^{\infty}(\tilde{M},\omega)$ (la condition requise d'invariance par $W^M(\tilde{M}_{0})$ r\'esulte de la construction). En composant avec l'hom\'eomorphisme $pw^{-1}:PW^{\infty}(\tilde{M},\omega)\to I(\tilde{M}({\mathbb R}),\omega)$ du th\'eor\`eme 1.4 de [IV], on obtient une application continue
$$(3) \qquad \begin{array}{ccc}C^{\infty}(\tilde{C})&\to&I(\tilde{M}({\mathbb R}),\omega)\\ f&\mapsto& \phi_{\tilde{M},b}(f).\\ \end{array}$$
Par construction, on a
$$\varphi_{f,b}(\tilde{\pi})=I^{\tilde{M}}(\tilde{\pi},\phi_{\tilde{M},b}(f))$$
pour toute $\omega$-repr\'esentation temp\'er\'ee et de longueur finie $\tilde{\pi}$ de $\tilde{M}({\mathbb R})$. Pour tout $Y\in {\cal A}_{\tilde{M}}$, on a alors
$$I^{\tilde{M}}(\tilde{\pi},Y,\phi_{\tilde{M},b}(f))=\int_{i{\cal A}_{\tilde{M}}^*}I^{\tilde{M}}(\tilde{\pi}_{\lambda},\phi_{\tilde{M},b}(f))e^{-<Y,\lambda>}\,d\lambda$$
$$=\int_{i{\cal A}_{\tilde{M}}^*}\varphi_{f,b}(\tilde{\pi}_{\lambda})e^{-<Y,\lambda>}\,d\lambda\int_{i{\cal A}_{\tilde{M}}^*}\int_{{\cal A}_{\tilde{M}}}J_{\tilde{M}}^{\tilde{G}}(\tilde{\pi}_{\lambda},X,f)b(X)\,dX\,e^{-<Y,\lambda>}\,d\lambda.$$
Or, par construction de $J_{\tilde{M}}^{\tilde{G}}(\tilde{\pi},X,f)$, on a l'\'egalit\'e 
$$J_{\tilde{M}}^{\tilde{G}}(\tilde{\pi}_{\lambda},X,f)=e^{<X,\lambda>}J_{\tilde{M}}^{\tilde{G}}(\tilde{\pi},X,f).$$
Par inversion de Fourier, on obtient
$$(4) \qquad I^{\tilde{M}}(\tilde{\pi},Y,\phi_{\tilde{M},b}(f))=J_{\tilde{M}}^{\tilde{G}}(\tilde{\pi},Y,f)b(Y).$$
Remarquons que ces relations  d\'eterminent enti\`erement $\phi_{\tilde{M},b}(f)$. Soit $b'$ un autre \'el\'ement de $C_{c}^{\infty}({\cal A}_{\tilde{M}})$. On sait que
$$I^{\tilde{M}}(\tilde{\pi},Y,\phi_{\tilde{M},b}(f)(b'\circ H_{\tilde{M}}))=b'(Y)I^{\tilde{M}}(\tilde{\pi},Y,\phi_{\tilde{M},b}(f)).$$
D'o\`u 
$$I^{\tilde{M}}(\tilde{\pi},Y,\phi_{\tilde{M},b}(f)(b'\circ H_{\tilde{M}}))=J_{\tilde{M}}^{\tilde{G}}(\tilde{\pi},Y,f)b(Y)b'(Y).$$
Cette relation est sym\'etrique en $b$ et $b'$. D'o\`u
$$I^{\tilde{M}}(\tilde{\pi},Y,\phi_{\tilde{M},b}(f)(b'\circ H_{\tilde{M}}))=I^{\tilde{M}}(\tilde{\pi},Y,\phi_{\tilde{M},b'}(f)(b\circ H_{\tilde{M}})).$$
Ces relations entra\^{\i}nent l'\'egalit\'e $\phi_{\tilde{M},b}(f)(b'\circ H_{\tilde{M}})=\phi_{\tilde{M},b'}(f)(b\circ H_{\tilde{M}})$. La famille $(\phi_{\tilde{M},b}(f))_{b\in C_{c}^{\infty}({\cal A}_{\tilde{M}})}$ appartient donc \`a l'analogue $P^{\tilde{M}}$ de l'espace $P$ pour $\tilde{M}$. Il en r\'esulte l'existence d'un unique \'el\'ement $\phi_{\tilde{M}}(f)\in I_{ac}(\tilde{M}({\mathbb R}),\omega)$ tel que $\phi_{\tilde{M},b}(f)=\phi_{\tilde{M}}(f)(b\circ H_{\tilde{M}})$ pour tout $b$. D'apr\`es la d\'efinition de la topologie sur $I_{ac}(\tilde{M}({\mathbb R}),\omega)$ et parce que les applications (3) sont continues pour tout $b$, l'application $f\mapsto \phi_{\tilde{M}}(f)$ est continue. Enfin, soit $\tilde{\pi}$ une $\omega$-repr\'esentation temp\'er\'ee et de longueur finie de $\tilde{M}({\mathbb R})$ et soit $X\in {\cal A}_{\tilde{M}}$. Choisissons une fonction $b\in C_{c}^{\infty}({\cal A}_{\tilde{M}})$ telle que $b(X)=1$. On a alors
$$I^{\tilde{M}}(\tilde{\pi},X,\phi_{\tilde{M}}(f))=I^{\tilde{M}}(\tilde{\pi},X,\phi_{\tilde{M}}(f))b(X)=I^{\tilde{M}}(\tilde{\pi},X,\phi_{\tilde{M}}(f)(b\circ H_{\tilde{M}}))$$
$$=I^{\tilde{M}}(\tilde{\pi},X,\phi_{\tilde{M},b}(f))=J_{\tilde{M}}^{\tilde{G}}(\tilde{\pi},X,f)b(X)$$
d'apr\`es (4), d'o\`u, puisque $b(X)=1$,
$$I^{\tilde{M}}(\tilde{\pi},X,\phi_{\tilde{M}}(f))=J_{\tilde{M}}^{\tilde{G}}(\tilde{\pi},X,f).$$
Notre application $\phi_{\tilde{M}}$, d\'efinie sur $C^{\infty}(\tilde{C})$, v\'erifie ainsi toutes les propri\'et\'es requises. 

 L'application $\phi_{\tilde{M}}$ \'etant caract\'eris\'ee par les \'egalit\'es ci-dessus, il est clair que, si $\tilde{C}'$ est un sous-ensemble compact de $\tilde{G}({\mathbb R})$ contenant $\tilde{C}$, l'application $\phi_{\tilde{M}}$ relative \`a $\tilde{C}$ est la restriction de celle relative \`a $\tilde{C}'$. Ces applications se recollent donc en une application $\phi_{\tilde{M}}:C_{c}^{\infty}(\tilde{G}({\mathbb R}))\to I_{ac}(\tilde{M}({\mathbb R}),\omega)$. Celle-ci est continue par d\'efinition de la topologie sur $C_{c}^{\infty}(\tilde{G}({\mathbb R}))$ et parce que les applications restreintes \`a $C^{\infty}(\tilde{C})$ le sont pour tout $\tilde{C}$.
 
 Soit maintenant $f\in C_{ac}^{\infty}(\tilde{G}({\mathbb R}))$. Soit $b\in C_{c}^{\infty}({\cal A}_{\tilde{M}})$. Choisissons une fonction $b^{\tilde{G}}\in C_{c}^{\infty}({\cal A}_{\tilde{G}})$ telle que $b^{\tilde{G}}$ vaille $1$ sur la projection dans ${\cal A}_{\tilde{G}}$ du support de $b$ dans ${\cal A}_{\tilde{M}}$. D\'efinissons une fonction $\phi_{\tilde{M},b}(f)\in I(\tilde{M}({\mathbb R}),\omega)$ par
 $$\phi_{\tilde{M},b}(f)=\phi_{\tilde{M}}(f(b^{\tilde{G}}\circ H_{\tilde{G}}))(b\circ H_{\tilde{M}}).$$
 De nouveau, des consid\'erations formelles montrent que l'\'egalit\'e (4) est satisfaite. En particulier, cette d\'efinition ne d\'epend pas du choix de $b^{\tilde{G}}$. L'application $f\mapsto \phi_{\tilde{M},b}(f)$ est continue car c'est la compos\'ee des trois applications continues
 $$f\mapsto f(b^{\tilde{G}}\circ H_{\tilde{G}})$$
 de  $C_{ac}^{\infty}(\tilde{G}({\mathbb R}))$ dans  $C_{c}^{\infty}(\tilde{G}({\mathbb R}))$,
 $$f\mapsto \phi_{\tilde{M}}(f)$$
de  $C_{c}^{\infty}(\tilde{G}({\mathbb R}))$  dans $I_{ac}(\tilde{M}({\mathbb R}),\omega)$ et
 $$f\mapsto f(b\circ H_{\tilde{M}})$$
de $I_{ac}(\tilde{M}({\mathbb R}),\omega)$ dans $I(\tilde{M}({\mathbb R}),\omega)$. Le m\^eme calcul que ci-dessus montre que de ces applications $f\mapsto \phi_{\tilde{M},b}(f)$ se d\'eduit une application continue $f\mapsto \phi_{\tilde{M}}(f)$ de $C_{ac}^{\infty}(\tilde{G}({\mathbb R}))$ dans $I_{ac}(\tilde{M}({\mathbb R}),\omega)$ qui v\'erifie les propri\'et\'es de l'\'enonc\'e. $\square$

 \bigskip

 \subsection{D\'efinition des int\'egrales orbitales pond\'er\'ees}
 Dans la premi\`ere section de [II], on a d\'efini des int\'egrales orbitales pond\'er\'ees $J_{\tilde{M}}^{\tilde{G}}(\gamma,\omega,f)$. Le corps de base \'etait alors non-archim\'edien. La d\'efinition reprenait \'evidemment celle d'Arthur, tout en la modifiant l\'eg\`erement. Il n'y a rien \`a changer sur le corps de base ${\mathbb R}$: on d\'efinit de m\^eme $J_{\tilde{M}}^{\tilde{G}}(\gamma,\omega,f)$ pour un espace de Levi $\tilde{M}\in {\cal L}(\tilde{M}_{0})$ et un \'el\'ement  $\gamma\in \tilde{M}({\mathbb R})$. Plus canoniquement, notons $D_{orb}(\tilde{M}({\mathbb R}),\omega)$ l'espace des distributions sur $C_{c}^{\infty}(\tilde{M}({\mathbb R}))$ engendr\'e par les int\'egrales orbitales relatives \`a des \'el\'ements $\gamma\in \tilde{M}({\mathbb R})$. Alors on peut d\'efinir $J_{\tilde{M}}^{\tilde{G}}(\boldsymbol{\gamma},{\bf f})$ pour $\boldsymbol{\gamma}\in D_{orb}(\tilde{M}({\mathbb R}),\omega)\otimes Mes(M({\mathbb R}))^*$ et ${\bf f}\in C_{c}^{\infty}(\tilde{G}({\mathbb R}))\otimes Mes(G({\mathbb R}))$.   La diff\'erence entre le cas archim\'edien et le cas non-archim\'edien est que l'on a l'\'egalit\'e $D_{orb}(\tilde{M}(F),\omega)=D_{g\acute{e}om}(\tilde{M}(F),\omega)$ quand $F$ est non-archim\'edien, tandis que l'on a seulement une inclusion $D_{orb}(\tilde{M}({\mathbb R}),\omega)\subset D_{g\acute{e}om}(\tilde{M}({\mathbb R}),\omega)$. Rappelons (cf. [I] 5.2) que ce dernier espace est celui des formes lin\'eaires continues sur   $I(\tilde{M}({\mathbb R}),\omega)$ qui, relev\'ees en des formes lin\'eaires sur $C_{c}^{\infty}(\tilde{M}({\mathbb R}))$, sont support\'ees par un nombre fini de classes de conjugaison par $M({\mathbb R})$. 
  Les constructions de [II] ne permettent pas de d\'efinir $J_{\tilde{M}}^{\tilde{G}}(\boldsymbol{\gamma},{\bf f})$  pour $\boldsymbol{\gamma}\in D_{g\acute{e}om}(\tilde{M}({\mathbb R}),\omega)\otimes Mes(M({\mathbb R}))^*$.

 Quoi qu'il en soit, on va d\'efinir les avatars $\omega$-\'equivariants $I_{\tilde{M}}^{\tilde{G}}(\boldsymbol{\gamma},{\bf f})$ pour $\boldsymbol{\gamma}\in D_{orb}(\tilde{M}({\mathbb R}),\omega)\otimes Mes(M({\mathbb R}))^*$ et ${\bf f}\in C_{c}^{\infty}(\tilde{G}({\mathbb R}))\otimes Mes(G({\mathbb R}))$, et m\^eme pour ${\bf f}\in C_{ac}^{\infty}(\tilde{G}({\mathbb R}))\otimes Mes(G({\mathbb R}))$. Les int\'egrales orbitales pond\'er\'ees $J_{\tilde{M}}^{\tilde{G}}(\boldsymbol{\gamma},{\bf f})$ s'\'etendent \`a ${\bf f}\in C_{ac}^{\infty}(\tilde{G}({\mathbb R}))\otimes Mes(G({\mathbb R}))$: on a
 $$J_{\tilde{M}}^{\tilde{G}}(\boldsymbol{\gamma},{\bf f})=J_{\tilde{M}}^{\tilde{G}}(\boldsymbol{\gamma},{\bf f}(b\circ H_{\tilde{G}}))$$
 pour toute fonction $b\in C_{c}^{\infty}({\cal A}_{\tilde{G}})$ telle que $b$ vaille $1$ dans un voisinage de l'image par $H_{\tilde{G}}$ du support de $\boldsymbol{\gamma}$. Quand $\boldsymbol{\gamma}$ est une int\'egrale orbitale associ\'ee \`a un \'el\'ement $\gamma\in \tilde{M}({\mathbb R})$,  Arthur d\'emontre en [A2] corollaire 6.2 que cette int\'egrale est donn\'ee par une mesure sur l'orbite de $\gamma$ qui est absolument continue relativement \`a la mesure de Haar. On en d\'eduit ais\'ement que 
 $${\bf f}\mapsto J_{\tilde{M}}^{\tilde{G}}(\boldsymbol{\gamma},{\bf f})$$
 est continue sur $C_{ac}^{\infty}(\tilde{G}({\mathbb R}))\otimes Mes(G({\mathbb R}))$. On d\'efinit une application lin\'eaire continue sur $C_{ac}^{\infty}(\tilde{G}({\mathbb R}))\otimes Mes(G({\mathbb R}))$, qui se factorise en une application lin\'eaire continue sur $I_{ac}(\tilde{G}({\mathbb R}),\omega)\otimes Mes(G({\mathbb R}))$, par la formule de r\'ecurrence
 $$I_{\tilde{M}}^{\tilde{G}}(\boldsymbol{\gamma},{\bf f})=J_{\tilde{M}}^{\tilde{G}}(\boldsymbol{\gamma},{\bf f})-\sum_{\tilde{L}\in {\cal L}(\tilde{M}),\tilde{L}\not=\tilde{G}}I_{\tilde{M}}^{\tilde{L}}(\boldsymbol{\gamma},\phi_{\tilde{L}}({\bf f})).$$
 La continuit\'e de cette application est assur\'ee par ce que l'on vient de dire ci-dessus, par celle des applications analogues relatives \`a $\tilde{L}\not=\tilde{G}$ que l'on suppose par r\'ecurrence, et par la continuit\'e de nos applications $\phi_{\tilde{L}}$. Une d\'emonstration assez formelle montre que l'application ainsi d\'efinie est $\omega$-\'equivariante (cf. [A3] proposition 4.1 dans le cas non tordu). Etant continue, elle se factorise donc en une application d\'efinie sur $I_{ac}(\tilde{G}({\mathbb R}),\omega)\otimes Mes(G({\mathbb R}))$ (cf. [I] 5.2).
 
 Cette d\'efinition \'etant pos\'ee, les int\'egrales orbitales $\omega$-\'equivariantes v\'erifient les m\^emes propri\'et\'es que dans le cas non-archim\'edien. 
 
 {\bf Variante.} Supposons $\tilde{G}=G$ et ${\bf a}=1$. Supposons fix\'ee une fonction $B$ comme en [II] 1.8. On note $D_{orb,unip}(G({\mathbb R}))$ le sous-espace des \'el\'ements de $D_{orb}(G({\mathbb R}))$ \`a support unipotent.  Pour $\boldsymbol{\gamma}\in D_{orb,unip}(M({\mathbb R}))\otimes Mes(M({\mathbb R}))^*$ et ${\bf f}\in I(G({\mathbb R}))\otimes Mes(G({\mathbb R}))$, on d\'efinit la variante $I_{M}^G(\boldsymbol{\gamma},B,{\bf f})$ comme en [II] 1.8.
 
 {\bf Variante.} Supposons $(G,\tilde{G},{\bf a})$ quasi-d\'eploy\'e et \`a torsion int\'erieure. Supposons fix\'e un syt\`eme de fonctions $B$ comme en [II] 1.9. Pour $\boldsymbol{\gamma}\in D_{orb}(\tilde{M}({\mathbb R}))\otimes Mes(M({\mathbb R}))^*$ et ${\bf f}\in I(\tilde{G}({\mathbb R}))\otimes Mes(G({\mathbb R}))$, on d\'efinit la variante $I_{\tilde{M}}^{\tilde{G}}(\boldsymbol{\gamma},B,{\bf f})$ comme en [II] 1.9.
 \bigskip

 Soit ${\cal O}$ une classe de conjugaison semi-simple dans $\tilde{M}({\mathbb R})$.  Rappelons que l'on note $D_{g\acute{e}om}({\cal O},\omega)$ le sous-espace des \'el\'ements de $D_{g\acute{e}om}(\tilde{M}({\mathbb R}),\omega)$ \`a support dans l'ensemble des \'el\'ements de $\tilde{M}({\mathbb R})$ dont la partie semi-simple appartient \`a ${\cal O}$.  On dit que que  ${\cal O}$ est $\tilde{G}$-\'equisinguli\`ere si elle v\'erifie la condition $M_{\gamma}=G_{\gamma}$ pour tout $\gamma\in {\cal O}$, cf. [II] 1.2(1). Plus g\'en\'eralement, une r\'eunion finie de classes de conjugaison semi-simples sera dite $\tilde{G}$-\'equisinguli\`ere si chacune de ces classes l'est.
 
 Supposons que ${\cal O}$ soit $\tilde{G}$-\'equisinguli\`ere.  Fixons $\gamma\in {\cal O}$.  Soit $f\in C_{c}^{\infty}(\tilde{G}({\mathbb R}))$. On sait alors qu'il existe une fonction $f'\in C_{c}^{\infty}(\tilde{M}({\mathbb R}))$ et un voisinage de $\gamma$ dans $\tilde{M}({\mathbb R})$ telle que, pour $\gamma'$ dans ce voisinage, on ait l'\'egalit\'e
  $$I^{\tilde{M}}(\gamma',\omega,f')=I_{\tilde{M}}^{\tilde{G}}(\gamma',\omega,f)$$
  cf. [II] 1.7(4). Fixons un tel $f'$. Pour $\boldsymbol{\gamma}\in D_{g\acute{e}om}( {\cal O},\omega)$, on peut d\'efinir $I^{\tilde{M}}(\boldsymbol{\gamma},f')$. La description que l'on a donn\'ee en [I] 5.2 de l'espace $D_{g\acute{e}om}( {\cal O},\omega)$ montre que ce terme
  ne d\'epend pas du choix de $f'$: il ne d\'epend que des  $I_{\tilde{M}}^{\tilde{G}}(\gamma',\omega,f)$ pour $\gamma'$ fortement r\'egulier voisin de $\gamma$. Plus pr\'ecis\'ement, avec les notations de [I] 5.2, il ne d\'epend que des fonctions $X\mapsto I_{\tilde{M}}^{\tilde{G}}(exp(X)\gamma,\omega,f)$ pour $X\in \mathfrak{t}^{\theta}({\mathbb R})$ tel que $exp(X)\gamma$ soit fortement r\'egulier, o\`u $\tilde{T}\in \tilde{{\cal T}}$. On peut donc poser
  $$I_{\tilde{M}}^{\tilde{G}}(\boldsymbol{\gamma},f)=I^{\tilde{M}}(\boldsymbol{\gamma},f').$$
  Il est  clair que, quand $\boldsymbol{\gamma}$ est une simple int\'egrale orbitale, cette d\'efinition co\"{\i}ncide avec la d\'efinition initiale.  
  
  Notons $D_{g\acute{e}om,\tilde{G}-\acute{e}qui}(\tilde{M}({\mathbb R}),\omega)$ le sous-espace de $D_{g\acute{e}om}(\tilde{M}({\mathbb R}),\omega)$ engendr\'e par les $D_{g\acute{e}om}({\cal O},\omega)$ pour les orbites ${\cal O}$ qui sont $\tilde{G}$-\'equisinguli\`eres.  En r\'eintroduisant les espaces de mesures de Haar, ce qui pr\'ec\`ede permet de d\'efinir   $I_{\tilde{M}}^{\tilde{G}}(\boldsymbol{\gamma},{\bf f})$ pour $\boldsymbol{\gamma}\in D_{g\acute{e}om,\tilde{G}-\acute{e}qui}( \tilde{M}({\mathbb R}),\omega)\otimes Mes(M({\mathbb R}))^*$ et ${\bf f}\in C_{c}^{\infty}(\tilde{G}({\mathbb R}))\otimes Mes(G({\mathbb R}))$ (ou ${\bf f}\in I(\tilde{G}({\mathbb R}),\omega)\otimes Mes(G({\mathbb R}))$). Dans les paragraphes suivants, on montrera que les constructions du cas non-archim\'edien s'\'etendent \`a ces termes. Montrons ici qu'ils se comportent de la fa\c{c}on attendue relativement \`a l'induction.
  
  \ass{Lemme}{Soit $\tilde{R}$ un espace de Levi contenu dans $\tilde{M}$.  Soit ${\cal O}$ une orbite semi-simple dans $\tilde{R}({\mathbb R})$. Notons ${\cal O}^{\tilde{M}}$ l'unique orbite pour l'action de $\tilde{M}({\mathbb R}) $  qui contient ${\cal O}$. On suppose  que ${\cal O}^{\tilde{M}}$ est $\tilde{G}$-\'equisinguli\`ere.
  
  (i) Pour tout $\tilde{L}\in \tilde{L} (\tilde{R})$ tel que $d_{\tilde{R}}^{\tilde{G}}(\tilde{M},\tilde{L})\not=0$, l'orbite ${\cal O}$ est $\tilde{L}$-\'equisinguli\`ere.
  
  (ii) Pour tout $\boldsymbol{\gamma}\in D_{g\acute{e}om}({\cal O},\omega)\otimes Mes(R({\mathbb R}))^*$ et tout ${\bf f}\in I(\tilde{G}({\mathbb R}),\omega)\otimes Mes(G({\mathbb R}))$, on a l'\'egalit\'e
  $$I_{\tilde{M}}^{\tilde{G}}(\boldsymbol{\gamma}^{\tilde{M}},{\bf f})=\sum_{\tilde{L}\in {\cal L}(\tilde{R})}d_{\tilde{R}}^{\tilde{G}}(\tilde{M},\tilde{L})I_{\tilde{R}}^{\tilde{L}}(\boldsymbol{\gamma},{\bf f}_{\tilde{L},\omega}).$$}
  
  {\bf Remarque.} L'assertion (i) donne un sens \`a la formule du (ii).
  
  \bigskip
  
  Preuve. L'assertion (i) est la remarque de [II] 2.11. Prouvons (ii). Oublions les mesures de Haar pour \^etre plus clair et notons $f$ plut\^ot que ${\bf f}$. Fixons encore $\gamma\in {\cal O}$. Introduisons  une fonction $g\in I(\tilde{M}({\mathbb R}),\omega)$ telle que, pour $\gamma'\in \tilde{M}({\mathbb R})$ assez r\'egulier et assez voisin de $\gamma$, on ait
 $$I_{\tilde{M}}^{\tilde{G}}(\gamma',\omega,f)=I^{\tilde{M}}(\gamma',\omega,g).$$
 De m\^eme, pour tout $\tilde{L}$ comme en (i), introduisons une fonction $g[\tilde{L}]\in I(\tilde{R}({\mathbb R}),\omega)$ telle que, pour $\gamma'\in \tilde{R}({\mathbb R})$ assez r\'egulier et assez voisin de $\gamma$, on ait
 $$I_{\tilde{R}}^{\tilde{L}}(\gamma',\omega,f_{\tilde{L},\omega})=I^{\tilde{R}}(\gamma',\omega,g[\tilde{L}]).$$
 On va montrer que
 
 (1) les deux \'el\'ements $g_{\tilde{R},\omega}$ et
 $$\sum_{\tilde{L}\in {\cal L}(\tilde{R})}d_{\tilde{R}}^{\tilde{G}}(\tilde{M},\tilde{L})g[\tilde{L}]$$
 de $I(\tilde{R}({\mathbb R}),\omega)$ co\"{\i}ncident au voisinage de $\gamma$.

 Pour $\gamma'\in \tilde{R}({\mathbb R})$ assez r\'egulier, on a par d\'efinition du terme constant
 $$I^{\tilde{R}}(\gamma',\omega,g_{\tilde{R},\omega})=I^{\tilde{M}}(\gamma',\omega,g),$$
 autrement dit, si $\gamma'$ est assez proche de $\gamma$,
$$I^{\tilde{R}}(\gamma',\omega,g_{\tilde{R},\omega})= I_{\tilde{M}}^{\tilde{G}}(\gamma',\omega,f).$$
La formule de l'\'enonc\'e que l'on cherche \`a prouver est valable dans le cas o\`u $\boldsymbol{\gamma}$ est une simple int\'egrale orbitale. En l'appliquant \`a l'int\'egrale orbitale associ\'ee \`a $\gamma'$, on obtient 
  $$I_{\tilde{M}}^{\tilde{G}}(\gamma',\omega,f)=\sum_{\tilde{L}\in {\cal L}(\tilde{R})}d_{\tilde{R}}^{\tilde{G}}(\tilde{M},\tilde{L})I_{\tilde{R}}^{\tilde{L}}(\gamma',\omega,f_{\tilde{L},\omega}),$$
  ou encore
  $$I^{\tilde{R}}(\gamma',\omega,g_{\tilde{R},\omega})=\sum_{\tilde{L}\in {\cal L}(\tilde{R})}d_{\tilde{R}}^{\tilde{G}}(\tilde{M},\tilde{L})I^{\tilde{R}}(\gamma',\omega,g[\tilde{L}])$$
  si $\gamma'$ est assez proche de $\gamma$. Cela prouve (1).

  Par d\'efinition
 $$I_{\tilde{M}}^{\tilde{G}}(\boldsymbol{\gamma}^{\tilde{M}},f)=I^{\tilde{M}}(\boldsymbol{\gamma}^{\tilde{M}},g).$$
 Par d\'efinition de l'induction des distributions,
 $$I^{\tilde{M}}(\boldsymbol{\gamma}^{\tilde{M}},g)=I^{\tilde{R}}(\boldsymbol{\gamma},g_{\tilde{R},\omega}).$$
 En appliquant (1), on obtient
 $$I_{\tilde{M}}^{\tilde{G}}(\boldsymbol{\gamma}^{\tilde{M}},f)=\sum_{\tilde{L}\in {\cal L}(\tilde{R})}d_{\tilde{R}}^{\tilde{G}}(\tilde{M},\tilde{L})I^{\tilde{R}}(\boldsymbol{\gamma},g[\tilde{L}]).$$
 Mais, pour tout $\tilde{L}$, on a par d\'efinition 
 $$I^{\tilde{R}}(\boldsymbol{\gamma},g[\tilde{L}])=I_{\tilde{R}}^{\tilde{L}}(\boldsymbol{\gamma},f_{\tilde{L},\omega})$$
 et la formule pr\'ec\'edente devient celle de l'\'enonc\'e. $\square$

 \bigskip
 
 \subsection{Int\'egrales orbitales pond\'er\'ees invariantes stables}
 On suppose dans ce paragraphe $(G,\tilde{G},\omega)$ quasi-d\'eploy\'e et \`a torsion int\'erieure. Soit $\tilde{M}$  un espace de Levi  de $\tilde{G}$.  Soit ${\cal O}$  une classe de conjugaison stable semi-simple dans $\tilde{M}({\mathbb R})$. D'apr\`es la d\'efinition du paragraphe pr\'ec\'edent, elle est $\tilde{G}$-\'equisinguli\`ere si toutes les classes de conjugaison par $\tilde{M}({\mathbb R})$ incluses dans ${\cal O}$ le sont. Remarquons que, puisque la d\'efinition de cette notion est g\'eom\'etrique,  dire que toutes ces classes sont $\tilde{G}$-\'equisinguli\`eres \'equivaut \`a dire que l'une d'elles l'est.  On note $D_{g\acute{e}om,\tilde{G}-\acute{e}qui}^{st}(\tilde{M}({\mathbb R}))$ la somme des  $D^{st}_{g\acute{e}om}({\cal O})$, o\`u ${\cal O}$ parcourt les classes de conjugaison stable semi-simples  $\tilde{G}$-\'equisinguli\`eres dans $\tilde{M}({\mathbb R})$. Pour $\boldsymbol{\delta}\in D_{g\acute{e}om,\tilde{G}-\acute{e}qui}^{st}(\tilde{M}({\mathbb R}))\otimes Mes(M({\mathbb R}))^*$ et pour ${\bf f}\in I(\tilde{G}({\mathbb R}))\otimes Mes(G({\mathbb R}))$, on d\'efinit l'int\'egrale orbitale pond\'er\'ee invariante stable par
 $$(1) \qquad S_{\tilde{M}}^{\tilde{G}}(\boldsymbol{\delta},{\bf f})=I_{\tilde{M}}^{\tilde{G}}(\boldsymbol{\delta},{\bf f})-\sum_{s\in Z(\hat{M})^{\Gamma_{{\mathbb R}}}/Z(\hat{G})^{\Gamma_{{\mathbb R}}}; s\not=1} i_{\tilde{M}}(\tilde{G},\tilde{G}'(s))S_{{\bf M}}^{{\bf G}'(s)}(\boldsymbol{\delta},{\bf f}^{{\bf G}'(s)}).$$
 On renvoie \`a [II] 1.10 pour diverses notations. Expliquons cette formule.  Puisque 
 $$D_{g\acute{e}om,\tilde{G}-\acute{e}qui}^{st}(\tilde{M}({\mathbb R}))\subset D_{g\acute{e}om,\tilde{G}-\acute{e}qui}(\tilde{M}({\mathbb R})),$$
  le premier terme $I_{\tilde{M}}^{\tilde{G}}(\boldsymbol{\delta},{\bf f})$ a \'et\'e  d\'efini au paragraphe pr\'ec\'edent. Soit $s\in Z(\hat{M})^{\Gamma_{{\mathbb R}}}/Z(\hat{G})^{\Gamma_{{\mathbb R}}}$. Le groupe $G'(s)$ n'est pas en g\'en\'eral un sous-groupe de $G$, mais son syst\`eme de racines est un sous-syst\`eme de celui de $G$. Il en r\'esulte ais\'ement que, pour $\gamma\in \tilde{M}({\mathbb R})$,  le syst\`eme de racines de $G'(s)_{\gamma}$ est un sous-syst\`eme de celui de $G_{\gamma}$. L'\'egalit\'e $G_{\gamma}=M_{\gamma}$ force $G'(s)_{\gamma}=M_{\gamma}$. Cela entra\^{\i}ne que $\boldsymbol{\delta}$ appartient \`a $D_{g\acute{e}om,\tilde{G}(s)-\acute{e}qui}^{st}(\tilde{M}({\mathbb R}))\otimes Mes(M({\mathbb R}))^*$. Pour $s\not=1$, la dimension de $G'(s)_{SC}$ est strictement inf\'erieure \`a celle de $G$. En raisonnant comme toujours par r\'ecurrence sur cette dimension, on peut supposer d\'efinies les int\'egrales orbitales pond\'er\'ees invariantes stables pour l'espace $\tilde{G}'(s)$. Pour d\'efinir le terme $S_{{\bf M}}^{{\bf G}'(s)}(\boldsymbol{\delta},{\bf f}^{{\bf G}'(s)})$, on a encore  besoin, d'une part, de propri\'et\'es   qui permettent d'\'etendre les d\'efinitions au cadre formel que l'on a introduit en [I] pour les donn\'ees endoscopiques. On reviendra ci-dessous sur ces propri\'et\'es. On a besoin d'autre part de supposer connu par r\'ecurrence le th\'eor\`eme ci-dessous.
 
 \ass{Th\'eor\`eme }{Pour  $\boldsymbol{\delta}\in D_{g\acute{e}om,\tilde{G}-\acute{e}qui}^{st}(\tilde{M}({\mathbb R}))\otimes Mes(M({\mathbb R}))^*$, la distribution ${\bf f}\mapsto S_{\tilde{M}}^{\tilde{G}}(\boldsymbol{\delta},{\bf f})$ est stable.}
 
 Ce th\'eor\`eme sera prouv\'e dans le paragraphe suivant.
 
 Abandonnons pour simplifier les espaces de mesures de Haar. Soit ${\cal O}$ une classe de conjugaison stable semi-simple et $\tilde{G}$-\'equisinguli\`ere dans $\tilde{M}({\mathbb R})$. Soit $f\in I(\tilde{G}({\mathbb R}))$. Notons $D_{orb}^{st}(\tilde{M}({\mathbb R}))$ l'espace des combinaisons lin\'eaires stables d'int\'egrales orbitales sur $\tilde{M}({\mathbb R})$. On montrera plus loin que
 
 (2) il existe une fonction $f'\in SI(\tilde{M}({\mathbb R}))$ et un voisinage de ${\cal O}$ dans $\tilde{M}({\mathbb R})$ tels que, pour tout  $\boldsymbol{\tau}\in D_{orb}^{st}(\tilde{M}({\mathbb R}))$ dont le support est form\'e  d'\'el\'ements $\tilde{G}$-r\'eguliers dans ce voisinage, on a l'\'egalit\'e
 $$S_{\tilde{M}}^{\tilde{G}}(\boldsymbol{\tau},f)=S^{\tilde{M}}(\boldsymbol{\tau},f');$$
 
 (3) pour $f'$ comme ci-dessus et pour tout $\boldsymbol{\delta}\in D_{g\acute{e}om}^{st}({\cal O})$, on a l'\'egalit\'e
 $$S_{\tilde{M}}^{\tilde{G}}(\boldsymbol{\delta},f)=S^{\tilde{M}}(\boldsymbol{\delta},f').$$
 
 Commen\c{c}ons \`a expliquer les propri\'et\'es formelles requises. Puisqu'on  va les appliquer par r\'ecurrence, il est l\'egitime de supposer connues toutes les propri\'et\'es n\'ecessaires des int\'egrales pond\'er\'ees invariantes stables telles qu'on les a d\'efinies ci-dessus. Consid\'erons des extensions compatibles
 $$1\to C_{\natural}\to G_{\natural}\to G\to 1\,\text{ et }\, \tilde{G}_{\natural}\to \tilde{G}$$
 o\`u $C_{\natural}$ est un tore central induit et o\`u $\tilde{G}_{\natural}$ est \`a torsion int\'erieure. Consid\'erons un caract\`ere $\lambda_{\natural}$ de $C_{\natural}({\mathbb R})$. 
 On note $\tilde{M}_{\natural}$ l'image r\'eciproque de $\tilde{M}$ dans $\tilde{G}_{\natural}$. Soit ${\cal O}$ une classe de conjugaison stable semi-simple et $\tilde{G}$-\'equisinguli\`ere dans $\tilde{M}({\mathbb R})$.  Fixons une classe de conjugaison stable ${\cal O}_{\natural}$ dans $\tilde{M}_{\natural}({\mathbb R})$ qui se projette sur ${\cal O}$. Soient $\boldsymbol{\delta}\in D_{g\acute{e}om,\lambda_{\natural}}^{st}(\tilde{M}_{\natural}({\mathbb R}),{\cal O})$ et $f\in I_{\lambda_{\natural}}(\tilde{G}_{\natural}({\mathbb R}))$.  Comme on l'a dit en [I] 5.6, il y a une application lin\'eaire surjective
 $$D_{g\acute{e}om}^{st}({\cal O}_{\natural})\to D_{g\acute{e}om,\lambda_{\natural}}^{st}(\tilde{M}_{\natural}({\mathbb R}),{\cal O}).$$
 Fixons une image r\'eciproque $\dot{\boldsymbol{\delta}}\in D_{g\acute{e}om}^{st}({\cal O}_{\natural})$ de $\boldsymbol{\delta}$. Il est facile de voir que la forme lin\'eaire $\varphi\mapsto S_{\tilde{M}_{\natural}}^{G_{\natural}}(\dot{\boldsymbol{\delta}},\varphi)$ sur $I(\tilde{G}_{\natural}({\mathbb R}))$ s'\'etend en une forme lin\'eaire sur $I_{ac}(\tilde{G}_{\natural}({\mathbb R}))$ par la formule 
 $$S_{\tilde{M}_{\natural}}^{\tilde{G}_{\natural}}(\dot{\boldsymbol{\delta}},\varphi)=S_{\tilde{M}_{\natural}}^{\tilde{G}_{\natural}}(\dot{\boldsymbol{\delta}},\varphi(b\circ H_{\tilde{G}_{\natural}})),$$
 o\`u $b\in C_{c}^{\infty}({\cal A}_{\tilde{G}_{\natural}})$ vaut $1$ dans un voisinage de l'image par $H_{\tilde{G}_{\natural}}$ de ${\cal O}_{\natural}$ (cf. [II] 1.10(2)).  La fonction $f$ appartient \`a $I_{ac}(\tilde{G}_{\natural}({\mathbb R}))$. On peut donc d\'efinir
 $S_{\tilde{M}_{\natural}}^{\tilde{G}_{\natural}}(\dot{\boldsymbol{\delta}},f)$. On doit prouver
 
 (4) $S_{\tilde{M}_{\natural}}^{\tilde{G}_{\natural}}(\dot{\boldsymbol{\delta}},f)$ ne d\'epend pas des choix de ${\cal O}_{\natural}$ et du rel\`evement $\dot{\boldsymbol{\delta}}$.
 
 Ainsi, on peut noter ce terme $S_{\tilde{M}_{\natural},\lambda_{\natural}}^{\tilde{G}_{\natural}}(\boldsymbol{\delta},f)$.
 
 Consid\'erons maintenant un autre couple d'extensions
 $$1\to C_{\flat}\to G_{\flat}\to G \to 1\text{ et }\tilde{G}_{\flat}\to \tilde{G}$$
 ainsi qu'un caract\`ere $\lambda_{\flat}$ de $C_{\flat}({\mathbb R})$  v\'erifiant les m\^emes hypoth\`eses que ci-dessus. On note $G_{\natural,\flat}$, resp.  $\tilde{G}_{\natural,\flat}$, les produits fibr\'es de $G_{\natural}$ et $G_{\flat}$ au-dessus de $G$, resp. de $\tilde{G}_{\natural}$ et $\tilde{G}_{\flat}$ au-dessus de $\tilde{G}$. Consid\'erons un caract\`ere $\lambda_{\natural,\flat}$ de $G_{\natural,\flat}({\mathbb R})$ dont la restriction \`a $C_{\natural}({\mathbb R})\times C_{\flat}({\mathbb R})$ soit $\lambda_{\natural}^{-1}\times \lambda_{\flat}$. Consid\'erons une fonction $\tilde{\lambda}_{\natural,\flat}$ sur $\tilde{G}_{\natural,\flat}({\mathbb R})$ qui se transforme selon le caract\`ere $\lambda_{\natural,\flat}$ (cf. [I] 2.6(i)). On d\'eduit de ces donn\'ees des isomorphismes 
 $$C_{c,\lambda_{\natural}}^{\infty}(\tilde{G}_{\natural}({\mathbb R}))\simeq C_{c,\lambda_{\flat}}^{\infty}(\tilde{G}_{\flat}({\mathbb R}))$$
 et
 $$D_{g\acute{e}om,\lambda_{\natural}}^{st}(\tilde{M}_{\natural}({\mathbb R}),{\cal O})\simeq D_{g\acute{e}om,\lambda_{\flat}}^{st}(\tilde{M}_{\flat}({\mathbb R}),{\cal O})$$
 cf. [II] 1.10. Pour $f_{\natural}$ et $f_{\flat}$, resp. $\boldsymbol{\delta}_{\natural}$ et $\boldsymbol{\delta}_{\flat}$, se correspondant par ces isomorphismes, on doit prouver que
 
 (5) $S_{\tilde{M}_{\natural},\lambda_{\natural}}^{\tilde{G}_{\natural}}(\boldsymbol{\delta}_{\natural},f_{\natural})=S_{\tilde{M}_{\flat},\lambda_{\flat}}^{\tilde{G}_{\flat}}(\boldsymbol{\delta}_{\flat},f_{\flat})$.
 
 En [II] 1.10, on a v\'erifi\'e les propri\'et\'es (4) et (5) lorsque le corps de base \'etait non-archim\'edien. On a vu que, gr\^ace \`a elles,  on d\'efinissait ais\'ement les termes $S_{{\bf M}'}^{{\bf G}'(s)}(\boldsymbol{\delta},{\bf f}^{{\bf G}'(s)})$ qui interviennent dans la formule (1). Les preuves de cette r\'ef\'erence restent correctes sur le corps de base ${\mathbb R}$, pourvu qu'on se limite \`a consid\'erer de v\'eritables int\'egrales orbitales. On ne les refait pas dans ce cas et on tient pour acquis que (4)  et (5) sont v\'erifi\'ees si les \'el\'ements $\dot{\boldsymbol{\delta}}$, $\boldsymbol{\delta}_{\natural}$ et $\boldsymbol{\delta}_{\flat}$ sont des int\'egrales orbitales stables. On va en d\'eduire que ces relations sont v\'erifi\'ees pour des \'el\'ements g\'en\'eraux.  
 
 Dans la situation de (4), on commence par montrer
 
 (6) il existe $f'_{\natural}\in C_{c,\lambda_{\natural}}^{\infty}(\tilde{M}_{\natural}({\mathbb R}))$ et un voisinage $\Omega$ de ${\cal O}$ dans $\tilde{M}({\mathbb R})$ tels que, pour tout $\dot{\boldsymbol{\tau}}\in D_{orb}^{st}(\tilde{M}_{\natural}({\mathbb R}))$ dont le support est form\'e d'\'el\'ements $\tilde{G}_{\natural}$-r\'eguliers se projetant dans $\Omega$, on ait l'\'egalit\'e
 $$S_{\tilde{M}_{\natural}}^{\tilde{G}_{\natural}}(\dot{\boldsymbol{\tau}},f)=S^{\tilde{M}_{\natural}}(\dot{\boldsymbol{\tau}},f'_{\natural}).$$
 
  Introduisons le groupe d\'eriv\'e $M_{\natural,der}$. On peut introduire un sous-tore $S$ de $Z(M_{\natural})^0$ de sorte que $\mathfrak{c}_{\natural}\oplus \mathfrak{s}=\mathfrak{z}_{M_{\natural}}$. 
   L'application produit
 $$(7) \qquad C_{\natural}({\mathbb R})\times S({\mathbb R}) \times M_{\natural,der}({\mathbb R})\to M_{\natural}({\mathbb R})$$
 est un homomorphisme de noyau fini et d'image un sous-groupe ouvert de $M_{\natural}({\mathbb R})$. Notons $\Delta$ le noyau et $\Sigma$ sa projection dans $C_{\natural}({\mathbb R})$. Le groupe $\Sigma$ contient le groupe $\Sigma'$ des $c\in C_{\natural}({\mathbb R})$ tels que $c{\cal O}_{\natural}={\cal O}_{\natural}$.
 D'apr\`es (2), on peut trouver $f''\in C_{c}^{\infty}(\tilde{M}_{\natural}({\mathbb R}))$ de sorte que, pour $\dot{\boldsymbol{\tau}}\in D_{orb}^{st}(\tilde{M}_{\natural}({\mathbb R}))$ dont le support est form\'e d'\'el\'ements $\tilde{G}_{\natural}$-r\'eguliers  dans un voisinage de ${\cal O}_{\natural}$, on ait l'\'egalit\'e
 $$S_{\tilde{M}_{\natural}}^{\tilde{G}_{\natural}}(\dot{\boldsymbol{\tau}},f)=S^{\tilde{M}_{\natural}}(\dot{\boldsymbol{\tau}},f'').$$
  On peut modifier $f''$ de sorte qu'elle soit nulle au voisinage de $c{\cal O}_{\natural}$ pour $c\in \Sigma-\Sigma'$.
  Le groupe $C_{\natural}({\mathbb R})$ agit sur les espaces de fonctions sur $\tilde{G}_{\natural}({\mathbb R})$ ou $\tilde{M}_{\natural}({\mathbb R})$ par $(\varphi,c)\mapsto \varphi^c$, o\`u $\varphi^c(\gamma)=\varphi(c\gamma)$. Il agit par dualit\'e sur les espaces de distributions. Soit $c\in \Sigma$ et  $\dot{\boldsymbol{\tau}}\in D_{orb}^{st}(\tilde{M}_{\natural}({\mathbb R}))$ dont le support est form\'e d'\'el\'ements $\tilde{G}_{\natural}$-r\'eguliers  dans un voisinage de ${\cal O}_{\natural}$. On a l'\'egalit\'e
  $$S^{\tilde{M}_{\natural}}(\dot{\boldsymbol{\tau}},(f'')^c)=S^{\tilde{M}_{\natural}}(\dot{\boldsymbol{\tau}}^{c^{-1}},f'').$$
  Si $c\not\in \Sigma'$, c'est nul . Si $c\in \Sigma'$, c'est 
  $$S_{\tilde{M}_{\natural}}^{\tilde{G}_{\natural}}(\dot{\boldsymbol{\tau}}^{c^{-1}},f)=S_{\tilde{M}_{\natural}}^{\tilde{G}_{\natural}}(\dot{\boldsymbol{\tau}},f^c)=\lambda_{\natural}(c)^{-1}S_{\tilde{M}_{\natural}}^{\tilde{G}_{\natural}}(\dot{\boldsymbol{\tau}},f),$$
  la deuxi\`eme \'egalit\'e ci-dessus se v\'erifiant  comme en [II] 1.10. On peut donc aussi bien remplacer $f''$ par
  $$\vert \Sigma'\vert ^{-1}\sum_{c\in \Sigma}\lambda_{\natural}(c)(f'')^c.$$
Fixons $\delta_{\natural}\in {\cal O}_{\natural}$.  D\'efinissons une fonction $f'_{\natural}$ sur $\tilde{M}_{\natural}({\mathbb R})$ de la fa\c{c}on suivante. Sur un \'el\'ement $m\delta_{\natural}$ o\`u $m\in M_{\natural}({\mathbb R})$ n'appartient pas \`a l'image de (7), on pose $f'_{\natural}(m\delta_{\natural})=0$. Pour    
  $$(c,s,x)\in C_{\natural}({\mathbb R})\times S({\mathbb R}) \times M_{\natural,der}({\mathbb R})$$
  On pose
  $$f'_{\natural}(csx\delta_{\natural})=\lambda_{\natural}(c)^{-1}f''(sx\delta_{\natural}).$$
  La d\'efinition est loisible: le membre de droite ne d\'epend que du produit $csx$. La fonction ainsi d\'efinie appartient \`a $C_{c,\lambda_{\natural}}^{\infty}(\tilde{M}_{\natural}({\mathbb R}))$. Soit $\dot{\boldsymbol{\tau}}\in D_{orb}^{st}(\tilde{M}_{\natural}({\mathbb R}))$ dont le support est une classe de conjugaison stable d'\'el\'ements $\tilde{G}_{\natural}$-r\'eguliers se projetant dans un voisinage assez petit de ${\cal O}$. On peut choisir  $c\in C_{\natural}({\mathbb R})$ tel que le support de $\dot{\boldsymbol{\tau}}^{c}$ soit contenu dans un petit voisinage de ${\cal O}_{\natural}$. On peut modifier $c$ de sorte que ce support soit aussi contenu dans l'ensemble $\{sx\delta_{\natural}; s\in S({\mathbb R}), x\in M_{\natural,der}({\mathbb R})\}$. On a alors
  $$ S^{\tilde{M}_{\natural}}(\dot{\boldsymbol{\tau}},f'_{\natural})=S^{\tilde{M}_{\natural}}(\dot{\boldsymbol{\tau}}^{c},(f'_{\natural})^c)=\lambda_{\natural}(c)^{-1}S^{\tilde{M}_{\natural}}(\dot{\boldsymbol{\tau}}^c,f'_{\natural}).$$
  D'apr\`es la propri\'et\'e du support de  $\dot{\boldsymbol{\tau}}^{c}$, on a
$$S^{\tilde{M}_{\natural}}(\dot{\boldsymbol{\tau}}^c,f'_{\natural})=S^{\tilde{M}_{\natural}}(\dot{\boldsymbol{\tau}}^c,f'')=  S_{\tilde{M}_{\natural}}^{\tilde{G}_{\natural}}(\dot{\boldsymbol{\tau}}^c,f)=S_{\tilde{M}_{\natural}}^{\tilde{G}_{\natural}}(\dot{\boldsymbol{\tau}},f^{c^{-1}})=\lambda_{\natural}(c)S_{\tilde{M}_{\natural}}^{\tilde{G}_{\natural}}(\dot{\boldsymbol{\tau}},f).$$
Cette suite d'\'egalit\'es montre que $f$ satisfait les conditions de (6).

Dans la situation de (4),  appliquons (3) \`a l'espace $\tilde{G}_{\natural}$, \`a la fonction $f$ et \`a  la fonction $f'_{\natural}$ fournie par (6) (on peut les multiplier par une fonction $b\circ H_{\tilde{G}_{\natural}}$ comme ci-dessus pour les remplacer par des fonctions \`a support compact). On obtient l'\'egalit\'e
$$S_{\tilde{M}_{\natural}}^{\tilde{G}_{\natural}}(\dot{\boldsymbol{\delta}},f)=S^{\tilde{M}_{\natural}}(\dot{\boldsymbol{\delta}},f'_{\natural}).$$
Par d\'efinition, ce dernier terme n'est autre que $S_{\lambda_{\natural}}^{\tilde{M}_{\natural}}(\boldsymbol{\delta},f'_{\natural})$. Il ne d\'epend pas du choix de $\dot{\boldsymbol{\delta}}$ donc le membre de gauche ci-dessus non plus. Cela  prouve (4). 

Dans la situation de (5), on construit de m\^eme des fonctions $f'_{\natural}$ et $f'_{\flat}$. La relation qui les d\'efinit et le fait que (5) soit v\'erifi\'e pour des \'el\'ements \`a support r\'egulier implique que ces deux fonctions se correspondent par l'isomorphisme
$$SI_{\lambda_{\natural}}(\tilde{M}_{\natural}({\mathbb R}))\simeq SI_{\lambda_{\flat}}(\tilde{M}_{\flat}({\mathbb R})).$$
Puisque $\boldsymbol{\delta}_{\natural}$ et $\boldsymbol{\delta}_{\flat}$ se correspondent, cela entra\^{\i}ne
$$S^{\tilde{M}_{\natural}}_{\lambda_{\natural}}(\boldsymbol{\delta}_{\natural},f'_{\natural})=S^{\tilde{M}_{\flat}}_{\lambda_{\flat}}(\boldsymbol{\delta}_{\flat},f'_{\flat}).$$
Mais on vient  de voir que
$$S_{\tilde{M}_{\natural}}^{G_{\natural}}(\boldsymbol{\delta}_{\natural},f_{\natural})=S^{\tilde{M}_{\natural}}_{\lambda_{\natural}}(\boldsymbol{\delta}_{\natural},f'_{\natural})$$
et
$$S_{\tilde{M}_{\flat}}^{G_{\flat}}(\boldsymbol{\delta}_{\flat},f_{\flat})=S^{\tilde{M}_{\flat}}_{\lambda_{\flat}}(\boldsymbol{\delta}_{\flat},f'_{\flat}).$$
L'\'egalit\'e (5) en r\'esulte.

Soit $s\in Z(\hat{M})^{\Gamma_{{\mathbb R}}}/Z(\hat{G})^{\Gamma_{{\mathbb R}}}$ avec $s\not=1$. La d\'emonstration ci-dessus, en particulier le fait que les fonctions $f'_{\natural}$ et $f'_{\flat}$ se correspondent, entra\^{\i}ne que les assertions (2) et (3) se g\'en\'eralisent sous la forme suivante. On fixe cette fois $f\in I({\bf G}'(s))$.

 (8) Il existe une fonction $f'\in SI({\bf M})$ et un voisinage de ${\cal O}$ dans $\tilde{M}({\mathbb R})$ tels que, pour tout  $\boldsymbol{\tau}\in D_{orb}^{st}({\bf M})$ dont le support est form\'e  d'\'el\'ements $\tilde{G}$-r\'eguliers dans ce voisinage, on a l'\'egalit\'e
 $$S_{{\bf M}}^{{\bf G}'(s)}(\boldsymbol{\tau},f)=S^{{\bf M}}(\boldsymbol{\tau},f');$$
 
 (9) pour $f'$ comme ci-dessus et pour tout $\boldsymbol{\delta}\in D_{g\acute{e}om}^{st}({\bf M},{\cal O})$ (avec une d\'efinition naturelle de cet espace), on a l'\'egalit\'e
 $$S_{{\bf M}}^{{\bf G}'(s)}(\boldsymbol{\delta},f)=S^{{\bf M}}(\boldsymbol{\delta},f').$$
 
 Venons-en \`a la preuve de (2) et (3). En fait, on peut identifier $SI({\bf M})$ et $SI(\tilde{M}({\mathbb R}))$ ainsi que $D_{g\acute{e}om}^{st}({\bf M},{\cal O})$ et $D_{g\acute{e}om}^{st}({\cal O})$. Pour $s\in Z(\hat{M})^{\Gamma_{{\mathbb R}}}/Z(\hat{G})^{\Gamma_{{\mathbb R}}}$ avec $s\not=1$, on applique les relations ci-dessus  \`a la fonction $f^{{\bf G}'(s)}$, on en d\'eduit une fonction $f'_{s}\in SI(\tilde{M}({\mathbb R}))$. En appliquant [II] 1.7(4) au terme $I_{\tilde{M}}^{\tilde{G}}(\boldsymbol{\delta},f)$, on a une fonction $f'_{0}\in I(\tilde{M}({\mathbb R}))$ qui v\'erifie des propri\'et\'es analogues. En fait, puisqu'on ne lui applique ici que des distributions $\boldsymbol{\delta}$ qui sont stables, on peut remplacer $f'_{0}$ par son image dans $ SI(\tilde{M}({\mathbb R}))$. On pose
 $$f'=f'_{0}-\sum_{s\in Z(\hat{M})^{\Gamma_{{\mathbb R}}}/Z(\hat{G})^{\Gamma_{{\mathbb R}}}; s\not=1} i_{\tilde{M}}(\tilde{G},\tilde{G}'(s))f'_{s}.$$
 Cette fonction v\'erifie (2) et (3). 
 
 \bigskip
 
 \subsection{Preuve du th\'eor\`eme 1.4}
 On consid\`ere d'abord un \'el\'ement $\boldsymbol{\delta}\in D_{orb}^{st}(\tilde{M}({\mathbb R}))\otimes Mes(M({\mathbb R}))^*$ dont le support est form\'e d'\'el\'ements fortement r\'eguliers dans $\tilde{G}$. S'il n'y a pas de torsion du tout, c'est-\`a-dire si $\tilde{G}=G$,  la stabilit\'e de la distribution ${\bf f}\mapsto S_{\tilde{M}}^{\tilde{G}}(\boldsymbol{\delta},{\bf f})$ a \'et\'e prouv\'ee par Arthur ([A4] theorem 1.1(b)) . Dans la section 2 de [III], on a effectu\'e une construction qui ram\`ene le cas \`a torsion int\'erieure au cas sans torsion.
  Dans les derniers paragraphes de cette section, on avait suppos\'e le corps de base non-archim\'edien. Comme on l'avait dit alors, ce n'\'etait que parce que certains objets n'avaient \'et\'e d\'efinis que dans ce cas. Maintenant que ces objets (\`a savoir l'application $\phi_{\tilde{M}}$ et les int\'egrales orbitales pond\'er\'ees stables) ont \'et\'e d\'efinis aussi dans le cas $F={\mathbb R}$, les r\'esultats de [III] 2 sont valables dans ce cas. En particulier la proposition [III] 2.8, qui afffirme que, pour $(G,\tilde{G},{\bf a})$ quasi-d\'eploy\'e et \`a torsion int\'erieure et pour $\boldsymbol{\delta}$ comme ci-dessus, la distribution ${\bf f}\mapsto S_{\tilde{M}}^{\tilde{G}}(\boldsymbol{\delta},{\bf f})$ est stable. Passons au cas d'un \'el\'ement $\boldsymbol{\delta}\in D_{g\acute{e}om,\tilde{G}-\acute{e}qui}^{st}(\tilde{M}({\mathbb R}))\otimes Mes(M({\mathbb R}))^*$. Soit ${\bf f}\in C_{c}^{\infty}(\tilde{G}({\mathbb R}))\otimes Mes(G({\mathbb R}))$ dont l'image dans $SI(\tilde{G}({\mathbb R}))\otimes Mes(G({\mathbb R}))$ est nulle. On  introduit une fonction $f'$ v\'erifiant 1.4(2) et (3). Le r\'esultat de stabilit\'e que l'on vient de prouver et la relation 1.4(2) entra\^{\i}nent que les int\'egrales orbitales stables assez r\'eguli\`eres de $f'$ sont nulles au voisinage du support de $\boldsymbol{\delta}$. La relation 1.4(3) entra\^{\i}ne alors que $S_{\tilde{M}}^{\tilde{G}}(\boldsymbol{\delta},{\bf f})=0$. Cela prouve le th\'eor\`eme.

 \bigskip
 
 \subsection{Une formule d'induction}
 
 Le triplet $(G,\tilde{G},{\bf a})$ est encore quasi-d\'eploy\'e et \`a torsion int\'erieure. Soient $\tilde{R}\subset \tilde{M}$ deux espaces de Levi. Soit ${\cal O}$ une classe de conjugaison stable semi-simple dans $\tilde{R}({\mathbb R})$, notons ${\cal O}^{\tilde{M}}$ l'unique classe de conjugaison stable dans $\tilde{M}({\mathbb R})$ qui contient ${\cal O}$.
 \ass{Proposition}{On suppose que ${\cal O}^{\tilde{M}}$ est $\tilde{G}$-\'equisinguli\`ere.
 
 (i) Pour tout espace de Levi $\tilde{L}\in {\cal L}(\tilde{R})$ tel que $e_{\tilde{R}}^{\tilde{G}}(\tilde{M},\tilde{L})\not=0$, la classe ${\cal O}$ est $\tilde{L}$-\'equisinguli\`ere.
 
 (ii) Pour $\boldsymbol{\delta}\in D^{st}_{g\acute{e}om}({\cal O})\otimes Mes(R({\mathbb R}))^*$ et ${\bf f}\in I(\tilde{G}({\mathbb R}))\otimes Mes(G({\mathbb R}))$, on a l'\'egalit\'e
 $$S_{\tilde{M}}^{\tilde{G}}(\boldsymbol{\delta}^{\tilde{M}},{\bf f})=\sum_{\tilde{L}\in {\cal L}(\tilde{R})}e_{\tilde{R}}^{\tilde{G}}(\tilde{M},\tilde{L})S_{\tilde{R}}^{\tilde{L}}(\boldsymbol{\delta},{\bf f}_{\tilde{L}}).$$}
 
 Le terme $e_{\tilde{R}}^{\tilde{G}}(\tilde{M},\tilde{L})$ est le m\^eme qu'en [II] 1.14. Le (i) est analogue \`a celui du lemme 1.3. Le (ii) se prouve comme dans le cas non-archim\'edien, cf. [II] 1.14.

 \bigskip
 
 \subsection{Int\'egrales orbitales pond\'er\'ees $\omega$-\'equivariantes et  endoscopie}
 Revenons au cas o\`u $(G,\tilde{G},{\bf a})$ est quelconque. Soit $\tilde{M}$ un espace de Levi de $\tilde{G}$  et soit ${\bf M}'=(M',{\cal M}',\tilde{\zeta})$ une donn\'ee endoscopique elliptique et relevante de $(M,\tilde{M},{\bf a}_{M})$. Soit ${\cal O}'$ une classe de conjugaison stable semi-simple dans $\tilde{M}'({\mathbb R})$. On note $D_{g\acute{e}om}^{st}({\bf M}',{\cal O}')$  le sous-espace des \'el\'ements de $D^{st}_{g\acute{e}om}({\bf M}')$ dont le support dans $\tilde{M}'({\mathbb R})$ est form\'e d'\'el\'ements dont la partie semi-simple appartient  \`a ${\cal O}'$. 
 Il ne  correspond pas toujours \`a ${\cal O}'$ de classe de conjugaison stable dans $\tilde{M}({\mathbb R})$, mais il lui correspond en tout cas une classe de conjugaison par $M({\mathbb C})$ dans $\tilde{M}({\mathbb C})$.   Nous dirons que ${\cal O}'$ est $\tilde{G}$-\'equisinguli\`ere si tout \'el\'ement $\gamma$ de cette classe dans $\tilde{M}({\mathbb C})$ v\'erifie $M_{\gamma}=G_{\gamma}$. On note $D^{st}_{g\acute{e}om,\tilde{G}-\acute{e}qui}({\bf M}')$ la somme des  $D_{g\acute{e}om}^{st}({\bf M}',{\cal O}')$ sur les classes de conjugaison stable semi-simples ${\cal O}'$ qui sont  $\tilde{G}$-\'equisinguli\`eres.
 
 Soient $\boldsymbol{\delta}\in  D^{st}_{g\acute{e}om,\tilde{G}-\acute{e}qui}({\bf M}')\otimes Mes(M'({\mathbb R}))^*$ et ${\bf f}\in C_{c}^{\infty}(\tilde{G}({\mathbb R}))\otimes Mes(G({\mathbb R}))$. Posons
 $$I_{\tilde{M}}^{\tilde{G},{\cal E}}({\bf M}',\boldsymbol{\delta},{\bf f})=\sum_{\tilde{s}\in \tilde{\zeta}Z(\hat{M})^{\Gamma_{{\mathbb R}},\hat{\theta}}/Z(\hat{G})^{\Gamma_{{\mathbb R}},\hat{\theta}}}i_{\tilde{M}'}(\tilde{G},\tilde{G}'(\tilde{s}))S_{{\bf M}'}^{{\bf G}'(\tilde{s})}(\boldsymbol{\delta},{\bf f}^{{\bf G}'(\tilde{s})}).$$
 On renvoie \`a [II] 1.12 pour les notations. Les propri\'et\'es formelles n\'ecessaires \`a la d\'efinition des termes ci-dessus ont \'et\'e vues dans le paragraphe 1.4.   Le point \`a souligner est que l'hypoth\`ese $\boldsymbol{\delta}\in  D^{st}_{g\acute{e}om,\tilde{G}-\acute{e}qui}({\bf M}')\otimes Mes(M'({\mathbb R}))^*$ entra\^{\i}ne $\boldsymbol{\delta}\in  D^{st}_{g\acute{e}om,\tilde{G}'(\tilde{s})-\acute{e}qui}({\bf M}')\otimes Mes(M'({\mathbb R}))^*$ pour tout $\tilde{s}$. En effet, soit  $\delta$ un \'el\'ement semi-simple de $\tilde{M}'({\mathbb R})$. Fixons un \'el\'ement $ \gamma\in \tilde{M}({\mathbb C})$   dans la classe de conjugaison g\'eom\'etrique correspondant \`a celle de $\delta$. Soit enfin $\tilde{s}$ comme ci-dessus. Alors il y a une injection du syst\`eme de racines du groupe $G'(\tilde{s})_{\delta}$ dans celui du groupe $G_{\gamma}$. Et l'ensemble des racines de $M'_{\delta}$ est celui des racines de $G'(\tilde{s})_{\delta}$ qui, par cette injection, deviennent des racines de $M_{\gamma}$. L'\'egalit\'e $M_{\gamma}=G_{\gamma}$ entra\^{\i}ne donc $M'_{\delta}=G'(\tilde{s})_{\delta}$. En notant ${\cal O}'$ la classe de conjugaison stable de $\delta$, on obtient que, si  ${\cal O}'$ est $\tilde{G}$-\'equisinguli\`ere, elle est aussi $\tilde{G}'(\tilde{s})$-\'equisinguli\`ere.  
  
   .

Au lieu d'un triplet $(G,\tilde{G},{\bf a})$, consid\'erons maintenant un triplet $(KG,K\tilde{G},{\bf a})$, o\`u $K\tilde{G}=(\tilde{G}_{p})_{p\in \Pi}$ est un $K$-espace, cf. [I] 1.11 dont on utilise les notations.  Soit $K\tilde{M}=(\tilde{M}_{p})_{p\in \Pi_{M}}\in {\cal L}(K\tilde{M}_{0})$.  Soit ${\bf M}'=(M',{\cal M}',\tilde{\zeta})$ une donn\'ee endoscopique elliptique et relevante de $(KM,K\tilde{M},{\bf a}_{M})$. Soit ${\cal O}'$ une classe de conjugaison stable semi-simple dans $\tilde{M}'({\mathbb R})$. Pour chaque composante connexe $\tilde{M}_{p}$ de $K\tilde{M}$, avec $p\in \Pi_{M}$, il correspond \`a ${\cal O}'$ une classe de conjugaison ${\cal O}_{p,{\mathbb C}}$ par $M_{p}({\mathbb C})$ dans $\tilde{M}_{p}({\mathbb C})$. Si l'on remplace $p$ par $q$,  ${\cal O}_{p,{\mathbb C}}$ est l'image par $\tilde{\phi}_{p,q}$ de ${\cal O}_{q,{\mathbb C}}$. Ainsi la condition $M_{p,\gamma}=G_{p,\gamma}$ pour tout $\gamma\in {\cal O}_{p,{\mathbb C}}$ est ind\'ependante de $p\in \Pi_{M}$. Si elle est v\'erifi\'ee, on dit que ${\cal O}'$ est $\tilde{G}$-\'equisinguli\`ere. On d\'efinit alors comme plus haut l'espace $D^{st}_{g\acute{e}om,\tilde{G}-\acute{e}qui}({\bf M}')$.  Soient $\boldsymbol{\delta}\in  D^{st}_{g\acute{e}om,\tilde{G}-\acute{e}qui}({\bf M}')\otimes Mes(M'({\mathbb R}))^*$ et ${\bf f}\in C_{c}^{\infty}(K\tilde{G}({\mathbb R}))\otimes Mes(G({\mathbb R}))$.  On pose 
$$I_{K\tilde{M}}^{K\tilde{G},{\cal E}}({\bf M}',\boldsymbol{\delta},{\bf f})=\sum_{\tilde{s}\in \tilde{\zeta}Z(\hat{M})^{\Gamma_{{\mathbb R}},\hat{\theta}}/Z(\hat{G})^{\Gamma_{{\mathbb R}},\hat{\theta}}}i_{\tilde{M}'}(\tilde{G},\tilde{G}'(\tilde{s}))S_{{\bf M}'}^{{\bf G}'(\tilde{s})}(\boldsymbol{\delta},{\bf f}^{{\bf G}'(\tilde{s})}).$$

 Les termes $I_{K\tilde{M}}^{K\tilde{G},{\cal E}}({\bf M}',\boldsymbol{\delta},{\bf f})$ que l'on a ainsi d\'efinis v\'erifient les m\^emes propri\'et\'es que dans le cas d'un corps de base non-archim\'edien. Ils sont invariants par le groupe d'automorphismes $Aut(K\tilde{M},{\bf M}')$, au sens pr\'ecis\'e en [II] 1.13. Ecrivons compl\`etement ce que devient la proposition 1.14 de [II].
  
 (a)  Soit $R'$ un groupe de Levi de $M'$ qui est relevant. On peut lui associer un $K$-espace de Levi $K\tilde{R}\in {\cal L}(K\tilde{M}_{0})$ tel que  $K\tilde{R}\subset K\tilde{M}$ et une donn\'ee endocopique  elliptique et relevante ${\bf R}'$ de $(KR,K\tilde{R},{\bf a}_{R})$. Soit ${\cal O}'$ une classe de conjugaison stable semi-simple dans $\tilde{R}'({\mathbb R})$. On note ${\cal O}^{\tilde{M}'}$ l'unique classe de conjugaison stable dans $\tilde{M}'({\mathbb R})$ qui contient ${\cal O}'$. On suppose que ${\cal O}^{\tilde{M}'}$ est $\tilde{G}$-\'equisinguli\`ere. Pour $\boldsymbol{\delta}\in D_{g\acute{e}om}^{st}({\bf R}',{\cal O}')\otimes Mes(R'({\mathbb R}))^*$ et ${\bf f}\in C_{c}^{\infty}(K\tilde{G}({\mathbb R}))\otimes Mes(G({\mathbb R}))$, on a alors l'\'egalit\'e
  
  (1)  $I_{K\tilde{M}}^{K\tilde{G},{\cal E}}({\bf M}',\boldsymbol{\delta}^{M'},{\bf f})=\sum_{K\tilde{L}\in {\cal L}(K\tilde{R})}d_{\tilde{R}}^{\tilde{G}}(\tilde{M},\tilde{L})I_{K\tilde{R}}^{K\tilde{L},{\cal E}}({\bf R}',\boldsymbol{\delta},{\bf f}_{\tilde{L},\omega})$.
  
  On a l'anologue du (i) du lemme 1.3: l'hypoth\`ese implique que ${\cal O}'$ est $\tilde{L}$-\'equisinguli\`ere pour tout $K\tilde{L}$ intervenant dans la somme, ce qui donne un sens aux termes de cette somme.

 (b)  Soit $R'$ un groupe de Levi de $M'$ qui n'est pas relevant. Soient ${\cal O}'$ et ${\cal O}^{\tilde{M}'}$ comme ci-dessus. On suppose encore que ${\cal O}^{\tilde{M}'}$ est $\tilde{G}$-\'equisinguli\`ere. Fixons des donn\'ees suppl\'ementaires $M'_{1}$,...,$\Delta_{1}$ pour ${\bf M}'$. On peut alors d\'efinir l'espace $D_{g\acute{e}om,\lambda_{1}}^{st}(\tilde{R}'_{1}({\mathbb R}),{\cal O}^{st}) $ de distributions sur $\tilde{R}'_{1}({\mathbb R})$. Pour $\boldsymbol{\delta}\in D_{g\acute{e}om,\lambda_{1}}^{st}(\tilde{R}'_{1}({\mathbb R}),{\cal O}')\otimes Mes(R'({\mathbb R}))^*$ et ${\bf f}\in C_{c}^{\infty}(K\tilde{G}({\mathbb R}))\otimes Mes(G({\mathbb R}))$, on a alors l'\'egalit\'e
  
(2) $I_{K\tilde{M}}^{K\tilde{G},{\cal E}}({\bf M}',\boldsymbol{\delta}^{M'},{\bf f})=0$. 

{\bf Remarque.} Dans la preuve de cette propri\'et\'e, l'usage du lemme 2.1 de [A5] fait dans le cas non-archim\'edien doit \^etre remplac\'e par celui du lemme 3.5 de [I].
\bigskip

Soit $p\in \Pi$ et soit ${\bf f}_{p}\in C_{c}^{\infty}(\tilde{G}_{p}({\mathbb R}))\otimes Mes(G({\mathbb R}))$. On peut identifier ${\bf f}_{p}$ \`a un \'el\'ement ${\bf f}\in C_{c}^{\infty}(K\tilde{G}({\mathbb R}))\otimes Mes(G({\mathbb R}))$ dont les composantes ${\bf f}_{q}$ sont nulles pour $q\in \Pi$, $q\not=p$. Si $p\in \Pi_{M}$, on a par d\'efinition 
$$I_{K\tilde{M}}^{K\tilde{G},{\cal E}}({\bf M}',\boldsymbol{\delta},{\bf f})=I_{\tilde{M}_{p}}^{K\tilde{G}_{p},{\cal E}}({\bf M}',\boldsymbol{\delta},{\bf f}).$$
Par contre, si $p\in \Pi-\Pi_{M}$, le membre de gauche est bien d\'efini tandis que celui de droite ne l'est pas. Il r\'esultera du th\'eor\`eme 1.10  (quand celui-ci sera prouv\'e) que le membre de gauche est nul. Mais ce n'est nullement \'evident a priori.
 
  \bigskip
  
  \subsection{Int\'egrales orbitales pond\'er\'ees $\omega$-\'equivariantes endoscopiques}
On consid\`ere un triplet $(KG,K\tilde{G},{\bf a})$ o\`u $K\tilde{G}$ est un $K$-espace. Soit $K\tilde{M} \in {\cal L}(K\tilde{M}_{0})$.  Soit ${\cal O}$ une classe de conjugaison stable semi-simple dans $K\tilde{M}({\mathbb R})$. On suppose que ${\cal O}$ est $\tilde{G}$-\'equisinguli\`ere, c'est-\`a-dire que pour tout $p\in \Pi_{M}$ tel que ${\cal O}\cap \tilde{M}_{p}({\mathbb R})\not=\emptyset$, cette intersection, qui est une classe de conjugaison stable dans $\tilde{M}_{p}({\mathbb R})$, est $\tilde{G}_{p}$-\'equisinguli\`ere.   On fixe un ensemble de repr\'esentants ${\cal E}(\tilde{M},{\bf a}_{M})$ des classes d'\'equivalence de donn\'ees endoscopiques de $(M,\tilde{M},{\bf a}_{M})$ qui sont elliptiques et relevantes. Pour tout ${\bf M}'\in {\cal E}(\tilde{M},{\bf a}_{M})$, il correspond \`a ${\cal O}$ une r\'eunion finie ${\cal O}_{\tilde{M}'}$ de classes de conjugaison stable dans $\tilde{M}'({\mathbb R})$.  D'apr\`es la d\'efinition du paragraphe pr\'ec\'edent, ces classes sont $\tilde{G}$-\'equisinguli\`eres.  Soit $\boldsymbol{\gamma}\in D_{g\acute{e}om}({\cal O})\otimes Mes(M({\mathbb R}))^*$. D'apr\`es la proposition 5.7 de [I], il existe une famille $(\boldsymbol{\delta}_{{\bf M}'})_{{\bf M}'\in {\cal E}(\tilde{M},{\bf a}_{M})}$, avec $\boldsymbol{\delta}_{{\bf M}'}\in D^{st}_{g\acute{e}om}({\bf M}',{\cal O}_{\tilde{M}'})\otimes Mes(M'({\mathbb R}))^*$ pour tout ${\bf M}'$, de sorte que
 $$(1) \qquad \boldsymbol{\gamma}=\sum_{{\bf M}'\in {\cal E}(\tilde{M},{\bf a}_{M})}transfert(\boldsymbol{\delta}_{{\bf M}'}).$$
 Fixons de tels objets. Soit ${\bf f}\in I(K\tilde{G}({\mathbb R}),\omega)\otimes Mes(G({\mathbb R}))$. La somme suivante est bien d\'efinie
 $$(2)\qquad \sum_{{\bf M}'\in {\cal E}(\tilde{M},{\bf a}_{M})}I_{K\tilde{M}}^{K\tilde{G}}({\bf M}',\boldsymbol{\delta}_{{\bf M}'},{\bf f}).$$
 De m\^eme qu'en [II] 1.15, les propri\'et\'es indiqu\'ees dans le paragraphe pr\'ec\'edent impliquent que cette somme ne d\'epend pas de la d\'ecomposition (1) choisie. On peut d\'efinir un terme $I_{K\tilde{M}}^{K\tilde{G},{\cal E}}(\boldsymbol{\gamma},{\bf f})$ comme \'etant cette somme (2). Par lin\'earit\'e, ce terme est donc d\'efini pour $\boldsymbol{\gamma}\in D_{g\acute{e}om,\tilde{G}-\acute{e}qui}(K\tilde{M}({\mathbb R}),\omega)\otimes Mes(M({\mathbb R}))^*$.
 
 Ce terme  a les m\^emes propri\'et\'es relativement \`a l'induction que les int\'egrales orbitales $\omega$-\'equivariantes. C'est-\`a-dire que soit $K\tilde{R} \in {\cal L}(K\tilde{M}_{0})$ tel que  $K\tilde{R}\subset K\tilde{M}$.  Soit ${\cal O}$ une classe de conjugaison stable semi-simple dans $K\tilde{R}({\mathbb R})$. Notons ${\cal O}^{\tilde{M}}$ l'unique classe de conjugaison stable dans $K\tilde{M}({\mathbb R}) $  qui contient ${\cal O}$. On suppose  que ${\cal O}^{\tilde{M}}$ est $\tilde{G}$-\'equisinguli\`ere. Alors, pour tout $\boldsymbol{\gamma}\in D_{g\acute{e}om}({\cal O},\omega)\otimes Mes(R({\mathbb R}))^*$ et tout ${\bf f}\in I(K\tilde{G}({\mathbb R}),\omega)\otimes Mes(G({\mathbb R}))$, on a l'\'egalit\'e
  $$I_{K\tilde{M}}^{K\tilde{G},{\cal E}}(\boldsymbol{\gamma}^{\tilde{M}},{\bf f})=\sum_{K\tilde{L}\in {\cal L}(K\tilde{R})}d_{\tilde{R}}^{\tilde{G}}(\tilde{M},\tilde{L})I_{K\tilde{R}}^{K\tilde{L}}(\boldsymbol{\gamma},{\bf f}_{K\tilde{L},\omega}).$$
  
  {\bf Remarque.} Si l'on ne travaillait pas avec un $K$-espace, mais seulement avec une composante, la somme (2) ne serait plus en g\'en\'eral ind\'ependante de la d\'ecomposition (1) choisie.
  
  \bigskip
  
  \subsection{Une propri\'et\'e locale des int\'egrales orbitales $\omega$-\'equivariantes endoscopiques}
  On conserve la m\^eme situation. Soit ${\cal O}$ une classe de conjugaison stable semi-simple dans $K\tilde{M}({\mathbb R})$.   On se d\'ebarrasse ici des espaces de mesures de Haar en fixant de telles mesures sur tous les groupes intervenant.
  
  \ass{Lemme}{On suppose que ${\cal O}$ est $\tilde{G}$-\'equisinguli\`ere. Soit $ f\in I(K\tilde{G}({\mathbb R}),\omega)$.
  
  (i) Il existe une fonction $f'\in I(K\tilde{M}({\mathbb R}),\omega)$ telle que, pour tout  $\boldsymbol{\tau}\in D_{orb}(K\tilde{M}({\mathbb R}),\omega)$ dont le support est form\'e  d'\'el\'ements $\tilde{G}$-r\'eguliers et assez voisins de ${\cal O}$, on ait l'\'egalit\'e
  $$I_{K\tilde{M}}^{K\tilde{G},{\cal E}}(\boldsymbol{\tau},f)=I^{K\tilde{M}}(\boldsymbol{\tau},f').$$
  
  (ii) Soit  $\boldsymbol{\gamma}\in D_{g\acute{e}om}({\cal O},\omega)$. Pour $f'$ comme en (i), on a l'\'egalit\'e
  $$I_{K\tilde{M}}^{K\tilde{G},{\cal E}}(\boldsymbol{\gamma},f)=I^{K\tilde{M}}(\boldsymbol{\gamma},f').$$}
  
  Preuve. Introduisons $\underline{la}$ paire de Borel \'epingl\'ee ${\cal E}^*$ de $M$, dont on note $T^*$ le tore. Posons
   $$V^{K\tilde{M}}=\left(((T^*/(1-\theta)(T^*))/W^{M,\theta})\times_{\mathfrak{Z}(M)}\mathfrak{Z}(\tilde{M})\right)^{\Gamma_{{\mathbb R}}},$$
   cf. [I] 1.8. Cet ensemble est une vari\'et\'e analytique r\'eelle. Les classes de conjugaison semi-simples g\'eom\'etriques dans $K\tilde{M}({\mathbb R})$ sont classifi\'ees par un sous-ensemble de $V^{K\tilde{M}}$. La classe de conjugaison stable ${\cal O}$ \'etant incluse dans une classe de conjugaison semi-simple g\'eom\'etrique, il lui correspond un point $v_{{\cal O}}\in V^{K\tilde{M}}$. 
  Soient $K\tilde{R} \in {\cal L}(K\tilde{M}_{0})$ tel que $K\tilde{R}\subset K\tilde{M}$ et ${\bf R}'=(R',{\cal R}',\tilde{\zeta})\in {\cal E}(K\tilde{R},{\bf a}_{R})$. Les classes de conjugaison semi-simples g\'eom\'etriques dans $\tilde{R}'({\mathbb R})$ sont de m\^eme classifi\'ees par une vari\'et\'e $V^{\tilde{R}'}$, qui s'envoie naturellement dans $V^{K\tilde{M}}$. On note ${\cal O}_{\tilde{R}'}$ la r\'eunion des classes de conjugaison stable semi-simples dans $\tilde{R}'({\mathbb R})$ dont l'image dans $V^{K\tilde{M}}$ est $v_{{\cal O}}$. 
  
   Soient $K\tilde{L}\in {\cal L}(K\tilde{R})$ tel que $d_{\tilde{R}}^{\tilde{G}}(\tilde{M},\tilde{L})\not=0$ et $\tilde{s}\in Z(\hat{R})^{\Gamma_{{\mathbb R}},\hat{\theta}}/Z(\hat{L})^{\Gamma_{{\mathbb R}},\hat{\theta}}$.  On a d\'ej\`a vu que ${\cal O}_{\tilde{R}'}$ \'etait $\tilde{L}'(\tilde{s})$-\'equisinguli\`ere (on veut dire par l\`a que cet ensemble est r\'eunion de classes de conjugaison stable qui sont $\tilde{L}'(\tilde{s})$-\'equisinguli\`eres). D'apr\`es une variante des relations (2) et (3) de 1.4 (voir aussi les relations (8) et (9) de ce paragraphe), il existe une fonction $g_{K\tilde{R},{\bf R}'}(K\tilde{L},\tilde{s})\in SI({\bf R}')$ telle que
  
  (1) pour tout $\boldsymbol{\tau}\in D^{st}_{orb}({\bf R}')$ dont le support est form\'e d'\'el\'ements assez r\'eguliers et assez voisins de ${\cal O}_{\tilde{R}'}$, on a l'\'egalit\'e
  $$S_{{\bf R}'}^{{\bf L}'(\tilde{s})}(\boldsymbol{\tau},(f_{K\tilde{L},\omega})^{{\bf L}'(\tilde{s})})=S^{{\bf R}'}(\boldsymbol{\tau},g_{K\tilde{R},{\bf R}'}(K\tilde{L},\tilde{s}));$$
  
  (2) pour tout $\boldsymbol{\delta}\in D^{st}_{g\acute{e}om}({\bf R}',{\cal O}_{\tilde{R}'})$, on a l'\'egalit\'e
  $$S_{{\bf R}'}^{{\bf L}'(\tilde{s})}(\boldsymbol{\delta},(f_{K\tilde{L},\omega})^{{\bf L}'(\tilde{s})})=S^{{\bf R}'}(\boldsymbol{\delta},g_{K\tilde{R},{\bf R}'}(K\tilde{L},\tilde{s})) .$$
  
 Introduisons le groupe $Aut^{K\tilde{M}}(K\tilde{R},{\bf R}')$ des automorphismes de $(K\tilde{R},{\bf R}')$, l'espace ambiant \'etant $K\tilde{M}$. Ce groupe agit sur $SI({\bf R}')$ et cette action se factorise par un quotient fini. Notons $\underline{Aut}^{K\tilde{M}}(K\tilde{R},{\bf R}')$ un tel quotient fini. Soit $x\in Aut^{K\tilde{M}}(K\tilde{R},{\bf R}')$. L'\'el\'ement $x$ permute les couples $(K\tilde{L},\tilde{s})$ intervenant   ci-dessus. On v\'erifie ais\'ement que la fonction $x(g_{K\tilde{R},{\bf R}'}(K\tilde{L},\tilde{s}))$ a les m\^emes propri\'et\'es que $g_{K\tilde{R},{\bf R}'}( x(K\tilde{L},\tilde{s}))$. On peut donc aussi bien remplacer chaque fonction   $g_{K\tilde{R},{\bf R}'}(K\tilde{L},\tilde{s})$ par
 $$\vert \underline{Aut}^{K\tilde{M}}(K\tilde{R},{\bf R}')\vert ^{-1}\sum_{x\in \underline{Aut}^{K\tilde{M}}(K\tilde{R},{\bf R}')}x(g_{K\tilde{R},{\bf R}'}(x^{-1}(K\tilde{L},\tilde{s}))).$$
  Cela fait, posons 
  $$(3) \qquad g_{K\tilde{R},{\bf R}'}=\sum_{K\tilde{L}\in {\cal L}(K\tilde{R})}d_{\tilde{R}}^{\tilde{G}}(\tilde{M},\tilde{L})\sum_{\tilde{s}\in Z(\hat{R})^{\Gamma_{{\mathbb R}},\hat{\theta}}/Z(\hat{L})^{\Gamma_{{\mathbb R}},\hat{\theta}}}i_{\tilde{R}'}(\tilde{L},\tilde{L}'(\tilde{s}))g_{K\tilde{R},{\bf R}'}(K\tilde{L},\tilde{s}).$$
  Alors cette fonction est invariante par $Aut^{K\tilde{M}}(K\tilde{R},{\bf R}')$. Fixons un ensemble de repr\'esentants ${\cal E}_{+}(K\tilde{M},{\bf a}_{M})$ des classes d'\'equivalence de couples $(K\tilde{R},{\bf R}')$ comme ci-dessus. Nos constructions   d\'efinissent une famille $(g_{K\tilde{R},{\bf R}'})_{(K\tilde{R},{\bf R}')\in {\cal E}_{+}(K\tilde{M},{\bf a}_{M})}$. Consid\'erons les conditions de la proposition 4.11 de [I]. La premi\`ere condition est l'invariance de $g_{K\tilde{R},{\bf R}'}$ par $Aut^{K\tilde{M}}(K\tilde{R},{\bf R}')$, qui est v\'erifi\'ee par construction. Pour la deuxi\`eme condition, soient ${\bf M}'\in {\cal E}(K\tilde{M},{\bf a}_{M})$, $\tilde{R}'$ un espace de Levi de $\tilde{M}'$ qui est relevant et soit $(K\tilde{R},{\bf R}')$ l'\'el\'ement de ${\cal E}_{+}(K\tilde{M},{\bf a}_{M})$ qui lui est associ\'e. La condition de [I] 4.11  impose l'\'egalit\'e 
  
  (4) $(g_{K\tilde{M},{\bf M}'})_{\tilde{R}'}=g_{K\tilde{R},{\bf R}'}$. 
  
  Celle-ci n'est pas v\'erifi\'ee en g\'en\'eral, mais elle l'est au voisinage de ${\cal O}_{\tilde{R}'}$. En effet, soit $\boldsymbol{\tau}\in D^{st}_{orb}({\bf R}')$ dont le support est form\'e d'\'el\'ements assez r\'eguliers et assez voisins de ${\cal O}_{\tilde{R}'}$. Les propri\'et\'es (1) et (3) nous disent que
  $$(5) \qquad S^{{\bf R}'}(\boldsymbol{\tau},g_{K\tilde{R},{\bf R}'})=\sum_{K\tilde{L}\in {\cal L}(K\tilde{R})}d_{\tilde{R}}^{\tilde{G}}(\tilde{M},\tilde{L})I_{K\tilde{R}}^{K\tilde{L},{\cal E}}({\bf R}',\boldsymbol{\tau},f_{K\tilde{L},\omega}).$$
  On a aussi 
  $$S^{{\bf R}'}(\boldsymbol{\tau},(g_{K\tilde{M},{\bf M}'})_{\tilde{R}'})=S^{{\bf M}'}(\boldsymbol{\tau}^{\tilde{M}'},g_{K\tilde{M},{\bf M}'}).$$
 La relation (5) appliqu\'ee au  couple $(K\tilde{M},{\bf M}')$ se simplifie en
 $$(6) \qquad S^{{\bf M}'}(\boldsymbol{\tau}',g_{K\tilde{M},{\bf M}'})=I_{K\tilde{M}}^{K\tilde{G},{\cal E}}({\bf M}',\boldsymbol{\tau}',f),$$
 pour tout $\boldsymbol{\tau}'\in D^{st}_{orb}({\bf M}')$ dont le support est form\'e d'\'el\'ements assez r\'eguliers et assez voisins de ${\cal O}_{\tilde{M}'}$.
En l'utilisant ,  on obtient
   $$S^{{\bf R}'}(\boldsymbol{\tau},(g_{K\tilde{M},{\bf M}'})_{\tilde{R}'})=I_{K\tilde{M}}^{K\tilde{G},{\cal E}}({\bf M}',\boldsymbol{\tau}^{\tilde{M}'},f).$$ 
   La relation 1.6(1) dit que ceci est \'egal au membre de droite de (5). Donc l'\'egalit\'e (4) est bien v\'erifi\'ee au voisinage de ${\cal O}_{\tilde{R}'}$.
   
   Nous ne d\'etaillerons pas la troisi\`eme condition de la proposition [I] 4.11: sa validit\'e locale r\'esulte comme ci-dessus de la relation 1.6(2). 
   
   Toute fonction $\xi$ sur $V^{K\tilde{M}}$ se rel\`eve en une fonction sur $K\tilde{M}({\mathbb R})$: c'est la fonction  $\gamma\mapsto \xi(v_{\gamma})$, o\`u, pour $\gamma\in K\tilde{M}({\mathbb R})$, $v_{\gamma}$ est le point de $V^{K\tilde{M}}$ qui classifie la partie semi-simple de $\gamma$. 
 Fixons un voisinage $\Omega_{1}$ de $v_{{\cal O}}$ et  un voisinage $\Omega_{2} $ de la cl\^oture de $\Omega_{1}$. Fixons une fonction  $\xi$ sur $V^{K\tilde{M}}$ qui est $C^{\infty}$, qui vaut $1$ sur $\Omega_{1}$ et vaut $0$ hors de $\Omega_{2}$. Pour tout $\tilde{R}'$ intervenant ci-dessus, $V^{\tilde{R}'}$ s'envoie dans $V^{K\tilde{M}}$. Par composition avec cette application, $\xi$ devient 
   une fonction sur  $V^{\tilde{R}'}$, qui se rel\`eve en une fonction  $\xi_{\tilde{R}'}$   sur $\tilde{R}'({\mathbb R})$. Rempla\c{c}ons $g_{K\tilde{R},{\bf R}'}$ par son produit avec $\xi_{\tilde{R}'}$. Cela ne retire rien aux propri\'et\'es de cette fonction. Mais il est plus ou moins clair que, si l'on choisit $\Omega_{2}$ assez petit, l'\'egalit\'e (4), qui n'\'etait v\'erifi\'ee que localement, est maintenant v\'erifi\'ee partout: on a annul\'e les fonctions l\`a o\`u l'\'egalit\'e n'\'etait pas v\'erifi\'ee. 
 
 On a maintenant obtenu une famille $(g_{K\tilde{R},{\bf R}'})_{(K\tilde{R},{\bf R}')\in {\cal E}_{+}(K\tilde{M},{\bf a}_{M})}$ qui v\'erifie les conditions de la proposition 4.11 de [I]. Cette proposition nous dit qu'il existe un unique $f'\in I(K\tilde{M}({\mathbb R}),\omega)$ tel que, pour tout ${\bf M}'\in {\cal E}(K\tilde{M},{\bf a}_{M})$, le terme $g_{K\tilde{M},{\bf M}'}$ soit \'egal au transfert $(f')^{{\bf M}'}$. Soit $\boldsymbol{\tau}\in D_{orb}(K\tilde{M}({\mathbb R}),\omega)$ dont le support est form\'e d'\'el\'ements $\tilde{G}$-r\'eguliers et assez voisins de ${\cal O}$. On peut \'ecrire
 $$\boldsymbol{\tau}=\sum_{{\bf M}'\in {\cal E}(K\tilde{M},{\bf a}_{M})}transfert(\boldsymbol{\tau}_{{\bf M}'}),$$
 o\`u $\boldsymbol{\tau}_{{\bf M}'}\in D^{st}_{orb}({\bf M}')$ a pour support des \'el\'ements assez r\'eguliers et voisins de ${\cal O}_{\tilde{M}'}$. Alors
 $$I^{K\tilde{M}}(\boldsymbol{\tau},f')=\sum_{{\bf M}'\in {\cal E}(K\tilde{M},{\bf a}_{M})}I^{K\tilde{M}}(transfert(\boldsymbol{\tau}_{{\bf M}'}),f')=\sum_{{\bf M}'\in {\cal E}(K\tilde{M},{\bf a}_{M})}S^{{\bf M}'}(\boldsymbol{\tau}_{{\bf M}'},(f')^{{\bf M}'})$$
 $$=\sum_{{\bf M}'\in {\cal E}(K\tilde{M},{\bf a}_{M})}S^{{\bf M}'}(\boldsymbol{\tau}_{{\bf M}'},g_{K\tilde{M},{\bf M}'})=\sum_{{\bf M}'\in {\cal E}(K\tilde{M},{\bf a}_{M})}I_{K\tilde{M}}^{K\tilde{G},{\cal E}}({\bf M}',\boldsymbol{\tau}_{{\bf M}'},f),$$
 cette derni\`ere \'egalit\'e provenant de (6). Par d\'efinition, la derni\`ere expression vaut $I_{K\tilde{M}}^{K\tilde{G},{\cal E}}(\boldsymbol{\tau},f)$. Cela prouve que la fonction $f'$ satisfait le (i) de l'\'enonc\'e. 
   
La preuve du (ii) est similaire, en utilisant la relation suivante pour ${\bf M}'\in {\cal E}(K\tilde{M},{\bf a}_{M})$:   
$$S^{{\bf M}'}(\boldsymbol{\delta},g_{K\tilde{M},{\bf M}'})=I_{K\tilde{M}}^{K\tilde{G},{\cal E}}({\bf M}',\boldsymbol{\delta},f),$$
pour tout $\boldsymbol{\delta}\in D^{st}_{g\acute{e}om}({\bf M}',{\cal O}_{\tilde{M}'})$.  Elle r\'esulte de la d\'efinition de $g_{K\tilde{M},{\bf M}'}$ et de (2) ci-dessus. Cela ach\`eve la preuve. $\square$

  \bigskip
  
  \subsection{Le th\'eor\`eme principal}
Evidemment, la d\'efinition des int\'egrales orbitales pond\'er\'ees $\omega$-\'equivariantes donn\'ee en 1.3 s'adapte aux $K$-espaces. Soit $K\tilde{M}$ un $K$-espace de Levi de $K\tilde{G}$.  
  
  \ass{Th\'eor\`eme (\`a prouver)}{Soient $\boldsymbol{\gamma}\in D_{g\acute{e}om,\tilde{G}-\acute{e}qui}(K\tilde{M}({\mathbb R}),\omega)\otimes Mes(M({\mathbb R}))^*$ et ${\bf f}\in I(K\tilde{G}({\mathbb R}))\otimes Mes(G({\mathbb R}))$. Alors on a l'\'egalit\'e
  $$I_{K\tilde{M}}^{K\tilde{G},{\cal E}}(\boldsymbol{\gamma},{\bf f})=I_{K\tilde{M}}^{K\tilde{G}}(\boldsymbol{\gamma},{\bf f}).$$}
  
  \bigskip

  \subsection{R\'eduction au cas des int\'egrales orbitales r\'eguli\`eres}

  \ass{Lemme}{Soit ${\bf f}\in I(K\tilde{G}({\mathbb R}))\otimes Mes(G({\mathbb R}))$. Supposons l'\'egalit\'e du th\'eor\`eme v\'erifi\'ee  pour tout  $\boldsymbol{\gamma}\in D_{orb}(K\tilde{M}({\mathbb R}),\omega)\otimes Mes(M({\mathbb R}))^*$ dont le support est form\'e d'\'el\'ements $\tilde{G}$-r\'eguliers. Alors elle l'est pour tout $\boldsymbol{\gamma}\in D_{g\acute{e}om,\tilde{G}-\acute{e}qui}(K\tilde{M}({\mathbb R}),\omega)\otimes Mes(M({\mathbb R}))^*$.} 
  
  Preuve. On  oublie comme toujours les questions de mesures. On fixe une classe de conjugaison stable semi-simple ${\cal O}$ dans $K\tilde{M}({\mathbb R})$ qui est $\tilde{G}$-\'equisinguli\`ere. On doit montrer que, sous l'hypoth\`ese de l'\'enonc\'e, le th\'eor\`eme est v\'erifi\'e pour $\boldsymbol{\gamma}\in D_{g\acute{e}om}({\cal O},\omega) $. D'apr\`es la d\'efinition de 1.3, on peut fixer une fonction $f''\in I(K\tilde{M}({\mathbb R}),\omega)$ telle que
  
   $$(1) \qquad I_{K\tilde{M}}^{K\tilde{G}}(\boldsymbol{\tau},f)=I^{K\tilde{M}}(\boldsymbol{\tau},f'')$$
  pour tout $ \boldsymbol{\tau}\in D_{orb}(K\tilde{M}({\mathbb R}),\omega)$ dont le support est form\'e d'\'el\'ements assez r\'eguliers et voisins de ${\cal O}$;
  
  $$(2) \qquad I_{K\tilde{M}}^{K\tilde{G}}(\boldsymbol{\gamma},f)=I^{K\tilde{M}}(\boldsymbol{\gamma},f'')$$
  pour tout $ \boldsymbol{\gamma}\in D_{g\acute{e}om}({\cal O},\omega)$.
  
  Le lemme 1.8 nous fournit une fonction $f'$ qui a les m\^emes propri\'et\'es relativement aux int\'egrales orbitales pond\'er\'ees endoscopiques:
  
  $$(3) \qquad I_{K\tilde{M}}^{K\tilde{G},{\cal E}}(\boldsymbol{\tau},f)=I^{K\tilde{M}}(\boldsymbol{\tau},f')$$
  pour tout $ \boldsymbol{\tau}\in D_{orb}(K\tilde{M}({\mathbb R}),\omega)$ dont le support est form\'e d'\'el\'ements assez r\'eguliers et voisins de ${\cal O}$;
  
  $$(4) \qquad I_{K\tilde{M}}^{K\tilde{G},{\cal E}}(\boldsymbol{\gamma},f)=I^{K\tilde{M}}(\boldsymbol{\gamma},f')$$
  pour tout $ \boldsymbol{\gamma}\in D_{g\acute{e}om}({\cal O},\omega)$.

L'hypoth\`ese de l'\'enonc\'e et les relations (1) et (3) impliquent  que $f'$ et $f''$ sont \'egales au voisinage de ${\cal O}$. Mais alors les relations  (2) et (4) impliquent  l'\'egalit\'e cherch\'ee. $\square$

\bigskip

\subsection{Elimination des $K$-espaces}
Revenons \`a un triplet $(G,\tilde{G},{\bf a})$ et \`a un espace de Levi $\tilde{M} $ de $\tilde{G}$. On peut inclure le triplet comme composante connexe d'un triplet $(KG,K\tilde{G},{\bf a})$. L'espace $\tilde{M}$ est alors une composante connexe d'un $K$-espace $K\tilde{M}$ et on ne perd rien \`a supposer que $K\tilde{M}\in {\cal L}(K\tilde{M}_{0})$. {\bf Supposons le th\'eor\`eme 1.10 prouv\'e pour ces donn\'ees}. Soient ${\bf M}'\in {\cal E}(K\tilde{M},{\bf a}_{M})$, $\boldsymbol{\delta}\in D_{g\acute{e}om,\tilde{G}-\acute{e}qui}^{st}({\bf M}')\otimes Mes(M'({\mathbb R}))^*$ et ${\bf f}\in I(\tilde{G}({\mathbb R}),\omega)\otimes Mes(G({\mathbb R}))$. On a d\'efini $I_{\tilde{M}}^{\tilde{G},{\cal E}}({\bf M}',\boldsymbol{\delta},{\bf f})$ en 1.6. On note ici $transfert(\boldsymbol{\delta})$ le transfert de $\boldsymbol{\delta}$ \`a $\tilde{M}({\mathbb R})$. On a alors l'\'egalit\'e
$$(1) \qquad I_{\tilde{M}}^{\tilde{G},{\cal E}}({\bf M}',\boldsymbol{\delta},{\bf f})=I_{\tilde{M}}^{\tilde{G}}(transfert(\boldsymbol{\delta}),{\bf f}).$$
En effet, introduisons la fonction ${\bf f}^{K\tilde{G}}\in I(K\tilde{G}({\mathbb R}),\omega)\otimes Mes(G({\mathbb R}))$ dont la composante sur $\tilde{G}({\mathbb R})$ est ${\bf f}$ et dont les autres composantes sont nulles. Introduisons le transfert $transfert^{K\tilde{G}}(\boldsymbol{\delta})$ de $\boldsymbol{\delta}$ \`a $K\tilde{M}({\mathbb R})$. Sa composante sur $\tilde{M}({\mathbb R})$ est $transfert(\boldsymbol{\delta})$. On ne conna\^{\i}t pas les autres composantes, mais peu importe. Par d\'efinition, on a
$$I_{\tilde{M}}^{\tilde{G},{\cal E}}({\bf M}',\boldsymbol{\delta},{\bf f})=I_{K\tilde{M}}^{K\tilde{G},{\cal E}}({\bf M}',\boldsymbol{\delta},{\bf f}^{K\tilde{G}})=I_{K\tilde{M}}^{K\tilde{G},{\cal E}}(transfert^{K\tilde{M}}(\boldsymbol{\delta}),{\bf f}^{K\tilde{G}}).$$
Gr\^ace au th\'eor\`eme, on obtient
$$I_{\tilde{M}}^{\tilde{G},{\cal E}}({\bf M}',\boldsymbol{\delta},{\bf f})=I_{K\tilde{M}}^{K\tilde{G}}(transfert^{K\tilde{M}}(\boldsymbol{\delta}),{\bf f}^{K\tilde{G}}).$$
Puisque ${\bf f}^{K\tilde{G}}$ est concentr\'e sur la composante $\tilde{G}({\mathbb R})$, ce dernier terme est \'egal par d\'efinition \`a $I_{\tilde{M}}^{\tilde{G}}(transfert(\boldsymbol{\delta}),{\bf f})$. Cela prouve (1). 

\bigskip

\subsection{Le cas quasi-d\'eploy\'e et \`a torsion int\'erieure}
Soit $(G,\tilde{G},{\bf a})$ un triplet quasi-d\'eploy\'e et \`a torsion int\'erieure. Soit $\tilde{M}$ un espace de Levi de $\tilde{G}$ et soit ${\bf M}'$ une donn\'ee endoscopique elliptique et relevante de $(M,\tilde{M},{\bf a}_{M})$.

\ass{Proposition}{Pour tout $\boldsymbol{\delta}\in D_{g\acute{e}om,\tilde{G}-\acute{e}qui}^{st}({\bf M}')\otimes Mes(M'({\mathbb R}))^*$ et  tout ${\bf f}\in I(\tilde{G}({\mathbb R}))\otimes Mes(G({\mathbb R}))$, on a l'\'egalit\'e
  $$I_{\tilde{M}}^{\tilde{G},{\cal E}}({\bf M}',\boldsymbol{\delta},{\bf f})=I_{\tilde{M}}^{\tilde{G}}(transfert(\boldsymbol{\delta}),{\bf f}).$$}
  
 Preuve. Un argument analogue \`a celui du paragraphe 1.11 nous ram\`ene au cas o\`u $\boldsymbol{\delta}$ appartient \`a $D_{orb}^{st}({\bf M}')\otimes Mes(M'({\mathbb R}))^*$ et le support de $\boldsymbol{\delta}$ est form\'e d'\'el\'ements $\tilde{G}$-r\'eguliers.  La preuve est alors  la m\^eme que celle de la proposition [III] 2.9(ii). C'est-\`a-dire que les constructions de la section 2 de [III] nous ram\`enent au cas d'un groupe sans torsion. Dans ce cas, l'assertion r\'esulte de [A5] th\'eor\`eme 1.1(a). $\square$

  \bigskip
  
  \section{Un nouvel espace de distributions}
  
   \bigskip
 
 \subsection{D\'efinition de l'espace $D_{tr-orb}(\tilde{G}({\mathbb R}))$}
  
 Soit $(G,\tilde{G},{\bf a})$ un triplet quasi-d\'eploy\'e et \`a torsion int\'erieure. On d\'efinit  un sous-espace $D_{tr-orb}(\tilde{G}({\mathbb R}))\subset D_{g\acute{e}om}(\tilde{G}({\mathbb R}))$ par r\'ecurrence sur $dim(G_{SC})$. Une fois cet espace d\'efini, on pose $D_{tr-orb}^{st}(\tilde{G}({\mathbb R}))=D_{tr-orb}(\tilde{G}({\mathbb R}))\cap D_{g\acute{e}om}^{st}(\tilde{G}({\mathbb R}))$. On a besoin d'\'etendre 
les d\'efinitions dans la situation habituelle suivante. On consid\`ere des extensions compatibles
$$1\to C_{\natural}\to G_{\natural}\to G\to 1 \,\,\text{et}\,\, \tilde{G}_{\natural}\to \tilde{G}$$
o\`u $C_{\natural}$ est un tore central induit et o\`u $\tilde{G}_{\natural}$ est encore \`a torsion int\'erieure. On se donne de plus un caract\`ere $\lambda_{\natural}$ de $C_{\natural}({\mathbb R})$. En supposant d\'efini l'espace $D_{tr-orb}(\tilde{G}_{\natural}({\mathbb R}))$, on note $D_{tr-orb,\lambda_{\natural}}(\tilde{G}_{\natural}({\mathbb R}))$ son image naturelle dans $D_{g\acute{e}om,\lambda_{\natural}}(\tilde{G}_{\natural}({\mathbb R}))$. De m\^eme, on note $D_{tr-orb,\lambda_{\natural}}^{st}(\tilde{G}_{\natural}({\mathbb R}))$  l'image naturelle de $D_{tr-orb}^{st}(\tilde{G}_{\natural}({\mathbb R}))$ dans $D_{g\acute{e}om,\lambda_{\natural}}^{st}(\tilde{G}_{\natural}({\mathbb R}))$. Supposons donn\'es d'autres objets
$$1\to C_{\flat}\to G_{\flat}\to G\to 1 \,\,, \tilde{G}_{\flat}\to \tilde{G}$$
et $\lambda_{\flat}$ v\'erifiant les m\^emes conditions. Supposons donn\'es comme en 1.4 un caract\`ere $\lambda_{\natural,\flat} $ de $G_{\natural,\flat}({\mathbb R})$ et une fonction $\tilde{\lambda}_{\natural,\flat}$ sur $\tilde{G}_{\natural,\flat}({\mathbb R})$ se transformant selon le caract\`ere $\lambda_{\natural,\flat}$. A l'aide de cette fonction, on construit les isomorphismes
habituels
$$D_{g\acute{e}om,\lambda_{\natural}}(\tilde{G}_{\natural}({\mathbb R}))\simeq D_{g\acute{e}om,\lambda_{\flat}}(\tilde{G}_{\flat}({\mathbb R}))$$
$$D_{g\acute{e}om,\lambda_{\natural}}^{st}(\tilde{G}_{\natural}({\mathbb R}))\simeq D_{g\acute{e}om,\lambda_{\flat}}^{st}(\tilde{G}_{\flat}({\mathbb R})).$$
On montrera plus loin que

(1) $D_{tr-orb,\lambda_{\natural}}(\tilde{G}_{\natural}({\mathbb R}))$ et $D_{tr-orb,\lambda_{\flat}}(\tilde{G}_{\flat}({\mathbb R}))$ se correspondent par le premier isomorphisme; $D_{tr-orb,\lambda_{\natural}}^{st}(\tilde{G}_{\natural}({\mathbb R}))$ et $D_{tr-orb,\lambda_{\flat}}^{st}(\tilde{G}_{\flat}({\mathbb R}))$ se correspondent par le second.

Si maintenant ${\bf G}'=(G',{\cal G}',s)$ est une donn\'ee endoscopique  relevante de $(G,\tilde{G},{\bf a})$, avec $G'\not=G$, on introduit des donn\'ees auxiliaires $G'_{1}$,...,$\Delta_{1}$. L'espace $D_{tr-orb,\lambda_{1}}^{st}(\tilde{G}_{1}({\mathbb R}))$ est bien d\'efini et on note $D_{tr-orb}^{st}({\bf G}')$ le sous-espace de $D_{g\acute{e}om}^{st}({\bf G}')$ auquel il s'identifie. D'apr\`es (1),  cette d\'efinition ne d\'epend pas du choix des donn\'ees auxiliaires. On peut montrer que, si ${\bf G}'_{0}$ est une donn\'ee endoscopique \'equivalente \`a ${\bf G}'$, l'isomorphisme d\'eduit d'une \'equivalence entre  $D_{g\acute{e}om}^{st}({\bf G}')$ et $D_{g\acute{e}om}^{st}({\bf G}'_{0})$ envoie $D_{tr-orb}^{st}({\bf G}')$ sur $D_{tr-orb}^{st}({\bf G}'_{0})$ (on ne donnera pas la preuve formelle similaire \`a celle de (1)).  Rappelons qu'en g\'en\'eral, on note ${\cal E}(\tilde{G},{\bf a})$ un ensemble de repr\'esentants des classes d'\'equivalence de donn\'ees endoscopiques elliptiques et relevantes de $(G,\tilde{G},{\bf a})$. On simplifie cette notation en ${\cal E}(\tilde{G})$ dans notre situation quasi-d\'eploy\'ee \`a torsion int\'erieure.  Cela \'etant, on d\'efinit $D_{tr-orb}(\tilde{G}({\mathbb R}))$ comme le sous-espace de $D_{g\acute{e}om}(\tilde{G}({\mathbb R}))$ engendr\'e par $D_{orb}(\tilde{G}({\mathbb R}))$ et par les  images par  transfert des espaces $D_{tr-orb}^{st}({\bf G}')$ pour ${\bf G}'\in {\cal E}(\tilde{G})$ tel que $G'\not=G$.

 Soit $(G,\tilde{G},{\bf a})$ un triplet qui n'est pas quasi-d\'eploy\'e et \`a torsion int\'erieure. On d\'efinit $D_{tr-orb}(\tilde{G}({\mathbb R}),\omega)$ comme le sous-espace de $D_{g\acute{e}om}(\tilde{G}({\mathbb R}),\omega)$ engendr\'e par $D_{orb}(\tilde{G}({\mathbb R}),\omega)$ et par les  images par  transfert des espaces $D_{tr-orb}^{st}({\bf G}')$ pour ${\bf G}'\in {\cal E}(\tilde{G},{\bf a})$.
 
 Soit $(KG,K\tilde{G},{\bf a})$ un $K$-triplet. On d\'efinit $D_{tr-orb}(K\tilde{G}({\mathbb R}),\omega)$ comme le sous-espace de $D_{g\acute{e}om}(K\tilde{G}({\mathbb R}),\omega)$ engendr\'e par $D_{orb}(K\tilde{G}({\mathbb R}),\omega)$ et par les  images par  transfert des espaces $D_{tr-orb}^{st}({\bf G}')$ pour ${\bf G}'\in {\cal E}(\tilde{G},{\bf a})$.
 
 {\bf Remarque.} Pour $p\in \Pi$, la projection naturelle $D_{g\acute{e}om}(K\tilde{G}({\mathbb R}),\omega)\to D_{g\acute{e}om}(\tilde{G}_{p}({\mathbb R}),\omega)$ envoie $D_{tr-orb}(K\tilde{G}({\mathbb R}),\omega)$ sur $D_{tr-orb}(\tilde{G}_{p}({\mathbb R}),\omega)$. J'ignore par contre si l'homomorphisme
 $$D_{tr-orb}(K\tilde{G}({\mathbb R}),\omega)\to \oplus_{p\in \Pi}D_{tr-orb}(\tilde{G}_{p}({\mathbb R}),\omega)$$
 est surjectif. 
 
 \bigskip
 
 Avant de prouver (1), on doit rappeler qu'en [I] 1.12, on a d\'efini   un groupe $G_{ab}({\mathbb R})$, un groupe r\'eductif connexe $G_{0}$ et des homomorphismes
 $$G({\mathbb R})\to G_{ab}({\mathbb R})\stackrel{N^G}{\to}G_{0,ab}({\mathbb R}).$$
 On a des objets similaires pour les espaces tordus, ainsi que des applications
 $$\tilde{G}({\mathbb R})\to \tilde{G}_{ab}({\mathbb R})\stackrel{N^{\tilde{G}}}{\to}\tilde{G}_{0,ab}({\mathbb R}).$$
 On a
 
 (2) $G_{0}=G$, $\tilde{G}_{0}=\tilde{G}$ et les applications $G({\mathbb R})\to G_{ab}({\mathbb R})$ et $\tilde{G}({\mathbb R})\to \tilde{G}_{ab}({\mathbb R})$ sont surjectives.
 
 Preuve. Les assertions pour les espaces tordus r\'esultent de celles pour les groupes. Le groupe $G_{0}$ est le groupe quasi-d\'eploy\'e tel que $\hat{G}_{0}=\hat{G}^{\hat{\theta},0}$. Ici, $\hat{\theta}=1$ et $G$ est quasi-d\'eploy\'e, d'o\`u $G_{0}=G$. Par d\'efinition, $G_{ab}({\mathbb R})=H^{1,0}(\Gamma_{{\mathbb R}};G_{SC}\to G)$. Fixons une paire de Borel $(B,T)$ de $G$ d\'efinie sur ${\mathbb R}$. L'application naturelle 
 $$H^{1,0}(\Gamma_{{\mathbb R}};T_{sc}\to T)\to H^{1,0}(\Gamma_{{\mathbb R}};G_{SC}\to G)$$ est bijective. On a une suite exacte
 $$T({\mathbb R})=H^0(\Gamma_{{\mathbb R}};T)\to H^{1,0}(\Gamma_{{\mathbb R}};T_{sc}\to T)\to H^1(\Gamma_{{\mathbb R}};T_{sc}).$$
 Mais $T_{sc}$ est un tore induit, donc le dernier groupe est nul. La premi\`ere fl\`eche est donc surjective et, a fortiori, l'homomorphisme $G({\mathbb R})\to G_{ab}({\mathbb R})$ est surjectif. $\square$

 Prouvons (1). Rappelons la suite exacte
 $$1\to \hat{G}\to \hat{G}_{\natural}\to \hat{C}_{\natural}\to 1.$$
 On fixe une paire de Borel \'epingl\'ee de $\hat{G}_{\natural}$ conserv\'ee par l'action galoisienne. Elle se restreint naturellement en une telle paire pour $\hat{G}$. En notant $\hat{T}$ et $\hat{T}_{\natural}$ leurs tores, on a la suite exacte
 $$1\to \hat{T}\to \hat{T}_{\natural}\to \hat{C}_{\natural}\to 1.$$
   Elle se restreint en une suite exacte
 $$(3) \qquad 1\to Z(\hat{G})\to Z(\hat{G}_{\natural})\to \hat{C}_{\natural}\to 1.$$
 Soit ${\bf G}'=(G',{\cal G}',s)$ une donn\'ee endoscopique de $(G,\tilde{G})$.  On ne perd rien \`a supposer que $s\in \hat{T}$. On a $\hat{G}'=\hat{G}_{s}$ (on rappelle que $\hat{G}_{s}$ est la composante neutre du centralisateur $Z_{\hat{G}}(s)$ de $s$ dans $\hat{G}$) . Posons $\hat{G}'_{\natural}=(\hat{G}_{\natural})_{s}$. On a $\hat{G}'_{\natural}=\hat{T}_{\natural}\hat{G}'=Z(\hat{G}_{\natural})\hat{G}'$ et la suite
 $$(4) \qquad 1\to \hat{G}'\to \hat{G}'_{\natural}\to \hat{C}_{\natural}\to 1$$
 est exacte. On pose ${\cal G}'_{\natural}=Z(\hat{G}_{\natural}){\cal G}'$. Ce sous-groupe de $^LG$ agit sur $\hat{G}'_{\natural}$, d'o\`u une action galoisienne sur ce groupe compatible avec la suite exacte ci-dessus. On introduit  un groupe $G'_{\natural}$ quasi-d\'eploy\'e sur ${\mathbb R}$ de sorte que $\hat{G}'_{\natural}$, muni de son action galoisienne,  en soit un groupe dual. On a une suite exacte duale de la pr\'ec\'edente
 $$1\to C_{\natural}\to G'_{\natural}\to G'\to 1.$$
 La donn\'ee ${\bf G}'_{\natural}=(G'_{\natural},{\cal G}'_{\natural},s)$ est une donn\'ee endoscopique pour $(G_{\natural},\tilde{G}_{\natural})$. On obtient une application ${\bf G}'\mapsto {\bf G}'_{\natural}$ qui, \`a une classe d'\'equivalence de donn\'ee endoscopique pour $(G,\tilde{G})$, associe une classe d'\'equivalence de donn\'ee endoscopique pour $(G_{\natural},\tilde{G}_{\natural})$.  Montrons qu'elle est bijective. Il suffit de d\'efinir son inverse. Pour cela, partons d'une donn\'ee endoscopique ${\bf G}'_{\natural}=(G'_{\natural},{\cal G}'_{\natural},s)$ pour $(G_{\natural},\tilde{G}_{\natural})$. On suppose $s\in \hat{T}_{\natural}$. Quitte \`a multiplier $s$ par un \'el\'ement de $Z(\hat{G}_{\natural})$, on peut gr\^ace \`a (3) supposer que $s\in \hat{T}$. On a comme ci-dessus $\hat{G}'_{\natural}=(\hat{G}_{\natural})_{s}=\hat{T}_{\natural}\hat{G}_{s}$, d'o\`u
 $\hat{G}'_{\natural}\cap \hat{G}=\hat{T}\hat{G}_{s}=\hat{G}_{s}$ puisque $\hat{T}\subset \hat{G}_{s}$.  En posant $\hat{G}'=\hat{G}_{s}$, on a encore la suite exacte (4). On introduit un groupe $G'$ quasi-d\'eploy\'e sur ${\mathbb R}$ de sorte que $\hat{G}'$ soit dual de $G'$. Posons ${\cal G}'={\cal G}'_{\natural}\cap {^LG}$. Montrons que le triplet ${\bf G}'=(G',{\cal G}',s)$ est une donn\'ee endoscopique pour $(G,\tilde{G})$.  La seule condition non \'evidente est la suivante. Il existe un cocycle $a:W_{{\mathbb R}}\to Z(\hat{G}_{\natural})$, qui est un cobord, de sorte que $sgw(s)^{-1}=a(w)g$ pour tout $w\in W_{{\mathbb R}}$ et  tout $(g,w)\in {\cal G}'_{\natural}$. Quand on se restreint \`a $(g,w)\in {\cal G}'$, tous les termes sauf \'eventuellement $a(w)$ sont dans $\hat{G}$. Ce dernier terme est donc lui-aussi dans $\hat{G}$, donc dans $Z(\hat{G}_{\natural})\cap \hat{G}=Z(\hat{G})$. Donc $a$ est  un cocycle de $W_{{\mathbb R}}$ dans $Z(\hat{G})$. Il faut voir que c'est un cobord. Cela r\'esulte de
 
 (5) l'homomorphisme naturel $H^1(W_{{\mathbb R}};Z(\hat{G}))\to H^1(W_{{\mathbb R}};Z(\hat{G}_{\natural}))$ est injectif.
 
 Il s'ins\`ere dans une suite exacte
 $$1\to Z(\hat{G})^{\Gamma_{{\mathbb R}}}\to Z(\hat{G}_{\natural})^{\Gamma_{{\mathbb R}}}\to \hat{C}_{\natural}^{\Gamma_{{\mathbb R}}} \to H^1(W_{{\mathbb R}};Z(\hat{G}))\to H^1(W_{{\mathbb R}};Z(\hat{G}_{\natural})).$$
 La suite (3) entra\^{\i}ne que l'homomorphisme $Z(\hat{G}_{\natural})^{\Gamma_{{\mathbb R}},0}\to \hat{C}_{\natural}^{\Gamma_{{\mathbb R}},0}$ est surjectif. Or $\hat{C}$ est induit donc $ \hat{C}_{\natural}^{\Gamma_{{\mathbb R}}}$ est connexe. Il en r\'esulte que le deuxi\`eme homomorphisme de la suite ci-dessus est surjectif, d'o\`u (5).
 
 Il est clair que l'application ${\bf G}'_{\natural}\mapsto {\bf G}'$ que l'on vient de construire est inverse de la pr\'ec\'edente, ce qui prouve la bijectivit\'e de ces deux applications.  
  
 Consid\'erons deux donn\'ees ${\bf G}'$ et ${\bf G}'_{\natural}$ qui se correspondent ainsi. Rappelons que, puisque $(G,\tilde{G},{\bf a})$ est quasi-d\'eploy\'e et \`a torsion int\'erieure, ${\bf G}'$ est relevante si et seulement si $\tilde{G}'({\mathbb R})\not=\emptyset$ ([I] lemme 1.9). De m\^eme pour $\tilde{G}'_{\natural}$. Puisque $C_{\natural}$ est induit, l'application $\tilde{G}'_{\natural}({\mathbb R})\to \tilde{G}'({\mathbb R})$ est surjective. Les deux ensembles sont vides ou non vides simultan\'em\'ent. Donc ${\bf G}'$ est relevante si et seulement si ${\bf G}'_{\natural}$ l'est.
 Supposons ces donn\'ees relevantes.  Fixons des donn\'ees suppl\'ementaires $G'_{1}$, $\tilde{G}'_{1}$, $C_{1}$, $\hat{\xi}_{1}$, $\Delta_{1}$ pour ${\bf G}'$. Notons $G'_{\natural,1}$ le produit fibr\'e de $G'_{\natural}$ et $G'_{1}$ au-dessus de $G'$. Le groupe dual $\hat{G}'_{\natural,1}$ s'identifie au quotient de $\hat{G}'\times Z(\hat{G}_{\natural})\times Z(\hat{G}_{1})$ par le sous-groupe des $(z_{1},z_{2},z_{3})\in Z(\hat{G}')$ tels que $z_{1}z_{2}z_{3}=1$. L'homomorphisme $\hat{\xi}_{1}:{\cal G}'\to {^LG}'_{1}$ se prolonge en un homomorphisme $\hat{\xi}_{\natural,1}:{\cal G}'_{\natural}\to {^LG}'_{\natural,1}$ de la fa\c{c}on suivante. Soit $(g_{\natural},w)\in {\cal G}'_{\natural}$. On \'ecrit $g_{\natural}=z_{\natural}g$, avec $z_{\natural}\in Z(\hat{G}_{\natural})$ et $(g,w)\in {\cal G}'$. Ecrivons $\hat{\xi}_{1}(g,w)=z_{1}h\times w$, avec $z_{1}\in Z(\hat{G}'_{1})$ et $h\in \hat{G}'$. On note $\hat{\xi}_{\natural,1}(g_{\natural},w)$ l'image de $(h,z_{\natural},z_{1})$ dans $\hat{G}'_{\natural,1}$. On v\'erifie que ce terme ne d\'epend pas des choix faits et que l'application $\hat{\xi}_{\natural,1}$ ainsi d\'efinie est un homomorphisme. On note $\tilde{G}'_{\natural,1}$ le produit fibr\'e de $\tilde{G}'_{\natural}$ (qui est l'espace de la donn\'ee endoscopique ${\bf G}'_{\natural}$) et de $\tilde{G}'_{1}$ au-dessus de $\tilde{G}'$. 
  Soit $(\delta_{\natural,1},\gamma_{\natural})$ un couples d'\'el\'ements de $\tilde{G}'_{\natural,1}({\mathbb R})\times \tilde{G}_{\natural}({\mathbb R})$ qui sont fortement r\'eguliers et qui se correspondent. Ce couple a une image naturelle $(\delta_{1},\gamma)\in \tilde{G}'_{1}({\mathbb R})\times \tilde{G}({\mathbb R})$ form\'ee d'\'el\'ements fortement r\'eguliers qui se correspondent. On pose $\Delta_{\natural,1}(\delta_{\natural,1},\gamma_{\natural})=\Delta_{1}(\delta_{1},\gamma)$. On v\'erifie que $G'_{\natural,1}$, $\tilde{G}'_{\natural,1}$, $C_{1}$, $\hat{\xi}_{\natural,1}$, $\Delta_{\natural,1}$ est un ensemble de donn\'ees auxiliaires pour ${\bf G}'_{\natural}$. 
  
  Consid\'erons maintenant d'autres donn\'ees index\'ees par $\flat$ comme en (1). Pour la m\^eme donn\'ee ${\bf G}'$ et les m\^emes donn\'ees auxiliaires, on construit de m\^eme ${\bf G}'_{\flat}$ et des donn\'ees auxiliaires $G'_{\flat,1}$,..., $\Delta_{\flat,1}$. On peut appliquer aussi ces constructions pour les produits fibr\'es $G'_{\natural,\flat}$ et $\tilde{G}'_{\natural,\flat}$. On a donc une donn\'ee endoscopique ${\bf G}'_{\natural,\flat}$ pour $(G'_{\natural,\flat},\tilde{G}'_{\natural,\flat})$ et des donn\'ees auxiliaires $G'_{\natural,\flat,1}$,...,$\Delta_{\natural,\flat,1}$. On v\'erifie que $G'_{\natural,\flat,1}$ est le produit fibr\'e de $G'_{\natural,1}$ et $G'_{\flat,1}$ au-dessus de $G'_{1}$ et que $\tilde{G}'_{\natural,\flat,1}$ est le produit fibr\'e de $\tilde{G}'_{\natural,1}$ et $\tilde{G}'_{\flat,1}$ au-dessus de $G'_{1}$. D'apr\`es [I] 1.12 et  (2) ci-dessus,  on a un homomorphisme naturel $G'_{ab}({\mathbb R})\to G_{ab}({\mathbb R})$ et une application compatible $\tilde{G}'_{ab}({\mathbb R})\to \tilde{G}_{ab}({\mathbb R})$.  De m\^eme, on a un homomorphisme $G'_{\natural,\flat,ab}({\mathbb R})\to G_{\natural,\flat,ab}({\mathbb R})$  et une application compatible  $\tilde{G}'_{\natural,\flat,ab}({\mathbb R})\to \tilde{G}_{\natural,\flat,ab}({\mathbb R})$. On les compose en un homomorphisme $G'_{\natural,\flat,1,ab}({\mathbb R})\to G_{\natural,\flat,ab}({\mathbb R})$  et une application compatible  $\tilde{G}'_{\natural,\flat,1,ab}({\mathbb R})\to \tilde{G}_{\natural,\flat,ab}({\mathbb R})$. Gr\^ace \`a (2), le caract\`ere $\lambda_{\natural,\flat}$ se factorise en un homomorphisme de $G_{\natural,\flat}({\mathbb R})$ et l'application $\tilde{\lambda}_{\natural,\flat}$ se factorise en une application d\'efinie sur $\tilde{G}_{\natural,\flat}({\mathbb R})$. Par composition avec les applications pr\'ec\'edentes, on obtient un caract\`ere $\lambda'_{\natural,\flat,1}$ de $G'_{\natural,\flat,1}({\mathbb R})$ et une application compatible $\tilde{\lambda}'_{\natural,\flat,1}$ sur $\tilde{G}'_{\natural,\flat,1}({\mathbb R})$. 
  
  Introduisons des notations pour nos applications
  $$\begin{array}{ccc}D_{tr-orb}(\tilde{G}_{\natural}({\mathbb R}))&&D_{tr-orb}(\tilde{G}_{\flat}({\mathbb R}))\\ p_{\lambda_{\natural}}\downarrow\,\,&&p_{\lambda_{\flat}}\downarrow\,\,\\ D_{tr-orb,\lambda_{\natural}}(\tilde{G}_{\natural}({\mathbb R}))&\stackrel{\iota}{\to}&D_{tr-orb,\lambda_{\flat}}(\tilde{G}_{\flat}({\mathbb R}))\\ \end{array}$$
  On veut prouver que, pour   $\boldsymbol{\gamma}_{\natural}\in D_{tr-orb}(\tilde{G}_{\natural}({\mathbb R}))$, il existe  $\boldsymbol{\gamma}_{\flat}\in D_{tr-orb}(\tilde{G}_{\flat}({\mathbb R}))$ tel que $\iota\circ p_{\lambda_{\natural}}(\boldsymbol{\gamma}_{\natural})=p_{\lambda_{\flat}}(\boldsymbol{\gamma}_{\flat})$. 
 Par d\'efinition, on peut \'ecrire $\boldsymbol{\gamma}_{\natural}$ comme somme d'un \'el\'ement $\boldsymbol{\gamma}_{\natural,orb}$ qui est une honn\^ete int\'egrale orbitale et d'une somme sur  ${\bf G}'\in {\cal E}(\tilde{G})$ tel que $G'\not=G$, d'\'el\'ements $transfert(\boldsymbol{\delta}_{\natural})$, o\`u $\boldsymbol{\delta}_{\natural}\in D^{st}_{tr-orb}({\bf G}'_{\natural})$. Pour les int\'egrales orbitales, il n'y a pas de probl\`eme. Il existe presque par d\'efinition un \'el\'ement $\boldsymbol{\gamma}_{\flat,orb}\in D_{orb}(\tilde{G}_{\flat}({\mathbb R}))$ tel que $\iota\circ p_{\lambda_{\natural}}(\boldsymbol{\gamma}_{\natural,orb})=p_{\lambda_{\flat}}(\boldsymbol{\gamma}_{\flat,orb})$. Fixons ${\bf G}'$. En fixant des donn\'ees auxiliaires comme ci-dessus, on identifie $\boldsymbol{\delta}_{\natural}$ \`a un \'el\'ement de $D^{st}_{tr-orb,\lambda_{1}}(\tilde{G}'_{\natural,1}({\mathbb R}))$ ($\lambda_{1}$ est le caract\`ere de $C_{1}({\mathbb R})$). On le rel\`eve en un \'el\'ement $\boldsymbol{\delta}_{\natural,1}\in D^{st}_{tr-orb}(\tilde{G}'_{\natural,1}({\mathbb R}))$. On  a le diagramme suivant
 $$(6) \qquad \begin{array}{ccccc}D^{st}_{g\acute{e}om}(\tilde{G}'_{\natural,1}({\mathbb R}))&\stackrel{p_{\lambda_{1}}}{\to}&D^{st}_{g\acute{e}om,\lambda_{1}}(\tilde{G}'_{\natural,1}({\mathbb R}))&\stackrel{transfert}{\to}&D_{g\acute{e}om}(G_{\natural}({\mathbb R}))\\ p_{\lambda_{\natural}}\downarrow\,\,&&p_{\lambda_{\natural}}\downarrow\,\,&& p_{\lambda_{\natural}}\downarrow\,\,\\ 
D^{st}_{g\acute{e}om,\lambda_{\natural}}(\tilde{G}'_{\natural,1}({\mathbb R}))&\stackrel{p_{\lambda_{1}}}{\to}&D^{st}_{g\acute{e}om,\lambda_{\natural}\times\lambda_{1}}(\tilde{G}'_{\natural,1}({\mathbb R}))& &D_{g\acute{e}om,\lambda_{\natural}}(G_{\natural}({\mathbb R}))\\ \end{array}$$ 
 On v\'erifie que la suite du bas se compl\`ete en un homomorphisme que l'on peut appeler $$transfert:D^{st}_{g\acute{e}om,\lambda_{\natural}\times\lambda_{1}}(\tilde{G}'_{\natural,1}({\mathbb R}))\to D_{g\acute{e}om,\lambda_{\natural}}(G_{\natural}({\mathbb R}))$$
 qui rend ce diagramme commutatif. On a alors
 $$p_{\lambda_{\natural}}\circ transfert(\boldsymbol{\delta}_{\natural})=p_{\lambda_{\natural}}\circ transfert\circ p_{\lambda_{1}}(\boldsymbol{\delta}_{\natural,1})=transfert\circ p_{\lambda_{1}}\circ (\boldsymbol{\delta}'_{\natural,1}),$$
 o\`u $\boldsymbol{\delta}'_{\natural,1}=p_{\lambda_{\natural}}(\boldsymbol{\delta}_{\natural,1})$. 
 De l'application $\tilde{\lambda}_{\natural,\flat,1}$ se d\'eduit un isomorphisme
 $$\iota^{st}_{1}:D^{st}_{g\acute{e}om,\lambda_{\natural}}(\tilde{G}'_{\natural,1}({\mathbb R}))\to D^{st}_{g\acute{e}om,\lambda_{\flat}}(\tilde{G}'_{\flat,1}({\mathbb R})).$$
 Il r\'esulte de sa construction qu'il se descend en un isomorphisme encore not\'e
$$\iota^{st}_{1}: D^{st}_{g\acute{e}om,\lambda_{\natural}\times\lambda_{1}}(\tilde{G}'_{\natural,1}({\mathbb R}))\to D^{st}_{g\acute{e}om,\lambda_{\flat}\times\lambda_{1}}(\tilde{G}'_{\flat,1}({\mathbb R})).$$
On v\'erifie que le diagramme suivant est commutatif
$$\begin{array}{ccccc} D^{st}_{g\acute{e}om,\lambda_{\natural}}(\tilde{G}'_{\natural,1}({\mathbb R}))&\stackrel{p_{\lambda_{1}}}{\to}&D^{st}_{g\acute{e}om,\lambda_{\natural}\times\lambda_{1}}(\tilde{G}'_{\natural,1}({\mathbb R}))& \stackrel{transfert}{\to}&D_{g\acute{e}om,\lambda_{\natural}}(G_{\natural}({\mathbb R}))\\ \iota^{st}_{1}\downarrow\,\,&&\iota^{st}_{1}\downarrow\,\,&&\iota \downarrow\,\,\\  D^{st}_{g\acute{e}om,\lambda_{\flat}}(\tilde{G}'_{\flat,1}({\mathbb R}))&\stackrel{p_{\lambda_{1}}}{\to}&D^{st}_{g\acute{e}om,\lambda_{\flat}\times\lambda_{1}}(\tilde{G}'_{\flat,1}({\mathbb R}))& \stackrel{transfert}{\to}&D_{g\acute{e}om,\lambda_{\flat}}(G_{\flat}({\mathbb R}))\\ \end{array}$$
En posant $\boldsymbol{\delta}'_{\flat,1}=\iota^{st}_{1}(\boldsymbol{\delta}'_{\natural,1})$, on obtient 
 $$\iota\circ p_{\lambda_{\natural}}\circ transfert(\boldsymbol{\delta}_{\natural})=\iota\circ transfert\circ p_{\lambda_{1}}(\boldsymbol{\delta}'_{\natural,1})= transfert\circ p_{\lambda_{1}}(\boldsymbol{\delta}'_{\flat,1}).$$
 Puisque $G'\not=G$, on peut appliquer (1) par r\'ecurrence en rempla\c{c}ant $\tilde{G}$, $\tilde{G}_{\natural}$ et $\tilde{G}_{\flat}$ par $\tilde{G}'_{1}$, $\tilde{G}'_{\natural,1}$ et $\tilde{G}'_{\flat,1}$. Puisque $\boldsymbol{\delta}'_{\natural,1}\in D^{st}_{tr-orb,\lambda_{\natural}}(\tilde{G}'_{\natural,1}({\mathbb R}))$, on obtient que $\boldsymbol{\delta}'_{\flat,1}\in D^{st}_{tr-orb,\lambda_{\flat}}(\tilde{G}'_{\flat,1}({\mathbb R}))$. On rel\`eve $\boldsymbol{\delta}'_{\flat,1}$ en un \'el\'ement $\boldsymbol{\delta}_{\flat,1}\in D^{st}_{tr-orb}(\tilde{G}'_{\flat,1}({\mathbb R}))$. Par un diagramme similaire \`a (6) pour les indices $\flat$, on a
 $$ transfert\circ p_{\lambda_{1}}(\boldsymbol{\delta}'_{\flat,1})=p_{\lambda_{\flat}}\circ transfert\circ p_{\lambda_{1}}(\boldsymbol{\delta}_{\flat,1}).$$
 L'\'el\'ement $p_{\lambda_{1}}(\boldsymbol{\delta}_{\flat,1})$ s'identifie \`a un \'el\'ement $\boldsymbol{\delta}_{\flat}\in D^{st}_{tr-orb}({\bf G}'_{\flat})$ et l'\'egalit\'e ci-dessus se r\'ecrit
 $$ transfert\circ p_{\lambda_{1}}(\boldsymbol{\delta}'_{\flat,1})=p_{\lambda_{\flat}}\circ transfert(\boldsymbol{\delta}_{\flat}).$$
 On a obtenu
 $$\iota\circ p_{\lambda_{\natural}}\circ transfert(\boldsymbol{\delta}_{\natural})=p_{\lambda_{\flat}}\circ transfert(\boldsymbol{\delta}_{\flat}).$$
 
 Notons $\boldsymbol{\gamma}_{\flat}$ la somme de $\boldsymbol{\gamma}_{\flat,orb}$ et de la somme pour tout ${\bf G}'\in {\cal E}(\tilde{G})$, avec $G'\not=G$, des \'el\'ements $transfert(\boldsymbol{\delta}_{\flat})$ que l'on vient de construire. Par d\'efinition, on a $\boldsymbol{\gamma}_{\flat}\in D_{tr-orb}(\tilde{G}_{\flat}({\mathbb R}))$ et on a prouv\'e l'\'egalit\'e $\iota\circ p_{\lambda_{\natural}}(\boldsymbol{\gamma}_{\natural})=p_{\lambda_{\flat}}(\boldsymbol{\gamma}_{\flat})$. Cela d\'emontre la premi\`ere assertion de (1).
 
 Prouvons la seconde assertion. On voit qu'elle r\'esulte de la propri\'et\'e suivante, qui porte sur une seule s\'erie de donn\'ees:
 
 (7) on a l'\'egalit\'e $D^{st}_{tr-orb,\lambda_{\natural}}(\tilde{G}_{\natural}({\mathbb R}))=D_{tr-orb,\lambda_{\natural}}(\tilde{G}_{\natural}({\mathbb R}))\cap D^{st}_{g\acute{e}om,\lambda_{\natural}}(\tilde{G}_{\natural}({\mathbb R}))$.  
 
 L'inclusion du membre de gauche dans celui de droite est \'evidente par d\'efinition. L'inclusion inverse signifie que, si $\boldsymbol{\gamma}$ est un \'el\'ement de $D_{tr-orb}(\tilde{G}_{\natural}({\mathbb R}))$ tel que $p_{\lambda_{\natural}}(\boldsymbol{\gamma})$ est stable, alors il existe $\boldsymbol{\delta}\in D^{st}_{tr-orb}(\tilde{G}_{\natural}({\mathbb R}))$ tel que $p_{\lambda_{\natural}}(\boldsymbol{\gamma})=p_{\lambda_{\natural}}(\boldsymbol{\delta})$. Avant de le prouver, \'enon\c{c}ons  deux propri\'et\'es de l'espace $D_{tr-orb}(\tilde{G}({\mathbb R}))$. 
  Le groupe $Z_{G}({\mathbb R})$ agit par translations sur $\tilde{G}({\mathbb R})$ et cons\'equemment sur l'espace de distributions $D_{g\acute{e}om}(\tilde{G}({\mathbb R}))$. Pour $z\in Z_{G}({\mathbb R})$, on note $\boldsymbol{\gamma}\mapsto \boldsymbol{\gamma}^z$ cette action. Alors
 
 (8) pour $z\in Z_{G}({\mathbb R})$ et $\boldsymbol{\gamma}\in D_{tr-orb}(\tilde{G}({\mathbb R}))$, on a $\boldsymbol{\gamma}^z\in D_{tr-orb}(\tilde{G}({\mathbb R}))$.
 
 Appelons caract\`ere affine de $\tilde{G}({\mathbb R})$ une fonction $\tilde{\chi}$ sur $\tilde{G}({\mathbb R})$ telle qu'il existe un caract\`ere $\chi$ de $G({\mathbb R})$ de sorte que $\tilde{\chi}(x\gamma)=\chi(x)\tilde{\chi}(\gamma)$ pour tous $x\in G({\mathbb R})$ et $\gamma\in \tilde{G}({\mathbb R})$. Un caract\`ere affine agit par multiplication sur $C_{c}^{\infty}(\tilde{G}({\mathbb R})$. Cette action se quotiente en une action sur $I(\tilde{G}({\mathbb R}))$ et on a aussi une action duale sur $D_{g\acute{e}om}(\tilde{G}({\mathbb R}))$ que l'on note $(\tilde{\chi},\boldsymbol{\gamma})\mapsto \tilde{\chi}\boldsymbol{\gamma}$. Alors
 
 (9) pour tout caract\`ere affine $\tilde{\chi}$ de $\tilde{G}({\mathbb R})$ et tout $\boldsymbol{\gamma}\in D_{tr-orb}(\tilde{G}({\mathbb R}))$, on a $\tilde{\chi}\boldsymbol{\gamma}\in D_{tr-orb}(\tilde{G}({\mathbb R}))$. 
 
 Les propri\'et\'es (8) et (9) sont claires si l'on remplace l'espace $D_{tr-orb}(\tilde{G}({\mathbb R}))$ par $D_{orb}(\tilde{G}({\mathbb R}))$. On les prouve alors par r\'ecurrence en \'etudiant comment se comportent les actions de $Z_{G}({\mathbb R})$ ou d'un caract\`ere affine relativement au transfert. 
 
 Revenons \`a la preuve de (7). Soit $\boldsymbol{\gamma}\in D_{tr-orb}(\tilde{G}_{\natural}({\mathbb R}))$. Introduisons le groupe d\'eriv\'e $G_{\natural,der}$ de $G_{\natural}$. L'image dans $G_{\natural,der}({\mathbb R})\backslash \tilde{G}_{\natural}({\mathbb R})$ du support de $\boldsymbol{\gamma}$ est finie. On peut donc \'ecrire $\boldsymbol{\gamma}=\sum_{i=1,...,n}\boldsymbol{\gamma}_{i}$,  de sorte que
 
- pour tout $i$,  $\boldsymbol{\gamma}_{i}$ appartient \`a $ D_{g\acute{e}om}(\tilde{G}_{\natural}({\mathbb R}))$;
 
 - l'image du support de $\boldsymbol{\gamma}_{i}$ dans $G_{\natural,der}({\mathbb R})\backslash \tilde{G}_{\natural}({\mathbb R})$  est un unique point $x_{i}$;
  
 - pour $i\not=j$, on a $x_{i}\not= x_{j}$.
  
  On peut r\'ecup\'erer chaque $\boldsymbol{\gamma}_{i}$ comme combinaison lin\'eaire de $\tilde{\chi}\boldsymbol{\gamma}$ pour des caract\`eres affines $\tilde{\chi}$ convenables de $\tilde{G}_{\natural}({\mathbb R})$. D'apr\`es (9), on a donc $\boldsymbol{\gamma}_{i}\in D_{tr-orb}(\tilde{G}_{\natural}({\mathbb R}))$  pour tout $i$. Soient $i\not=j$, supposons que $x_{i}=cx_{j}$ avec $c\in C_{\natural}({\mathbb R})$.  Posons $\boldsymbol{\gamma}'=\boldsymbol{\gamma}+\lambda_{\natural}(c)^{-1}\boldsymbol{\gamma}_{j}^c-\boldsymbol{\gamma}_{j}$. Cet \'el\'ement est  encore dans $D_{tr-orb}(\tilde{G}_{\natural}({\mathbb R}))$ d'apr\`es (8)  et v\'erifie $p_{\lambda_{\natural}}(\boldsymbol{\gamma}')=p_{\lambda_{\natural}}(\boldsymbol{\gamma})$. Mais, pour $\boldsymbol{\gamma}'$, la composante $\boldsymbol{\gamma}_{j}$ est remplac\'ee par $\lambda_{\natural}(c)^{-1}\boldsymbol{\gamma}_{j}^c$ et la projection de son support est $x_{i}$. Cette composante s'ajoute \`a $\boldsymbol{\gamma}_{i}$ et  on a diminu\'e le nombre des composantes. En poursuivant, on arrive \`a un \'el\'ement que l'on note encore $\boldsymbol{\gamma}$, qui appartient \`a $D_{tr-orb}(\tilde{G}_{\natural}({\mathbb R}))$ et a m\^eme image par $p_{\lambda_{\natural}}$ que le $\boldsymbol{\gamma}$ initial, mais dont l'ensemble $\{x_{1},...,x_{n}\}$ associ\'e v\'erifie la condition: si $i\not=j$, on a $x_{i}\not \in C_{\natural}({\mathbb R})x_{j}$. Posons $\Delta=C_{\natural}({\mathbb R})\cap G_{\natural,der}({\mathbb R})$. C'est un groupe fini. Gr\^ace \`a (8), on peut moyenner $\boldsymbol{\gamma}$ par ce groupe sans changer les propri\'et\'es ci-dessus de cet \'el\'ement et supposer que $\boldsymbol{\gamma}^c=\lambda_{\natural}(c)\boldsymbol{\gamma}$ pour tout $c\in \Delta$. Supposons alors que $p_{\lambda_{\natural}}(\boldsymbol{\gamma})$ soit stable. On a prouv\'e en [II] 1.10(5) qu'un $\boldsymbol{\gamma}$ v\'erifiant toutes les hypoth\`eses ci-dessus \'etait stable. Donc $\boldsymbol{\gamma}\in D_{tr-orb}^{st}(\tilde{G}_{\natural}({\mathbb R}))$. Cela prouve (7) et la seconde assertion de (1). $\square$

 \bigskip
 
 \subsection{Premi\`eres propri\'et\'es de l'espace $D_{tr-orb}(\tilde{G}({\mathbb R}),\omega)$}
 
 Consid\'erons un triplet $(G,\tilde{G},{\bf a})$ quelconque. Soit ${\cal O}$ une classe de conjugaison stable d'\'el\'ements semi-simples dans $\tilde{G}({\mathbb R})$. Posons $D_{tr-orb}({\cal O},\omega)=D_{tr-orb}(\tilde{G}({\mathbb R}),\omega)\cap D_{g\acute{e}om}({\cal O},\omega)$. Dans le cas o\`u $(G,\tilde{G},{\bf a})$ est quasi-d\'eploy\'e et \`a torsion int\'erieure, on pose $D_{tr-orb}^{st}({\cal O},\omega)=D_{tr-orb}^{st}(\tilde{G}({\mathbb R}),\omega)\cap D^{st}_{g\acute{e}om}({\cal O},\omega)$. On a
 
 (1) $D_{tr-orb}(\tilde{G}({\mathbb R}),\omega)=\oplus_{{\cal O}}D_{tr-orb}({\cal O},\omega)$ o\`u ${\cal O}$ parcourt toutes les classes de conjugaison stable d'\'el\'ements semi-simples dans $\tilde{G}({\mathbb R})$;
 
 (2) si $(G,\tilde{G},{\bf a})$ est quasi-d\'eploy\'e et \`a torsion int\'erieure,  $D^{st}_{tr-orb}(\tilde{G}({\mathbb R}),\omega)=\oplus_{{\cal O}}D^{st}_{tr-orb}({\cal O},\omega)$.
 
 Preuve. L'espace $D_{tr-orb}(\tilde{G}({\mathbb R}),\omega)$ est la somme de $D_{orb}(\tilde{G}({\mathbb R}),\omega)$ et des espaces $transfert(D_{tr-orb}({\bf G}'))$, o\`u ${\bf G}'$ parcourt ${\cal E}(\tilde{G},{\bf a})$, avec la restriction $G'\not=G$ si $(G,\tilde{G},{\bf a})$ est quasi-d\'eploy\'e et \`a torsion int\'erieure. Pour prouver (1), il suffit de montrer que chacun de ces espaces v\'erifie une d\'ecomposition analogue. C'est clair pour l'espace $D_{orb}(\tilde{G}({\mathbb R}),\omega)$. Pour un espace $transfert(D_{tr-orb}({\bf G}'))$, cela r\'esulte par r\'ecurrence de l'assertion (2) appliqu\'ee \`a ${\bf G}'$. Cela prouve (1). Supposons  $(G,\tilde{G},{\bf a})$ quasi-d\'eploy\'e et \`a torsion int\'erieure. Pour un \'el\'ement $\boldsymbol{\gamma}\in D_{g\acute{e}om}(\tilde{G}({\mathbb R}))$ s'\'ecrivant $\sum_{{\cal O}}\boldsymbol{\gamma}_{{\cal O}}$, o\`u $\boldsymbol{\gamma}_{{\cal O}}\in D_{g\acute{e}om}({\cal O})$, l'\'el\'ement $\boldsymbol{\gamma}$ est stable si et seulement si $\boldsymbol{\gamma}_{{\cal O}}$ est stable pour tout ${\cal O}$. Alors (2) r\'esulte de (1). $\square$
 
 Pour un $K$-triplet $(KG,K\tilde{G},{\bf a})$ et une classe de conjugaison stable ${\cal O}$ d'\'el\'ements semi-simples dans $K\tilde{G}({\mathbb R})$, on d\'efinit de m\^eme $D_{tr-orb}({\cal O},\omega)$ et on a encore
 
 (3) $D_{tr-orb}(K\tilde{G}({\mathbb R}),\omega)=\oplus_{{\cal O}}D_{tr-orb}({\cal O},\omega)$ .
 
 Revenons \`a un triplet $(G,\tilde{G},{\bf a})$. 
 Soit $\tilde{M}$ un espace de Levi de $\tilde{G}$. On a

(4) l'homomorphisme d'induction envoie $D_{tr-orb}(\tilde{M}({\mathbb R}))\otimes Mes(M({\mathbb R}))^*$ dans $D_{tr-orb}(\tilde{G}({\mathbb R}))\otimes Mes(G({\mathbb R}))^*$; si $(G,\tilde{G},{\bf a})$ est quasi-d\'eploy\'e et \`a torsion int\'erieure,  l'homomorphisme d'induction envoie $D_{tr-orb}^{st}(\tilde{M}({\mathbb R}))\otimes Mes(M({\mathbb R}))^*$ dans $D^{st}_{tr-orb}(\tilde{G}({\mathbb R}))\otimes Mes(G({\mathbb R}))^*$.

Preuve. On oublie les espaces de mesures. La deuxi\`eme assertion r\'esulte de la premi\`ere puisque l'induction conserve la stabilit\'e. La premi\`ere assertion est vraie si on remplace les espaces $D_{tr-orb}$ par les espaces d'int\'egrales orbitales $D_{orb}$. Il nous suffit donc de fixer une donn\'ee ${\bf M}'\in {\cal E}(\tilde{M})$, avec $M'\not=M$ dans le cas o\`u $(G,\tilde{G},{\bf a})$ est quasi-d\'eploy\'e et \`a torsion int\'erieure, de fixer $\boldsymbol{\delta}\in D_{tr-orb}^{st}({\bf M}')$ et de prouver que $(transfert(\boldsymbol{\delta}))^{\tilde{G}}$ appartient \`a $D_{tr-orb}(\tilde{G}({\mathbb R})))$. Il existe un \'el\'ement ${\bf G}'\in {\cal E}(\tilde{G})$ dont ${\bf M}'$ soit une "donn\'ee de Levi".   Puisque le transfert commute \`a l'induction, on a
$$(transfert(\boldsymbol{\delta}))^{\tilde{G}}=transfert(\boldsymbol{\delta}^{{\bf G}'}).$$
En raisonnant par r\'ecurrence,  on peut supposer que $\boldsymbol{\delta}^{{\bf G}'}$ appartient \`a $D_{tr-orb}^{st}({\bf G}')$. Alors l'\'el\'ement ci-dessus appartient par d\'efinition \`a $D_{tr-orb}(\tilde{G}({\mathbb R}))$. $\square$

 Pour un $K$-triplet $(KG,K\tilde{G},{\bf a})$ et un $K$-espace de Levi $K\tilde{M}\in {\cal L}(K\tilde{M}_{0})$, on a de m\^eme
 
 (5)  l'homomorphisme d'induction envoie $D_{tr-orb}(K\tilde{M}({\mathbb R}))\otimes Mes(M({\mathbb R}))^*$ dans \,\,\,\,\,\,
 $D_{tr-orb}(K\tilde{G}({\mathbb R}))\otimes Mes(G({\mathbb R}))^*$.

 \bigskip
 \subsection{Un lemme de s\'eparation}
 On consid\`ere un triplet $(G,\tilde{G},{\bf a})$  quelconque et  un espace de Levi $\tilde{M}$ de $\tilde{G}$.  On a d\'efini en [II] 3.1 un ensemble ${\cal J}_{\tilde{M}}^{\tilde{G}}$ et, pour tout $J$ dans cet ensemble, un espace $U_{J}$ de germes au point $1$ de fonctions d\'efinies presque partout sur $A_{\tilde{M}}(F)$. Le corps de base $F$ \'etait non-archim\'edien dans cette r\'ef\'erence. Les m\^emes d\'efinitions valent sur le corps de base ${\mathbb R}$. On ne les reprend pas en se contentant de renvoyer \`a [II] 3.1. Toutefois, il nous faut donner une d\'emonstration de la propri\'et\'e essentielle [II] 3.1(3). Celle donn\'ee dans cette r\'ef\'erence ne s'adapte pas au corps de base ${\mathbb R}$. Signalons en passant que la propri\'et\'e [II] 3.1(2) devient fausse sur ${\mathbb R}$. Mais elle ne nous servait qu'\`a d\'emontrer (3).
 
 \ass{Lemme}{Soient $J\in {\cal J}_{\tilde{M}}^{\tilde{G}}$  et $u\in U_{J}$. Supposons que $\tilde{M}\not=\tilde{G}$ et que $u$ soit \'equivalent \`a une constante. Alors $u=0$ et cette constante est nulle.}
 
 Preuve. L'\'el\'ement $J$ est form\'e de familles $\underline{\alpha}=(\alpha_{i})_{i=1,...,n}$, o\`u 

- $n=a_{\tilde{M}}-a_{\tilde{G}}$;

- les $\alpha_{i}$ sont des racines de $A_{\tilde{M}}$ dans $G$ lin\'eairement ind\'ependantes;

- le r\'eseau $\oplus_{i=1,...,n}{\mathbb Z}\alpha_{i}$ qu'elles engendrent dans $\mathfrak{a}_{\tilde{M}}^*({\mathbb R})$ est un r\'eseau fix\'e $R_{J}$. 

A une telle famille, on associe la fonction
$$a\mapsto u_{\underline{\alpha}}(a)=\prod_{i=1,...,n}log(\vert \alpha_{i}(a)-\alpha_{i}(a)^{-1}\vert _{{\mathbb R}})$$
d\'efinie presque partout sur $ A_{\tilde{M}}({\mathbb R})$. Introduisons la relation de $\pm $-\'equivalence dans $J$: $\underline{\alpha}=(\alpha_{i})_{i=1,...,n}$ est $\pm $-\'equivalent \`a $\underline{\alpha}'=(\alpha'_{i})_{i=1,...,n}$ si et seulement si on a l'\'egalit\'e ensembliste $\{\pm \alpha_{i}; i=1,...,n\}=\{\pm \alpha'_{i}; i=1,...,n\}$ (en adoptant une notation additive pour les racines). La fonction $u_{\underline{\alpha}}$ ne d\'epend que de la classe de $\pm $-\'equivalence de $\underline{\alpha}$. Fixons un sous-ensemble $\underline{J}\subset J$ de repr\'esentants des classes de $\pm $-\'equivalence.
L'\'el\'ement $u$ est une combinaison lin\'eaire
$$u=\sum_{\underline{\alpha}\in \underline{J}}c_{\underline{\alpha}}u_{\underline{\alpha}},$$
avec des coefficients complexes $c_{\underline{\alpha}}$. On se limite \`a un voisinage de $1$ dans $A_{\tilde{M}}({\mathbb R})$. Tout \'el\'ement dans un tel voisinage s'\'ecrit de fa\c{c}on unique $a=exp(H)$ avec $H$ proche de $0$. On pose $d(a)=\vert \vert H\vert \vert $, o\`u $\vert \vert .\vert \vert $ est la norme euclidienne fix\'ee sur $\mathfrak{a}_{\tilde{M}}({\mathbb R})\simeq {\cal A}_{\tilde{M}}$. 
Dire que $u$ est \'equivalent \`a une constante $c$ signifie qu'il existe $r>0$ de sorte que, si l'on se restreint \`a un domaine d\'efini par les relations
$$\vert \alpha(a)-1\vert _{{\mathbb R}}>C  d(a),$$
o\`u $C$ est une constante positive fix\'ee, il existe $C'>0$ tel que l'on ait une minoration
$$\vert u(a)-c\vert \leq C' d(a)^r$$
pour $a$ dans le domaine assez proche de $1$. 
Cette notion d'\'equivalence se descend \`a l'alg\`ebre de Lie. Pour chaque $\underline{\alpha}\in \underline{J}$, introduisons la fonction
$$H\mapsto v_{\underline{\alpha}}(H)=\prod_{i=1,...,n}log(\vert 2\alpha_{i}(H)\vert _{{\mathbb R}})$$
d\'efinie presque partout sur $\mathfrak{a}_{\tilde{M}}({\mathbb R})$.  Pour tout $i$, 
on a l'\'egalit\'e
$$log(\vert exp(\alpha_{i}(H))-exp(-\alpha_{i}(H))\vert _{{\mathbb R}})=log(\vert 2\alpha_{i}(H)\vert _{{\mathbb R}})+ log(\vert \frac{exp(\alpha_{i}(H))-exp(-\alpha_{i}(H))}{2\alpha_{i}(H)}\vert _{{\mathbb R}}).$$
Le second terme est analytique au voisinage de $H=0$ et nul en ce point. Cela entra\^{\i}ne que les fonctions $v_{\underline{\alpha}}$ et $H\mapsto u_{\underline{\alpha}}(exp(H))$ sont \'equivalentes. Posons
$$v=\sum_{\underline{\alpha}\in \underline{J}}c_{\underline{\alpha}}v_{\underline{\alpha}}.$$
Alors $v$ est \'equivalent \`a la constante $c$. Cela entra\^{\i}ne

(1) $v(H)=c$ pour tout $H\in \mathfrak{a}_{\tilde{M}}({\mathbb R})$. 

En effet, fixons un point $H$ en position g\'en\'erale. Consid\'erons l'ensemble $\{tH; t\in {\mathbb R}, 0<t<2\}$. Il est contenu dans un domaine comme ci-dessus. En cons\'equence, la limite de $v(tH)-c$ est nulle quand $t$ tend vers $0$. Pour tout $\underline{\alpha}\in J$, la fonction $t\mapsto v_{\underline{\alpha}}(tH)$ est polynomiale en $log(t)$ pour $t>0$. Donc aussi $v(tH)-c$. Quand $t$ tend vers $0$, $log(t)$ tend vers $+\infty$. Un  polyn\^ome en $log(t)$ ne peut tendre vers $0$ que s'il est identiquement nul. Donc $v(tH)-c=0$ pour tout $t$. Appliqu\'ee \`a $t=1$, cette relation donne (1).

Les fonctions $v_{\underline{\alpha}}$ se quotientent en des fonctions sur $\mathfrak{a}_{\tilde{M}}({\mathbb R})/\mathfrak{a}_{\tilde{G}}({\mathbb R})$. On ne perd rien \`a supposer, pour simplifier les notations, que $\mathfrak{a}_{\tilde{G}}=\{0\}$. 
Consid\'erons une fonction
$$f=\sum_{\underline{\alpha}\in \underline{J}}x_{\underline{\alpha}}v_{\underline{\alpha}},$$
avec des coefficients $x_{\underline{\alpha}}\in {\mathbb C}$. 
Notons $J(f)$ l'ensemble des $\underline{\alpha}\in \underline{J}$ tels que $x_{\underline{\alpha}}\not=0$. Consid\'erons la r\'eunion des $\underline{\alpha}\in J(f)$, vus comme des ensembles de formes lin\'eaires sur $\mathfrak{a}_{\tilde{M}}({\mathbb R})$. C'est un ensemble fini de racines, notons-le $\Sigma(f)$. Chacune d'elles d\'etermine l'hyperplan de $ \mathfrak{a}_{\tilde{M}}({\mathbb R})$ sur lequel elle s'annule. D'o\`u un ensemble fini d'hyperplans. Le compl\'ementaire dans $ \mathfrak{a}_{\tilde{M}}({\mathbb R})$ de la r\'eunion de ces hyperplans est r\'eunion finie de c\^ones. La fonction $f$ est clairement analytique sur chacun d'eux. Soit ${\cal C}$ l'un de ces c\^ones et soit $d\in {\mathbb C}$. On va montrer

(2) si $f(H)-d$ est identiquement nul sur ${\cal C}$, alors $J(f)$ est vide.

On raisonne par r\'ecurrence sur un entier $N\geq1$: on montre que les relations "$f(H)-d$ identiquement nul sur ${\cal C}$" et "$J(f)$ a $N$ \'el\'ements" sont contradictoires. L'assertion est  \'evidente si  $N=1$: une fonction $v_{\underline{\alpha}}$ n'est certainement pas constante sur un c\^one ouvert (le nombre de racines $n=a_{\tilde{M}}-a_{\tilde{G}}=a_{\tilde{M}}$ \'etant strictement positif d'apr\`es l'hypoth\`ese $\tilde{M}\not=\tilde{G}$). Soit  $N\geq2$, supposons que $f(H)-d$ soit identiquement nul sur ${\cal C}$ et que $J(f)$ ait $N$ \'el\'ements. Notons que, si l'on fixe $\underline{\alpha}=(\alpha_{i})_{i=1,...,n}\in J(f)$, ${\cal C}$ est contenu dans l'une des composantes connexes du compl\'ementaire des hyperplans noyaux des $\alpha_{i}$. Puisque les $\alpha_{i}$ sont lin\'eairement ind\'ependants, il en r\'esulte que ${\cal C}$ a au moins $n$ murs. C'est-\`a-dire qu'il  
  y a au moins $n$ hyperplans  ${\cal H}_{j}$ pour $j=1,...,n$, d'\'equations $\beta_{j}(H)=0$, avec $\beta_{j}\in \Sigma(f)$, de sorte .que l'intersection de ${\cal H}_{j}$ avec l'adh\'erence de ${\cal C}$ contienne un ouvert de ${\cal H}_{j}$.  On choisit de tels hyperplans. Il y a certainement un $j$  et un $\underline{\alpha}=(\alpha_{i})_{i=1,...,n}\in J(f)$ tel que $\beta_{j}\not=\pm \alpha_{i}$ pour tout $i$. Sinon, tout $\underline{\alpha}\in J(f)$ contiendrait $\pm \beta_{1},...,\pm \beta_{n}$, sa classe de $\pm $-\'equivalence serait uniquement d\'etermin\'ee et $J(f)$ ne contiendrait qu'un \'el\'ement. Fixons un $j$ comme ci-dessus et posons simplement $\beta=\beta_{j}$. Il existe un $\underline{\alpha}=(\alpha_{i})_{i=1,...,n}\in J(f)$ et un $i$  tel que $\beta= \alpha_{i}$. Cela traduit    simplement l'appartenance de $\beta_{j}$ \`a $\Sigma(f)$. Notons $J_{1}(f)$ le sous-ensemble des $\underline{\alpha}=(\alpha_{i})_{i=1,...,n}\in J(f)$ tels que $ \beta$ soit l'un des $\pm\alpha_{i}$ et  $J_{2}(f)$ son compl\'ementaire dans $J(f)$. Ce que l'on vient de dire signifie que $J_{1}(f)$ et $J_{2}(f)$ sont tous deux non vides. Notons ${\cal H}$ l'hyperplan d\'efini par $\beta(H)=0$, soit $\varpi$  un \'el\'ement de $\mathfrak{a}_{\tilde{M}}({\mathbb R})$ orthogonal \`a ${\cal H}$ et tel  que ${\cal C}$ soit contenu dans ${\cal H}+{\mathbb R}_{>0}\varpi$. Montrons que l'on peut trouver un ensemble ouvert non vide $U\subset {\cal H}$ et un r\'eel $\epsilon>0$ de sorte que

(3) $\{H+t\varpi; H\in U, 0<t<\epsilon\}\subset {\cal C}$;

(4) pour $\alpha\in \Sigma(f)$ avec $\alpha\not=\pm \beta$, on a $\alpha(H)\not=0$ pour tout $H\in U$. 

On peut certainement assurer (3) par l'hypoth\`ese que ${\cal H}$ est un bord de ${\cal C}$. Pour assurer (4), il suffit de retirer de $U$ le sous-ensemble des $H$ qui sont annul\'es par une racine $\alpha\in \Sigma(f)$, $\alpha\not=\pm \beta$. Il faut v\'erifier que l'ensemble obtenu reste non vide. Il suffit pour cela que les $\alpha$ en question ne s'annulent pas identiquement sur ${\cal H}$, ou encore qu'ils ne soient pas proportionnels \`a $\beta$. Supposons $\alpha\in \Sigma(f)$ et $\alpha=e \beta$, avec $e\in {\mathbb Q}$ et $e\not=\pm 1$. Par d\'efinition, $\alpha$ intervient dans une famille $\underline{\alpha}\in J(f)$ tandis que $\beta$ intervient dans une famille $\underline{\beta}\in J(f)$. Ces deux familles ne sont pas les m\^emes: une famille ne peut pas contenir $\beta$ et $e\beta$ car ces deux \'el\'ements ne sont pas lin\'eairement ind\'ependants. L'hypoth\`ese $\alpha=e\beta$ avec $e\not=\pm 1$ interdit aux \'el\'ements de la famille $\underline{\alpha}$ d'engendrer le m\^eme r\'eseau que les \'el\'ements de la famille $\underline{\beta}$. Cela contredit la d\'efinition de $J$. D'o\`u les assertions ci-dessus.

On fixe $U$ et $\epsilon$ comme ci-dessus. Fixons $H\in U$, soit  $t\in ]0,\epsilon[$.  Pour $\underline{\alpha}=(\alpha_{i})_{i=1,...,n}\in J_{1}(f)$, on peut supposer $\alpha_{1}=\pm \beta$. On pose $\underline{\alpha}'=(\alpha_{i})_{i=2,...,n}$. On a
$$v_{\underline{\alpha}}(H+t\varpi)=log(b t)v_{\underline{\alpha}'}(H+t\varpi),$$
o\`u $b=\vert 2\beta(\varpi)\vert _{{\mathbb R}}$. Dans $\underline{\alpha}'$ n'interviennent que des racines $\alpha$ v\'erifiant l'hypoth\`ese de (4), donc pour lesquelles $\alpha(H+t\varpi)$ ne s'annule pas en $t=0$. Il en r\'esulte que $t\mapsto v_{\underline{\alpha}'}(H+t\varpi)$ se prolonge en une fonction analytique en $t$ au voisinage de $t=0$.  Pour la m\^eme raison, si $\underline{\alpha}\in J_{2}(f)$, la fonction $t\mapsto v_{\underline{\alpha}}(H+t\varpi)$ se prolonge en une fonction analytique en $t$ au voisinage de $t=0$.  Posons
$$f'_{1}=\sum_{\underline{\alpha}\in J_{1}(f)}x_{\underline{\alpha}}v_{\underline{\alpha}'},$$
$$f_{2}=\sum_{\underline{\alpha}\in J_{2}(f)}x_{\underline{\alpha}}v_{\underline{\alpha}}.$$
Alors les fonctions $t\mapsto f'_{1}(H+t\varpi)$ et $t\mapsto f_{2}(H+t\varpi)$ se prolongent en des fonctions analytiques en $t$ au voisinage de $t=0$. De plus
$$(5)\qquad   log(bt)f_{1}'(H+t\varpi)+f_{2}(H+t\varpi)-d=f(H+t\varpi)-d=0$$
 d'apr\`es l'hypoth\`ese et (3). Si $f'_{1}(H+t\varpi)$ n'est pas identiquement nul, on en d\'eduit 
$$log(bt)=\frac{d-f_{2}(H+t\varpi)}{f'_{1}(H+t\varpi)},$$
et $log(bt)$ se prolonge en une fonction m\'eromorphe au voisinage de $t=0$. C'est impossible. Donc $f_{1}'(H+t\varpi)$ est identiquement nul. D'apr\`es (5), on a aussi $f_{2}(H+t\varpi)=d$.  Cela est vrai pour $H\in U$ et $t\in ]0,\epsilon[$. Donc $f_{2}(H)=d$ pour $H$  dans un ouvert non vide de $\mathfrak{a}_{\tilde{M}}({\mathbb R})$. La fonction $f_{2}$ est du m\^eme type que $f$. Il lui est associ\'e un ensemble fini de c\^ones dans lesquels elle est analytique. L'assertion pr\'ec\'edente entra\^{\i}ne que $f_{2}(H)=d$ pour $H$ dans l'un de ces c\^ones.  Le nombre d'\'el\'ements de $J_{2}(f)$ est compris entre $1$ et $N-1$. L'hypoth\`ese de r\'ecurrence dit que ces deux propri\'et\'es sont contradictoires. Cela ach\`eve la preuve de (2). 

Achevons la preuve du lemme. Les deux assertions (1) et (2) entra\^{\i}nent que les coefficients $c_{\underline{\alpha}}$ de $v$ sont nuls. Donc la fonction $u$ initiale est nulle et alors aussi la constante $c$. $\square$

 {\bf Variante.} Supposons $G=\tilde{G}$ et ${\bf a}=1$. On fixe une fonction $B$ comme en [II] 1.8.  On a d\'efini dans cette r\'ef\'erence l'ensemble $\Sigma(A_{M},B)$. On en d\'eduit en ensemble ${\cal J}_{M}^G(B)$ similaire au ${\cal J}_{\tilde{M}}^{\tilde{G}}$ pr\'ec\'edent. Le lemme reste valable pour cet ensemble.
 
 {\bf Variante.} Supposons $(G,\tilde{G},{\bf a})$ quasi-d\'eploy\'e et \`a torsion int\'erieure. Fixons un syst\`eme de fonctions $B$ comme en [II] 1.9. Pour un \'el\'ement semi-simple $\eta\in \tilde{M}({\mathbb R})$, on a d\'efini dans cette r\'ef\'erence l'ensemble $\Sigma(A_{M},B_{\eta})$. On en d\'eduit un ensemble ${\cal J}_{\tilde{M}}^{\tilde{G}}(B_{\eta})$ pour lequel le lemme reste valable.

 \bigskip
 
 \subsection{Programme d'extension des d\'efinitions}  
   On consid\`ere les trois situations suivantes.
   
   (A) On se donne un $K$-triplet $(KG,K\tilde{G},{\bf a})$, un $K$-espace de Levi $K\tilde{M}\in {\cal L}(K\tilde{M}_{0})$ et une classe de conjugaison stable semi-simple ${\cal O}$ dans $K\tilde{M}({\mathbb R})$. On note ${\cal O}^{K\tilde{G}}$ la classe de conjugaison stable dans $K\tilde{G}({\mathbb R})$ qui la contient. 
   
   (B) On se donne un triplet $(G,\tilde{G},{\bf a})$ quasi-d\'eploy\'e et \`a torsion int\'erieure, un espace de Levi $\tilde{M}$ de $\tilde{G}$ et une classe de conjugaison stable semi-simple ${\cal O}$ dans $\tilde{M}({\mathbb R})$. On note ${\cal O}^{\tilde{G}}$ la classe de conjugaison stable dans $\tilde{G}({\mathbb R})$ qui la contient. On fixe un syst\`eme de fonctions $B$ sur $\tilde{G}({\mathbb R})$ comme en [II] 1.9.
   
   (C) On se donne un groupe $G$, un Levi $M$ de $G$ et une fonction $B$ sur $G({\mathbb R})$ comme en [II] 1.8. On consid\`ere la classe de conjugaison stable dans $M({\mathbb R})$ r\'eduite \`a $\{1\}$.  
   
   Dans le cas (A), 
  \'ecrivons $K\tilde{M}=(\tilde{M}_{p})_{p\in \Pi^M}$. Pour $p,q\in \Pi^M$, les tores $A_{\tilde{M}_{p}}$ et $A_{\tilde{M}_{q}}$ s'identifient.  On note $A_{K\tilde{M}}$ ce tore commun. Il s'en d\'eduit une identification des ensembles ${\cal J}_{\tilde{M}_{p}}^{\tilde{G}_{p}}$ et ${\cal J}_{\tilde{M}_{q}}^{\tilde{G}_{q}}$. On note ${\cal J}_{K\tilde{M}}^{K\tilde{G}}$ cet ensemble commun. On se propose de d\'efinir
 
  - pour tout $J\in {\cal J}_{K\tilde{M}}^{K\tilde{G}}$, une application lin\'eaire
  $$\rho_{J}^{K\tilde{G}}:D_{tr-orb}( {\cal O},\omega)\otimes Mes(M({\mathbb R}))^*\to U_{J}\otimes (D_{g\acute{e}om}({\cal O},\omega)\otimes Mes(M({\mathbb R}))^*)/Ann_{{\cal O}}^{K\tilde{G}};$$

 - pour tout $\boldsymbol{\gamma}\in D_{tr-orb}({\cal O},\omega)\otimes Mes(M({\mathbb R}))^*$, une application lin\'eaire ${\bf f}\mapsto I_{K\tilde{M}}^{K\tilde{G}}(\boldsymbol{\gamma},{\bf f})$ sur $I(K\tilde{G}({\mathbb R}),\omega)\otimes Mes(G({\mathbb R}))$.
 
 Dans le cas (B), on se propose de d\'efinir
 
   - pour tout $J\in {\cal J}_{\tilde{M}}^{\tilde{G}}(B_{{\cal O}})$, deux applications lin\'eaires
  $$\rho_{J}^{\tilde{G}}:D_{tr-orb}( {\cal O})\otimes Mes(M({\mathbb R}))^*\to U_{J}\otimes (D_{g\acute{e}om}({\cal O})\otimes Mes(M({\mathbb R}))^*)/Ann_{{\cal O}}^{\tilde{G}},$$
  et
  $$\sigma_{J}^{\tilde{G}}:D_{tr-orb}^{st}( {\cal O})\otimes Mes(M({\mathbb R}))^*\to U_{J}\otimes (D_{g\acute{e}om}^{st}({\cal O})\otimes Mes(M({\mathbb R}))^*)/Ann_{{\cal O}}^{st,\tilde{G}};$$

 - pour tout $\boldsymbol{\gamma}\in D_{tr-orb}({\cal O})\otimes Mes(M({\mathbb R}))^*$, une application lin\'eaire ${\bf f}\mapsto I_{\tilde{M}}^{\tilde{G}}(\boldsymbol{\gamma},B, {\bf f})$ sur $I(\tilde{G}({\mathbb R}))\otimes Mes(G({\mathbb R}))$;
 
  - pour tout $\boldsymbol{\delta}\in D_{tr-orb}^{st}({\cal O})\otimes Mes(M({\mathbb R}))^*$, une application lin\'eaire ${\bf f}\mapsto S_{\tilde{M}}^{\tilde{G}}(\boldsymbol{\delta},B,{\bf f})$ sur $SI(\tilde{G}({\mathbb R}))\otimes Mes(G({\mathbb R}))$.
  
  Dans le cas (C), on note par un indice $unip$ les objets relatifs \`a la classe $\{1\}\subset M({\mathbb R})$. On pose par exemple  $D_{tr-unip}(M({\mathbb R}))=D_{tr-orb}(\{1\})$. On se propose de d\'efinir
  
  - pour tout $J\in {\cal J}_{M}^G(B)$, une application lin\'eaire
  $$\rho_{J}^G:D_{tr-unip}(M({\mathbb R}))\otimes Mes(M({\mathbb R}))^*\to U_{J}\otimes (D_{unip}(M({\mathbb R}))\otimes Mes(M({\mathbb R}))^*)/Ann_{unip}^G;$$
  
  - pour tout $\boldsymbol{\gamma}\in D_{tr-unip}(M({\mathbb R}))\otimes Mes(M({\mathbb R}))^*$, une application lin\'eaire ${\bf f}\mapsto I_{M}^G(\boldsymbol{\gamma},B,{\bf f})$ sur $I(G({\mathbb R}))\otimes Mes(G({\mathbb R}))$.
  
  Si de plus $G$ est quasi-d\'eploy\'e, on se propose de d\'efinir
  
   - pour tout $J\in {\cal J}_{M}^G(B)$, une application lin\'eaire
  $$\sigma_{J}^G:D_{tr-unip}^{st}(M({\mathbb R}))\otimes Mes(M({\mathbb R}))^*\to U_{J}\otimes (D_{unip}^{st}(M({\mathbb R}))\otimes Mes(M({\mathbb R}))^*)/Ann_{unip}^{st,G};$$
  
  - pour tout $\boldsymbol{\delta}\in D_{tr-unip}^{st}(M({\mathbb R}))\otimes Mes(M({\mathbb R}))^*$, une application lin\'eaire ${\bf f}\mapsto S_{M}^G(\boldsymbol{\delta},B,{\bf f})$ sur $SI(G({\mathbb R}))\otimes Mes(G({\mathbb R}))$.  
  
  Dans le cas (A), consid\'erons une donn\'ee endoscopique ${\bf M}'=(M',{\cal M}',\tilde{\zeta})$ de $(KM,K\tilde{M},{\bf a})$, elliptique et relevante. Consid\'erons une classe de conjugaison stable semi-simple ${\cal O}'$ dans $\tilde{M}'({\mathbb R})$ qui correspond \`a ${\cal O}$. Soient $\boldsymbol{\delta}\in D_{tr-orb}^{st}({\cal O}')\otimes Mes(M'({\mathbb R}))^*$ et $a\in A_{K\tilde{M}}({\mathbb R})$ en position g\'en\'erale et proche de $1$. Nos hypoth\`eses de r\'ecurrence et quelques formalit\'es que nous passerons nous autorisent \`a d\'efinir
  
  - pour $J\in {\cal J}_{K\tilde{M}}^{K\tilde{G}}$, le terme
  $$\rho_{J}^{K\tilde{G},{\cal E}}({\bf M}',\boldsymbol{\delta},a)=\sum_{\tilde{s}\in \tilde{\zeta}Z(\hat{M})^{\Gamma_{{\mathbb R}},\hat{\theta}}/Z(\hat{G})^{\Gamma_{{\mathbb R}},\hat{\theta}}}i_{\tilde{M}'}(\tilde{G},\tilde{G}'(\tilde{s}))\sum_{J'\in {\cal J}_{\tilde{M}'}^{\tilde{G}'(\tilde{s})}(B_{{\cal O}'}^{\tilde{G}});J'\mapsto J}transfert(\sigma_{J'}^{{\bf G}'(\tilde{s})}(\boldsymbol{\delta},\xi(a))),$$
  cf. [II] 3.8; c'est la valeur en $a$ d'un \'el\'ement de $U_{J}\otimes (D_{g\acute{e}om}({\cal O},\omega)\otimes Mes(M({\mathbb R}))^*)/Ann_{{\cal O}}^{K\tilde{G}}$.
  
  -  pour ${\bf f}\in I(K\tilde{G}({\mathbb R}),\omega)\otimes Mes(G({\mathbb R}))$, l'int\'egrale endoscopique
  $$I_{K\tilde{M}}^{K\tilde{G},{\cal E}}({\bf M}',\boldsymbol{\delta},{\bf f})=\sum_{\tilde{s}\in \tilde{\zeta}Z(\hat{M})^{\Gamma_{{\mathbb R}},\hat{\theta}}/Z(\hat{G})^{\Gamma_{{\mathbb R}},\hat{\theta}}}i_{\tilde{M}'}(\tilde{G},\tilde{G}'(\tilde{s}))S_{{\bf M}'}^{{\bf G}'(\tilde{s})}(\boldsymbol{\delta},B^{\tilde{G}},{\bf f}^{{\bf G}'(\tilde{s})}).$$
  
  Il y a des variantes de ces d\'efinitions dans les cas (B) et (C). Les notations doivent \^etre adapt\'ees de fa\c{c}on \'evidente. L'unique diff\'erence est que, dans le cas (B) et dans le cas (C) avec $G$ quasi-d\'eploy\'e, les hypoth\`eses de r\'ecurrence ne permettent de d\'efinir ces termes que si $M'\not=M$. 
  
  Venons-en aux conditions impos\'ees aux applications que l'on se propose de d\'efinir. Consid\'erons le cas (A). Pour $\boldsymbol{\gamma}\in D_{orb}({\cal O},\omega)\otimes Mes(M({\mathbb R}))^*$ on a d\'ej\`a d\'efini en 1.3 une application lin\'eaire ${\bf f}\mapsto I_{K\tilde{M}}^{K\tilde{G}}(\boldsymbol{\gamma},{\bf f})$. D'autre part, les d\'efinitions du cas non-archim\'edien de [II] 3.4 s'appliquent et fournissent pour tout $J\in {\cal J}_{K\tilde{M}}^{K\tilde{G}}$ un \'el\'ement 
  $$\rho_{J}^{K\tilde{G}}(\boldsymbol{\gamma})\in U_{J}\otimes (D_{orb}({\cal O},\omega)\otimes Mes(M({\mathbb R}))^*)/Ann_{{\cal O}}^{K\tilde{G}}$$
  On impose
  
  (1) les  d\'efinitions co\"{\i}ncident pour $\boldsymbol{\gamma}\in D_{orb}({\cal O},\omega)\otimes Mes(M({\mathbb R}))^*$.
  
  Pour ${\bf M}'$, ${\cal O}'$, $\boldsymbol{\delta}$ comme ci-dessus, on impose
  
  (2) on a l'\'egalit\'e
  $$I_{K\tilde{M}}^{K\tilde{G}}(transfert(\boldsymbol{\delta}),{\bf f})=I_{K\tilde{M}}^{K\tilde{G},{\cal E}}({\bf M}',\boldsymbol{\delta},{\bf f})$$
  pour tout ${\bf f}\in I(K\tilde{G}({\mathbb R}),\omega)\otimes Mes(G({\mathbb R}))$;
  
  (3) on a l'\'egalit\'e
  $$\rho_{J}^{K\tilde{G}}(transfert(\boldsymbol{\delta}),a)=\rho_{J}^{K\tilde{G},{\cal E}}({\bf M}',\boldsymbol{\delta},a)$$
  pour tout $J\in {\cal J}_{K\tilde{M}}^{K\tilde{G}}$ et tout $a\in A_{K\tilde{M}}({\mathbb R})$ en position g\'en\'erale et proche de $1$. 
  
  On impose
  
  (4) pour tout $\boldsymbol{\gamma}\in D_{tr-orb}({\cal O},\omega)\otimes Mes(M({\mathbb R}))^*$,  tout $J\in {\cal J}_{K\tilde{M}}^{K\tilde{G}}$ et tout $a\in A_{K\tilde{M}}({\mathbb R})$ en position g\'en\'erale et proche de $1$, la distribution induite $\rho_{J}^{K\tilde{G}}(\boldsymbol{\gamma},a)^{K\tilde{G}}$ appartient \`a $D_{tr-orb}({\cal O}^{K\tilde{G}},\omega)\otimes Mes(G({\mathbb R}))^*$;
  
  (5) pour tout $\boldsymbol{\gamma}\in D_{tr-orb}({\cal O},\omega)\otimes Mes(M({\mathbb R}))^*$ et tout ${\bf f}\in I(K\tilde{G}({\mathbb R}),\omega)\otimes Mes(G({\mathbb R}))$, le germe en $1$ de la fonction $a\mapsto I_{K\tilde{M}}^{K\tilde{G}}(a\boldsymbol{\gamma},{\bf f})$, qui est d\'efinie pour $a\in A_{K\tilde{M}}({\mathbb R})$ en position g\'en\'erale et proche de $1$, est \'equivalent \`a
  $$\sum_{K\tilde{L}\in {\cal L}(K\tilde{M})} \sum_{J\in {\cal J}_{K\tilde{M}}^{K\tilde{L}}}I_{K\tilde{L}}^{K\tilde{G}}(\rho_{J}^{K\tilde{L}}(\boldsymbol{\gamma},a)^{K\tilde{L}},{\bf f}).$$
  
  Notons que $a\boldsymbol{\gamma}$ est $\tilde{G}$-\'equisingulier, donc $I_{K\tilde{M}}^{K\tilde{G}}(a\boldsymbol{\gamma},{\bf f})$ est d\'efini d'apr\`es 1.3. D'autre part, d'apr\`es (4), les termes $I_{K\tilde{L}}^{K\tilde{G}}(\rho_{J}^{K\tilde{L}}(\boldsymbol{\gamma},a)^{K\tilde{L}},{\bf f})$ sont d\'efinis (ou plus exactement, le seront quand notre programme sera rempli). 
  
  Nos termes doivent v\'erifier les propri\'etes habituelles de compatibilit\'e \`a l'induction. Soit $K\tilde{R}\in {\cal L}(K\tilde{M}_{0})$ tel que $K\tilde{R}\subset K\tilde{M}$. Soit ${\cal O}_{K\tilde{R}}$ une classe de conjugaison stable semi-simple dans $K\tilde{R}({\mathbb R})$, supposons que ${\cal O}$ soit la classe de conjugaison stable dans $K\tilde{M}({\mathbb R})$ qui la contient. Soit $\boldsymbol{\gamma}\in D_{tr-orb}({\cal O}_{K\tilde{R}},\omega)\otimes Mes(R({\mathbb R}))^*$. On impose
  
  (6) pour tout ${\bf f}\in I(K\tilde{G}({\mathbb R}),\omega)\otimes Mes(G({\mathbb R}))$, on a l'\'egalit\'e
  $$I_{K\tilde{M}}^{K\tilde{G}}(\boldsymbol{\gamma}^{K\tilde{M}},{\bf f})=\sum_{K\tilde{L}\in {\cal L}(K\tilde{R})}d_{\tilde{R}}^{\tilde{G}}(\tilde{M},\tilde{L})I_{K\tilde{R}}^{K\tilde{L}}(\boldsymbol{\gamma},{\bf f}_{K\tilde{L},\omega});$$
  
  (7) pour tout $J\in {\cal J}_{K\tilde{M}}^{K\tilde{G}}$ et tout $a\in A_{K\tilde{M}}({\mathbb R})$ en position g\'en\'erale et proche de $1$, on a l'\'egalit\'e
  $$\rho_{J}^{K\tilde{G}}(\boldsymbol{\gamma}^{K\tilde{M}},a)=\sum_{K\tilde{L}\in {\cal L}(K\tilde{R}), J\in {\cal J}_{K\tilde{R}}^{K\tilde{L}}}d_{\tilde{R}}^{\tilde{G}}(\tilde{M},\tilde{L})\rho_{J}^{K\tilde{L}}(\boldsymbol{\gamma},a)^{K\tilde{M}},$$
  cf. [II] 3.10.
  
  Les propri\'et\'es (4) \`a (7) ont  des analogues dans les cas (B) et (C), qui n'en diff\`erent que par la notation.  Dans le cas (B) et dans le cas (C) avec $G$ quasi-d\'eploy\'e, on a des propri\'et\'es similaires pour les termes stables. Ecrivons-les dans le cas (B). On impose
  
   (8) pour tout $\boldsymbol{\delta}\in D_{tr-orb}^{st}({\cal O})\otimes Mes(M({\mathbb R}))^*$,  tout $J\in {\cal J}_{\tilde{M}}^{\tilde{G}}(B_{{\cal O}})$ et tout $a\in A_{M}({\mathbb R})$ en position g\'en\'erale et proche de $1$, la distribution induite $\sigma_{J}^{\tilde{G}}(\boldsymbol{\delta},a)^{\tilde{G}}$ appartient \`a $D_{tr-orb}^{st}({\cal O}^{\tilde{G}})\otimes Mes(G({\mathbb R}))^*$;
  
  (9) pour tout $\boldsymbol{\delta}\in D_{tr-orb}^{st}({\cal O})\otimes Mes(M({\mathbb R}))^*$ et tout ${\bf f}\in SI(\tilde{G}({\mathbb R}),\omega)\otimes Mes(G({\mathbb R}))$, le germe en $1$ de la fonction $a\mapsto S_{\tilde{M}}^{\tilde{G}}(a\boldsymbol{\delta},{\bf f})$, qui est d\'efinie pour $a\in A_{M}({\mathbb R})$ en position g\'en\'erale et proche de $1$, est \'equivalent \`a
  $$\sum_{\tilde{L}\in {\cal L}(\tilde{M})} \sum_{J\in {\cal J}_{\tilde{M}}^{\tilde{L}}(B_{{\cal O}})}S_{\tilde{L}}^{\tilde{G}}(\sigma_{J}^{\tilde{L}}(\boldsymbol{\delta},a)^{\tilde{L}},B,{\bf f}).$$
  
  Soit $\tilde{R} $ un espace de Levi tel que $\tilde{R}\subset \tilde{M}$. Soit ${\cal O}_{\tilde{R}}$ une classe de conjugaison stable semi-simple dans $\tilde{R}({\mathbb R})$, supposons que ${\cal O}$ soit la classe de conjugaison stable dans $\tilde{M}({\mathbb R})$ qui la contient. Soit $\boldsymbol{\delta}\in D_{tr-orb}^{st}({\cal O}_{\tilde{R}})\otimes Mes(R({\mathbb R}))^*$. On impose
  
  (10) pour tout ${\bf f}\in S(\tilde{G}({\mathbb R}),\omega)\otimes Mes(G({\mathbb R}))$, on a l'\'egalit\'e
  $$S_{\tilde{M}}^{\tilde{G}}(\boldsymbol{\delta}^{\tilde{M}},B,{\bf f})=\sum_{\tilde{L}\in {\cal L}(\tilde{R})}e_{\tilde{R}}^{\tilde{G}}(\tilde{M},\tilde{L})S_{\tilde{R}}^{\tilde{L}}(\boldsymbol{\delta},B,{\bf f}_{\tilde{L}});$$
  
  (11) pour tout $J\in {\cal J}_{\tilde{M}}^{\tilde{G}}(B_{{\cal O}})$ et tout $a\in A_{M}({\mathbb R})$ en position g\'en\'erale et proche de $1$, on a l'\'egalit\'e
  $$\sigma_{J}^{\tilde{G}}(\boldsymbol{\delta}^{\tilde{M}},a)=\sum_{\tilde{L}\in {\cal L}(\tilde{R}), J\in {\cal J}_{\tilde{R}}^{\tilde{L}}(B_{{\cal O}_{\tilde{R}}})}e_{\tilde{R}}^{\tilde{G}}(\tilde{M},\tilde{L})\sigma_{J}^{\tilde{L}}(\boldsymbol{\delta},a)^{\tilde{M}}.$$
  
  On veut aussi que les applications stables  soient d\'eduites des non-stables par les formules habituelles. C'est-\`a-dire que, soit $\boldsymbol{\delta}\in  D_{tr-orb}^{st}({\cal O})\otimes Mes(M({\mathbb R}))^*$. On impose
  
  (12) pour  tout ${\bf f}\in I(\tilde{G}({\mathbb R}),\omega)\otimes Mes(G({\mathbb R}))$, on a l'\'egalit\'e
  $$S_{\tilde{M}}^{\tilde{G}}(\boldsymbol{\delta},B,{\bf f})=I_{\tilde{M}}^{\tilde{G}}(\boldsymbol{\delta},B,{\bf f})-\sum_{s\in Z(\hat{M})^{\Gamma_{{\mathbb R}}}/Z(\hat{G})^{\Gamma_{{\mathbb R}}}, s\not=1}i_{\tilde{M}}(\tilde{G},\tilde{G}'(s))S_{{\bf M}}^{{\bf G}'(s)}(\boldsymbol{\delta},B,{\bf f}^{{\bf G}'(s)});$$
  
  (13)  pour tout $J\in {\cal J}_{\tilde{M}}^{\tilde{G}}(B_{{\cal O}})$ et tout $a\in A_{M}({\mathbb R})$ en position g\'en\'erale et proche de $1$, on a l'\'egalit\'e
  $$\sigma_{J}^{\tilde{G}}(\boldsymbol{\delta},a)=\rho_{J}^{\tilde{G}}(\boldsymbol{\delta},a)-\sum_{s\in Z(\hat{M})^{\Gamma_{{\mathbb R}}}/Z(\hat{G})^{\Gamma_{{\mathbb R}}}, s\not=1,J\in {\cal J}_{\tilde{M}}^{\tilde{G}'(s)}(B_{{\cal O}})}i_{\tilde{M}}(\tilde{G},\tilde{G}'(s))\sigma_{J}^{\tilde{G}'(s)}(\boldsymbol{\delta},a).$$
  
  Ce programme sera r\'ealis\'e inconditionnellement dans les cas (B) et (C) dans les sections 3 et 4. Dans le cas (A), il sera r\'ealis\'e   dans la section 5 sous une hypoth\`ese qui sera expliqu\'ee dans le paragraphe suivant.
  
  \bigskip
  
  \subsection{R\'eduction des conditions impos\'ees dans le cas (A)}
  
 Consid\'erons le cas (A) du paragraphe pr\'ec\'edent. Dans le cas o\`u $K\tilde{M}=K\tilde{G}$, on voit que notre probl\`eme admet une solution unique. Pour $\boldsymbol{\gamma}\in D_{tr-orb}({\cal O},\omega)\otimes Mes(G({\mathbb R}))^*$ et ${\bf f}\in I(K\tilde{G}({\mathbb R}),\omega)\otimes Mes(G({\mathbb R}))$, on a  $I_{K\tilde{G}}^{K\tilde{G}}(\boldsymbol{\gamma},{\bf f})=I^{K\tilde{G}}(\boldsymbol{\gamma},{\bf f})$. L'ensemble ${\cal J}_{K\tilde{G}}^{K\tilde{G}}$ est r\'eduit \`a l'\'el\'ement vide. L'application $\rho_{\emptyset}^{K\tilde{G}}$ est l'identit\'e, modulo l'identification $U_{\emptyset}={\mathbb C}$. On suppose d\'esormais $K\tilde{M}\not=K\tilde{G}$. 
 
 On impose l'hypoth\`ese suivante
 
 (Hyp) pour tout $\boldsymbol{\gamma}\in D_{ orb}(K\tilde{M}({\mathbb R}),\omega)\otimes Mes(M({\mathbb R}))^*$ dont le support est form\'e d'\'el\'ements fortement $\tilde{G}$-r\'eguliers et pour tout ${\bf f}\in I(K\tilde{G}({\mathbb R}),\omega)\otimes Mes(G({\mathbb R}))$, on a l'\'egalit\'e
 $I_{K\tilde{M}}^{K\tilde{G},{\cal E}}(\boldsymbol{\gamma},{\bf f})= I_{K\tilde{M}}^{K\tilde{G}}(\boldsymbol{\gamma},{\bf f})$.
 
 D'apr\`es le lemme 1.11, cette hypoth\`ese implique que la m\^eme \'egalit\'e vaut sous l'hypoth\`ese plus faible $\boldsymbol{\gamma}\in D_{ g\acute{e}om,\tilde{G}-\acute{e}qui}(K\tilde{M}({\mathbb R}),\omega)\otimes Mes(M({\mathbb R}))^*$.
 
 Soit $\boldsymbol{\gamma}\in D_{tr-orb}({\cal O},\omega)\otimes Mes(M({\mathbb R}))^*$. Par d\'efinition, il existe
 
- $\boldsymbol{\gamma}_{orb}\in D_{orb}({\cal O},\omega)\otimes Mes(M({\mathbb R}))^*$,

- une famille finie de donn\'ees endoscopiques $({\bf M}'_{i})_{i=1,...,n}$ de $(M,\tilde{M},{\bf a}_{M})$, elliptiques et relevantes;

- pour tout $i=1,...,n$, une classe de conjugaison stable semi-simple ${\cal O}'_{i}$ dans $\tilde{M}'_{i}({\mathbb R})$ correspondant \`a ${\cal O}$ et un \'el\'ement $\boldsymbol{\delta}_{i}\in D^{st}_{tr-orb}({\bf M}'_{i},{\cal O}'_{i})\otimes Mes(M'({\mathbb R}))^*$,

de sorte que
 $$(1) \qquad \boldsymbol{\gamma}=\boldsymbol{\gamma}_{orb}+\sum_{i=1,...,n}transfert(\boldsymbol{\delta}_{i}).$$
 Cette d\'ecomposition n'est toutefois pas uniquement d\'etermin\'ee.
 
 Les conditions (1), (2), (3) de 2.4 imposent les \'egalit\'es
 $$(2) \qquad I_{K\tilde{M}}^{K\tilde{G}}(\boldsymbol{\gamma},{\bf f})=I_{K\tilde{M}}^{K\tilde{G}}(\boldsymbol{\gamma}_{orb},{\bf f})+\sum_{i=1,...,n}I_{K\tilde{M}}^{K\tilde{G},{\cal E}}({\bf M}'_{i},\boldsymbol{\delta}_{i},{\bf f})$$
 pour tout ${\bf f}\in I(K\tilde{G}({\mathbb R}),\omega)\otimes Mes(G({\mathbb R}))$ 
 et
 $$(3) \qquad \rho_{J}^{K\tilde{G}}(\boldsymbol{\gamma},a)=\rho_{J}^{K\tilde{G}}(\boldsymbol{\gamma}_{orb},a)+\sum_{i=1,...,n}\rho_{J}^{K\tilde{G},{\cal E}}({\bf M}'_{i},\boldsymbol{\delta}_{i},a)$$
 pour tout $J\in {\cal J}_{K\tilde{M}}^{K\tilde{G}}$ et tout $a\in A_{K\tilde{M}}({\mathbb R})$ en position g\'en\'erale et proche de $1$. L'assertion d'existence d'applications v\'erifiant les propri\'et\'es (1), (2) et (3) de 2.4 revient \`a dire que les membres de droite de (2) et (3) ci-dessus ne d\'ependent pas de la d\'ecomposition (1). 
 
 L'ensemble ${\cal J}_{K\tilde{M}}^{K\tilde{G}}$ admet un unique \'el\'ement maximal, cf. [II] 3.1. C'est l'\'el\'ement $J_{max}$ tel que, pour $\underline{\alpha}=\{\alpha_{1},...,\alpha_{n}\}\in J_{max}$, le ${\mathbb Z}$-module $R_{J_{max}}$ engendr\'e par les $\alpha_{i}$ contient toute racine de $A_{K\tilde{M}}({\mathbb R})$ (c'est-\`a-dire toute racine de $A_{\tilde{M}_{p}}({\mathbb R})$ dans $\mathfrak{g}_{p}({\mathbb R})$, pour un quelconque $p\in \Pi^M$). Supposons
 
 (4) pour tout $J\in {\cal J}_{K\tilde{M}}^{K\tilde{G}}$, $J\not=J_{max}$, il existe une application lin\'eaire $\rho_{J}^{K\tilde{G}}$ v\'erifiant les propri\'et\'es (1) et (3) de 2.4.

  Comme on vient de le montrer, cette application est uniquement d\'etermin\'ee. On va montrer qu'en admettant cette propri\'et\'e (4), et sous l'hypoth\`ese (Hyp), on peut r\'ealiser le programme de 2.4. Consid\'erons une donn\'ee endoscopique ${\bf M}'=(M',{\cal M}',\tilde{\zeta})$ de $(KM,K\tilde{M},{\bf a})$, elliptique et relevante. Consid\'erons une classe de conjugaison stable semi-simple ${\cal O}'$ dans $\tilde{M}'({\mathbb R})$ qui correspond \`a ${\cal O}$. Soient $\boldsymbol{\delta}\in D_{tr-orb}^{st}({\cal O}')\otimes Mes(M'({\mathbb R}))^*$ et $ {\bf f}\in I(K\tilde{G}({\mathbb R}),\omega)\otimes Mes(G({\mathbb R}))$. Montrons que
  
  (5) le germe en $1$ de la fonction
  $$a\mapsto I_{K\tilde{M}}^{K\tilde{G},{\cal E}}({\bf M}',a\boldsymbol{\delta},{\bf f}),$$
  qui est d\'efinie pour $a\in A_{K\tilde{M}}({\mathbb R})$ en position g\'en\'erale et proche de $1$,
  est \'equivalent \`a
  $$I_{K\tilde{M}}^{K\tilde{G},{\cal E}}({\bf M}',\boldsymbol{\delta},{\bf f})+\sum_{J\in {\cal J}_{K\tilde{M}}^{K\tilde{G}}}I^{K\tilde{G}}(\rho_{J}^{K\tilde{G},{\cal E}}({\bf M}',\boldsymbol{\delta},a)^{K\tilde{G}},{\bf f})$$
  $$+\sum_{K\tilde{L}\in {\cal L}(K\tilde{M}),K\tilde{L}\not=K\tilde{M},K\tilde{G}}\sum_{J\in {\cal J}_{K\tilde{M}}^{K\tilde{L}}}I_{K\tilde{L}}^{K\tilde{G}}(\rho_{J}^{K\tilde{L}}(transfert(\boldsymbol{\delta}),a)^{K\tilde{L}},{\bf f}).$$
  
  Preuve. On reprend la preuve de la proposition [II] 3.9. Elle montre que le germe en $1$ de la fonction $a\mapsto I_{K\tilde{M}}^{K\tilde{G},{\cal E}}({\bf M}',a\boldsymbol{\delta},{\bf f})$ est \'equivalent \`a celui de la fonction qui, \`a $a$, associe
  $$(6) \qquad\sum_{K\tilde{L}\in {\cal L}(\tilde{M})}\sum_{J\in {\cal J}_{K\tilde{M}}^{K\tilde{L}}}\sum_{\tilde{s}\in \tilde{\zeta}Z(\hat{M})^{\Gamma_{{\mathbb R}},\hat{\theta}}/Z(\hat{L})^{\Gamma_{{\mathbb R}},\hat{\theta}}}i_{\tilde{M}'}(\tilde{L},\tilde{L}'(\tilde{s}))$$
  $$\sum_{J'\in {\cal J}_{\tilde{M}'}^{\tilde{L}'(\tilde{s})}(B^{\tilde{G}}_{{\cal O}'}),J'\mapsto J}I_{K\tilde{L}}^{K\tilde{G},{\cal E}}({\bf L}'(\tilde{s}),\sigma_{J'}^{{\bf L}'(\tilde{s})}(\boldsymbol{\delta},\xi(a))^{{\bf L}'(\tilde{s})},{\bf f}).$$
  Consid\'erons un $K\tilde{L}$ tel que $K\tilde{L}\not=K\tilde{M}$, $K\tilde{L}\not=K\tilde{G}$. Alors on conna\^{\i}t par r\'ecurrence les propri\'et\'es de tous les termes. En particulier
  $$I_{K\tilde{L}}^{K\tilde{G},{\cal E}}({\bf L}'(\tilde{s}),\sigma_{J'}^{{\bf L}'(\tilde{s})}(\boldsymbol{\delta},\xi(a))^{{\bf L}'(\tilde{s})},{\bf f})=I_{K\tilde{L}}^{K\tilde{G}}( transfert(\sigma_{J'}^{{\bf L}'(\tilde{s})}(\boldsymbol{\delta},\xi(a))^{{\bf L}'(\tilde{s})}),{\bf f})$$
  $$=I_{K\tilde{L}}^{K\tilde{G}}( transfert(\sigma_{J'}^{{\bf L}'(\tilde{s})}(\boldsymbol{\delta},\xi(a)))^{K\tilde{L}},{\bf f}).$$
 D'apr\`es la d\'efinition de $\rho_{J}^{K\tilde{L},{\cal E}}({\bf M}',\boldsymbol{\delta},a)$, la somme en $\tilde{s}$ et $J'$ devient
 $$I_{K\tilde{L}}^{K\tilde{G}}(\rho_{J}^{K\tilde{L},{\cal E}}({\bf M}',\boldsymbol{\delta},a)^{K\tilde{L}},{\bf f}).$$
 Puisque $\rho_{J}^{K\tilde{L},{\cal E}}({\bf M}',\boldsymbol{\delta},a)=\rho_{J}^{K\tilde{L}}(transfert(\boldsymbol{\delta}),a)$, la sous-somme index\'ee par $K\tilde{L}$ de la formule (6) est \'egale \`a celle de la formule (5). Consid\'erons maintenant l'espace de Levi $K\tilde{G}$. On a trivialement
 $$I_{K\tilde{G}}^{K\tilde{G},{\cal E}}({\bf G}'(\tilde{s}),\sigma_{J'}^{{\bf G}'(\tilde{s})}(\boldsymbol{\delta},\xi(a))^{{\bf G}'(\tilde{s})},{\bf f})=I^{K\tilde{G}}( transfert(\sigma_{J'}^{{\bf G}'(\tilde{s})}(\boldsymbol{\delta},\xi(a))^{{\bf G}'(\tilde{s})}),{\bf f})$$
 $$=I^{K\tilde{G}}( transfert(\sigma_{J'}^{{\bf G}'(\tilde{s})}(\boldsymbol{\delta},\xi(a)))^{K\tilde{G}},{\bf f}).$$ 
 De nouveau, la sous-somme index\'ee par $K\tilde{G}$ dans la formule (6)   co\"{\i}ncide avec la deuxi\`eme somme de la formule (5). Enfin, pour $K\tilde{L}=K\tilde{M}$, la somme en $J$ se r\'eduit \`a l'unique terme $J=\emptyset$, la somme en $\tilde{s}$ se r\'eduit \`a l'unique terme $\tilde{s}=\tilde{\zeta}$. La sous-somme index\'ee par $K\tilde{M}$ dans la formule (6) co\"{\i}ncide avec le premier terme de la formule (5). Cela d\'emontre (5).

 Conservons les donn\'ees ${\bf M}'$ et ${\cal O}'$. 
  Soit $\tilde{R}'$ un espace de Levi de $\tilde{M}'$. Soit ${\cal O}'_{\tilde{R}'}$ une classe de conjugaison stable dans $\tilde{R}'({\mathbb R})$, supposons que ${\cal O}'$ soit la classe de conjugaison stable dans $\tilde{M}'({\mathbb R})$ qui la contient. Supposons $\tilde{R}'$ relevant. On construit comme en [I] 3.4 un $K$-espace de Levi $K\tilde{R}\in {\cal L}(K\tilde{M}_{0})$ tel que $K\tilde{R}\subset K\tilde{M}$, de sorte que $\tilde{R}'$ soit l'espace d'une donn\'ee endoscopique ${\bf R}'$ de $(KR,K\tilde{R},{\bf a})$ qui est elliptique et relevante. Soit $\boldsymbol{\delta}\in D_{tr-orb}^{st}({\bf R}',{\cal O}'_{\tilde{R}'})\otimes Mes(R'({\mathbb R}))^*$. Alors
  
  (7) pour tout ${\bf f}\in I(K\tilde{G}({\mathbb R}),\omega)\otimes Mes(G({\mathbb R}))$, on a l'\'egalit\'e
  $$I_{K\tilde{M}}^{K\tilde{G},{\cal E}}({\bf M}',\boldsymbol{\delta}^{{\bf M}'},{\bf f})=\sum_{K\tilde{L}\in {\cal L}(K\tilde{R})}d_{\tilde{R}}^{\tilde{G}}(\tilde{M},\tilde{L})I_{K\tilde{R}}^{K\tilde{L},{\cal E}}({\bf R}',\boldsymbol{\delta},{\bf f}_{K\tilde{L},\omega});$$
  
  (8) pour tout $J\in {\cal J}_{K\tilde{M}}^{K\tilde{G}}$ et tout $a\in A_{K\tilde{M}}({\mathbb R})$ en position g\'en\'erale et proche de $1$, on a l'\'egalit\'e
  $$\rho_{J}^{K\tilde{G},{\cal E}}({\bf M}',\boldsymbol{\delta}^{{\bf M}'},a)=\sum_{K\tilde{L}\in {\cal L}(K\tilde{R}), J\in {\cal J}_{K\tilde{R}}^{K\tilde{L}}}d_{\tilde{R}}^{\tilde{G}}(\tilde{M},\tilde{L})\rho_{J}^{K\tilde{L},{\cal E}}({\bf R}',\boldsymbol{\delta},a)^{K\tilde{M}}.$$
  
  La preuve de ces propri\'et\'es  est analogue \`a celle de la proposition [II] 1.14(i). 
  
  Consid\'erons maintenant un \'el\'ement $ \boldsymbol{\gamma}\in D_{tr-orb}({\cal O},\omega)\otimes Mes(M({\mathbb R}))^*$. Choisissons une d\'ecomposition (1). Etudions le germe en $1$ de la fonction qui, \`a $a$, associe 
 $$(9) \qquad  I_{K\tilde{M}}^{K\tilde{G}}(a\boldsymbol{\gamma}_{orb},{\bf f})+\sum_{i=1,...,n}I_{K\tilde{M}}^{K\tilde{G},{\cal E}}({\bf M}'_{i},\xi_{i}(a)\boldsymbol{\delta}_{i},{\bf f}).$$
 On applique la propri\'et\'e (5) \`a chaque terme de la somme en $i$. Le premier terme $I_{K\tilde{M}}^{K\tilde{G}}(a\boldsymbol{\gamma}_{orb},{\bf f})$ satisfait la propri\'et\'e 2.4(5), les termes $\rho_{J}^{K\tilde{G}}(\boldsymbol{\gamma}_{orb},a)$ \'etant d\'efinis comme en [II] 3.2. En effet, la preuve de cette r\'ef\'erence s'applique. D'autre part, pour $J\in {\cal J}_{K\tilde{M}}^{K\tilde{G}}$, $J\not=J_{max}$,  l'hypoth\`ese (4) assure la validit\'e de l'\'egalit\'e (3). 
 Notons $\underline{I}_{K\tilde{M}}^{K\tilde{G}}(\boldsymbol{\gamma},{\bf f})$ le membre de droite de (2).   Notons $\underline{\rho}_{J_{max}}^{K\tilde{G}}(\boldsymbol{\gamma},a)$ le membre de droite de (3) pour $J=J_{max}$.   On obtient alors que le germe en $1$ de la fonction qui \`a $a$ associe (9), 
 est \'equivalent \`a celui de la fonction qui, \`a $a$, associe
  $$(10)\qquad \underline{I}_{K\tilde{M}}^{K\tilde{G}}(\boldsymbol{\gamma},{\bf f})+I^{K\tilde{G}}(\underline{\rho}_{J_{max}}(\boldsymbol{\gamma},a)^{K\tilde{G}},{\bf f})+\sum_{J\in {\cal J}_{K\tilde{M}}^{K\tilde{G}},J\not=J_{max}}I^{K\tilde{G}}(\rho_{J}^{K\tilde{G}}(\boldsymbol{\gamma},a)^{K\tilde{G}},{\bf f})$$
  $$+\sum_{K\tilde{L}\in {\cal L}(K\tilde{M}),K\tilde{L}\not=K\tilde{M},K\tilde{G}}\sum_{J\in {\cal J}_{K\tilde{M}}^{K\tilde{L}}}I_{K\tilde{L}}^{K\tilde{G}}(\rho_{J}^{K\tilde{L}}( \boldsymbol{\gamma},a)^{K\tilde{L}},{\bf f}).$$
  L'hypoth\`ese (Hyp) assure que (9) est \'egal \`a $ I_{K\tilde{M}}^{K\tilde{G}}(a\boldsymbol{\gamma},{\bf f})$. Ce terme est ind\'ependant de la d\'ecomposition (1). Dans (10), tous les termes sauf deux sont aussi ind\'ependants de cette d\'ecomposition. On obtient que le germe en $1$ de la somme de ces deux termes restants est ind\'ependant de cette d\'ecomposition, \`a \'equivalence pr\`es. Pr\'ecis\'ement, le germe en $1$ de
  $$ \underline{I}_{K\tilde{M}}^{K\tilde{G}}(\boldsymbol{\gamma},{\bf f})+I^{K\tilde{G}}(\underline{\rho}_{J_{max}}(\boldsymbol{\gamma},a)^{K\tilde{G}},{\bf f})$$
  est bien d\'etermin\'e \`a \'equivalence pr\`es. Comme fonction de $a$, le premier terme est constant, tandis que le second appartient \`a $U_{J_{max}}$. Le lemme 2.3 assure que chaque terme est bien d\'etermin\'e. Que le premier soit bien d\'etermin\'e pour tout ${\bf f}$ signifie  que le membre de droite de (2)  est ind\'ependant de la d\'ecomposition (1). Que le deuxi\`eme terme soit bien d\'etermin\'e pour tout ${\bf f}$ signifie que $\underline{\rho}_{J_{max}}(\boldsymbol{\gamma},a)^{K\tilde{G}}$ est bien d\'etermin\'e, ou encore que $\underline{\rho}_{J_{max}}(\boldsymbol{\gamma},a)$ l'est, modulo $Ann_{{\cal O}}^{K\tilde{G}}$.  C'est-\`a-dire que le membre de droite de (3) est ind\'ependant de la d\'ecomposition (1). Comme on l'a dit, cela assure l'existence de termes v\'erifiant les conditions (1), (2) et (3) de 2.4. Ces termes \'etant maintenant bien d\'efinis, l'expression (10) n'est autre que la somme figurant dans 2.4(5) et le calcul pr\'ec\'edent prouve cette relation. 
  
Pour v\'erifier 2.4(4), on peut supposer soit que $\boldsymbol{\gamma}\in D_{orb}({\cal O},\omega)\otimes Mes(M({\mathbb R}))^*$, soit qu'il existe ${\bf M}'$, ${\cal O}'$ et $\boldsymbol{\delta}$ comme plus haut tels que $\boldsymbol{\gamma}=transfert(\boldsymbol{\delta})$. Dans le premier cas, la propri\'et\'e 2.4(4) est claire: on a m\^eme $\rho_{J}^{K\tilde{G}}(\boldsymbol{\gamma},a)^{K\tilde{G}}\in D_{orb}({\cal O}^{K\tilde{G}},\omega)\otimes Mes(G({\mathbb R}))^*$. Dans le second cas, on applique la d\'efinition 
$$\rho_{J}^{K\tilde{G}}(\boldsymbol{\gamma},a)^{K\tilde{G}}=\rho_{J}^{K\tilde{G},{\cal E}}({\bf M}',\boldsymbol{\gamma},a)^{K\tilde{G}} =\sum_{\tilde{s}\in \tilde{\zeta}Z(\hat{M})^{\Gamma_{{\mathbb R}},\hat{\theta}}/Z(\hat{G})^{\Gamma_{{\mathbb R}},\hat{\theta}}}i_{\tilde{M}'}(\tilde{G},\tilde{G}'(\tilde{s}))$$
$$\sum_{J'\in {\cal J}_{\tilde{M}'}^{\tilde{G}'(\tilde{s})}(B^{\tilde{G}}_{{\cal O}'}),J'\mapsto J}\left(transfert(\sigma_{J'}^{{\bf G}'(\tilde{s})}(\boldsymbol{\delta},\xi(a)))\right)^{K\tilde{G}}.$$
Par commutation du transfert \`a l'induction, on r\'ecrit le dernier terme
$$transfert\left(\sigma_{J'}^{{\bf G}'(\tilde{s})}(\boldsymbol{\delta},\xi(a))^{{\bf G}'(\tilde{s})}\right).$$
Par r\'ecurrence, on sait que le terme entre parenth\`ese appartient \`a $D_{tr-orb}^{st}({\bf G}'(\tilde{s}),{\cal O}^{_{'}\tilde{G}'(\tilde{s})})\otimes Mes(G'(\tilde{s};{\mathbb R}))^*$, avec une notation \'evidente. Son transfert appartient donc \`a $D_{tr-orb}({\cal O}^{K\tilde{G}},\omega)\otimes Mes(G({\mathbb R}))^*$. Donc aussi $\rho_{J}^{K\tilde{G}}(\boldsymbol{\gamma},a)^{K\tilde{G}}$. 
  
Pour d\'emontrer les propri\'et\'es 2.4 (6) et (7), on peut de nouveau supposer soit que $\boldsymbol{\gamma}\in D_{orb}({\cal O}_{K\tilde{R}},\omega)\otimes Mes(R({\mathbb R}))^*$, soit qu'il existe ${\bf R}'$, ${\cal O}'_{\tilde{R}'}$ et $\boldsymbol{\delta}$ comme plus haut tels que $\boldsymbol{\gamma}=transfert(\boldsymbol{\delta})$. Dans le premier cas, les propri\'et\'es se d\'emontrent comme dans le cas non-archim\'edien, cf. [II] 1.7 et [II] 3.10. Dans le second cas, elles r\'esultent de (7) et (8) ci-dessus. Cela r\'ealise enti\`erement notre programme. 

  \bigskip
  
  \subsection{R\'eduction des conditions impos\'ees dans le cas (C)}
  On consid\`ere le cas (C) de 2.4. Comme en [III] 1.1, on peut affaiblir nos hypoth\`eses de r\'ecurrence car il est clair qu'en partant de notre groupe non tordu, toutes nos constructions ne font appara\^{\i}tre que de tels groupes. Nos hypoth\`eses sont donc les suivantes. Si $G$ est quasi-d\'eploy\'e, on suppose toutes les assertions connues pour des groupes quasi-d\'eploy\'es $G'$ tels que $dim(G'_{SC})< dim(G_{SC})$. Si $G$ n'est pas quasi-d\'eploy\'e, on suppose toutes les assertions connues pour des groupes quasi-d\'eploy\'es $G'$ tels que $dim(G'_{SC})\leq dim(G_{SC})$ et toutes les assertions connues pour des groupes non-quasi-d\'eploy\'es $G'$ tels que $dim(G'_{SC})< dim(G_{SC})$. Pour une assertion relative \`a un Levi $M$ de $G$, on suppose toutes les assertions connues pour le m\^eme groupe $G$ et pour tout Levi $L\in {\cal L}(M)$ avec $L\not=M$.
 
  De nouveau, le probl\`eme pos\'e en 2.4 est \`a peu pr\`es tautologique si $M=G$. On suppose $M\not=G$.  
   
 L'ensemble ${\cal J}_{M}^G(B)$ a encore un \'el\'ement maximal $J_{max}$. On suppose

 (1) pour tout $J\in {\cal J}_{M}^G(B)$, $J\not=J_{max}$, il existe une application lin\'eaire $\rho_{J}^{G}$ v\'erifiant les propri\'et\'es (1) et (3) de 2.4.
 
 Si $G$ n'est pas quasi-d\'eploy\'e,   on impose  une hypoth\`ese similaire \`a celle du paragraphe pr\'ec\'edent, que l'on formule diff\'eremment. Soit ${\bf M}'=(M',{\cal M}',\zeta)$ une donn\'ee endoscopique elliptique et relevante de $M$. On suppose
 
 (Hyp) pour tout $\boldsymbol{\delta}\in D_{orb}^{st}(M'({\mathbb R}))\otimes Mes(M'({\mathbb R}))^*$ dont le support est form\'e d'\'el\'ements fortement $G$-r\'eguliers et pour tout ${\bf f}\in I(G({\mathbb R}))\otimes Mes(G({\mathbb R}))$, on a l'\'egalit\'e
 $$I_{M}^{G,{\cal E}}({\bf M}',\boldsymbol{\delta},{\bf f})= I_{M}^{G}(transfert(\boldsymbol{\delta}),{\bf f}).$$
 Alors la  m\^eme d\'emonstration que dans le paragraphe pr\'ec\'edent s'applique et notre programme est r\'ealis\'e. Supposons maintenant $G$ quasi-d\'eploy\'e. On n'a plus besoin de l'hypoth\`ese (Hyp): elle est v\'erifi\'ee d'apr\`es la proposition 1.13. Dans un premier temps, les hypoth\`eses de r\'ecurrence ne permettent de d\'emontrer les propri\'et\'es (5),  (7) et (8) du paragraphe pr\'ec\'edent que pour des donn\'ees ${\bf M}'$ telles que $M'\not=M$. Mais, dans la d\'ecomposition (1) de ce paragraphe, il n'appara\^{\i}t par d\'efinition que de telles donn\'ees. La d\'emonstration s'applique et d\'emontre l'existence d'applications v\'erifiant les propri\'et\'es (1), (2) et (3) de 2.4. On obtient aussi les propri\'et\'es (4)  \`a (7) de ce paragraphe.
 
 Passons aux variantes stables des applications. Soit $\boldsymbol{\delta}\in D_{tr-unip}^{st}(M({\mathbb R}))\otimes Mes(M({\mathbb R}))^*$. Pour ${\bf f}\in I(G({\mathbb R}))\otimes Mes(G({\mathbb R}))$, on d\'efinit $S_{M}^G(\boldsymbol{\delta},B,{\bf f})$ par la formule habituelle 2.4(12). Pour $J\in {\cal J}_{M}^G(B)$ et $a\in A_{M}({\mathbb R})$ en position g\'en\'erale et proche de $1$, on d\'efinit $\sigma_{J}^G(\boldsymbol{\delta},a)$ par la formule 2.4(13). On suppose
 
 (2)  pour tout $J\in {\cal J}_{M}^G(B)$, $J\not=J_{max}$, $\sigma_{J}^G(\boldsymbol{\delta},a)$ appartient \`a 
 $$(D_{unip}^{st}(M({\mathbb R}))\otimes Mes(M({\mathbb R}))^*)/Ann_{unip}^{st,G}.$$
 
 La m\^eme preuve qu'en [II] 3.7 montre que le germe en $1$ de la fonction 
 $$a\mapsto S_{M}^G(a\boldsymbol{\delta},{\bf f})$$
 est \'equivalent \`a
 $$\sum_{J\in {\cal J}_{M}^G(B)}I^G(\sigma_{J}^G(\boldsymbol{\delta},a)^G,{\bf f})+\sum_{L\in {\cal L}(M),L\not=G}\sum_{J\in {\cal J}_{M}^L(B)}S_{L}^G(\sigma_{J}^L(\boldsymbol{\delta},a)^L,B,{\bf f}).$$
 Supposons que ${\bf f}$ soit instable, c'est-\`a-dire que son image dans $SI(G({\mathbb R}))\otimes Mes(G({\mathbb R}))$ soit nulle. D'apr\`es le th\'eor\`eme 1.4,  $S_{M}^G(a\boldsymbol{\delta},{\bf f})=0$. Dans la somme ci-dessus, tous les termes sauf deux sont nuls, soit par hypoth\`ese de r\'ecurrence (pour $L\not=M$, $L\not=G$), soit d'apr\`es l'hypoth\`ese (2) (pour les $J\in {\cal J}_{M}^G(B)$, $J\not=J_{max}$). On obtient que le germe en $1$ de
$$I^G(\sigma_{J_{max}}^G(\boldsymbol{\delta},a)^G,{\bf f})+S_{M}^G(\boldsymbol{\delta},B,{\bf f})$$
est \'equivalent \`a $0$. Comme fonction de $a$, le premier terme appartient \`a $U_{J_{max}}$, le second est constant. Le lemme 2.3 implique que ces deux termes sont nuls. Cela \'etant vrai pour tout ${\bf f}$ instable, cela signifie que la distribution ${\bf f}\mapsto S_{M}^G(\boldsymbol{\delta},B,{\bf f})$ et $\sigma_{J_{max}}^G(\boldsymbol{\delta},a)^G$ sont stables. D'apr\`es le lemme [I] 5.13,  cette derni\`ere propri\'et\'e \'equivaut \`a ce que $\sigma_{J_{max}}^G(\boldsymbol{\delta},a)$ appartienne \`a  $(D_{unip}^{st}(M({\mathbb R}))\otimes Mes(M({\mathbb R}))^*)/Ann_{unip}^{st,G}$. 

La propri\'et\'e 2.4(8) r\'esulte de la d\'efinition 2.4(13) par le m\^eme argument qui nous a permis de prouver 2.4(4) dans le paragraphe pr\'ec\'edent. La propri\'et\'e 2.4(9) r\'esulte du calcul ci-dessus. Enfin, les propri\'et\'es 2.4(10) et 2.4(11) se prouvent comme en [II] 1.14 (ii). Notons que ces trois derni\`eres propri\'et\'es impliquent \`a leur tour les propri\'et\'es (5), (7) et (8) du paragraphe pr\'ec\'edent pour la donn\'ee ${\bf M}'={\bf M}$, cas que l'on avait laiss\'e en suspens. On a ainsi r\'ealis\'e notre programme.

  \bigskip
  
  \subsection{R\'eduction des conditions impos\'ees dans le cas (B)}
  On consid\`ere le cas (B) de 2.4. On n'a pas besoin de l'analogue des hypoth\`eses (Hyp) de 2.5 et 2.6: elle est v\'erifi\'ee d'apr\`es la proposition  1.13. Par contre, il s'av\`ere que la distinction d'un \'el\'ement maximal de ${\cal J}_{\tilde{M}}^{\tilde{G}}(B_{{\cal O}})$ n'est plus pertinente. On impose
  
  (1) pour tout $J\in {\cal J}_{\tilde{M}}^{\tilde{G}}(B_{{\cal O}})$,  il existe une application lin\'eaire $\rho_{J}^{\tilde{G}}$ v\'erifiant les propri\'et\'es (1) et (3) de 2.4.
  
  On d\'efinit l'application $\sigma_{J}^{\tilde{G}}$ par la formule 2.4(13). On impose

 (2)  pour tout $J\in {\cal J}_{\tilde{M}}^{\tilde{G}}(B_{{\cal O}})$, tout $\boldsymbol{\delta}\in D_{tr-orb}^{st}({\cal O})\otimes Mes(M({\mathbb R}))^*$ et tout $a\in A_{M}({\mathbb R})$ en position g\'en\'erale et proche de $1$, $\sigma_{J}^{\tilde{G}}(\boldsymbol{\delta},a)$ appartient \`a $(D_{g\acute{e}om}^{st}({\cal O})\otimes Mes(M({\mathbb R}))^*)/Ann_{{\cal O}}^{st,\tilde{G}}$.
 
 A l'aide de ces deux hypoth\`eses, la m\^eme preuve que dans le paragraphe pr\'ec\'edent r\'ealise notre programme. 

\bigskip

\section{Extension des d\'efinitions, cas des groupes non tordus}

\bigskip

\subsection{Rappel des r\'esultats d'Arthur}
Dans cette section, on consid\`ere un groupe non tordu, c'est-\`a-dire un triplet $(G,\tilde{G},{\bf a})$ tel que $\tilde{G}=G$ et ${\bf a}=1$. On fixe une fonction $B$ comme en [II] 1.8.  On affaiblit nos hypoth\`eses de r\'ecurrence comme on l'a expliqu\'e en 2.6. 
Soient $M$ un Levi de $G$ et ${\bf M}'$ une donn\'ee endoscopique de $M$ elliptique et relevante. Le r\'esultat suivant a \'et\'e prouv\'e par Arthur ([A4] th\'eor\`eme 1.1).

(1) Soit $\boldsymbol{\delta}\in D_{ orb}^{st}({\bf M}')\otimes Mes(M'({\mathbb R}))^*$ dont le support est form\'e d'\'el\'ements fortement $G$-r\'eguliers et soit ${\bf f}\in I(G({\mathbb R}))\otimes Mes(G({\mathbb R}))$. Alors on a l'\'egalit\'e
$$I_{M}^G(transfert(\boldsymbol{\delta}),{\bf f})=I_{M}^{G,{\cal E}}({\bf M}',\boldsymbol{\delta},{\bf f}).$$

Autrement dit, l'hypoth\`ese (Hyp) de 2.6 est v\'erifi\'ee. Rappelons qu'une variante du lemme 1.11 montre que cet \'enonc\'e s'\'etend aux \'el\'ements $\boldsymbol{\delta}\in D_{g\acute{e}om,\tilde{G}-\acute{e}qui}^{st}({\bf M}')\otimes Mes(M'({\mathbb R}))^*$.

\bigskip

\subsection{R\'ealisation du programme de 2.4}
Soient $M$ un Levi de $G$ et $J\in {\cal J}_{M}^G(B)$. En [II] 3.1, on a associ\'e \`a $J$ un groupe $G_{J}$ contenant $M$ comme Levi. On a $dim(G_{J,SC})\leq dim(G_{SC})$ et cette in\'egalit\'e est stricte si $J$ n'est pas l'\'el\'ement maximal. Soit $\boldsymbol{\gamma}\in D_{tr-unip}(M({\mathbb R}))\otimes Mes(M({\mathbb R}))^*$. On l'\'ecrit
$$(1) \qquad \boldsymbol{\gamma}=\boldsymbol{\gamma}_{orb}+\sum_{i=1,...,n}transfert(\boldsymbol{\delta}_{i})$$
comme en 2.5(1). Les donn\'ees endoscopiques ${\bf M}'_{i}$ qui apparaissent ici v\'erifient $M'_{i}\not=M$ si $G$ est quasi-d\'eploy\'e. Pour $J\in {\cal J}_{M}^G(B)$ et $a\in A_{M}({\mathbb R})$ en position g\'en\'erale et proche de $1$, posons
$$\underline{\rho}_{J}^G(\boldsymbol{\gamma},a)=\rho_{J}^G(\boldsymbol{\gamma}_{orb},a)+\sum_{i=1,...,n}\rho_{J}^{G,{\cal E}}({\bf M}'_{i},\boldsymbol{\delta}_{i},a).$$
Comme dans le cas non-archim\'edien (cf. [II] 3.2), $\rho_{J}^G(\boldsymbol{\gamma}_{orb},a)$ est l'image de $\rho_{J}^{G_{J}}(\boldsymbol{\gamma}_{orb},a)$ par la projection
$$(2) \qquad( D_{unip}(M({\mathbb R}))\otimes Mes(M({\mathbb R}))^*)/Ann_{unip}^{G_{J}}\to ( D_{unip}(M({\mathbb R}))\otimes Mes(M({\mathbb R}))^*)/Ann_{unip}^G.$$
La m\^eme preuve qu'en [III] proposition 1.4(i) montre que, pour tout $i=1,...,n$, $\rho_{J}^{G,{\cal E}}({\bf M}'_{i},\boldsymbol{\delta}_{i},a)$ est l'image de $\rho_{J}^{G_{J},{\cal E}}({\bf M}'_{i},\boldsymbol{\delta}_{i},a)$ par la projection (2). Il en r\'esulte que $\underline{\rho}_{J}^G(\boldsymbol{\gamma},a)$ est l'image par cette projection de l'\'el\'ement
$$\rho_{J}^{G_{J}}(\boldsymbol{\gamma}_{orb},a)+\sum_{i=1,...,n}\rho_{J}^{G_{J},{\cal E}}({\bf M}'_{i},\boldsymbol{\delta}_{i},a).$$
Supposons $J$ non maximal. Alors $dim(G_{J,SC})<dim(G_{SC})$ et, par r\'ecurrence, le terme ci-dessus est ind\'ependant de la d\'ecomposition (1): il vaut $\rho_{J}^{G_{J}}(\boldsymbol{\gamma},a)$.  Il en r\'esulte que $\underline{\rho}_{J}^G(\boldsymbol{\gamma},a)$ est ind\'ependant de cette d\'ecomposition. Cela prouve l'assertion (1) de 2.6.

Supposons $G$ quasi-d\'eploy\'e. Soient $\boldsymbol{\delta}\in D_{tr-unip}^{st}(M({\mathbb R}))\otimes Mes(M({\mathbb R}))^*$, $J\in {\cal J}_{M}^G(B)$  et $a\in A_{M}({\mathbb R})$ en position g\'en\'erale et proche de $1$. Comme on l'a dit en 2.6, on d\'efinit $\sigma_{J}^G(\boldsymbol{\delta},a)$ par l'\'egalit\'e 2.4(13). La m\^eme preuve qu'en [III] proposition 1.2(i) montre que

(3)  $\sigma_{J}^G(\boldsymbol{\delta},a)$ est l'image de $i_{J}^{G_{J}}\sigma_{J}^{G_{J}}(\boldsymbol{\delta},a)$ par la projection (2).

 Si $J$ n'est pas maximal, on sait par r\'ecurrence que $\sigma_{J}^{G_{J}}(\boldsymbol{\delta},a)$  est stable. Il en r\'esulte que $\sigma_{J}^{G_{J}}(\boldsymbol{\delta},a)$  l'est aussi. Cela prouve l'assertion (2) de 2.6. 

On a ainsi v\'erifi\'e les trois assertions dont on a vu en 2.6 qu'elles suffisaient \`a r\'ealiser le programme de 2.4.

 \bigskip
 
 \subsection{Passage \`a un rev\^etement}
 Consid\'erons deux  groupes $G$ et $G_{\sharp}$ r\'eductifs et connexes et un sous-tore  $Z\subset Z(G)$. On pose $G_{\flat}=Z\times G_{\sharp}$ et on suppose donn\'ee une suite exacte
 $$(1) \qquad 1\to \Xi_{b}\to G_{\flat}\stackrel{q}{\to} G\to 1$$
 o\`u $\Xi_{b}$ est un sous-groupe fini central de $G_{\flat}$ et o\`u $q$ se restreint \`a $Z$ en le plongement naturel. Consid\'erons un voisinage ouvert $V_{\sharp}$ de $1$ dans $G_{\sharp}({\mathbb R})$, invariant par l'action naturelle de $G_{ad}({\mathbb R})\simeq G_{\sharp,ad}({\mathbb R})$ et assez petit pour que $q$ se restreigne en un isomorphisme de $Z({\mathbb R})\times V_{\sharp}$ sur $V=q(Z({\mathbb R})\times V_{\sharp})$. On a d\'efini en [III] 3.1 des applications lin\'eaires $\iota^*_{G_{\sharp},G}:D_{g\acute{e}om}(V_{\sharp})\to D_{g\acute{e}om}(V)$ et $\iota^*_{G,G_{\sharp}}:D_{g\acute{e}om}(V)\to D_{g\acute{e}om}(V_{\sharp})$. Leur d\'efinition conserve un sens sur le corps de base ${\mathbb R}$. Leur description est un peu plus compliqu\'ee que dans le cas non-archim\'edien, car les distributions sur $Z({\mathbb R})$ sont un peu plus compliqu\'ees. Notons $D_{g\acute{e}om}(V)_{\sharp}$ le sous-espace des \'el\'ements de $D_{g\acute{e}om}(V)$ qui s'annulent sur toute fonction $f\in C_{c}^{\infty}(V)$ dont la restriction \`a $q(V_{\sharp})$ est nulle. Notons $D_{g\acute{e}om}(V_{\sharp})^{G({\mathbb R})}$ le sous-espace des \'el\'ements de $D_{g\acute{e}om}(V_{\sharp})$ qui sont invariants par l'action de $G({\mathbb R})$. On a une projection naturelle $p:D_{g\acute{e}om}(V_{\sharp})\to D_{g\acute{e}om}(V_{\sharp})^{G({\mathbb R})}$. L'application $\iota^*_{G_{\sharp},G}$ se factorise en
 $$(2) \qquad D_{g\acute{e}om}(V_{\sharp})\stackrel{p}{\to} D_{g\acute{e}om}(V_{\sharp})^{G({\mathbb R})}\stackrel{\iota^*_{G_{\sharp},G}}{\simeq}D_{g\acute{e}om}(V)_{\sharp}\subset D_{g\acute{e}om}(V).$$
 En sens inverse, rappelons que tout \'el\'ement  $H\in Sym(\mathfrak{ z})$ d\'efinit un op\'erateur diff\'erentiel $\partial_{H}$ sur $C_{c}^{\infty}(Z({\mathbb R}))$. Pour $z\in Z({\mathbb R})$ et $H\in Sym(\mathfrak{ z})$, notons $\partial_{z,H}$ l'\'el\'ement de $D_{g\acute{e}om}(Z({\mathbb R}))$ d\'efini par $I^Z(\partial_{z,H},f)=(\partial_{H}f)(z)$. Alors tout \'el\'ement $\boldsymbol{\gamma}\in D_{g\acute{e}om}(Z({\mathbb R}))$ s'\'ecrit de fa\c{c}on unique $\boldsymbol{\gamma}=\sum_{z\in Z({\mathbb R})}\partial_{z,H_{z}}$, avec $H_{z}=0$ pour presque tout $z$. Fixons une base $(H_{i})_{i\in {\mathbb N}}$ de $Sym(\mathfrak{ z})$ form\'ee d'\'el\'ements homog\`enes. On suppose $H_{0}=1$. Alors tout \'el\'ement $\boldsymbol{\gamma}\in D_{g\acute{e}om}(V)$ s'\'ecrit de fa\c{c}on unique sous la forme
 $$\boldsymbol{\gamma}=\sum_{z\in Z({\mathbb R})}\sum_{i\in {\mathbb N}}\partial_{z,H_{i}}\otimes \boldsymbol{\gamma}_{z,i},$$
 o\`u $\boldsymbol{\gamma}_{z,i}$ est un \'el\'ement de $D_{g\acute{e}om}(V)_{\sharp}$. Notons $\boldsymbol{\gamma}_{z,i,\sharp}$ l'unique \'el\'ement de $D_{g\acute{e}om}(V_{\sharp})^{G({\mathbb R})}$ tel que $\iota^*_{G_{\sharp},G}(\boldsymbol{\gamma}_{z,i,\sharp})=\boldsymbol{\gamma}_{z,i}$. On a alors
 $$(3) \qquad \iota^*_{G,G_{\sharp}}(\boldsymbol{\gamma})=\sum_{z\in Z({\mathbb R})}\boldsymbol{\gamma}_{z,0,\sharp}.$$
 L'application $\iota_{G_{\sharp},G}^*$ se restreint en une application
 $$  D_{unip}(G_{\sharp}({\mathbb R}))\to D_{unip}(G({\mathbb R})).$$
 Contrairement au cas non-archim\'edien, elle n'est pas surjective. Par exemple, pour un \'el\'ement non nul  $H\in \mathfrak{z}({\mathbb R})$, la distribution
 $f\mapsto \frac{d}{dt}f(exp(tH))_{\vert t=0}$ sur $C_{c}^{\infty}(V)$ est \`a support unipotent mais n'appartient pas \`a l'image. Il est toutefois clair que l'application se restreint en une application surjective
 $$D_{orb,unip}(G_{\sharp}({\mathbb R}))\to D_{orb,unip}(G({\mathbb R})).$$
 
 Supposons $G$ quasi-d\'eploy\'e. Les constructions s'adaptent aux distributions stables (en supposant $V_{\sharp}$ invariant par conjugaison stable) et nos applications se restreignent en des applications entre les espaces $D_{g\acute{e}om}^{st}(V_{\sharp})$ et $D_{g\acute{e}om}^{st}(V)$. Les groupes $G$ et $G_{\sharp}$ ont m\^eme groupe adjoint. L'action par conjugaison sur $G_{\sharp}$ d'un \'el\'ement de $G({\mathbb R})$ s'identifie \`a celle d'un \'el\'ement de $G_{\sharp,AD}({\mathbb R})$. Or toute distribution stable sur $G_{\sharp}({\mathbb R})$ est invariante par conjugaison par $G_{\sharp,AD}({\mathbb R})$. Il en r\'esulte que la projection $p$ de la relation (2) est l'identit\'e sur $D^{st}_{g\acute{e}om}(V_{\sharp})$. A fortiori
 
 (4) l'application $\iota_{G_{\sharp},G}^*:D^{st}_{unip}(G_{\sharp}({\mathbb R}))\to D^{st}_{unip}(G({\mathbb R}))$ est injective. 
 
 \ass{Lemme}{(i) Les homomorphismes $\iota_{G_{\sharp},G}^*$ et $\iota_{G,G_{\sharp}}^*$ se restreignent en des isomorphismes inverses l'un de l'autre entre $D_{tr-unip}(G({\mathbb R}))$ et $\iota^*_{G,G_{\sharp}}(D_{tr-unip}(G({\mathbb R})))$. On a les inclusions
 $$\iota^*_{G,G_{\sharp}}(D_{tr-unip}(G({\mathbb R})))\subset p(D_{tr-unip}(G_{\sharp}({\mathbb R})))\subset D_{tr-unip}(G_{\sharp}({\mathbb R})).$$

 (ii) Si $G$ est quasi-d\'eploy\'e, on a les relations
  $$\iota^*_{G,G_{\sharp}}(D_{tr-unip}(G({\mathbb R})))=p(D_{tr-unip}(G_{\sharp}({\mathbb R})))\subset D_{tr-unip}(G_{\sharp}({\mathbb R}))$$
  et 
  $$\iota^*_{G,G_{\sharp}}(D^{st}_{tr-unip}(G({\mathbb R})))=p(D^{st}_{tr-unip}(G_{\sharp}({\mathbb R})))= D^{st}_{tr-unip}(G_{\sharp}({\mathbb R})).$$}
 
 Preuve. On fixe des mesures de Haar sur chacun des groupes. On va prouver l'inclusion
 $$(5) \qquad D_{tr-unip}(G({\mathbb R}))\subset \iota^*_{G_{\sharp},G}(D_{tr-unip}(G_{\sharp}({\mathbb R}))).$$
 Il suffit de d\'emontrer qu'un ensemble de g\'en\'erateurs de $D_{tr-unip}(G({\mathbb R}))$ est contenu dans le membre de droite. D'apr\`es 2.2(3), l'espace $D_{tr-unip}(G({\mathbb R}))$ est engendr\'e par $D_{orb,unip}(G({\mathbb R}))$ et les espaces $transfert(D_{tr-unip}^{st}({\bf G}'))$, o\`u ${\bf G}'$ parcourt les donn\'ees endoscopiques elliptiques et relevantes de $G$, avec $G'\not=G$ si $G$ est quasi-d\'eploy\'e. Un \'el\'ement de $D_{orb,unip}(G({\mathbb R}))$ appartient l'espace de droite de (5): tout \'el\'ement unipotent de $G({\mathbb R})$ est l'image naturelle d'un \'el\'ement unipotent de $G_{\sharp}({\mathbb R})$. Soit ${\bf G}'=(G',{\cal G}',s)$ une donn\'ee endoscopique elliptique et relevante de $G$, avec $G'\not=G$ si $G$ est quasi-d\'eploy\'e. Consid\'erons un \'el\'ement $\boldsymbol{\gamma}\in D_{tr-unip}(G({\mathbb R})) $ de la forme $transfert(\boldsymbol{\delta})$, o\`u $\boldsymbol{\delta}\in D_{tr-unip}^{st}({\bf G}')$.  
 Dualement \`a la suite exacte (1), on a une suite exacte
 $$(6) \qquad 1\to \hat{\Xi}_{b}\to \hat{G}\to \hat{Z}\times \hat{G}_{\sharp}\to 1,$$
 o\`u $\hat{\Xi}_{b}$ est un sous-groupe fini central de $\hat{G}$.
 L'\'el\'ement $s\in \hat{G}$ s'envoie sur  un \'el\'ement $(z,s_{\sharp})$ de $\hat{Z}\times \hat{G}_{\sharp}$. En notant $\hat{G}'_{\sharp}$ la composante neutre de $Z_{\hat{G}_{\sharp}}(s_{\sharp})$, on a la suite exacte
 $$1\to \hat{\Xi}_{b}\to \hat{G}'\to \hat{Z}\times \hat{G}'_{\sharp}\to 1.$$
 Le groupe ${\cal G}'/\hat{\Xi}_{b}$ contient $\hat{Z}$ donc est de la forme $\hat{Z}\times {\cal G}'_{\sharp}$, o\`u ${\cal G}'_{\sharp}$ est un sous-groupe de ${^LG}_{\sharp}$. Ce groupe d\'efinit une action galoisienne sur $\hat{G}'_{\sharp}$. 
 On introduit un groupe quasi-d\'eploy\'e $G'_{\sharp}$ sur ${\mathbb R}$ dont le groupe dual soit $\hat{G}'_{\sharp}$. Alors ${\bf G}'_{\sharp}=(G'_{\sharp},{\cal G}'_{\sharp},s_{\sharp})$ est une donn\'ee endoscopique pour $G_{\sharp}$. On v\'erifie ais\'ement qu'elle est elliptique et relevante. Remarquons que les groupes $G'_{\sharp}$ et $G'$ sont dans la m\^eme situation que $G_{\sharp}$ et $G$. C'est-\`a-dire que l'on a une suite exacte
 $$1\to \Xi_{b}\to Z\times G'_{\sharp}\to G'\to 1.$$
 On peut appliquer par r\'ecurrence la relation (5) \`a $G'$ et $G'_{\sharp}$. Modulo quelques formalit\'es, on peut aussi bien l'appliquer \`a ${\bf G}'$ et ${\bf G}'_{\sharp}$. On obtient que $\boldsymbol{\delta}$ appartient \`a l'espace $\iota_{{\bf G}'_{\sharp},{\bf G}'}^*(D_{tr-unip}({\bf G}'_{\sharp}))$. Mais, d'apr\`es (3), l'application $\iota_{{\bf G}'_{\sharp},{\bf G}'}^*\circ\iota_{{\bf G}',{\bf G}'_{\sharp}}^*$ est l'identit\'e sur cet espace. Posons $\boldsymbol{\delta}_{\sharp}=\iota_{{\bf G}',{\bf G}'_{\sharp}}^*(\boldsymbol{\delta})$. Alors $\boldsymbol{\delta}=\iota_{{\bf G}'_{\sharp},{\bf G}'}^*(\boldsymbol{\delta}_{\sharp})$. D'apr\`es le (ii) du lemme appliqu\'e par r\'ecurrence, on a
  $\boldsymbol{\delta}_{\sharp}\in D_{tr-unip}({\bf G}'_{\sharp})$.  On a l'\'egalit\'e
 $$(7) \qquad transfert \circ \iota^*_{{\bf G}'_{\sharp},{\bf G}'}=\iota_{G_{\sharp},G}^*\circ transfert$$
 cf. [III] 3.7(4). On en d\'eduit $\boldsymbol{\gamma}=\iota_{G_{\sharp},G}^*\circ transfert(\boldsymbol{\delta}_{\sharp})$. Mais $transfert(\boldsymbol{\delta}_{\sharp})$ appartient \`a $D_{tr-unip}(G_{\sharp}({\mathbb R}))$ par d\'efinition de cet espace. Cela d\'emontre que $\boldsymbol{\gamma}$ appartient au membre de droite de (5). Cela prouve cette assertion. 
 
 Le membre de droite de (5) est contenu dans $D_{g\acute{e}om}(V)_{\sharp}$ d'apr\`es (2). La relation (3) entra\^{\i}ne que, pour tout sous-espace $X$ de $D_{g\acute{e}om}(V)_{\sharp}$, les homomorphismes $\iota_{G_{\sharp},G}^*$ et $\iota_{G,G_{\sharp}}^*$ se restreignent en des isomorphismes inverses l'un de l'autre entre $X$ et $\iota_{G,G_{\sharp}}^*(X)$.
 Donc (5) entra\^{\i}ne la premi\`ere assertion de l'\'enonc\'e. Soit $\boldsymbol{\gamma}\in D_{tr-unip}(G({\mathbb R}))$, \'ecrivons $\boldsymbol{\gamma}=\iota^*_{G_{\sharp},G}(\boldsymbol{\gamma}_{\sharp})$, avec $\boldsymbol{\gamma}_{\sharp}\in D_{tr-unip}(G_{\sharp}({\mathbb R}))$. Puisque $\boldsymbol{\gamma}\in D_{g\acute{e}om}(V)_{\sharp}$  comme on vient de le voir, la relation (3) entra\^{\i}ne que $\iota^*_{G,G_{\sharp}}(\boldsymbol{\gamma})=p(\boldsymbol{\gamma}_{\sharp})$. D'o\`u la premi\`ere inclusion de (i). Par simple transport de structure, l'espace $D_{tr-unip}(G_{\sharp}({\mathbb R}))$ est stable par tout automorphisme de $G_{\sharp}$. En particulier, il est stable par l'action de $G({\mathbb R})$. On en d\'eduit $p(D_{tr-unip}(G_{\sharp}({\mathbb R})))\subset D_{tr-unip}(G_{\sharp}({\mathbb R}))$, ce qui ach\`eve la preuve de (i).
 
 Supposons $G$ quasi-d\'eploy\'e. On va prouver l'inclusion
 $$(8) \qquad \iota^*_{G_{\sharp},G}(D_{tr-unip}(G_{\sharp}({\mathbb R})))\subset D_{tr-unip}(G({\mathbb R})).$$
 De nouveau, il suffit de prouver que, pour $\boldsymbol{\gamma}_{\sharp}$ dans un ensemble de g\'en\'erateurs de $D_{tr-unip}(G_{\sharp}({\mathbb R}))$, on a $\iota^*_{G_{\sharp},G}(\boldsymbol{\gamma}_{\sharp})\in D_{tr-unip}(G({\mathbb R}))$. C'est clair si $\boldsymbol{\gamma}_{\sharp}\in D_{orb,unip}(G_{\sharp}({\mathbb R}))$. Soit ${\bf G}'_{\sharp}=(G'_{\sharp},{\cal G}'_{\sharp},s_{\sharp})$ une donn\'ee endoscopique elliptique et relevante pour $G_{\sharp}$, telle que $G'_{\sharp}\not=G$. Soit $\boldsymbol{\gamma}_{\sharp}\in D_{tr-unip}(G_{\sharp}({\mathbb R}))$ de la forme $transfert(\boldsymbol{\delta}_{\sharp})$, avec $\boldsymbol{\delta}_{\sharp}\in D_{tr-unip}({\bf G}'_{\sharp})$.  
 Rappelons qu'\`a la donn\'ee endoscopique $ {\bf G}'_{\sharp}$ est associ\'e un caract\`ere $\omega_{\sharp}$ de  $G_{\sharp,ad}({\mathbb R})$ selon lequel se transforme le facteur de transfert, cf. [I] 2.7. Il se restreint en un caract\`ere de $G({\mathbb R})$. Si cette restriction est non triviale, l'image de $transfert(\boldsymbol{\delta}_{\sharp})$ par la projection $p$ est nulle. Autrement dit $p(\boldsymbol{\gamma}_{\sharp})=0$. Mais $\iota^*_{G_{\sharp},G}\circ p=\iota^*_{G_{\sharp},G}$, donc $\iota^*_{G_{\sharp},G}(\boldsymbol{\gamma}_{\sharp})=0$. A fortiori, $\iota^*_{G_{\sharp},G}(\boldsymbol{\gamma}_{\sharp})\in D_{tr-unip}(G({\mathbb R}))$. On est ramen\'e au cas o\`u la restriction de $\omega_{\sharp}$ \`a $G({\mathbb R})$ est triviale. Montrons  que
 
 (9) pour $G$ quasi-d\'eploy\'e, l'homomorphisme de $H^1(W_{{\mathbb R}};Z(\hat{G}))$ dans le groupe des caract\`eres de $G({\mathbb R})$ est bijective.
 
 Que $G$ soit quasi-d\'eploy\'e ou pas,  $H^1(W_{{\mathbb R}};Z(\hat{G}))$ s'identifie au groupe des caract\`eres de $G_{ab}({\mathbb R})$. En effet, en fixant un sous-tore maximal $T$ de $G$ d\'efini sur ${\mathbb R}$, ces deux groupes s'identifient respectivement \`a $H^{1,0}(W_{{\mathbb R}};\hat{T}\to\hat{T}_{ad})$ et $H^{1,0}(\Gamma_{{\mathbb R}};T_{sc}\to T)$ et on conna\^{\i}t l'assertion pour ces derniers groupes. L'assertion (9) \'equivaut \`a dire que, si $G$ est quasi-d\'eploy\'e, l'homomorphisme $G({\mathbb R})\to G_{ab}({\mathbb R})$ est surjectif. On a d\'emontr\'e   cela en 2.1(2).

 Fixons $s\in \hat{G}$ qui se projette sur $(1,s_{\sharp})$ par l'homomorphisme de la suite (6). Notons ${\cal G}'$ l'image r\'eciproque de $\hat{Z}\times {\cal G}'_{\sharp}$ dans ${^LG}$.  Il agit sur $\hat{G}_{s}$ et munit ce groupe d'une action galoisienne. On introduit un groupe $G'$ quasi-d\'eploy\'e sur ${\mathbb R}$ dont $\hat{G}_{s}$ est le groupe dual. Posons ${\bf G}'=(G',{\cal G}',s)$. Montrons que c'est une donn\'ee endoscopique pour $G$. La seule condition non \'evidente est la suivante. Pour $w\in W_{{\mathbb R}}$ et $(g,w)\in {\cal G}'$, on a une \'egalit\'e
 $sgw(s)^{-1}=a(w)g$, o\`u $a$ est un cocycle de $W_{{\mathbb R}}$ dans $Z(\hat{G})$. Il faut voir que ce cocycle est un cobord. En reprenant les d\'efinitions de [I] 2.7, on voit que la restriction de $\omega_{\sharp}$ \`a $G({\mathbb R})$ est pr\'ecis\'ement associ\'ee \`a la classe de $a$. L'hypoth\`ese que cette restriction est triviale jointe \`a (9) entra\^{\i}ne que $a$ est un cobord.  On poursuit alors la d\'emonstration comme ci-dessus. En appliquant (8) par r\'ecurrence \`a ${\bf G}'$ et ${\bf G}'_{\sharp}$, on a $\iota^*_{{\bf G}'_{\sharp},{\bf G}'}(\boldsymbol{\delta}_{\sharp})\in D_{tr-unip}({\bf G}')$. Puisque $\iota^*_{{\bf G}'_{\sharp},{\bf G}'}$ pr\'eserve la stabilit\'e, on a m\^eme $\iota^*_{{\bf G}'_{\sharp},{\bf G}'}(\boldsymbol{\delta}_{\sharp})\in D^{st}_{tr-unip}({\bf G}')$.  Alors
 $$\iota_{G_{\sharp},G}^*(\boldsymbol{\gamma}_{\sharp})=\iota^*_{G_{\sharp},G}\circ transfert(\boldsymbol{\delta}_{\sharp})=transfert\circ \iota^*_{{\bf G}'_{\sharp},{\bf G}'}(\boldsymbol{\delta}_{\sharp}).$$
 Ce dernier \'el\'ement appartient \`a $ transfert(D_{tr-unip}^{st}({\bf G}'))$ qui est inclus dans $ D_{tr-unip}(G({\mathbb R}))$. Cela prouve (8). 

 Les premi\`eres assertions de (ii) r\'esultent de (8) comme (i) r\'esultait de (5). 
 
Comme on l'a dit, cf. (4), $p$ est l'identit\'e sur les distributions stables.  A fortiori, $p(D^{st}_{tr-unip}(G_{\sharp}({\mathbb R})))=D^{st}_{tr-unip}(G_{\sharp}({\mathbb R}))$, ce qui est la derni\`ere \'egalit\'e de l'\'enonc\'e. On a aussi
  $$D_{tr-unip}^{st}(G_{\sharp}({\mathbb R}))=p(D_{tr-unip}^{st}(G_{\sharp}({\mathbb R})))\subset p(D_{tr-unip}(G_{\sharp}({\mathbb R})))\cap D_{g\acute{e}om}^{st}(G_{\sharp}({\mathbb R}))$$
 $$\subset D_{tr-unip}(G_{\sharp}({\mathbb R}))\cap D_{g\acute{e}om}^{st}(G_{\sharp}({\mathbb R}))=D_{tr-unip}^{st}(G_{\sharp}({\mathbb R})).$$ 
 Les inclusions de cette suite sont forc\'ement des \'egalit\'es, donc
  $$D_{tr-unip}^{st}(G_{\sharp}({\mathbb R}))= p(D_{tr-unip}(G_{\sharp}({\mathbb R})))\cap D_{g\acute{e}om}^{st}(G_{\sharp}({\mathbb R})).$$
 Puisque $\iota^*_{G_{\sharp},G}$ et $\iota^*_{G,G_{\sharp}}$ sont des isomorphismes inverses entre $D_{tr-unip}(G({\mathbb R}))$ et $p(D_{tr-unip}(G_{\sharp}({\mathbb R})))$ et qu'ils pr\'eservent la stabilit\'e, on a
 $$\iota^*_{G,G_{\sharp}}(D^{st}_{tr-unip}(G({\mathbb R})))=p(D_{tr-unip}(G_{\sharp}({\mathbb R})))\cap D_{g\acute{e}om}^{st}(G_{\sharp}({\mathbb R}))=D_{tr-unip}^{st}(G_{\sharp}({\mathbb R})).$$
  C'est l'avant-derni\`ere \'egalit\'e de (ii) qu'il restait \`a prouver. $\square$

 \bigskip
 
 \subsection{Rev\^etement et applications $\rho_{J}$ et $\sigma_{J}$}
 On conserve la situation pr\'ec\'edente. Soit $M$ un Levi de $G$. On note $M_{\flat}$ son image r\'eciproque. On a $M_{\flat}=Z\times M_{\sharp}$, o\`u $M_{\sharp}$ est un Levi de $G_{\sharp}$. La fonction $B$ sur $G$ se rel\`eve en une fonction encore not\'ee $B$ sur $G_{\sharp}$. Les ensembles ${\cal J}_{M}^G(B)$ et ${\cal J}_{M_{\sharp}}^{G_{\sharp}}(B)$ s'identifient et, pour $J$ dans cet ensemble, on peut identifier les deux espaces de germes $U_{J}$ relatifs aux groupes ambiants $G$ et $G_{\sharp}$, cf. [III] 3.2. On fixe une mesure de Haar sur $Z({\mathbb R})$, qui permet d'identifier $Mes(G({\mathbb R}))$ \`a $Mes(G_{\sharp}({\mathbb R}))$ et $Mes(M({\mathbb R}))$ \`a $Mes(M_{\sharp}({\mathbb R}))$, cf. [III] 3.1. 
 
 \ass{Lemme}{(i) Soient $J\in {\cal J}_{M}^G(B)$ et $\boldsymbol{\gamma}_{\sharp}\in D_{tr-unip}(M_{\sharp}({\mathbb R}))\otimes Mes(M_{\sharp}({\mathbb R}))^*$. Supposons $\iota^*_{G_{\sharp},G}(\boldsymbol{\gamma}_{\sharp})\in D_{tr-unip}(M({\mathbb R}))\otimes Mes(M({\mathbb R}))^*$. Alors on a l'\'egalit\'e
 $$\rho_{J}^G\circ \iota^*_{M_{\sharp},M}(\boldsymbol{\gamma}_{\sharp})=\iota^*_{M_{\sharp},M}\circ \rho_{J}^{G_{\sharp}}(\boldsymbol{\gamma}_{\sharp}).$$
 
 (ii) Supposons $G $ quasi-d\'eploy\'e. Soient $J\in {\cal J}_{M}^G(B)$ et $\boldsymbol{\delta}_{\sharp}\in D_{tr-unip}^{st}(M_{\sharp}({\mathbb R}))\otimes Mes(M_{\sharp}({\mathbb R}))^*$. Alors on a l'\'egalit\'e
  $$\sigma_{J}^G\circ \iota^*_{M_{\sharp},M}(\boldsymbol{\delta}_{\sharp})=\iota^*_{M_{\sharp},M}\circ \sigma_{J}^{G_{\sharp}}(\boldsymbol{\delta}_{\sharp}).$$}
  
  {\bf Remarque.} Les $\iota^*_{M_{\sharp},M}$ intervenant dans les membres de droite sont en fait tensoris\'es avec l'identit\'e de $U_{J}$.
  
  Preuve. Pour simplifier, on abandonne les espaces de mesures. L'action de $M({\mathbb R})$ sur $M_{\sharp}({\mathbb R})$  se prolonge en une action $G({\mathbb R})$ sur $G_{\sharp}({\mathbb R})$. Par simple transport de structure, on voit alors que
  $$p\circ \rho_{J}^{G_{\sharp}}(\boldsymbol{\gamma}_{\sharp})=\rho_{J}^{G_{\sharp}}\circ p(\boldsymbol{\gamma}_{\sharp}).$$
  Parce que $\iota^*_{M_{\sharp},M}\circ p=\iota^*_{M_{\sharp},M}$, le membre de gauche de l'\'egalit\'e du (i) ne d\'epend que de $p(\boldsymbol{\gamma}_{\sharp})$. Pour la m\^eme raison et gr\^ace \`a l'\'egalit\'e ci-dessus, on voit que c'est aussi le cas du membre de droite. On peut alors reformuler l'assertion (i) de la fa\c{c}on suivante. Soit $\boldsymbol{\gamma}\in  D_{tr-unip}(M({\mathbb R})) $. Alors il existe $\boldsymbol{\gamma}_{\sharp}\in D_{tr-unip}(M_{\sharp}({\mathbb R})) $ tel que  $\iota^*_{G_{\sharp},G}(\boldsymbol{\gamma}_{\sharp})=\boldsymbol{\gamma}$ et que l'\'egalit\'e du (i) soit v\'erifi\'ee. Il suffit de d\'emontrer cette assertion quand $\boldsymbol{\gamma}$ parcourt un ensemble de g\'en\'erateurs de $D_{tr-unip}(M({\mathbb R})) $. Pour $\boldsymbol{\gamma}\in D_{orb,unip}(M({\mathbb R})) $, on choisit $\boldsymbol{\gamma}_{\sharp}\in D_{orb,unip}(M_{\sharp}({\mathbb R})) $ et l'\'egalit\'e r\'esulte des d\'efinitions comme dans le cas non-archim\'edien. 
 Soit maintenant ${\bf M}'=(M',{\cal M}',\zeta)$ une donn\'ee endoscopique elliptique et relevante de $M$, avec $M'\not=M$ si $G$ est quasi-d\'eploy\'e. Consid\'erons un \'el\'ement $\boldsymbol{\gamma}\in D_{orb,unip}(M({\mathbb R})) $ de la forme $transfert(\boldsymbol{\delta})$, o\`u $\boldsymbol{\delta}\in D^{st}_{tr-unip}({\bf M}')$. Comme dans la preuve du lemme pr\'ec\'edent (o\`u l'on remplace le groupe $G$ par $M$), on construit une donn\'ee endoscopique ${\bf M}'_{\sharp}=(M'_{\sharp},{\cal M}'_{\sharp},\zeta_{\sharp})$ de $M_{\sharp}$ qui est elliptique et relevante. On choisit $\boldsymbol{\delta}_{\sharp}\in D^{st}_{tr-unip}({\bf M}'_{\sharp})$ tel que $\iota^*_{{\bf M}'_{\sharp},{\bf M}'}(\boldsymbol{\delta}_{\sharp})=\boldsymbol{\delta}$. On peut choisir $\boldsymbol{\gamma}_{\sharp}=transfert(\boldsymbol{\delta}_{\sharp})$. D'apr\`es 2.4(3), on a
 $$\rho_{J}^G(\boldsymbol{\gamma},a)=\rho_{J}^{G,{\cal E}}({\bf M}',\boldsymbol{\delta},a)\,\,\text{et}\,\, \rho_{J}^{G_{\sharp}}(\boldsymbol{\gamma}_{\sharp},a_{\sharp})=\rho_{J}^{G_{\sharp},{\cal E}}({\bf M}'_{\sharp},\boldsymbol{\delta}_{\sharp},a_{\sharp})$$
 pour $a\in A_{M}({\mathbb R})$ et $a_{\sharp}\in A_{M_{\sharp}}({\mathbb R})$. C'est-\`a-dire
 $$\rho_{J}^G(\boldsymbol{\gamma},a)=\sum_{s\in \zeta Z(\hat{M})^{\Gamma_{{\mathbb R}}}/Z(\hat{G})^{\Gamma_{{\mathbb R}}}}i_{M'}(G,G'(s))\sum_{J'\in {\cal J}_{M'}^{G'(s)}(B); J'\mapsto J}transfert(\sigma_{J'}^{{\bf G}'(s)}(\boldsymbol{\delta},\xi(a))),$$
 $$\rho_{J}^{G_{\sharp}}(\boldsymbol{\gamma}_{\sharp},a_{\sharp})=\sum_{s\in \zeta_{\sharp} Z(\hat{M}_{\sharp})^{\Gamma_{{\mathbb R}}}/Z(\hat{G}_{\sharp})^{\Gamma_{{\mathbb R}}}}i_{M'_{\sharp}}(G_{\sharp},G_{\sharp}'(s))\sum_{J'\in {\cal J}_{M'_{\sharp}}^{G'_{\sharp}(s)}(B); J'\mapsto J}transfert(\sigma_{J'}^{{\bf G}'_{\sharp}(s)}(\boldsymbol{\delta}_{\sharp},\xi(a_{\sharp}))).$$
On a
$$Z(\hat{M})^{\Gamma_{{\mathbb R}}}/Z(\hat{G})^{\Gamma_{{\mathbb R}}}=Z(\hat{M}_{ad})^{\Gamma_{{\mathbb R}}}= Z(\hat{M}_{\sharp,ad})^{\Gamma_{{\mathbb R}}}=Z(\hat{M}_{\sharp})^{\Gamma_{{\mathbb R}}}/Z(\hat{G}_{\sharp})^{\Gamma_{{\mathbb R}}}.$$
Les sommations en $s$ des deux formules ci-dessus s'identifient. Pour $s$ fix\'e, on v\'erifie ais\'ement l'\'egalit\'e
 $$i_{M'}(G,G'(s))=i_{M'_{\sharp}}(G_{\sharp},G_{\sharp}'(s)).$$
 Les sommations en $J'$ d'identifient. Supposons que $a$ et $a_{\sharp}$ se correspondent par la relation $a=q(a_{\sharp})$. Les termes $\xi(a)$ et $\xi(a_{\sharp})$ se correspondent par une relation similaire.  Pour $J'$ fix\'e, on peut appliquer par r\'ecurrence l'assertion (ii) du lemme:
 $$\sigma_{J'}^{{\bf G}'(s)}(\boldsymbol{\delta},\xi(a))=\iota^*_{{\bf M}'_{\sharp},{\bf M}'}(\sigma_{J'}^{{\bf G}'_{\sharp}(s)}(\boldsymbol{\delta}_{\sharp},\xi(a_{\sharp}))).$$
 En utilisant 3.7(6), on obtient
$$ transfert(\sigma_{J'}^{{\bf G}'(s)}(\boldsymbol{\delta},\xi(a)))=\iota^*_{M_{\sharp},M}\circ  transfert(\sigma_{J'}^{{\bf G}'_{\sharp}(s)}(\boldsymbol{\delta}_{\sharp},\xi(a_{\sharp}))).$$
On en d\'eduit 
$$\rho_{J}^G(\boldsymbol{\gamma},a)=\iota^*_{M_{\sharp},M}(\rho_{J}^{G_{\sharp}}(\boldsymbol{\gamma}_{\sharp},a_{\sharp})).$$
Cela prouve (i). Le (ii) s'en d\'eduit comme dans le lemme [III] 3.6. $\square$

\bigskip

\subsection{Un r\'esultat d'induction}
On suppose $G$ quasi-d\'eploy\'e. Soit $M$ un Levi de $G$. 

\ass{Lemme}{Soit $\boldsymbol{\gamma}\in D_{unip}(M_{SC}({\mathbb R}))\otimes Mes(M_{SC}({\mathbb R}))^*$. Supposons que la distribution induite $(\iota^*_{M_{SC},M}(\boldsymbol{\gamma}))^{G}$ soit stable. Alors il existe $\boldsymbol{\delta}\in D_{unip}^{st}(M_{SC}({\mathbb R}))\otimes Mes(M_{SC}({\mathbb R}))^*$ tel que l'on ait l'\'egalit\'e
$$(\iota^*_{M_{SC},M}(\boldsymbol{\gamma}))^{G}=(\iota^*_{M_{SC},M}(\boldsymbol{\delta}))^{G}.$$}

{\bf Remarque.} On rappelle que $M_{SC}$ est le rev\^etement simplement connexe du groupe d\'eriv\'e de $M$. Le lemme serait plus facile mais plus faible si l'on rempla\c{c}ait ce groupe par l'image r\'eciproque $M_{sc}$ de $M$ dans $G_{SC}$. 

\bigskip

Preuve. On fixe des mesures de Haar sur tous les groupes intervenant.  Introduisons les espaces $I(G({\mathbb R}))_{unip,loc}$ et $SI(G({\mathbb R}))_{unip,loc}$ de [I] 5.2 et [I] 5.5. L'indice $unip$ remplace comme \`a notre habitude l'indice ${\cal O}$ de ce paragraphe, cette classe ${\cal O}$ \'etant r\'eduite \`a $\{1\}$. Rappelons la description que l'on a donn\'ee en [I] 5.2 de l'espace $I(G({\mathbb R}))_{unip,loc}$. On fixe un ensemble ${\cal T}$ de repr\'esentants des classes de conjugaison par $G({\mathbb R})$ de sous-tores maximaux de $G$. Pour $T\in {\cal T}$, notons $M_{T}$ le plus petit Levi de $G$ contenant $T$. Il est d\'etermin\'e par la condition $A_{M_{T}}=A_{T}$. Notons $\Sigma^{M_{T}}(T)$ l'ensemble des racines de $T$ dans $\mathfrak{m}_{T}$. Ce sont les racines imaginaires de $T$. Notons $\mathfrak{t}_{\star}$ le sous-ensemble des $X\in \mathfrak{t}$ tels que $\alpha(X)\not=0$ pour tout $\alpha\in \Sigma^{M_{T}}(T)$. On note $\underline{\Omega}_{T}$ l'ensemble des composantes connexes de $\mathfrak{t}_{\star}({\mathbb R})$ et on fixe un \'el\'ement $\Omega_{T}\in \underline{\Omega}_{T}$.  Posons $W_{{\mathbb R }}(T)=Norm_{G({\mathbb R})}(T)/T({\mathbb R})$ et $W(T)=Norm_{G}(T)/T$. On a $W_{{\mathbb R}}(T)\subset W(T)^{\Gamma_{{\mathbb R}}}$.  Pour $w\in W(T)^{\Gamma_{{\mathbb R}}}$, notons $n(w)$ le nombre de racines $\alpha\in \Sigma^{M_{T}}(T)$ telles que l'hyperplan noyau de $\alpha$ s\'epare $\Omega_{T}$ de $w(\Omega_{T})$. On pose $\epsilon(w)=(-1)^{n(w)}$. C'est un caract\`ere de $ W(T)^{\Gamma_{{\mathbb R}}}$. En rempla\c{c}ant $G$ par $M_{T}$, on d\'efinit le sous-groupe $W^{M_{T}}(T)\subset W(T)$. Remarquons que $\Gamma_{{\mathbb R}}$ agit trivialement sur ce sous-groupe  et que $W^{M_{T}}(T)$ agit transitivement sur $\underline{\Omega}_{T}$. 
On note ${\mathbb C}[[\mathfrak{t}({\mathbb R})]]$ l'espace des s\'eries formelles sur $\mathfrak{t}({\mathbb R})$. Soient $f\in I(G({\mathbb R}))$, $T\in {\cal T}$ et $\Omega\in \underline{\Omega}_{T}$. On d\'efinit une fonction $\phi_{f,T,\Omega}$ sur $\Omega$ par $\phi_{f,T,\Omega}(X)=I^G(exp(X),f)$. Cette fonction se prolonge en une fonction $C^{\infty}$ dans un voisinage de $\Omega$. On note $\varphi_{f,T,\Omega}$ son d\'eveloppement en s\'erie formelle en $0$. Pour $w\in W_{{\mathbb R}}(T)$ et $X\in \mathfrak{t}({\mathbb R})$, on a l'\'egalit\'e
$$(1) \qquad \varphi_{f,T,w(\Omega)}(w(X))=\varphi_{f,T,\Omega}(X).$$
L'application $f\mapsto (\varphi_{f,T,\Omega})_{T\in {\cal T},\Omega\in \underline{\Omega}_{T}}$ se quotiente en une injection
$$I(G({\mathbb R}))_{unip,loc}\to \oplus_{T\in {\cal T},\Omega\in \underline{\Omega}_{T}}{\mathbb C}[[\mathfrak{t}({\mathbb R})]].$$
Cette injection est un hom\'eomorphisme de l'espace de d\'epart sur son image, laquelle est ferm\'ee dans l'espace d'arriv\'ee.
L'espace $I(G({\mathbb R}))_{unip,loc}$ est muni d'une filtration $({\cal F}^nI(G({\mathbb R}))_{unip,loc})_{n=a_{G}-1,...,a_{M_{0}}}$, o\`u $M_{0}$ est un Levi minimal. L'espace ${\cal F}^nI(G({\mathbb R}))_{unip,loc}$ est celui des $f\in I(G({\mathbb R}))_{unip,loc}$ tels que $\varphi_{f,T,\Omega}=0$ pour tous $T$, $\Omega$ tels que $a_{T}>n$. Pour d\'ecrire le gradu\'e $Gr^nI(G({\mathbb R}))_{unip,loc}$, introduisons les notations suivantes. Pout tout $T\in {\cal T}$, le groupe   $W_{{\mathbb R}}(T)$ agit naturellement sur ${\mathbb C}[[\mathfrak{t}({\mathbb R})]]$. On note ${\mathbb C}[[\mathfrak{t}({\mathbb R})]]^{W_{{\mathbb R}}(T),\epsilon}$ le sous-espace  isotypique pour l'action du groupe $W_{{\mathbb R}}(T)$, de type  $\epsilon$. Notons ${\cal T}^n$ le sous-ensemble des $T\in {\cal T}$ tels que $a_{T}=n$. Pour $f\in {\cal F}^nI(G({\mathbb R}))_{unip,loc}$ et $T\in {\cal T}^n$, la s\'erie $\varphi_{f,T,\Omega_{T}}$ appartient \`a ${\mathbb C}[[\mathfrak{t}({\mathbb R})]]^{W_{{\mathbb R}}(T),\epsilon}$. Pour une autre composante $\Omega\in \underline{\Omega}_{T}$, soit $w\in W ^{M_{T}}(T)$ tel que $\Omega=w(\Omega_{T})$. On a alors 
$$(2) \qquad \varphi_{f,T,\Omega}=\epsilon(w)\varphi_{f,T,\Omega_{T}}.$$
 L'application $f\mapsto (\varphi_{f,T,\Omega_{T}})_{T\in {\cal T}^n}$ se quotiente en un hom\'eomorphisme
$$Gr^nI(G({\mathbb R}))_{unip,loc}\simeq \oplus_{T\in {\cal T}^n}{\mathbb C}[[\mathfrak{t}({\mathbb R})]]^{W_{{\mathbb R}}(T),\epsilon}.$$

Des descriptions analogues valent pour l'espace $SI(G({\mathbb R}))_{unip,loc}$. On doit seulement remplacer pour tout $T$ le groupe $W_{{\mathbb R}}(T)$ par $W(T)^{\Gamma_{{\mathbb R}}}$. Notons $s^G:I(G({\mathbb R}))_{unip,loc}\to SI(G({\mathbb R}))_{unip,loc}$ la projection naturelle et $s^{G,n}$ l'application d\'eduite entre les gradu\'es de degr\'e $n$. On voit que $s^{G,n}$ est la somme sur $T\in {\cal T}$ des applications  
$$\begin{array}{ccc}{\mathbb C}[[\mathfrak{t}({\mathbb R})]]^{W_{{\mathbb R}}(T),\epsilon}&\to& {\mathbb C}[[\mathfrak{t}({\mathbb R})]]^{W(T)^{\Gamma_{{\mathbb R}}},\epsilon}\\ \varphi&\mapsto&\sum_{w\in W(T)^{\Gamma_{{\mathbb R}}}/W_{{\mathbb R}}(T)}\epsilon(w)w(\varphi).\\  \end{array}$$
 Cette application admet une section naturelle
 $$\begin{array}{ccc}{\mathbb C}[[\mathfrak{t}({\mathbb R})]]^{W(T)^{\Gamma_{{\mathbb R}},\epsilon}}&\to&{\mathbb C}[[\mathfrak{t}({\mathbb R})]]^{W_{{\mathbb R}}(T),\epsilon}\\ \varphi&\mapsto& [W(T)^{\Gamma_{{\mathbb R}}}:W_{{\mathbb R}}(T)]^{-1}\varphi.\\ \end{array}$$
 On obtient ainsi une identification de $Gr^nSI(G({\mathbb R}))_{unip,loc}$  \`a un sous-espace de $Gr^nI(G({\mathbb R}))_{unip,loc}$.

Il est facile de reprendre les preuves des paragraphes 4.15 et 4.16 de [I] (qui concernaient les espaces $I(G({\mathbb R}))$ et $SI(G({\mathbb R}))$) et de montrer que les r\'esultats de ces paragraphes valent pour nos espaces $I(G({\mathbb R}))_{unip,loc}$ et $SI(G({\mathbb R}))_{unip,loc}$.

On d\'ecrit de m\^eme les espaces $I(M({\mathbb R}))_{unip,loc}$, $I(M_{SC}({\mathbb R}))_{unip,loc}$ etc... et leurs gradu\'es. Selon notre habitude, on  ajoute  si n\'ecessaire des exposants $M$ ou $M_{SC}$ dans les notations pour les objets  obtenus en rempla\c{c}ant l'espace ambiant $G$ par $M$ ou $M_{SC}$. 

On dispose de l'application de restriction $res_{M}:I(G({\mathbb R}))_{unip,loc}\to I(M({\mathbb R}))_{unip,loc}$. De la projection $M_{SC}({\mathbb R})\to M({\mathbb R})$ se d\'eduit une application naturelle $\iota_{M_{SC},M}:I(M({\mathbb R}))_{unip,loc}\to I(M_{SC}({\mathbb R}))_{unip,loc}$. Notons $res_{M_{SC}}=\iota_{M_{SC},M}\circ res_{M}:I(G({\mathbb R}))_{unip,loc}\to I(M_{SC}({\mathbb R}))_{unip,loc}$. D\'ecrivons cette application. Soit $f\in I(G({\mathbb R}))_{unip,loc}$, posons $f'=res_{M_{SC}}(f)$. Soient $T'\in {\cal T}^{M}$ et $\Omega'\in \underline{\Omega}_{T'}$. La fonction $\varphi_{f',T',\Omega'}$ appartient \`a ${\mathbb C}[[\mathfrak{t}'_{sc}({\mathbb R})]]$, o\`u $T'_{sc}$ est l'image r\'eciproque de $T'$ dans $M_{SC}$. Fixons $g\in G({\mathbb R})$ tel que $ad_{g}(T')\in {\cal T}$. Posons $T=ad_{g}(T')$ et $\Omega=ad_{g}(\Omega')$. Alors, pour $X'\in \mathfrak{t}'_{sc}({\mathbb R})$, on a l'\'egalit\'e
$$\varphi_{f',T',\Omega'}(X')=\varphi_{f,T,\Omega}(ad_{g}(X')).$$
Pour $f'\in I(M_{SC}({\mathbb R}))_{unip,loc}$, $T\in {\cal T}$ et $\Omega\in \underline{\Omega}_{T}$, consid\'erons la condition

$(Hyp)_{f',T,\Omega}$ il existe $\varphi\in {\mathbb C}[[\mathfrak{t}({\mathbb R})]]$ telle que, pour tout $T'\in {\cal T}^{M}$, tout $g\in G({\mathbb R})$ tel que $ad_{g}(T')=T$ et tout $X'\in \mathfrak{t}'_{sc}({\mathbb R})$, on ait l'\'egalit\'e
$$\varphi(ad_{g}(X')=\varphi_{f',T',ad_{g^{-1}}(\Omega)}(X').$$

Notons $V$ le  sous-espace des $f'\in I(M_{SC}({\mathbb R}))_{unip,loc}$ tels que, pour tous $T\in {\cal T}$ et $\Omega\in \underline{\Omega}_{T}$, la condition $(Hyp)_{f',T,\Omega}$ soit v\'erifi\'ee.
D'apr\`es la description ci-dessus,  l'image de $res_{M_{SC}}$ est contenue dans $V$. 
L'application $res_{M_{SC}}$ est compatible aux filtrations. On note $res_{M_{SC}}^n$ l'application d\'eduite entre gradu\'es de degr\'e $n$. On note $p^n:{\cal F}^nI(G({\mathbb R}))_{unip,loc}\to Gr^nI(G({\mathbb R}))_{unip,loc}$ et $p^{M_{SC},n}:{\cal F}^nI(M_{SC}({\mathbb R}))_{unip,loc}\to Gr^nI(M_{SC}({\mathbb R}))_{unip,loc}$ les projections naturelles. Nous allons prouver que

(3) l'image de $res_{M_{SC}}$ est \'egale \`a $V$ ;

(4) l'intersection avec ${\cal F}^nI(M_{SC}({\mathbb R}))_{unip,loc}$ de l'image de $res_{M_{SC}}$ est l'image par $res_{M_{SC}}$ de ${\cal F}^nI(G({\mathbb R}))_{unip,loc}$;

(5) soit $f^n\in Ker(res_{M_{SC}}^n)$; alors il existe $f\in Ker(res_{M_{SC}})\cap {\cal F}^nI(G({\mathbb R}))_{unip,loc}$ telle que $p^n(f)=f^n$. 

  On a vu que l'image de $res_{M_{SC}}$ \'etait contenue dans $V$. On va d'abord prouver

(6) soit $f'\in V\cap {\cal F}^nI(M_{SC}({\mathbb R}))_{unip,loc}$; alors il existe $f^n\in Gr^nI(G({\mathbb R}))_{unip,loc}$ tel que $res^n_{M_{SC}}(f^n)=p^{M_{SC},n}(f')$. 

Soit $T\in {\cal T}^n$. Consid\'erons une fonction $\varphi\in {\mathbb C}[[\mathfrak{t}({\mathbb R})]]$ satisfaisant l'hypoth\`ese $(Hyp)_{f',T,\Omega_{T}}$.  Soit $w\in W_{{\mathbb R}}(T)$. Montrons  que la fonction $\epsilon(w)w(\varphi)$ satisfait encore cette hypoth\`ese. On fixe $x\in Norm_{G({\mathbb R})}(T)$ tel que $w$ soit la restriction de $ad_{x}$ \`a $T$. Soient $T'\in {\cal T}^M$, $g\in G({\mathbb R})$ tel que $ad_{g}(T')=T$ et $X'\in \mathfrak{t}'_{sc}({\mathbb R})$. On a
$$\epsilon(w)w(\varphi)(ad_{g}(X'))=\epsilon(w)\varphi(ad_{x^{-1}g}(X')=\epsilon(w)\varphi_{f',T',ad_{g^{-1}x}(\Omega_{T})}(X'),$$
d'apr\`es l'hypoth\`ese $(Hyp)_{f',T,\Omega_{T}}$ appliqu\'ee \`a $T'$ et $x^{-1}g$. Soit $w'\in W^{M_{T'}}(T')$ tel que $ad_{g^{-1}x}(\Omega_{T})=w'(ad_{g^{-1}}(\Omega_{T})$. 
On a $T'\in {\cal T}^{M,n}$. L'hypoth\`ese $f'\in  {\cal F}^nI(M_{SC}({\mathbb R}))_{unip,loc}$ entra\^{\i}ne que l'analogue de la condition (2) est v\'erifi\'ee. C'est-\`a-dire que
$$\varphi_{f',T',ad_{g^{-1}x}(\Omega_{T})}(X')=\epsilon'(w')\varphi_{f',T',ad_{g^{-1}}(\Omega_{T})}(X'),$$
en notant pour plus de pr\'ecision $\epsilon'$ le caract\`ere de $W^{M_{T'}}(T')$. Mais l'application $\alpha\mapsto ad_{g^{-1}}(\alpha)$ se restreint en une bijection entre l'ensemble des racines $\alpha\in \Sigma^{M_{T}}(T)$ telles que l'hyperplan noyau de $\alpha$ s\'epare $\Omega_{T}$ de $w(\Omega_{T})$ et l'ensemble des racines $\beta\in \Sigma^{M_{T'}}(T')$ telles que l'hyperplan noyau de $\beta$ s\'epare $ad_{g^{-1}}(\Omega_{T})$ de $ad_{g^{-1}x}(\Omega_{T})$. Il en r\'esulte que $\epsilon(w)=\epsilon'(w')$. On obtient alors
$$\epsilon(w)w(\varphi)(ad_{g}(X'))=\varphi_{f',T',ad_{g^{-1}}(\Omega_{T})}(X').$$
 Donc $\epsilon(w)w(\varphi)$ v\'erifie $(Hyp)_{f',T,\Omega_{T}}$. On peut alors remplacer $\varphi$ par
 $$\vert W_{{\mathbb R}}(T)\vert ^{-1}\sum_{w\in W_{{\mathbb R}}(T)}\epsilon(w)w(\varphi).$$
 Cette fonction satisfait encore $(Hyp)_{f',T,\Omega_{T}}$ et appartient \`a ${\mathbb C}[[\mathfrak{t}({\mathbb R})]]^{W_{{\mathbb R}}(T),\epsilon}$.  Notons plus pr\'ecis\'ement $\varphi_{T}$ cette fonction. La famille $(\varphi_{T})_{T\in {\cal T}^n}$ s'identifie \`a un \'el\'ement $f^n\in Gr^nI(G({\mathbb R}))_{unip,loc}$. On v\'erifie que $res^n_{M_{SC}}(f^n)=p^{M_{SC},n}(f')$. Cela prouve (6). 
 
 Montrons maintenant que
 
(7) $V\cap {\cal F}^nI(M_{SC}({\mathbb R}))_{unip,loc}$ est l'image par $res_{M_{SC}}$ de ${\cal F}^nI(G({\mathbb R}))_{unip,loc}$.

On raisonne par r\'ecurrence sur $n$. C'est clair si $n<a_{G}$ puisqu'alors tous les espaces sont nuls. Soit $n\geq a_{G}$, supposons l'assertion d\'emontr\'ee pour $n-1$. Soit $f'\in V\cap {\cal F}^nI(M_{SC}({\mathbb R}))_{unip,loc}$. D'apr\`es (6), il existe $f\in {\cal F}^nI(G({\mathbb R}))_{unip,loc}$ tel que $p^{M_{SC},n}(f'-res_{M_{SC}}(f))=0$. Cela entra\^{\i}ne $f'-res_{M_{SC}}(f)\in {\cal F}^{n-1}I(M_{SC}({\mathbb R}))_{unip,loc}$.  Puisque l'image de $res_{M_{SC}}$ est contenu dans $V$, on a aussi $f'-res_{M_{SC}}(f)\in V$. En appliquant l'hypoth\`ese de r\'ecurrence, $f'-res_{M_{SC}}(f)$ appartient \`a l'image  par $res_{M_{SC}}$ de ${\cal F}^{n-1}I(G({\mathbb R}))_{unip,loc}$. Cela entra\^{\i}ne que $f'$ appartient \`a l'image  par $res_{M_{SC}}$ de ${\cal F}^nI(G({\mathbb R}))_{unip,loc}$, d'o\`u (7).

Pour $n$ maximal, (7) implique (3). Une fois (3) connue, (4) est \'equivalent \`a (7). Enfin, soit $f^n\in Ker(res_{M_{SC}}^n)$. Choisissons $f_{0}\in {\cal F}^nI(G({\mathbb R}))_{unip,loc}$ tel que $p^n(f_{0})=f^n$. L'hypoth\`ese $f^n\in Ker(res_{M_{SC}}^n)$ implique  $res_{M_{SC}}(f_{0})\in {\cal F}^{n-1}I(M_{SC}({\mathbb R}))_{unip,loc}$. En appliquant (4), il existe $f_{1}\in {\cal F}^{n-1}I(G({\mathbb R}))_{unip,loc}$ tel que $res_{M_{SC}}(f_{1})=res_{M_{SC}}(f_{0})$. L'\'el\'ement $f=f_{0}-f_{1}$ appartient \`a $Ker(res_{M_{SC}})\cap {\cal F}^nI(G({\mathbb R}))_{unip,loc}$ et v\'erifie $p^n(f)=f^n$. Cela prouve (5). 

Montrons que

(8) l'image de l'application $res_{M_{SC}}$ est ferm\'ee dans $I(M_{SC}({\mathbb R}))_{unip,loc}$.

 Pour tout $T'\in {\cal T}^{M}$,    fixons un ensemble de repr\'esentants ${\cal W}_{T'}$ du quotient $\{g\in G({\mathbb R}); ad_{g}(T')\in {\cal T}\}/T'({\mathbb R})$.  Consid\'erons l'espace
 $$X=\oplus_{T'\in {\cal T}^M,\Omega'\in \underline{\Omega}_{T'},g\in {\cal W}_{T'}}{\mathbb C}[[\mathfrak{t}_{sc}'({\mathbb R})]].$$
 L'espace 
 $$Y=\oplus_{T'\in {\cal T}^M,\Omega'\in \underline{\Omega}_{T'} }{\mathbb C}[[\mathfrak{t}_{sc}'({\mathbb R})]]$$
 s'identifie au sous-espace ferm\'e des \'el\'ements de $X$ dont les composantes sont ind\'ependantes de $g\in {\cal W}_{T'}$. On a dit que $I(M_{SC}({\mathbb R}))_{unip,loc}$ s'identifiait \`a un sous-espace ferm\'e de $Y$. Il s'identifie donc aussi \`a un sous-espace ferm\'e de $X$. Posons
 $$Z=\oplus_{T\in {\cal T},\Omega\in \underline{\Omega}_{T}}{\mathbb C}[[\mathfrak{t}({\mathbb R})]].$$
 On d\'efinit une application $\psi:Z\to X$ de la fa\c{c}on suivante. Soit $(\varphi_{T,\Omega})_{T\in {\cal T},\Omega\in \underline{\Omega}_{T}}\in Z$. Soient $T'\in {\cal T}^M$, $\Omega'\in \underline{\Omega}_{T'}$ et $g\in {\cal W}_{T'}$. Pour $X'\in \mathfrak{t}'_{sc}({\mathbb R})$, posons
 $$\varphi_{T',\Omega',g}(X')=\varphi_{ad_{g}(T),ad_{g}(\Omega')}(ad_{g}(X')).$$
 L'image de $(\varphi_{T,\Omega})_{T\in {\cal T},\Omega\in \underline{\Omega}_{T}}$ par $\psi$ est alors la famille $(\varphi_{T',\Omega',g})_{T'\in {\cal T}^M,\Omega'\in \underline{\Omega}_{T'},g\in {\cal W}_{T'}}$. L'assertion (3) revient \`a dire que l'image de $res_{M_{SC}}$ est \'egale \`a l'intersection dans $X$ de $I(M_{SC}({\mathbb R}))_{unip,loc}$ et de l'image de $\psi$. Pour prouver (8), il suffit donc de prouver que l'image de $\psi$ est ferm\'ee dans $X$. L'espace des series formelles, disons sur $\mathfrak{t}({\mathbb R})$, est le produit sur $i\in {\mathbb N}$ des espaces de polyn\^omes homog\`enes de degr\'e $i$ sur $\mathfrak{t}({\mathbb R})$. Ce r\'esultat s'\'etend bien s\^ur \`a nos espaces $X$ et $Z$: on a $X=\prod_{i\in {\mathbb N}}X_{i}$ et $Z=\prod_{i\in {\mathbb N}}Z_{i}$. L'application $\psi$ est le produit d'applications $\psi_{i}:Z_{i}\to X_{i}$. Soit alors $x\in X$ qui est dans l'adh\'erence de l'image de $\psi$. Pour tout $i$, la composante $x_{i}\in X_{i}$ est dans l'adh\'erence de l'image de $\psi_{i}$. Or les espaces $X_{i}$ et $Z_{i}$ sont de dimensions finies donc l'image de $\psi_{i}$ est ferm\'ee. On peut donc choisir $z_{i}\in Z_{i}$ tel que $\psi_{i}(z_{i})=x_{i}$. Mais alors, l'\'el\'ement $z=\prod_{i\in {\mathbb N}}z_{i}\in Z$ v\'erifie $\psi(z)=x$. Donc $x$ appartient \`a l'image de $\psi$. Cela d\'emontre (8).

On a de m\^eme une application $res_{M_{SC}}^{st}=\iota_{M_{SC},M}\circ res_{M}^{st}:SI(G({\mathbb R}))_{unip,loc}\to SI(M_{SC}({\mathbb R}))_{unip,loc}$. Elle se d\'ecrit  de la m\^eme fa\c{c}on que $res_{M_{SC}}$.   Des propri\'et\'es analogues \`a (3), (4), (5) et (8) valent pour cette application.

Nos applications sont compatibles aux filtrations. Pour tout $n$, on en d\'eduit un diagramme commutatif
$$\begin{array}{ccc}Gr^nI(G({\mathbb R}))_{unip,loc}&\stackrel{res^n_{M_{SC}}}{\to}&Gr^nI(M_{SC}({\mathbb R}))_{unip,loc}\\ s^{G,n}\downarrow\,&&s^{M_{SC},n}\downarrow\,\\ Gr^nSI(G({\mathbb R}))_{unip,loc}&\stackrel{res^{st,n}_{M_{SC}}}{\to}&Gr^nSI(M_{SC}({\mathbb R}))_{unip,loc}\\ \end{array}$$

Montrons que

(9) $Im(res^n_{M_{SC}})\cap Ker(s^{M_{SC},n})=res^n_{M_{SC}}(Ker(s^{G,n}))$.

Comme on l'a vu ci-dessus, on peut identifier $Gr^nSI(G({\mathbb R}))_{unip,loc}$ \`a un sous-espace de $Gr^nI(G({\mathbb R}))_{unip,loc}$. On a alors
$$Gr^nI(G({\mathbb R}))_{unip,loc}=Ker(s^{G,n})\oplus Gr^nSI(G({\mathbb R}))_{unip,loc}.$$
De m\^eme
$$Gr^nI(M_{SC}({\mathbb R}))_{unip,loc}=Ker(s^{M_{SC},n})\oplus Gr^nSI(M_{SC}({\mathbb R}))_{unip,loc} .$$
Les deux sections sont compatibles aux applications $res^n_{M_{SC}}$ et $res^{st,n}_{M_{SC}}$. Pr\'ecis\'ement, $res^n_{M_{SC}}$ envoie $Gr^nSI(G({\mathbb R}))_{unip,loc}$ dans $Gr^nSI(M_{SC}({\mathbb R}))_{unip,loc}$  et co\"{\i}ncide avec $res^{st,n}_{M_{SC}}$ sur $Gr^nSI(G({\mathbb R}))_{unip,loc}$. Cela r\'esulte des descriptions ci-dessus et de la propri\'et\'e suivante

(10) pour tout $T'\in {\cal T}^M$, on a l'\'egalit\'e 
$$[W(T')^{\Gamma_{{\mathbb R}}}:W_{{\mathbb R}}(T')]=[W^M(T')^{\Gamma_{{\mathbb R}}}:W^M_{{\mathbb R}}(T')].$$

On v\'erifie que les suites suivantes sont exactes
$$1\to W^{M_{T'}}(T')^{\Gamma_{{\mathbb R}}}\to W(T')^{\Gamma_{{\mathbb R}}}\to W(M_{T'})\to 1,$$
$$1\to W^{M_{T'}}_{{\mathbb R}}(T')\to W_{{\mathbb R}}(T')\to W(M_{T'})\to 1.$$
On en d\'eduit l'\'egalit\'e
$$[W(T')^{\Gamma_{{\mathbb R}}}:W_{{\mathbb R}}(T')]=[W^{M_{T'}}(T')^{\Gamma_{{\mathbb R}}}:W^{M_{T'}}_{{\mathbb R}}(T')].$$
On a de m\^eme
$$[W^M(T')^{\Gamma_{{\mathbb R}}}:W^M_{{\mathbb R}}(T')]=[W^{M_{T'}}(T')^{\Gamma_{{\mathbb R}}}:W^{M_{T'}}_{{\mathbb R}}(T')].$$
D'o\`u (10). 

Soit alors $f\in Gr^nI(G({\mathbb R}))_{unip,loc}$, supposons $res^n_{M_{SC}}(f)\in Ker(s^{M_{SC},n})$. Ecrivons $f=f_{0}+f^{st}$, avec $f_{0}\in Ker(s^{G,n})$ et $f^{st}\in Gr^nSI(G({\mathbb R}))_{unip,loc}$. On a $s^{M_{SC},n}\circ res^n_{M_{SC}}(f)=res^n_{M_{SC}}(f^{st})$. Ce terme est nul par hypoth\`ese. Donc $res^n_{M_{SC}}(f)=res^n_{M_{SC}}(f_{0})$, c'est-\`a-dire $res^n_{M_{SC}}(f)$ appartient \`a $res^n_{M_{SC}}(Ker(s^{G,n}))$. Cela prouve (9). 

Montrons que l'assertion (9) se rel\`eve des gradu\'es aux espaces eux-m\^emes, c'est-\`a-dire

(11) $Im(res_{M_{SC}})\cap Ker(s^{M_{SC}})=res_{M_{SC}}(Ker(s^G))$. 

\noindent On d\'emontre   par r\'ecurrence sur $n$ que

(12) $res_{M_{SC}}({\cal F}^nI(G({\mathbb R}))_{unip,loc})\cap Ker(s^{M_{SC}})\subset res_{M_{SC}}(Ker(s^G))$. 

Soit $f'\in res_{M_{SC}}({\cal F}^nI(G({\mathbb R}))_{unip,loc})\cap Ker(s^{M_{SC}})$. Alors $p^{M_{SC},n}(f')\in Im(res^n_{M_{SC}})\cap Ker(s^{M_{SC},n})$. D'apr\`es (9), il existe $f^n\in Ker(s^{G,n})$ tel que $ res^n_{M_{SC}}(f^n)=p^{M_{SC},n}(f')$. Comme on l'a vu en [I] 4.17 (on a dit que les r\'esultats de ce paragraphe s'appliquaient \`a nos espaces localis\'es), il existe $f\in {\cal F}^nI(G({\mathbb R}))_{unip,loc}\cap Ker(s^G)$ tel que $p^n(f)=f^n$.   Posons $f'_{0}=f'-res_{M_{SC}}(f)$. Les \'egalit\'es ci-dessus entra\^{\i}nent $p^{M_{SC},n}(f'_{0})=0$, donc $f'_{0}\in {\cal F}^{n-1}I(M_{SC}({\mathbb R}))_{unip,loc}$. L'\'el\'ement $f'_{0}$ appartient encore \`a l'image de $res_{M_{SC}}$. En appliquant (4), on a $f'_{0}\in res_{M_{SC}}({\cal F}^{n-1}I(G({\mathbb R}))_{unip,loc})$. Puisque $f\in Ker(s^G)$, on a $res_{M_{SC}}(f)\in Ker(s^{M_{SC}})$ donc aussi $f'_{0}\in Ker(s^{M_{SC}})$. En appliquant l'hypoth\`ese de r\'ecurrence, on obtient $f'_{0}\in res_{M_{SC}}(Ker(s^G))$. Alors $f'=f'_{0}+ res_{M_{SC}}(f)$ appartient aussi \`a $res_{M_{SC}}(Ker(s^G))$. Cela prouve (12).

Pour $n$ maximal, (12) entra\^{\i}ne que le membre de gauche de (11) est inclus dans celui de droite. L'inclusion oppos\'ee est \'evidente. D'o\`u (11).

 Notons $I^{inst}(G({\mathbb R}))$ le noyau de l'application naturelle $I(G({\mathbb R}))\to SI(G({\mathbb R}))$. Rappelons que
 
 (13) l'image naturelle de $I^{inst}(G({\mathbb R}))$ dans $I(G({\mathbb R}))_{unip,loc}$ est dense dans $Ker(s^G)$. 
 
 En effet, c'est l'assertion [I] 5.15 (4) dans le cas particulier $k=1$ et $\tilde{M}_{1}=\tilde{G}=G$. 

Venons-en \`a la preuve du lemme.  Soit $\boldsymbol{\gamma}$ comme dans l'\'enonc\'e. C'est une forme lin\'eaire continue sur $I(M_{SC}({\mathbb R}))_{unip,loc}$. On cherche un \'el\'ement $\boldsymbol{\delta}\in D^{st}_{unip}(M_{SC}({\mathbb R}))$, c'est-\`a-dire une forme lin\'eaire continue sur $SI(M_{SC}({\mathbb R}))_{unip,loc}$, ou encore une forme lin\'eaire continue sur $I(M_{SC}({\mathbb R}))_{unip,loc}$ nulle sur $Ker(s^{M_{SC}})$. La condition $(\iota^*_{M_{SC},M}(\boldsymbol{\gamma}))^G=(\iota^*_{M_{SC},M}(\boldsymbol{\delta}))^G$ revient \`a dire que $\boldsymbol{\gamma}$ et $\boldsymbol{\delta}$ co\"{\i}ncident sur l'image de $res_{M_{SC}}$. Cette image est ferm\'ee d'apr\`es (8) et la somme de cette image avec $Ker(s^{M_{SC}})$ est aussi ferm\'ee: c'est l'image r\'eciproque par $s^{M_{SC}}$ de l'image de $res^{st}_{M_{SC}}$ qui est ferm\'ee par l'analogue de (8). La condition n\'ecessaire et suffisante pour que $\boldsymbol{\delta}$ existe est donc que $\boldsymbol{\gamma}$ annule $Im(res_{M_{SC}})\cap Ker(s^{M_{SC}})$. Ou encore, d'apr\`es (11), que $\boldsymbol{\gamma}$ annule $res_{M_{SC}}(Ker(s^G))$. D'apr\`es (13),  cet espace est l'adh\'erence de l'image par $res_{M_{SC}}$ de l'image naturelle de $I^{inst}(G({\mathbb R}))$ dans $I(G({\mathbb R}))_{unip,loc}$. L'hypoth\`ese que  $(\iota^*_{M_{SC},M}(\boldsymbol{\gamma}))^G$ est stable signifie que $\boldsymbol{\gamma}$ annule cette image. Etant continue, $\boldsymbol{\gamma}$ annule aussi son adh\'erence, c'est-\`a-dire $res_{M_{SC}}(Ker(s^G))$. Cela ach\`eve la d\'emonstration. $\square$

\bigskip

\subsection{Un corollaire}
 Soient $M$ un Levi de $G$ et $J\in {\cal J}_{M}^G(B)$. 
 \ass{Corollaire}{(i) Soient $\boldsymbol{\gamma}\in D_{tr-unip}(M({\mathbb R}))\otimes Mes(M({\mathbb R}))^*$ et $a\in A_{M}({\mathbb R})$ en position g\'en\'erale et proche de $1$. Alors il existe $\boldsymbol{\tau}\in D_{unip}(M_{SC}({\mathbb R}))\otimes Mes(M_{SC}({\mathbb R}))^*$ tel que $\rho_{J}^G(\boldsymbol{\gamma},a)$ soit l'image de $\iota^*_{M_{SC},M}(\boldsymbol{\tau})$ modulo $Ann_{unip}^G$.
 
 (ii) Supposons $G$ quasi-d\'eploy\'e.  Soient $\boldsymbol{\delta}\in D_{tr-unip}^{st}(M({\mathbb R}))\otimes Mes(M({\mathbb R}))^*$ et $a\in A_{M}({\mathbb R})$ en position g\'en\'erale et proche de $1$. Alors il existe $\boldsymbol{\tau}\in D_{unip}^{st}(M_{SC}({\mathbb R}))\otimes Mes(M_{SC}({\mathbb R}))^*$ tel que $\sigma_{J}^G(\boldsymbol{\delta},a)$ soit l'image de $\iota^*_{M_{SC},M}(\boldsymbol{\tau})$ modulo $Ann_{unip}^{st,G}$.}
 
 Preuve. Pour $\boldsymbol{\gamma}\in D_{orb,unip}(M({\mathbb R}))\otimes Mes(M({\mathbb R}))^*$,  l'\'el\'ement $\rho_{J}^G(\boldsymbol{\gamma},a)$ appartient par d\'efinition \`a $D_{orb,unip}(M({\mathbb R}))\otimes Mes(M({\mathbb R}))^*$ (plus exactement, est l'image modulo $Ann_{unip}^G$ d'un \'el\'ement de cet espace). Puisque $D_{orb,unip}(M({\mathbb R}))\otimes Mes(M({\mathbb R}))^*$ est l'image par  $\iota^*_{M_{SC},M}$ de $D_{orb,unip}(M_{SC}({\mathbb R}))\otimes Mes(M_{SC}({\mathbb R}))^*$, l'assertion (i) vaut pour $\boldsymbol{\gamma}$. Soit ${\bf M}'=(M',{\cal M}',\zeta)$ une donn\'ee endoscopique elliptique et relevante de $M$, avec $M'\not=M$ si $G$ est quasi-d\'eploy\'e. Soit $\boldsymbol{\delta}\in D_{tr-unip}^{st}({\bf M}')\otimes Mes(M'({\mathbb R}))^*$, consid\'erons l'\'el\'ement $\boldsymbol{\gamma}=transfert(\boldsymbol{\delta})$. Par d\'efinition, on a
 $$\rho_{J}^G(\boldsymbol{\gamma},a)=\rho_{J}^{G,{\cal E}}({\bf M}',\boldsymbol{\delta},a) =\sum_{s\in \zeta Z(\hat{M})^{\Gamma_{{\mathbb R}}}/Z(\hat{G})^{\Gamma_{{\mathbb R}}}} i_{M'}(G,G'(s))$$
 $$\sum_{J'\in {\cal J}_{M'}^{G'(s)}(B); J'\mapsto J}transfert(\sigma_{J'}^{{\bf G}'(s)}(\boldsymbol{\delta},\xi(a))).$$
 Fixons $s$ et $J'$ apparaissant dans cette somme. En appliquant le (ii) de l'\'enonc\'e par r\'ecurrence, on peut supposer que $\sigma_{J'}^{{\bf G}'(s)}(\boldsymbol{\delta},\xi(a))$ est l'image  par $\iota^*_{M'_{SC},M'}$ d'un \'el\'ement $\boldsymbol{\tau}_{s,J'}\in D_{unip}^{st}(M'_{SC}({\mathbb R}))\otimes Mes(M'_{SC}({\mathbb R}))^*$. On sait que la donn\'ee ${\bf M}'$ d\'etermine une donn\'ee endoscopique ${\bf M}'_{sc}$ de  $M_{SC}$. L'application  $\iota^*_{M'_{SC},M'}$ se factorise en $\iota^*_{M'_{sc},M'}\circ \iota^*_{M'_{SC},M'_{sc}}$. Posons $\boldsymbol{\tau}'_{s,J'}=\iota^*_{M'_{SC},M'_{sc}}(\boldsymbol{\tau}_{s,J'})$. On a alors $\sigma _{J'}^{{\bf G}'(s)}(\boldsymbol{\delta},\xi(a))=\iota^*_{M'_{sc},M'}(\boldsymbol{\tau}'_{s,J'})$, d'o\`u
 $$transfert(\sigma_{J'}^{{\bf G}'(s)}(\boldsymbol{\delta},\xi(a)))=\iota^*_{M_{SC},M}\circ transfert(\boldsymbol{\tau}'_{s,J'}).$$
 On a alors $\rho_{J}^G(\boldsymbol{\gamma},a)=\iota^*_{M_{SC},M}(\boldsymbol{\tau})$, o\`u
 $$\boldsymbol{\tau}=\sum_{s\in \zeta Z(\hat{M})^{\Gamma_{{\mathbb R}}}/Z(\hat{G})^{\Gamma_{{\mathbb R}}}} i_{M'}(G,G'(s))\sum_{J'\in {\cal J}_{M'}^{G'(s)}(B); J'\mapsto J}transfert(\boldsymbol{\tau}'_{s,J'}).$$
 Cela d\'emontre (i) pour $\boldsymbol{\gamma}$ et cela ach\`eve la preuve de cette assertion (i).
 
 Supposons $G$ quasi-d\'eploy\'e. Pour  $\boldsymbol{\delta}\in D_{tr-unip}^{st}(M({\mathbb R}))\otimes Mes(M({\mathbb R}))^*$, on voit de m\^eme, en appliquant (i) et les hypoth\`eses de r\'ecurrence, que $\sigma_{J}^G(\boldsymbol{\delta},a)$ est l'image par $\iota^*_{M_{SC},M}$ d'un \'el\'ement $\boldsymbol{\tau}'\in D_{unip}(M_{SC}({\mathbb R}))\otimes Mes(M_{SC}({\mathbb R}))^*$. On sait que l'induite $\sigma_{J}^G(\boldsymbol{\delta},a)^G$ est stable. En appliquant le lemme  3.5, il existe $\boldsymbol{\tau}\in  D^{st}_{unip}(M_{SC}({\mathbb R}))\otimes Mes(M_{SC}({\mathbb R}))^*$ tel que les induites \`a $G$ des distributions $\iota^*_{M_{SC},M}(\boldsymbol{\tau}')$ et $\iota^*_{M_{SC},M}(\boldsymbol{\tau})$  soient \'egales. Mais alors $\sigma_{J}^G(\boldsymbol{\delta},a)$ est l'image de $\iota^*_{M_{SC},M}(\boldsymbol{\tau})$ modulo $Ann_{unip}^{st,G}$. Cela prouve (ii). $\square$
 
 \bigskip
 
 \section{Extension des d\'efinitions, cas quasi-d\'eploy\'e et \`a torsion int\'erieure}
 
  \bigskip
 
 \subsection{Descente et endoscopie}
 On a rappel\'e dans la section 5 de [III] les liens entre descente et endoscopie. Le corps de base y \'etait non-archim\'edien.   Presque tout reste valable sur notre corps de base r\'eel. Il faut toutefois  modifier l\'eg\`erement l'identit\'e cruciale [III] 5.1(3). Rappelons  bri\`evement la situation. Dans ce paragraphe, le triplet $(G,\tilde{G},{\bf a})$ est quelconque. On consid\`ere une donn\'ee endoscopique ${\bf G}'=(G',{\cal G}',\tilde{s})$ de $(G,\tilde{G},{\bf a})$ elliptique et relevante. On fixe des donn\'ees auxiliaires $G'_{1}$,...,$\Delta_{1}$. Consid\'erons un diagramme $(\epsilon,B',T',B,T,\eta)$ joignant deux \'el\'ements semi-simples $\epsilon\in \tilde{G}'({\mathbb R})$ et  $\eta\in \tilde{G}({\mathbb R})$. On suppose $G'_{\epsilon}$ quasi-d\'eploy\'e. On  fixe un rel\`evement $\epsilon_{1}$ de $\epsilon$ dans $\tilde{G}'_{1}({\mathbb R})$. On fixe une d\'ecomposition d'alg\`ebres de Lie  $\mathfrak{g}'_{1}=\mathfrak{c}_{1}\oplus \mathfrak{g}'$. On construit  comme en [I] 1.2 une action galoisienne quasi-d\'eploy\'ee sur  $G_{\eta}$ qui conserve une paire de Borel \'epingl\'ee de ce groupe compl\'etant la paire $(B\cap G_{\eta},T\cap G_{\eta})$. On note $\bar{G}$ ce groupe muni de cette action. Celle-ci est telle que l'application identit\'e $\psi:G_{\eta}\to \bar{G}$ est un torseur int\'erieur. On a construit en [W2] 3.5 une donn\'ee endoscopique $\bar{{\bf G}}'=(\bar{G}',\bar{{\cal G}}',\bar{s})$ de $\bar{G}_{SC}$. Le couple $(\bar{G}'_{SC},G'_{\epsilon,SC})$ se compl\`ete en un triplet endoscopique non standard. Soit $y\in G$ tel que $y\sigma(y)^{-1}\in I_{\eta}=Z(G)^{\theta}G_{\eta}$. On pose $\eta[y]=y^{-1}\eta y$. Alors $\psi\circ ad_{y}$ est un torseur int\'erieur de $G_{\eta[y]}$ sur $\bar{G}$. Ainsi, la donn\'ee $\bar{{\bf G}}'$ est aussi une donn\'ee endoscopique pour  $G_{\eta[y]}$. Supposons qu'elle soit relevante. Puisque $\bar{G}_{SC}$ est simplement connexe, il n'est pas besoin de donn\'ees auxiliaires pour cette donn\'ee et on peut fixer pour celle-ci un facteur de transfert $\Delta(y)$. Soit $Y\in \mathfrak{g}'_{\epsilon}({\mathbb R})$ en position g\'en\'erale et proche de $0$.  Modulo les isomorphismes de [III] 5.1, on le d\'ecompose en $Y_{sc}+Z$ avec $Y_{sc}\in \mathfrak{g}'_{\epsilon,SC}({\mathbb R})$ et $Z\in \mathfrak{z}(G'_{\epsilon};{\mathbb R})$. On d\'ecompose $Z$ en $Z_{1}+Z_{2}$, avec $Z_{1}\in \mathfrak{z}(\bar{G};{\mathbb R})$ et $Z_{2}\in \mathfrak{z}(\bar{G}';{\mathbb R})$. On transf\`ere $Y_{sc}$ en un \'el\'ement $\bar{Y}_{sc}\in \bar{\mathfrak{g}}'_{SC}({\mathbb R})$.  On pose $\bar{Y}=\bar{Y}_{sc}+Z_{2}$. Supposons que $\bar{Y}$ se transf\`ere en un \'el\'ement $X[y]_{sc}\in \mathfrak{g}_{\eta[y],SC}({\mathbb R})$. On pose $X[y]=X[y]_{sc}+Z_{1}$. On a
 
 (1) il existe $b\in \mathfrak{z}(G'_{\epsilon};{\mathbb R})^*$ et, pour tout $y$ comme ci-dessus, il existe $d(y)\in {\mathbb C}^{\times}$ de sorte que, pour toutes donn\'ees $Y$, $X[y]$ comme ci-dessus, on ait l'\'egalit\'e
 $$d(y)\Delta(y)(exp(\bar{Y}),exp(X[y]_{sc}))=e^{<b,Z>}\Delta_{1}(exp(Y)\epsilon_{1},exp(X[y])\eta[y]).$$
 
 Preuve. Commen\c{c}ons par fixer des tores maximaux dans nos diff\'erents groupes et supposons que les \'el\'ements $Y$, $\bar{Y}$ et $X[y]$ appartiennent aux alg\`ebres de Lie de ces tores. Le lemme [I] 2.8 entra\^{\i}ne qu'il existe $(b_{1},b_{2})\in \mathfrak{z}(G'_{1};{\mathbb R})^*\oplus (1-\theta)(\mathfrak{z}(G,{\mathbb R})^*)$ de sorte que la fonction
 $$e^{-<b_{1},Y>-<b_{2},X[y]>}\Delta_{1}(exp(Y)\epsilon_{1},exp(X[y])\eta[y])$$
 soit localement constante. On a prolong\'e par exemple $b_{1}$ en une forme lin\'eaire sur $\mathfrak{g}'_{1}({\mathbb R})$ nulle sur $\mathfrak{g}'_{1,SC}({\mathbb R})$. L'\'egalit\'e se simplifie puisque la projection de $X[y]$ dans $\mathfrak{z}(G;{\mathbb R})$  est invariante par $\theta$: on a $<b_{2},X[y]>=0$. On a aussi $<b_{1},Y>=<b_{0},Z>$, o\`u $b_{0}$ est la restriction de $b_{1}$ \`a $\mathfrak{z}(G'_{\epsilon};{\mathbb R})$. Un m\^eme r\'esultat vaut pour le facteur $\Delta(y)(exp(\bar{Y}),exp(X[y]_{sc}))$: il existe $b'_{1}\in \mathfrak{z}(\bar{G}';{\mathbb R})^*$ tel que la fonction
 $$e^{-<b'_{1},Z_{2}>}\Delta(y)(exp(\bar{Y}),exp(X[y]_{sc}))$$
 soit localement constante. On peut identifier $b'_{1}$ \`a un \'el\'ement de $\mathfrak{z}(G'_{\epsilon};{\mathbb R})^*$, cf. [III] 5.1(1). En posant $b=b'_{1}-b_{0}$, on obtient que le rapport
 $$(2) \qquad e^{<b,Z>}\Delta_{1}(exp(Y)\epsilon_{1},exp(X[y])\eta[y]) \Delta(y)(exp(\bar{Y}),exp(X[y]_{sc}))^{-1}$$
 est localement constant.
  Il r\'esulte des d\'efinitions de [I] 2.8 que les termes $b_{1}$ et $b'_{1}$ ne d\'ependent ni des tores fix\'es, ni de l'\'el\'ement $y$ (les constructions se situent dans les groupes duaux et ne voient pas $y$).  L'\'el\'ement $b$ non plus. On peut calculer la fonction localement constante (2) de la m\^eme fa\c{c}on qu'en [W2] chapitre 10. On obtient qu'elle est constante, les tores \'etant fix\'es. Il reste \`a prouver que cette constante ne d\'epend pas des tores. Dans la notation de [W2], il s'agit de prouver l'\'egalit\'e $\Delta(T,\underline{T})=1$.  La m\'ethode de [W2] paragraphe 12 consistait \`a plonger la situation locale dans une situation d\'efinie sur un corps de nombres et \`a utiliser une formule de produit pour d\'eduire la valeur de $\Delta(T,\underline{T})$ de valeurs analogues dans une situation non ramifi\'ee. La m\^eme m\'ethode s'applique et ram\`ene l'\'egalit\'e \`a prouver \`a l'\'egalit\'e analogue dans une situation non ramifi\'ee sur un corps non-archim\'edien. La preuve dans ce cas a \'et\'e faite dans [W2]. Cela prouve (1). $\square$

  Le facteur exponentiel qui intervient en (1) affecte la description du transfert des distributions si celles-ci font intervenir des d\'erivations dans la direction centrale de $G'_{\epsilon}$. Mais ce n'est pas le cas si on se limite aux distributions provenant de $G'_{\epsilon,SC}$ via l'application $\iota^*_{G'_{\epsilon,SC},G_{\epsilon}}$ de  la section 3. En particulier, la relation [III] 5.2(7) reste valide.

 \bigskip
 
 \subsection{Localisation}
 Dans la suite de la  section, on consid\`ere un triplet $(G,\tilde{G},{\bf a})$ quasi-d\'eploy\'e et \`a torsion int\'erieure. Soit $\eta$ un \'el\'ement semi-simple de $ \tilde{G}({\mathbb R})$. Conform\'ement \`a la convention de [I] 5.10, nous notons   $$desc_{\eta}^{\tilde{G},*}:D_{g\acute{e}om}(G_{\eta}({\mathbb R}))\otimes Mes(G_{\eta}({\mathbb R}))^*\to D_{g\acute{e}om}(\tilde{G}({\mathbb R}))\otimes Mes(G({\mathbb R}))^*$$
 l'application lin\'eaire d\'eduite de la descente d'Harish-Chandra. Comme on l'a dit dans cette r\'ef\'erence, elle n'est d\'efinie que sur les \'el\'ements de  $D_{g\acute{e}om}(G_{\eta}({\mathbb R}))\otimes Mes(G_{\eta}({\mathbb R}))^*$ dont le support est assez voisin de l'origine. 
 Le groupe  $ Z_{G}(\eta;{\mathbb R})$ agit naturellement dans $D_{unip}(G_{\eta}({\mathbb R}))\otimes Mes(G_{\eta}({\mathbb R}))^*$ via son quotient fini  $ Z_{G}(\eta:{\mathbb R})/G_{\eta}({\mathbb R})$.  Notons ${\cal O}$ la classe de conjugaison de $\eta$ par $G({\mathbb R})$. Alors la restriction de $desc_{\eta}^{\tilde{G},*}$ \`a $D_{unip}(G_{\eta}({\mathbb R}))\otimes Mes(G_{\eta}({\mathbb R}))^*$ se factorise en
 $$(1) \qquad D_{unip}(G_{\eta}({\mathbb R}))\otimes Mes(G_{\eta}({\mathbb R}))^*\stackrel{p}{\to}(D_{unip}(G_{\eta}({\mathbb R}))\otimes Mes(G_{\eta}({\mathbb R}))^*)^{Z_{G}(\eta;{\mathbb R})}$$
 $$\stackrel{desc_{\eta}^{\tilde{G},*}}{\simeq}D_{g\acute{e}om}({\cal O})\otimes Mes(G({\mathbb R}))^*\subset D_{g\acute{e}om}(\tilde{G}({\mathbb R}))\otimes Mes(G({\mathbb R}))^*,$$
 o\`u l'exposant $Z_{G}(\eta;{\mathbb R})$ d\'esigne selon l'usage le sous-espace des invariants et o\`u $p$ est la projection naturelle sur ce sous-espace.
 
 Supposons $G_{\eta}$ quasi-d\'eploy\'e. On a de m\^eme une application lin\'eaire  
  $$desc_{\eta}^{st,\tilde{G},*}:D^{st}_{g\acute{e}om}(G_{\eta}({\mathbb R}))\otimes Mes(G_{\eta}({\mathbb R}))^*\to D^{st}_{g\acute{e}om}(\tilde{G}({\mathbb R}))\otimes Mes(G({\mathbb R}))^*,$$
  cf. [I] 5.10.
 Posons $\Xi_{\eta}=Z_{G}(\eta)/G_{\eta}$. Le groupe $\Xi_{\eta}^{\Gamma_{{\mathbb R}}}$ agit naturellement sur  $D^{st}_{unip}(G_{\eta}({\mathbb R}))\otimes Mes(G_{\eta}({\mathbb R}))^*$. Notons ${\cal O}^{st}$ la classe de conjugaison stable de $\eta$. Alors la restriction $desc_{\eta}^{st,\tilde{G},*}$ \`a $D^{st}_{unip}(G_{\eta}({\mathbb R}))\otimes Mes(G_{\eta}({\mathbb R}))^*$ se factorise en
 $$(2) \qquad D^{st}_{unip}(G_{\eta}({\mathbb R}))\otimes Mes(G_{\eta}({\mathbb R}))^*\stackrel{p^{st}}{\to}(D^{st}_{unip}(G_{\eta}({\mathbb R}))\otimes Mes(G_{\eta}({\mathbb R}))^*)^{\Xi_{\eta}^{\Gamma_{{\mathbb R}}}}$$
 $$\stackrel{desc_{\eta}^{st,\tilde{G},*}}{\simeq}D^{st}_{g\acute{e}om}({\cal O}^{st})\otimes Mes(G({\mathbb R}))^*\subset D^{st}_{g\acute{e}om}(\tilde{G}({\mathbb R}))\otimes Mes(G({\mathbb R}))^*,$$
 avec des notations similaires \`a celles ci-dessus. Notons ${\cal Y}(\eta)$ l'ensemble des $y\in G$ tels que  $y\sigma(y)^{-1}\in G_{\eta}$ pour tout $\sigma\in \Gamma_{{\mathbb R}}$. Fixons un ensemble de repr\'esentants $\dot{{\cal Y}}(\eta)$ de l'ensemble de doubles classes $G_{\eta}\backslash {\cal Y}(\eta)/G({\mathbb R})$. Pour $y\in {\cal Y}(\eta)$, on pose $\eta[y]=y^{-1}\eta y$. L'application $ad_{y}$ se restreint en un torseur int\'erieur de $G_{\eta[y]}$ sur sa forme quasi-d\'eploy\'ee $G_{\eta}$, gr\^ace auquel on peut transf\'erer une distribution stable sur $G_{\eta}({\mathbb R})$ en une distribution sur $G_{\eta[y]}({\mathbb R})$. On note $transfert_{y}$ cette application.  Les applications $desc^{\tilde{G},*}_{\eta}$ et $desc^{st,\tilde{G},*}_{\eta}$ sont reli\'ees par l'\'egalit\'e
 $$(3) \qquad desc^{st,\tilde{G},*}_{\eta}=\sum_{y\in \dot{{\cal Y}}(\eta)}desc_{\eta[y]}^{\tilde{G},*}\circ transfert_{y},$$
 cf. [I] 5.10.
 
 \bigskip
 
 \subsection{Localisation des espaces $D_{tr-orb}({\cal O})$}
 Soit ${\cal O}$ une classe de conjugaison stable d'\'el\'ements semi-simples dans $\tilde{G}({\mathbb R})$. Fixons $\eta\in {\cal O}$ tel que $G_{\eta}$ soit quasi-d\'eploy\'e. On introduit un ensemble $\dot{{\cal Y}}(\eta)$ comme dans le paragraphe pr\'ec\'edent. En g\'en\'eral, l'ensemble $\{\eta[y]; y\in \dot{{\cal Y}}(\eta)\}$ n'est pas un ensemble de repr\'esentants des classes de conjugaison par $G({\mathbb R})$ dans ${\cal O}$. Il est plus gros. On peut toutefois supposer que si $\eta[y]$ est conjugu\'e \`a $\eta[y']$, alors ces deux points sont \'egaux. On peut alors fixer un sous-ensemble $\dot{{\cal X}}(\eta)\subset \dot{{\cal Y}}(\eta)$ de sorte que $\{\eta[y]; y\in \dot{{\cal X}}(\eta)\}$ soit un ensemble de repr\'esentants des classes de conjugaison par $G({\mathbb R})$ dans ${\cal O}$.

 \ass{Lemme}{(i)   L'espace $D_{tr-orb}({\cal O})\otimes Mes(G({\mathbb R}))^*$ est contenu dans la somme sur $y\in \dot{{\cal X}}(\eta)$ des images par $desc_{\eta[y]}^{\tilde{G},*}$ de $D_{tr-unip}(G_{\eta[y]}({\mathbb R}))\otimes Mes(G_{\eta}({\mathbb R}))^*$.
 
 (ii) L'espace $D^{st}_{tr-orb}({\cal O})\otimes Mes(G({\mathbb R}))^*$ est contenu dans l'image par $desc_{\eta}^{st,\tilde{G},*}$ de $D^{st}_{tr-unip}(G_{\eta}({\mathbb R}))\otimes Mes(G_{\eta}({\mathbb R}))^*$.}

 Preuve. Oublions les espaces de mesures. Consid\'erons (i). Il est clair que $D_{orb}({\cal O})$ est contenu dans la somme sur $y\in \dot{{\cal X}}(\eta)$ des images par $desc_{\eta[y]}^{\tilde{G},*}$ de $D_{orb,unip}(G_{\eta[y]}({\mathbb R}))$. Consid\'erons une donn\'ee endoscopique ${\bf G}'=(G',{\cal G}',s)$ de $(G,\tilde{G})$ avec $G'\not=G$. Soit $\boldsymbol{\gamma}\in D_{tr-orb}({\cal O})$, supposons qu'il existe $\boldsymbol{\delta}\in D^{st}_{tr-orb}({\bf G}')$ tel que $\boldsymbol{\gamma}=transfert(\boldsymbol{\delta})$. On veut prouver que $\boldsymbol{\gamma}$ appartient \`a la somme des images par $desc_{\eta[y]}^{\tilde{G},*}$ de $D_{tr-unip}(G_{\eta[y]}({\mathbb R}))$. Par lin\'earit\'e, on peut fixer une classe de conjugaison stable ${\cal O}'$ dans $G'({\mathbb R})$ qui se transf\`ere en ${\cal O}$ et supposer que $\boldsymbol{\delta}\in D^{st}_{tr-orb}({\bf G}',{\cal O}')$. On fixe $\epsilon\in {\cal O}'$ tel que $G'_{\epsilon}$ soit quasi-d\'eploy\'e. On fixe des donn\'ees auxiliaires $G'_{1}$,...$\Delta_{1}$ pour ${\bf G}'$ et un point $\epsilon_{1}\in G'_{1}({\mathbb R})$ au-dessus de $\epsilon$. On applique (ii) par r\'ecurrence (plus exactement, on applique une forme de (ii) adapt\'ee comme toujours \`a la situation des donn\'ees auxiliaires; on laisse cette adaptation au lecteur). On peut donc fixer $\boldsymbol{\delta}_{\epsilon_{1}}\in D_{tr-unip}^{st}(G'_{1,\epsilon_{1}}({\mathbb R}))$  tel que $\boldsymbol{\delta}=desc_{\epsilon_{1}}^{st,\tilde{G}'_{1},*}(\boldsymbol{\delta}_{\epsilon_{1}})$. Appliquons le lemme 3.3(ii): on peut fixer $\boldsymbol{\delta}_{SC}\in D_{tr-unip}^{st}(G'_{\epsilon,SC}({\mathbb R}))$ tel que $\boldsymbol{\delta}_{\epsilon_{1}}=\iota^*_{G'_{\epsilon,SC},G'_{1,\epsilon_{1}}}(\boldsymbol{\delta}_{SC})$. Alors $\boldsymbol{\gamma}$ est calcul\'e par la formule [III] 5.2(7). Cette formule se simplifie dans notre situation quasi-d\'eploy\'ee et \`a torsion int\'erieure. Les constantes $c(y)$ valent $1$. Surtout, le groupe $\bar{G}'_{SC}$, qui en g\'en\'eral est en situation d'endoscopie non standard avec $G'_{\epsilon,SC}$, est ici \'egal \`a ce groupe. On obtient le r\'esultat suivant. Pour tout $y\in \dot{{\cal Y}}(\eta)$, on note $\boldsymbol{\delta}[y]$ l'image de $\iota^*_{\bar{G}'_{SC},\bar{G}'}(\boldsymbol{\delta}_{SC})$ par l'application  $transfert_{y}$ \`a $G_{\eta[y],SC}({\mathbb R})$, avec la convention que ce transfert est nul si ${\bf \bar{G}}'$ n'est pas relevante pour $G_{\eta[y],SC}$. Alors
 $$(1) \qquad \boldsymbol{\gamma}=\sum_{y\in \dot{{\cal Y}}(\eta)}d(y)desc_{\eta[y]}^{\tilde{G},*}\circ \iota^*_{G_{\eta[y],SC},G_{\eta[y]}}(\boldsymbol{\delta}[y]).$$
 Le lemme 3.3(ii) assure que $\iota^*_{\bar{G}'_{SC},\bar{G}'}(\boldsymbol{\delta}_{SC})\in D_{tr-unip}^{st}(\bar{G}'({\mathbb R}))$. On a alors $\boldsymbol{\delta}[y]\in D_{tr-unip}(G_{\eta[y],SC}({\mathbb R}))$. Mais, parce que $G_{\eta[y]}$ n'est pas quasi-d\'eploy\'e,  le lemme 3.3 n'assure pas que $\iota^*_{G_{\eta[y],SC},G_{\eta[y]}}(\boldsymbol{\delta}[y])$ appartient \`a $D_{tr-unip}(G_{\eta[y]}({\mathbb R}))$. Pour le voir, il faut utiliser encore une simplification due \`a notre situation quasi-d\'eploy\'ee et \`a torsion int\'erieure. Non seulement $G'_{\epsilon,SC}$ s'identifie au groupe $\bar{G}'_{SC}$, mais $G'_{\epsilon}$ lui-m\^eme s'identifie au groupe d'une donn\'ee endoscopique de $\bar{G}$, ou encore de $G_{\eta[y]}$. Notons ${\bf G}''$ cette donn\'ee. La donn\'ee pr\'ec\'edente ${\bf G}'$ est d\'eduite de ${\bf G}''$ par la construction de la preuve du lemme 3.3. Posons $\boldsymbol{\delta}'=\iota^*_{\bar{G}'_{SC},\bar{G}'}(\boldsymbol{\delta}_{SC})$. On peut consid\'erer que c'est un \'el\'ement de $D_{tr-unip}^{st}({\bf G}')$. On a alors
$$\iota^*_{G_{\eta[y],SC},G_{\eta[y]}}(\boldsymbol{\delta}[y])=\iota^*_{G_{\eta[y],SC},G_{\eta[y]}}\circ transfert_{y}(\boldsymbol{\delta}')=transfert_{y}\circ\iota^*_{{\bf G}',{\bf G}''}(\boldsymbol{\delta}').$$
Parce que $G'_{\epsilon}$ est quasi-d\'eploy\'e, le lemme 3.3(ii) assure que  $\iota^*_{{\bf G}',{\bf G}''}(\boldsymbol{\delta}')$ appartient \`a $D_{tr-unip}^{st}({\bf G}'')$. Alors  $\iota^*_{G_{\eta[y],SC},G_{\eta[y]}}(\boldsymbol{\delta}[y])$ appartient \`a $D_{tr-unip}(G_{\eta[y]}({\mathbb R}))$. La formule (1) se r\'ecrit
$$ \boldsymbol{\gamma}=\sum_{y\in \dot{{\cal X}}(\eta)}desc_{\eta[y]}^{\tilde{G},*}(\boldsymbol{\gamma}[y]),$$
o\`u
$$\boldsymbol{\gamma}[y]=\sum_{y'\in \dot{{\cal Y}}(\eta); \eta[y']=\eta[y]}d(y') \iota^*_{G_{\eta[y'],SC},G_{\eta[y']}}(\boldsymbol{\delta}[y']).$$
Cela assure que $\boldsymbol{\gamma}$ appartient \`a l'espace indiqu\'e en (i) et cela ach\`eve la preuve de cette assertion.

 Soit maintenant $\boldsymbol{\delta}\in  D^{st}_{tr-orb}({\cal O})$. D'apr\`es ce que l'on vient de prouver, on peut \'ecrire 
$$(2) \qquad  \boldsymbol{\delta}=\sum_{y\in \dot{{\cal X}}(\eta)}desc_{\eta[y]}^{\tilde{G},*}(\boldsymbol{\gamma}[y]),$$ 
avec $\boldsymbol{\gamma}[y]\in D_{tr-unip}(G_{\eta[y]}({\mathbb R}))$ pour tout $y$. Appliquons 4.2(1) en rempla\c{c}ant $\eta$ par $\eta[y]$. On note $p_{y}$ la projection qui intervient dans cette relation. Cette application consiste \`a moyenner sur un groupe d'automorphismes de $G_{\eta[y]}$. Comme on l'a d\'ej\`a dit, l'espace $D_{tr-unip}(G_{\eta[y]}({\mathbb R}))$ est invariant par automorphismes, donc aussi par $p_{y}$. 
On peut remplacer $\boldsymbol{\gamma}[y]$ par $p_{y}(\boldsymbol{\gamma}[y])$, qui appartient encore \`a  $D_{tr-unip}(G_{\eta[y]}({\mathbb R}))$. D'autre part, d'apr\`es 4.2(2), on peut \'ecrire
$$(3) \qquad \boldsymbol{\delta}=desc_{\eta}^{st,\tilde{G},*}(\boldsymbol{\delta}_{\eta}),$$
avec $\boldsymbol{\delta}_{\eta}\in D_{unip}^{st}(G_{\eta}({\mathbb R}))^{\Xi_{\eta}^{\Gamma_{{\mathbb R}}}}$. Appliquons 4.2(3). On obtient
$$(4) \qquad \boldsymbol{\delta}= \sum_{y\in \dot{{\cal X}}(\eta)}desc_{\eta[y]}^{\tilde{G},*}(\boldsymbol{\gamma}'[y]),$$
o\`u
$$\boldsymbol{\gamma}'[y]=\sum_{y'\in \dot{{\cal Y}}(\eta); \eta[y']=\eta[y]}transfert_{y'}(\boldsymbol{\delta}_{\eta}).$$
Fixons $y'\in \dot{{\cal Y}}(\eta)$. Pour $g\in Z_{G}(\eta[y'],{\mathbb R})$, on v\'erifie que $ad_{g}\circ transfert_{y'}=transfert_{y'}\circ ad_{y'g(y')^{-1}}$. On a $y'g(y')^{-1}\in Z_{G}(\eta)$ et on v\'erifie que l'image de cet \'el\'ement dans $\Xi_{\eta}$ est fixe par $\Gamma_{{\mathbb R}}$. Donc $\boldsymbol{\delta}_{\eta}$ est  fixe par $ad_{y'g(y')^{-1}}$, donc $transfert_{y'}(\boldsymbol{\delta}_{\eta})$ est fixe par $Z_{G}(\eta[y'],{\mathbb R})$. Il en r\'esulte que, pour tout $y\in \dot{{\cal X}}(\eta)$, $\boldsymbol{\gamma}'[y]$ est invariant par $Z_{G}(\eta[y],{\mathbb R})$. Puisque $desc_{\eta[y]}^{\tilde{G},*}$ est injectif sur les distributions fix\'ees par ce groupe, les \'egalit\'es (2) et (4) entra\^{\i}nent que $\boldsymbol{\gamma}[y]=\boldsymbol{\gamma}'[y]$ pour tout $y$. Appliquons cela \`a $y=1$ (on peut supposer que $1$ appartient \`a notre syst\`eme de repr\'esentants $\dot{{\cal X}}(\eta)$). Cela entra\^{\i}ne $\boldsymbol{\gamma}'[1] \in D_{tr-unip}(G_{\eta}({\mathbb R}))$, c'est-\`a-dire
$$\sum_{y'\in \dot{{\cal Y}}(\eta); \eta[y']=\eta}transfert_{y'}(\boldsymbol{\delta}_{\eta})\in D_{tr-unip}(G_{\eta}({\mathbb R})).$$
Un \'el\'ement $y'\in {\cal Y}(\eta)$ tel que $\eta[y']=\eta$ appartient \`a $Z_{G}(\eta)$ et son image dans $\Xi_{\eta}$ est fixe par $\Gamma_{{\mathbb R}}$. On voit que $transfert_{y'}$ n'est autre que l'action de l'\'el\'ement de $\Xi_{\eta}^{\Gamma_{{\mathbb R}}}$ ainsi d\'efini. Puisque $\boldsymbol{\delta}_{\eta}$ est fixe par cette action, on obtient  $transfert_{y'}(\boldsymbol{\delta}_{\eta})=\boldsymbol{\delta}_{\eta}$. Le membre de gauche de la relation pr\'ec\'edente n'est autre que $\boldsymbol{\delta}_{\eta}$ multipli\'e par le nombre d'\'el\'ements de l'ensemble de sommation. Donc $\boldsymbol{\delta}_{\eta}\in D_{tr-unip}(G_{\eta}({\mathbb R}))$. Puisque c'est une distribution stable, on a $\boldsymbol{\delta}_{\eta}\in D^{st}_{tr-unip}(G_{\eta}({\mathbb R}))$. Alors l'assertion (ii) r\'esulte de (3). $\square$

\bigskip

\subsection{Un r\'esultat d'induction}
Soit ${\cal O}^{\tilde{G}}$ une classe de conjugaison stable semi-simple dans $\tilde{G}({\mathbb R})$. Fixons $\eta\in {\cal O}^{\tilde{G}}$ tel que $G_{\eta}$ soit quasi-d\'eploy\'e. Notons $\tilde{M}$ le commutant de $A_{G_{\eta}}$ dans $\tilde{M}$. C'est un espace de Levi de $\tilde{G}$, on a $\eta\in\tilde{M}({\mathbb R})$ et $M_{\eta}=G_{\eta}$. On note ${\cal O}$ la classe de conjugaison stable de $\eta$ dans $\tilde{M}({\mathbb R})$.

\ass{Lemme}{Les applications d'induction
$$D_{tr-orb}({\cal O})\otimes Mes(M({\mathbb R}))^*\to D_{tr-orb}({\cal O}^{\tilde{G}})\otimes Mes(G({\mathbb R}))^*$$
et 
$$D^{st}_{tr-orb}({\cal O})\otimes Mes(M({\mathbb R}))^*\to D^{st}_{tr-orb}({\cal O}^{\tilde{G}})\otimes Mes(G({\mathbb R}))^*$$
sont surjectives. }

Preuve. On oublie les espaces de mesures.  L'application 
$$(1) \qquad M_{\eta}\backslash {\cal Y}^M(\eta)/M({\mathbb R})\to G_{\eta}\backslash {\cal Y}(\eta)/G({\mathbb R})$$
est toujours injective et l'hypoth\`ese $M_{\eta}=G_{\eta}$ entra\^{\i}ne qu'elle est surjective, cf. [I] lemme 5.12.  Un syst\`eme de repr\'esentants $\dot{{\cal Y}}^M(\eta)$ du premier  ensemble de doubles classes est donc aussi un tel syst\`eme pour le second ensemble. Un \'el\'ement $\boldsymbol{\gamma}\in D_{orb}({\cal O}^{\tilde{G}})$ est combinaison lin\'eaire d'int\'egrales orbitales dans $\tilde{G}({\mathbb R})$ associ\'ees \`a des \'el\'ements $u \eta[y]$ o\`u $y\in \dot{{\cal Y}}^M(\eta)$ et $u$ est un \'el\'ement unipotent de $ G_{\eta[y]}({\mathbb R})=M_{\eta[y]}({\mathbb R})$. Il est induit de la m\^eme  combinaison lin\'eaire d'int\'egrales orbitales dans $\tilde{M}({\mathbb R})$ associ\'ees aux m\^emes \'el\'ements. Cette combinaison lin\'eaire appartient \`a $D_{orb}({\cal O})$.  Soit maintenant ${\bf G}'=(G',{\cal G}',s)$ une donn\'ee endoscopique elliptique et relevante de $(G,\tilde{G})$ avec $G'\not=G$. Soit $\boldsymbol{\gamma}\in D_{tr-orb}({\cal O})$ tel qu'il existe $\boldsymbol{\delta}\in D^{st}_{tr-orb}({\bf G}')$ de sorte que $\boldsymbol{\gamma}=transfert(\boldsymbol{\delta})$. On veut prouver que $\boldsymbol{\gamma}$ est induit d'un \'el\'ement de $D_{tr-orb}({\cal O}^{\tilde{M}})$. Par lin\'earit\'e, on peut fixer une classe de conjugaison stable semi-simple ${\cal O}^{\tilde{G}'}$ dans $\tilde{G}'({\mathbb R})$ et supposer $\boldsymbol{\delta}\in D_{tr-orb}({\bf G}',{\cal O}^{\tilde{G}'})$. On fixe $\epsilon\in {\cal O}^{\tilde{G}'}$ avec $G'_{\epsilon}$ quasi-d\'eploy\'e et on fixe un diagramme $(\epsilon,B',T',B,T,\eta)$. Notons $\tilde{R}'$ le commutant de $A_{G'_{\epsilon}}$ dans $\tilde{G}'$. Puisque $(G,\tilde{G})$ est quasi-d\'eploy\'e et \`a torsion int\'erieure, $\tilde{R}'$ d\'etermine   un espace de Levi $\tilde{R}$ de $\tilde{G}$ de sorte que $T\subset R$ et que  l'isomorphisme $T'\simeq T$ se restreigne en un isomorphisme  $A_{R'}\simeq A_{R}$. Alors $\tilde{R}'$ est l'espace d'une donn\'ee endoscopique elliptique ${\bf R}'$ de $(R,\tilde{R})$.  L'\'el\'ement $\eta$ appartient \`a $\tilde{R}({\mathbb R})$. On a les relations $A_{M}=A_{G_{\eta}}\subset A_{G'_{\epsilon}}=A_{R'}=A_{R}$, donc $\tilde{R}\subset \tilde{M}$. Notons ${\cal O}'$ la classe de conjugaison stable de $\epsilon$ dans $\tilde{R}'({\mathbb R})$.  En appliquant par r\'ecurrence la deuxi\`eme assertion du lemme (adapt\'ee \`a la situation endoscopique), il existe $\boldsymbol{\delta}_{{\bf R}'}\in D^{st}_{tr-orb}({\bf R}',{\cal O}')$ tel que $\boldsymbol{\delta}=(\boldsymbol{\delta}_{{\bf R}'})^{{\bf G}'}$. Posons $\boldsymbol{\gamma}_{\tilde{R}}=transfert(\boldsymbol{\delta}_{{\bf R}'})$. Cette distribution appartient \`a $D_{tr-orb}({\cal O}_{\tilde{R}})$, o\`u ${\cal O}_{\tilde{R}}$ est la classe de conjugaison stable de $\eta$ dans $\tilde{R}({\mathbb R})$. Par commutation du transfert \`a l'induction, $\boldsymbol{\gamma}=(\boldsymbol{\gamma}_{\tilde{R}})^{\tilde{G}}=(\boldsymbol{\gamma}_{\tilde{M}})^{\tilde{G}}$, o\`u $\boldsymbol{\gamma}_{\tilde{M}}=(\boldsymbol{\gamma}_{\tilde{R}})^{\tilde{M}}$. Puisque $\boldsymbol{\gamma}_{\tilde{M}}\in D_{tr-orb}({\cal O})$, cela d\'emontre l'assertion cherch\'ee, d'o\`u la premi\`ere assertion du lemme. 

Soit maintenant $\boldsymbol{\delta}\in D^{st}_{tr-orb}({\cal O}^{\tilde{G}})$. D'apr\`es ce que l'on vient de prouver, il existe $\boldsymbol{\gamma}_{\tilde{M}}\in D_{tr-orb}({\cal O})$ tel que $\boldsymbol{\delta}=(\boldsymbol{\gamma}_{\tilde{M}})^{\tilde{G}}$. D'apr\`es le lemme [I] 5.12(iii), il existe aussi $\boldsymbol{\delta}_{\tilde{M}}\in D^{st}_{g\acute{e}om}({\cal O})$ tel que $\boldsymbol{\delta}=(\boldsymbol{\delta}_{\tilde{M}})^{\tilde{G}}$. Introduisons le groupe $N$ des \'el\'ements $n\in G({\mathbb R})$ tels que $ad_{n}$ conserve $\tilde{M}$ et ${\cal O}$. Il agit via son quotient fini $N/(N\cap M({\mathbb R}))$ sur $D_{g\acute{e}om}({\cal O})$. Cette action conserve (par transport de structure) les sous-espaces $D_{tr-orb}({\cal O})$ et $D^{st}_{g\acute{e}om}({\cal O})$. Evidemment, l'induction est insensible \`a l'action de ce groupe, c'est-\`a-dire que $(ad_{n}(\boldsymbol{\gamma}_{1}))^{\tilde{G}}=\boldsymbol{\gamma}_{1}^{\tilde{G}}$ pour tout $\boldsymbol{\gamma}_{1}\in D_{g\acute{e}om}({\cal O})$ et tout $n\in N$. On peut donc remplacer les \'el\'ements $\boldsymbol{\gamma}_{\tilde{M}}$ et $\boldsymbol{\delta}_{\tilde{M}}$ par leur moyenne sous l'action de $N$ sans changer les propri\'et\'es pr\'ec\'edentes. Autrement dit, on peut supposer $\boldsymbol{\gamma}_{\tilde{M}}$ et $\boldsymbol{\delta}_{\tilde{M}}$ invariantes par $N$. D'apr\`es le lemme [I] 5.12(iii), 
 l'induction induit un isomorphisme de l'espace des invariants $D_{g\acute{e}om}({\cal O})^N$ sur $D_{g\acute{e}om}({\cal O}^{\tilde{G}})$.
Cela entra\^{\i}ne  $\boldsymbol{\gamma}_{\tilde{M}}=\boldsymbol{\delta}_{\tilde{M}}$. Cet  \'el\'ement appartient donc \`a  $D^{st}_{tr-orb}({\cal O})$. Cela prouve que $\boldsymbol{\delta}$ est induit d'un \'el\'ement de cet espace, ce qui d\'emontre la seconde assertion du lemme.

\bigskip

\subsection{D\'efinition des termes $\rho_{J}^{\tilde{G}}$ et $\sigma_{J}^{\tilde{G}}$, premier cas}
On fixe  pour la suite de la section un syst\`eme de fonctions $B$ comme en [II] 1.9. Soient $\tilde{M}$ un espace de Levi de $\tilde{G}$ et ${\cal O}$ une classe de conjugaison stable semi-simple dans $\tilde{M}({\mathbb R})$. On note ${\cal O}^{\tilde{G}}$ la classe de conjugaison stable dans $\tilde{G}({\mathbb R})$ qui la contient.  
  Fixons $\eta\in {\cal O}$ avec $M_{\eta}$ quasi-d\'eploy\'e, introduisons le sous-espace de Levi $\tilde{R}$ de $\tilde{M}$ tel que $A_{R}=A_{M_{\eta}}$. On note ${\cal O}_{\tilde{R}}$ la classe de conjugaison stable de $\eta$ dans $\tilde{R}({\mathbb R})$. On suppose dans ce paragraphe

(1)  $\tilde{R}\not=\tilde{M}$.

Soit $J\in {\cal J}_{\tilde{M}}^{\tilde{G}}(B_{{\cal O}})$. Pour $\boldsymbol{\gamma}\in D_{tr-orb}({\cal O})\otimes Mes(M({\mathbb R}))^*$, fixons gr\^ace au lemme 4.4 (appliqu\'e avec $\tilde{G}$ et $\tilde{M}$ remplac\'es par $\tilde{M}$ et $\tilde{R}$) un \'el\'ement $\boldsymbol{\gamma}_{\tilde{R}}\in D_{tr-orb}({\cal O}_{\tilde{R}})\otimes Mes(R({\mathbb R}))^*$ tel que $\boldsymbol{\gamma}=(\boldsymbol{\gamma}_{\tilde{R}})^{\tilde{M}}$. Pour $a\in A_{M}({\mathbb R})$ en position g\'en\'erale et proche de $1$, posons
$$(2) \qquad \rho_{J}^{\tilde{G}}(\boldsymbol{\gamma},a)=\sum_{\tilde{L}\in {\cal L}(\tilde{R}); J\in {\cal J}_{\tilde{R}}^{\tilde{L}}(B_{{\cal O}_{\tilde{R}}})}d_{\tilde{R}}^{\tilde{G}}(\tilde{M},\tilde{L})\rho_{J}^{\tilde{L}}(\boldsymbol{\gamma}_{\tilde{R}},a)^{\tilde{M}}.$$
Tous les termes sont d\'efinis par r\'ecurrence, en vertu de l'hypoth\`ese (1).
On a

(3) ce terme  ne d\'epend pas du choix de $\boldsymbol{\gamma}_{\tilde{R}}$. 

Preuve. On a introduit un groupe $N\subset M({\mathbb R})$ dans la preuve du lemme 4.4. Ce groupe agit sur $D_{tr-orb}({\cal O}_{\tilde{R}})\otimes Mes(R({\mathbb R}))^*$ et aussi sur ${\cal L}(\tilde{R})$. On v\'erifie que, pour $n\in N$, $ad_{n}$ permute les espaces de Levi $\tilde{L}$ tels que $d_{\tilde{R}}^{\tilde{G}}(\tilde{M},\tilde{L})\not=0$ et $J\in {\cal J}_{\tilde{R}}^{\tilde{L}}(B_{{\cal O}_{\tilde{R}}})$. Par transport de structure, il est plus ou moins clair que
 $$\rho_{J}^{ad_{n}(\tilde{L})}(ad_{n}(\boldsymbol{\gamma}_{\tilde{R}}),a)=ad_{n}(\rho_{J}^{\tilde{L}}(\boldsymbol{\gamma}_{\tilde{R}},a)).$$
 Puisque l'induction de $\tilde{R}$ \`a $\tilde{M}$ est insensible \`a la composition avec $ad_{n}$, on obtient que remplacer $\boldsymbol{\gamma}_{\tilde{R}}$ par $ad_{n}(\boldsymbol{\gamma}_{\tilde{R}})$ ne change pas le membre de droite de (2). On peut alors aussi bien remplacer $\boldsymbol{\gamma}_{\tilde{R}}$ par sa projection naturelle sur l'espace des invariants par $N$. Mais, d'apr\`es le lemme [I] 5.12(iii), cette projection est uniquement d\'etermin\'ee par $\boldsymbol{\gamma}$. L'assertion (3) s'ensuit. $\square$
 
 La formule  (2) d\'efinit une application lin\'eaire
 $$\rho_{J}^{\tilde{G}}:D_{tr-orb}({\cal O})\otimes Mes(M({\mathbb R}))^*\to U_{J}\otimes (D_{g\acute{e}om}({\cal O})\otimes Mes(M({\mathbb R}))^*)/Ann_{{\cal O}}^{\tilde{G}}.$$
 On va prouver qu'elle v\'erifie les conditions (1) et (3) de 2.4, c'est-\`a-dire:  
 
 (4)  cette d\'efinition co\"{\i}ncide avec celle de [II] 3.4 pour $\boldsymbol{\gamma}\in D_{orb}({\cal O})\otimes Mes(M({\mathbb R}))^*$;
 
 (5) soit ${\bf M}'$ une donn\'ee endoscopique elliptique et relevante de $(M,\tilde{M})$ avec $M'\not=M$; soit ${\cal O}'$ une classe de conjugaison stable semi-simple dans $\tilde{M}'({\mathbb R})$ correspondant \`a ${\cal O}$;  soit $\boldsymbol{\delta}\in D_{tr-orb}^{st}({\bf M}',{\cal O}')\otimes Mes(M'({\mathbb R}))^*$; alors $\rho_{J}^{\tilde{G},{\cal E}}({\bf M}',\boldsymbol{\delta})=\rho_{J}^{\tilde{G}}(transfert(\boldsymbol{\delta}))$. 
 
 Dans la situation de (4), on peut choisir $\boldsymbol{\gamma}_{\tilde{R}}\in D_{orb}({\cal O}_{\tilde{R}})\otimes Mes(R({\mathbb R}))^*$. L'\'egalit\'e cherch\'ee r\'esulte alors du lemme [II] 3.10. Dans la situation de (5), on voit en reprenant la preuve du lemme 4.4 que l'on peut introduire une donn\'ee endoscopique ${\bf R}'$ de $(R,\tilde{R})$, qui est une donn\'ee de Levi de ${\bf M}'$, une classe de conjugaison stable ${\cal O}'_{\tilde{R}'}$ dans $\tilde{R}'({\mathbb R})$ correspondant \`a ${\cal O}'$ et \`a ${\cal O}_{\tilde{R}}$, et un \'el\'ement $\boldsymbol{\delta}_{{\bf R}'}\in D_{tr-orb}^{st}({\bf R}',{\cal O}'_{\tilde{R}'})$, de sorte que $\boldsymbol{\delta}=(\boldsymbol{\delta}_{{\bf R}'})^{{\bf M}'}$. D'apr\`es l'analogue de la relation 2.5(8) (qui est valide d'apr\`es l'hypoth\`ese $M'\not=M$), on a alors
 $$(6) \qquad \rho_{J}^{\tilde{G},{\cal E}}({\bf M}',\boldsymbol{\delta},a)=\sum_{\tilde{L}\in {\cal L}(\tilde{R}); J\in {\cal J}_{\tilde{R}}^{\tilde{L}}(B_{{\cal O}_{\tilde{R}}})}d_{\tilde{R}}^{\tilde{G}}(\tilde{M},\tilde{L})\rho_{J}^{\tilde{L},{\cal E}}({\bf R}',\boldsymbol{\delta}_{{\bf R}'},a)^{\tilde{M}}.$$
 En utilisant les hypoth\`eses de r\'ecurrence, on a
 $$\rho_{J}^{\tilde{L},{\cal E}}({\bf R}',\boldsymbol{\delta}_{{\bf R}'},a)=\rho_{J}^{\tilde{L}}(\boldsymbol{\gamma}_{\tilde{R}},a)$$
 o\`u $\boldsymbol{\gamma}_{\tilde{R}}=transfert(\boldsymbol{\delta}_{{\bf R}'})$. Mais $\boldsymbol{\gamma}_{\tilde{R}}$ est un \'el\'ement de $D_{tr-orb}({\cal O}_{\tilde{R}})\otimes Mes(R({\mathbb R}))^*$ tel que $(\boldsymbol{\gamma}_{\tilde{R}})^{\tilde{M}}=transfert(\boldsymbol{\delta})$. Alors le membre de droite de (6) co\"{\i}ncide avec celui de (2) pour $\boldsymbol{\gamma}=transfert(\boldsymbol{\delta})$. Cela prouve (5).

 On d\'efinit une application lin\'eaire 
 $$\sigma_{J}^{\tilde{G}}:D_{tr-orb}^{st}({\cal O})\otimes Mes(M({\mathbb R}))^*\to U_{J}\otimes (D_{g\acute{e}om}({\cal O})\otimes Mes(M({\mathbb R}))^*)/Ann_{{\cal O}}^{\tilde{G}}$$
 par la formule habituelle  2.4(13). On doit prouver qu'elle prend ses valeurs dans $U_{J}\otimes (D_{g\acute{e}om}^{st}({\cal O})\otimes Mes(M({\mathbb R}))^*)/Ann_{{\cal O}}^{st,\tilde{G}}$. 
  Soit $\boldsymbol{\delta}\in D_{tr-orb}^{st}({\cal O})\otimes Mes(M({\mathbb R}))^*$. D'apr\`es le lemme 4.4, on peut fixer $\boldsymbol{\delta}_{\tilde{R}}\in D_{tr-orb}^{st}({\cal O}_{\tilde{R}})\otimes Mes(R({\mathbb R}))^*$ de sorte que $\boldsymbol{\delta}=(\boldsymbol{\delta}_{\tilde{R}})^{\tilde{M}}$. A partir des relations d'induction d\'ej\`a connues, en particulier la formule (2), on prouve formellement que
  $$(7) \qquad \sigma_{J}^{\tilde{G}}(\boldsymbol{\delta},a)=\sum_{\tilde{L}\in {\cal L}(\tilde{R}); J\in {\cal J}_{\tilde{R}}^{\tilde{L}}(B_{{\cal O}_{\tilde{R}}})}e_{\tilde{R}}^{\tilde{G}}(\tilde{M},\tilde{L})\sigma_{J}^{\tilde{L}}(\boldsymbol{\delta}_{\tilde{R}},a)^{\tilde{M}}$$
  pour $a\in A_{M}({\mathbb R})$ en position g\'en\'erale et proche de $1$. Les distributions du membre de droite sont stables modulo $Ann_{{\cal O}}^{\tilde{G}}$, d'o\`u l'assertion cherch\'ee. 
  
  On a ainsi v\'erifi\'e les conditions (1) et (2) de 2.7 et on a vu dans ce paragraphe qu'elles suffisaient \`a r\'ealiser le programme de 2.4. Cela  valide ce programme sous l'hypoth\`ese (1).
  
  \bigskip
  
  \subsection{ D\'efinition des termes $\rho_{J}^{\tilde{G}}$ et $\sigma_{J}^{\tilde{G}}$, deuxi\`eme cas}
  On fixe encore $\eta\in {\cal O}$ avec $M_{\eta}$ quasi-d\'eploy\'e et on suppose maintenant 
  que  $A_{M_{\eta}}=A_{M}$. 
  
  Du syst\`eme de fonctions $B$ se d\'eduit une fonction $B_{\eta}$ sur le syst\`eme de racines de $G_{\eta}$. Pour simplifier, on note encore $B$ cette fonction.   L'\'egalit\'e $A_{M}= A_{M_{\eta}}$ identifie les racines de $A_{M_{\eta}}$ dans $G_{\eta}$ \`a des racines de $A_{M}$ dans $G$. Si $A_{G}\subsetneq A_{G_{\eta}}$, l'ensemble ${\cal J}_{\tilde{M}}^{\tilde{G}}(B_{{\cal O}})$ est vide et on n'a rien \`a d\'emontrer. Supposons $A_{G}=A_{G_{\eta}}$. Alors les deux ensembles 
  ${\cal J}_{\tilde{M}}^{\tilde{G}}(B_{{\cal O}})$  et  ${\cal J}_{M_{\eta}}^{G_{\eta}}(B)$ s'identifient. Les deux espaces $U_{J}$ possibles associ\'es \`a un \'el\'ement $J$ de cet ensemble s'identifient par l'\'egalit\'e    $A_{M}= A_{M_{\eta}}$.
Soient  $\boldsymbol{\gamma}\in D_{tr-orb}({\cal O})\otimes Mes(M({\mathbb R}))^*$ et $J\in {\cal J}_{\tilde{M}}^{\tilde{G}}(B_{{\cal O}})$. Gr\^ace au lemme 4.3(i), on \'ecrit
$$\boldsymbol{\gamma}=\sum_{y\in \dot{{\cal X}}(\eta)}desc_{\eta[y]}^{\tilde{M},*}(\boldsymbol{\gamma}[y]),$$
avec $\boldsymbol{\gamma}[y]\in D_{tr-unip}(M_{\eta[y]}({\mathbb R}))\otimes Mes(M_{\eta[y]}({\mathbb R}))^*$. On pose
$$(1) \qquad \rho_{J}^{\tilde{G}}(\boldsymbol{\gamma})=\sum_{y\in \dot{{\cal X}}(\eta)}desc_{\eta[y]}^{\tilde{M},*}(\rho_{J}^{G_{\eta[y]}}(\boldsymbol{\gamma}[y])).$$
Les \'el\'ements $\rho_{J}^{G_{\eta[y]}}(\boldsymbol{\gamma}[y])$ ne sont d\'efinis que modulo l'annulateur $Ann_{unip}^{G_{\eta[y]}}$ de l'application d'induction de $M_{\eta(y]}$ \`a $G_{\eta[y]}$. Il est clair que $desc_{\eta[y]}^{\tilde{M},*}$ envoie cet annulateur dans $Ann_{{\cal O}}^{\tilde{G}}$ de l'application d'induction de $\tilde{M}$ \`a $\tilde{G}$. L'\'el\'ement $\rho_{J}^{\tilde{G}}(\boldsymbol{\gamma})$ est donc bien d\'efini dans 
$$U_{J}\otimes (D_{g\acute{e}om}({\cal O})\otimes Mes(M({\mathbb R}))^*)/Ann_{{\cal O}}^{\tilde{G}}.$$  
La d\'efinition  ne d\'epend pas des choix faits. En effet, pour $y\in \dot{{\cal X}}(\eta)$, notons $p_{y}$ la projection de 4.2(1) relative \`a $\eta[y]$. Puisque $desc_{\eta[y]}^{\tilde{M},*}\circ p_{y}=desc_{\eta[y]}^{\tilde{M},*}$, on peut remplacer $desc_{\eta[y]}^{\tilde{M},*}$ par $desc_{\eta[y]}^{\tilde{M},*}\circ p_{y}$ dans les \'egalit\'es pr\'ec\'edentes. On a d\'ej\`a dit que $\rho_{J}^{G_{\eta[y]}}$ \'etait \'equivariant par $p_{y}$, c'est-\`a-dire que 
$$p_{y}(\rho_{J}^{G_{\eta[y]}}(\boldsymbol{\gamma}[y]))=\rho_{J}^{G_{\eta[y]}}(p_{y}(\boldsymbol{\gamma}[y])).$$
L'\'el\'ement $p_{y}(\boldsymbol{\gamma}[y])$ est uniquement d\'etermin\'e et appartient \`a $D_{tr-unip}(M_{\eta[y]}({\mathbb R}))\otimes Mes(M_{\eta[y]}({\mathbb R}))^*$. Cela prouve que la formule (1) ne d\'epend pas des choix de $\boldsymbol{\gamma}[y]$. De m\^eme, changer l'ensemble de repr\'esentants $\dot{{\cal X}}(\eta)$ revient \`a transporter les distributions $\boldsymbol{\gamma}[y]$ par des isomorphismes, ce qui ne modifie pas le r\'esultat.

 On va montrer que l'application $\rho_{J}^{\tilde{G}}$ ainsi d\'efinie v\'erifie les conditions (1) et (3) de 2.4, c'est-\`a-dire

 (2)  cette d\'efinition co\"{\i}ncide avec celle de [II] 3.4 pour $\boldsymbol{\gamma}\in D_{orb}({\cal O})\otimes Mes(M({\mathbb R}))^*$;
 
 (3) soit ${\bf M}'$ une donn\'ee endoscopique elliptique et relevante de $(M,\tilde{M})$ avec $M'\not=M$; soit ${\cal O}'$ une classe de conjugaison stable semi-simple dans $\tilde{M}'({\mathbb R})$ correspondant \`a ${\cal O}$; soit $\boldsymbol{\delta}\in D_{tr-orb}^{st}({\bf M}',{\cal O}')\otimes Mes(M'({\mathbb R}))^*$; alors $\rho_{J}^{\tilde{G},{\cal E}}({\bf M}',\boldsymbol{\delta})=\rho_{J}^{\tilde{G}}(transfert(\boldsymbol{\delta}))$. 
 
 Dans la situation de (2), on peut choisir des $\boldsymbol{\gamma}[y]\in D_{orb,unip}(M_{\eta[y]}({\mathbb R}))\otimes Mes(M_{\eta[y]}({\mathbb R}))^*$. La preuve est alors la m\^eme qu'en [III] 4.1. 
 
 Dans la situation de (3), fixons $\epsilon\in {\cal O}'$ avec $M'_{\epsilon}$ quasi-d\'eploy\'e. On fixe un diagramme $(\epsilon,B^{M'},T',B^M,T,\eta)$ d'espaces ambiants $\tilde{M}'$ et $\tilde{M}$. On a $A_{M}=A_{M'}$ par ellipticit\'e et on suppose d'abord 
 
 (4)  $A_{M'_{\epsilon}}=A_{M'}$. 
 
 On fixe des donn\'ees auxiliaires $M'_{1}$,...,$\Delta_{1}$ pour ${\bf M}'$ et un \'el\'ement $\epsilon_{1}\in \tilde{M}'_{1}({\mathbb R})$ au-dessus de $\epsilon$. Gr\^ace au lemme 4.3(ii), on peut fixer $\boldsymbol{\delta}_{\epsilon}\in D^{st}_{tr-unip}(M'_{1,\epsilon_{1}}({\mathbb R}))\otimes Mes(M'_{\epsilon}({\mathbb R}))^*$ de sorte que
 $$\boldsymbol{\delta}=desc_{\epsilon_{1}}^{st,\tilde{M}'_{1},*}(\boldsymbol{\delta}_{\epsilon}).$$
 Remarquons qu'on peut identifier un voisinage stablement invariant de $1$ dans $M'_{1,\epsilon_{1}}({\mathbb R})$ \`a un tel voisinage de $1$ dans le produit  $C_{1}({\mathbb R})\times M'_{\epsilon}({\mathbb R})$. Ainsi $ D^{st}_{tr-unip}(M'_{1,\epsilon_{1}}({\mathbb R})\otimes Mes(M'_{\epsilon}({\mathbb R}))^*$ s'identifie \`a $ D^{st}_{tr-unip}(M'_{\epsilon}({\mathbb R})\otimes Mes(M'_{\epsilon}({\mathbb R}))^*$. On a expliqu\'e en 4.1 que les r\'esultats de la section 5 de [III] valaient, mutatis mutandis, sur notre corps de base r\'eel. On peut alors reprendre la preuve de [III] 7.1. Dans notre situation quasi-d\'eploy\'ee et \`a torsion int\'erieure, les constructions se simplifient. Le groupe $M'_{\epsilon}$ appara\^{\i}t comme le groupe de la donn\'ee endoscopique $\bar{{\bf M}}'$ de cette r\'ef\'erence. On n'a plus besoin de passer aux rev\^etements simplement connexes des groupes d\'eriv\'es. On obtient les formules parall\`eles
 $$ transfert(\boldsymbol{\delta})=\sum_{y\in \dot{{\cal Y}}^M(\eta)} desc_{\eta[y]}^{\tilde{M},*}\circ transfert_{y}(\boldsymbol{\delta}_{\epsilon}),$$
 $$ \rho_{J}^{\tilde{G},{\cal E}}({\bf M}',\boldsymbol{\delta},a)=\sum_{y\in \dot{{\cal Y}}^M(\eta)} desc_{\eta[y]}^{\tilde{M},*}(\rho_{J}^{G_{\eta[y]},{\cal E}}(\bar{{\bf M}}',\boldsymbol{\delta}_{\epsilon},a))$$
 pour tout $a\in A_{M}({\mathbb R})$ en position g\'en\'erale et proche de $1$.  D'apr\`es 2.4(3)  qui est d\'ej\`a prouv\'e pour le groupe $G_{\eta[y]}$ (puisqu'il n'est pas tordu), on a  
$$\rho_{J}^{G_{\eta[y]},{\cal E}}(\bar{{\bf M}}',\boldsymbol{\delta}_{\epsilon},a)=\rho_{J}^{G_{\eta[y]}}(transfert_{y}(\boldsymbol{\delta}_{\epsilon}),a)$$
pour tout $y$. On a une application surjective $p:\dot{{\cal Y}}^M(\eta)\to \dot{{\cal X}}^M(\eta)$ et on peut supposer nos syst\`emes de repr\'esentants choisis de sorte que $\eta[y]$ soit constant sur les fibres. Pour $y\in \dot{{\cal X}}^M(\eta)$, posons
$$\boldsymbol{\gamma}[y]=\sum_{y'\in p^{-1}(y)}transfert_{y'}(\boldsymbol{\delta}_{\epsilon}).$$
Les formules ci--dessus se r\'ecrivent
 $$ transfert(\boldsymbol{\delta})=\sum_{y\in \dot{{\cal X}}^M(\eta)} desc_{\eta[y]}^{\tilde{M},*}(\boldsymbol{\gamma}[y])$$
 et
 $$\rho_{J}^{\tilde{G},{\cal E}}({\bf M}',\boldsymbol{\delta},a)=\sum_{y\in \dot{{\cal X}}^M(\eta)} desc_{\eta[y]}^{\tilde{M},*}(\rho_{J}^{G_{\eta[y]}}(\boldsymbol{\gamma}[y],a)).$$
 Il suffit d'appliquer la d\'efinition (1) pour conclure \`a l'\'egalit\'e  $\rho_{J}^{\tilde{G},{\cal E}}({\bf M}',\boldsymbol{\delta})=\rho_{J}^{\tilde{G}}(transfert(\boldsymbol{\delta}))$. Cela prouve (3) sous l'hypoth\`ese (4). 
 
 Supposons maintenant $A_{M'_{\epsilon}}\not= A_{M'}$. On introduit l'espace de Levi $\tilde{R}'$ de $\tilde{M}'$ tel que $A_{M'_{\epsilon}}=A_{R'}$. C'est l'espace d'une donn\'ee endoscopique elliptique d'un espace de Levi $\tilde{R}$ de $\tilde{M}$ contenant $\eta$. On note ${\cal O}'_{\tilde{R}'}$ la classe de conjugaison stable de $\epsilon$ dans $\tilde{R}'({\mathbb R})$ et ${\cal O}_{\tilde{R}}$ celle de $\eta$ dans $\tilde{R}({\mathbb R})$. Appliquant le lemme 4.4, on peut fixer $\boldsymbol{\delta}_{{\bf R}'}\in D_{tr-orb}^{st}({\bf R}',{\cal O}'_{\tilde{R}'})\otimes Mes(R'({\mathbb R}))^*$ de sorte que $\boldsymbol{\delta}=(\boldsymbol{\delta}_{{\bf R}'})^{{\bf M}'}$.  
 On applique l'analogue de 2.5(8) qui est valide puisque $M'\not=M$:
 $$(5) \qquad \rho_{J}^{\tilde{G},{\cal E}}({\bf M}',\boldsymbol{\delta},a)=\sum_{\tilde{L}\in {\cal L}(\tilde{R}), J\in {\cal J}_{\tilde{R}}^{\tilde{L}}(B_{{\cal O}_{\tilde{R}}})}d_{\tilde{R}}^{\tilde{G}}(\tilde{M},\tilde{L})\rho_{J}^{\tilde{L},{\cal E}}({\bf R}',\boldsymbol{\delta}_{{\bf R}'},a)^{\tilde{M}}.$$
 Fixons $\tilde{L}$ apparaissant ci-dessus.  Par r\'ecurrence, on peut supposer que
 $$\rho_{J}^{\tilde{L},{\cal E}}({\bf R}',\boldsymbol{\delta}_{{\bf R}'},a)=\rho_{J}^{\tilde{L}}(transfert(\boldsymbol{\delta}_{{\bf R}'}),a)).$$
  On a les relations $A_{R}\subset A_{R_{\eta}}\subset A_{R'_{\epsilon}}=A_{R'}=A_{R}$. Donc $A_{R}=A_{R_{\eta}}$. On est donc dans la situation de d\'epart, avec $\tilde{G}$ et $\tilde{M}$  remplac\'es par $\tilde{L}$ et $\tilde{R}$. Ecrivons
  $$(6) \qquad transfert(\boldsymbol{\delta}_{{\bf R}'})=\sum_{y\in \dot{{\cal X}}^R(\eta)}desc_{\eta[y]}^{\tilde{R},*}(\boldsymbol{\gamma}[y]).$$
  On peut supposer que $\rho_{J}^{\tilde{L}}(transfert(\boldsymbol{\delta}_{{\bf R}'},a))$ 
  est donn\'e par l'analogue de la formule (1), \`a savoir
 $$\rho_{J}^{\tilde{L}}(transfert(\boldsymbol{\delta}_{{\bf R}'}),a)=\sum_{y\in \dot{{\cal X}}^R(\eta)}desc_{\eta[y]}^{\tilde{R},*}(\rho_{J}^{L_{\eta[y]}}(\boldsymbol{\gamma}[y],a)).$$
 D'o\`u par induction
 $$\rho_{J}^{\tilde{L},{\cal E}}({\bf R}',\boldsymbol{\delta}_{{\bf R}'},a)^{\tilde{M}}=\rho_{J}^{\tilde{L}}(transfert(\boldsymbol{\delta}_{{\bf R}'}),a)^{\tilde{M}}=\sum_{y\in \dot{{\cal X}}^R(\eta)}desc_{\eta[y]}^{\tilde{M},*}(\rho_{J}^{L_{\eta[y]}}(\boldsymbol{\gamma}[y],a)^{M_{\eta[y]}}).$$
 Revenant \`a la formule (5), on  obtient
$$ \rho_{J}^{\tilde{G},{\cal E}}({\bf M}',\boldsymbol{\delta},a)= \sum_{y\in \dot{{\cal X}}^R(\eta)}desc_{\eta[y]}^{\tilde{M},*}(X_{J}[y]),$$
o\`u
$$X_{J}[y]=\sum_{\tilde{L}\in {\cal L}(\tilde{R}), J\in {\cal J}_{\tilde{R}}^{\tilde{L}}(B_{{\cal O}})}d_{\tilde{R}}^{\tilde{G}}(\tilde{M},\tilde{L})\rho_{J}^{L_{\eta[y]}}(\boldsymbol{\gamma}[y],a)^{M_{\eta[y]}}.$$
Notons $E$ l'ensemble des $\tilde{L}\in {\cal L}(\tilde{R})$ tels que $ J\in {\cal J}_{\tilde{R}}^{\tilde{L}}(B_{{\cal O}_{\tilde{R}}})$ et $d_{\tilde{R}}^{\tilde{G}}(\tilde{M},\tilde{L})\not=0$. Fixons $y\in \dot{{\cal X}}^R(\eta)$.  Notons $E_{\eta[y]}$ l'ensemble des $L'\in {\cal L}^{G_{\eta[y]}}(R_{\eta[y]})$ tels que $J\in {\cal J}_{R_{\eta[y]}}^{L'}(B)$ et $d_{R_{\eta[y]}}^{G_{\eta[y]}}(M_{\eta[y]},L')\not=0$. Montrons que

(7) l'application $\tilde{L}\mapsto L_{\eta[y]}$ se restreint en une bijection de $E$ sur $E_{\eta[y]}$; pour $\tilde{L}\in E$, on a l'\'egalit\'e $d_{\tilde{R}}^{\tilde{G}}(\tilde{M},\tilde{L})=d_{R_{\eta[y]}}^{G_{\eta[y]}}(M_{\eta[y]},L_{\eta[y]})$.

On ne perd rien \`a supposer $y=1$ (on n'utilisera pas ici le fait que $G_{\eta}$ est quasi-d\'eploy\'e). Comme on l'a dit plus haut, le fait que ${\cal J}_{\tilde{R}}^{\tilde{L}}(B_{{\cal O}_{\tilde{R}}})$ soit non-vide implique que $A_{L_{\eta}}=A_{L}$ et que ${\cal J}_{\tilde{R}}^{\tilde{L}}(B_{{\cal O}_{\tilde{R}}})={\cal J}_{R_{\eta}}^{L_{\eta}}(B)$. Si $\tilde{L} \in E$, on a donc $J\in  {\cal J}_{R_{\eta}}^{L_{\eta}}(B)$ et les \'egalit\'es ${\cal A}_{S}={\cal A}_{S_{\eta}}$ pour $S=G,M,L,R$. Il en r\'esulte que $d_{R_{\eta}}^{G_{\eta}}(M_{\eta},L_{\eta})=d_{\tilde{R}}^{\tilde{G}}(\tilde{M},\tilde{L})\not=0$. Donc $L_{\eta}\in E_{\eta}$.  Notons que l'\'egalit\'e pr\'ec\'edente est la derni\`ere assertion de (7). L'\'egalit\'e $A_{L_{\eta}}=A_{L}$ implique que $\tilde{L}$ est uniquement d\'etermin\'e par $L_{\eta}$.  R\'eciproquement, pour $L'\in E_{\eta}$, on d\'efinit $\tilde{L}$ par l'\'egalit\'e $A_{L}=A_{L'}$ et des arguments analogues montrent que $\tilde{L}\in E$. Cela prouve (7).

Gr\^ace \`a (8), on peut r\'ecrire
$$X_{J}[y]=\sum_{L'\in {\cal L}^{G_{\eta[y]}}(R_{\eta[y]}) , J\in {\cal J}_{R_{\eta[y]}}^{L'}(B)}d_{R_{\eta[y]}}^{G_{\eta[y]}}(M_{\eta[y]},L')\rho_{J}^{L'}(\boldsymbol{\gamma}[y],a)^{M_{\eta[y]}}.$$
En vertu de la formule de descente 2.4(7) d\'ej\`a prouv\'ee pour le groupe $G_{\eta[y]}$, on obtient
$$X_{J}[y]=\rho_{J}^{G_{\eta[y]}}(\boldsymbol{\gamma}[y]^{M_{\eta[y]}},a),$$
d'o\`u
$$ \rho_{J}^{\tilde{G},{\cal E}}({\bf M}',\boldsymbol{\delta},a)= \sum_{y\in \dot{{\cal X}}^R(\eta)}desc_{\eta[y]}^{\tilde{M},*}(\rho_{J}^{G_{\eta[y]}}(\boldsymbol{\gamma}[y]^{M_{\eta[y]}},a)).$$
Par induction, la formule (6) donne
$$transfert(\boldsymbol{\delta})= (transfert(\boldsymbol{\delta}_{{\bf R}'}))^{\tilde{M}}=\sum_{y\in \dot{{\cal X}}^R(\eta)}desc_{\eta[y]}^{\tilde{M},*}(\boldsymbol{\gamma}[y]^{M_{\eta[y]}}).$$
Il y a une application naturelle de $\dot{{\cal X}}^R(\eta)$ dans $\dot{{\cal X}}^M(\eta)$: \`a un \'el\'ement $y'$ du premier ensemble, on associe l'unique \'el\'ement $y$ du second tel que les points $\eta[y']$ et $\eta[y]$ soient conjugu\'es par un \'el\'ement de $M({\mathbb R})$. En fixant une telle conjugaison, on peut identifier $\boldsymbol{\gamma}[y']^{M_{\eta[y']}}$ \`a un \'el\'ement de $D_{tr-orb}(M_{\eta[y]}({\mathbb R}))\otimes Mes(M_{\eta[y]}({\mathbb R}))^*$. En notant $\boldsymbol{\gamma}[y]$ la somme des $\boldsymbol{\gamma}[y']^{M_{\eta[y']}}$ sur les \'el\'ements $y'$ s'envoyant sur $y$, on obtient 
$$transfert(\boldsymbol{\delta})=\sum_{y\in \dot{{\cal X}}^M(\eta)}desc_{\eta[y]}^{\tilde{M},*}(\boldsymbol{\gamma}[y])$$
et
$$ \rho_{J}^{\tilde{G},{\cal E}}({\bf M}',\boldsymbol{\delta},a)= \sum_{y\in \dot{{\cal X}}^M(\eta)}desc_{\eta[y]}^{\tilde{M},*}(\rho_{J}^{G_{\eta[y]}}(\boldsymbol{\gamma}[y],a)).$$
Mais alors, la d\'efinition (1) conduit \`a l'\'egalit\'e
$$\rho_{J}^{\tilde{G},{\cal E}}({\bf M}',\boldsymbol{\delta},a)=\rho_{J}^{\tilde{G}}(transfert(\boldsymbol{\delta}),a).$$
Cela ach\`eve la preuve de (3).

 On d\'efinit l'application $\sigma_{J}^{\tilde{G}}$ par la formule  2.4(13). Soit $\boldsymbol{\delta}\in D_{tr-orb}^{st}({\cal O})\otimes Mes(M({\mathbb R}))^*$. D'apr\`es le lemme  4.3 (ii), on peut choisir  $\boldsymbol{\delta}_{\eta}\in D_{tr-unip}^{st}(M_{\eta}({\mathbb R}))\otimes Mes(M_{\eta}({\mathbb R}))^*$ de sorte que $\boldsymbol{\delta}=desc_{\eta}^{st,\tilde{M},*}(\boldsymbol{\delta}_{\eta})$. Pour $a\in A_{M}({\mathbb R})$ en position g\'en\'erale et proche de $1$, on a
 $$(9)\qquad \sigma_{J}^{\tilde{G}}(\boldsymbol{\delta},a)=e_{\tilde{M}}^{\tilde{G}}(\eta)desc_{\eta}^{st,\tilde{M},*}(\sigma_{J}^{G_{\eta}}(\boldsymbol{\delta}_{\eta},a)),$$
 cf. [III] 4.3 pour la d\'efinition de $e_{\tilde{M}}^{\tilde{G}}(\eta)$. Cela se prouve comme en [III] 7.3, \`a partir de la formule de descente d\'ej\`a prouv\'ee pour le terme $\rho_{J}^{\tilde{G}}(\boldsymbol{\delta},a)$. La formule (9) implique que $\sigma_{J}^{\tilde{G}}(\boldsymbol{\delta},a)$ est stable modulo $Ann_{{\cal O}}^{\tilde{G}}$.  
 
 On a ainsi v\'erifi\'e les conditions (1) et (2) de 2.7. Elles impliquent la validit\'e du programme de 2.4.

 \bigskip
 
 \section{Extension des d\'efinitions, cas g\'en\'eral}
 
 \bigskip
 
 \subsection{Un r\'esultat compl\'ementaire pour l'endoscopie non standard}
 On consid\`ere ici un triplet endoscopique non standard $(G_{1},G_{2},j_{*})$, cf. [III] 6.1 dont on reprend les notations.  Rappelons qu'il y a une correspondance bijective entre classes de conjugaison stable semi-simples dans $\mathfrak{g}_{1}({\mathbb R})$ et  classes de conjugaison stable semi-simples dans $\mathfrak{g}_{2}({\mathbb R})$. De cette correspondance r\'esulte un isomorphisme
 $$ SI(\mathfrak{g}_{1}({\mathbb R}))\otimes Mes(G_{1}({\mathbb R}))\simeq
 SI(\mathfrak{g}_{2}({\mathbb R}))\otimes Mes(G_{2}({\mathbb R}))^*.$$
 En effet, les espaces $SI(\mathfrak{g}_{i}({\mathbb R}))$ sont d\'ecrits par Shelstad (celle-ci traite le cas des groupes mais la description vaut a fortiori pour les alg\`ebres de Lie). Le point essentiel de cette descrition  sont les conditions de saut.  Mais celles-ci sont insensibles au remplacement d'une racine par un multiple r\'eel (ici rationnel) de cette racine. L'assertion ci-dessus s'ensuit. 
 
 Soient $M_{1}$ et $M_{2}$ des  Levi de $G_{1}$ et $G_{2}$ qui se correspondent.  En rempla\c{c}ant ci-dessus $G_{1}$ et $G_{2}$ par $M_{1}$ et $M_{2}$, puis en dualisant, on obtient un isomorphisme
 $$(1) \qquad D^{st}_{g\acute{e}om}(\mathfrak{m}_{1}({\mathbb R}))\otimes Mes(M_{1}({\mathbb R}))^*\simeq
 D^{st}_{g\acute{e}om}(\mathfrak{m}_{2}({\mathbb R}))\otimes Mes(M_{2}({\mathbb R}))^*.$$
 En le restreignant aux distributions \`a support nilpotent, puis en passant aux groupes par l'exponentielle, on obtient un isomorphisme
 $$(2) \qquad D^{st}_{unip}(M_{1}({\mathbb R}))\otimes Mes(M_{1}({\mathbb R}))^*\simeq D^{st}_{unip}(M_{2}({\mathbb R}))\otimes Mes(M_{2}({\mathbb R}))^*.$$
 En reprenant la preuve du lemme [II] 3.1, on voit que cet isomorphisme envoie $Ann_{unip}^{G_{1},st}$ sur $Ann_{unip}^{G_{2},st}$.
 
  On impose
 
 (3) la fonction $b$ (cf. [III] 6.1) est constante sur l'ensemble de racines $\Sigma^{M_{2}}(T_{2})$. 
 
 Montrons que
 
 (4) sous l'hypoth\`ese (3), l'isomorphisme (2) se restreint en un isomorphisme
 $$D^{st}_{tr-unip}(M_{1}({\mathbb R}))\otimes Mes(M_{1}({\mathbb R}))^*\simeq D^{st}_{tr-unip}(M_{2}({\mathbb R}))\otimes Mes(M_{2}({\mathbb R}))^*.$$
 
 L'hypoth\`ese (3) entra\^{\i}ne que les rev\^etements simplement connexes des groupes d\'eriv\'es de $M_{1}$ et $M_{2}$ sont isomorphes. Notons $M_{SC}$ ce rev\^etement commun et notons $\underline{b}$ la valeur constante de $b$ sur $\Sigma^{M_{2}}(T_{2})$. L'automorphisme $X\mapsto \underline{b}X$ de $\mathfrak{m}_{SC}$ induit un automorphisme $\boldsymbol{\gamma}\mapsto \boldsymbol{\gamma}[\underline{b}]$ de $D^{st}_{ nil}(\mathfrak{m}_{SC}({\mathbb R}))\otimes Mes(M_{SC}({\mathbb R}))^*$, qui se rel\`eve en un automorphisme de $D^{st}_{unip}(M_{SC}({\mathbb R}))\otimes Mes(M_{SC}({\mathbb R}))^*$ not\'e de la m\^eme fa\c{c}on.  On a un diagramme commutatif
 $$\begin{array}{ccc}D^{st}_{unip}(M_{1}({\mathbb R}))\otimes Mes(M_{1}({\mathbb R}))^*&\simeq&D^{st}_{unip}(M_{2}({\mathbb R}))\otimes Mes(M_{2}({\mathbb R}))^*\\ \iota^*_{M_{SC},M_{1}}\uparrow\,\,&&\,\,\uparrow \iota^*_{M_{SC},M_{2}}\\ D^{st}_{unip}(M_{SC}({\mathbb R}))\otimes Mes(M_{SC}({\mathbb R}))^*&\stackrel{\boldsymbol{\gamma}\mapsto \boldsymbol{\gamma}[\underline{b}]}{\to}&D^{st}_{unip}(M_{SC}({\mathbb R}))\otimes Mes(M_{SC}({\mathbb R}))^*.\\ \end{array}$$
 D'apr\`es le lemme 3.3, les applications $\iota^*_{M_{SC},M_{i}}$ pour $i=1,2$ deviennent des isomorphismes quand on remplace les espaces $D^{st}_{unip}$ par les espaces $D^{st}_{tr-unip}$. Il suffit donc de prouver que l'espace $D^{st}_{tr-unip}(M_{SC}({\mathbb R}))\otimes Mes(M_{SC}({\mathbb R}))^*$ est stable par l'application $\boldsymbol{\gamma}\mapsto \boldsymbol{\gamma}[\underline{b}]$.  Plus g\'en\'eralement, consid\'erons un groupe r\'eductif connexe $G$ sur ${\mathbb R}$ et un r\'eel non nul $r$. On d\'efinit comme ci-dessus l'automorphisme $\boldsymbol{\gamma}\mapsto \boldsymbol{\gamma}[r]$ de $D_{unip}(G({\mathbb R}))\otimes Mes(G({\mathbb R}))^*$. On montre que
 
 (5)  l'espace $D_{tr-unip}(G({\mathbb R}))\otimes Mes(G({\mathbb R}))^*$ est invariant par cet automorphisme;
 
 (6) si $G$ est quasi-d\'eploy\'e, $D_{tr-unip}^{st}(G({\mathbb R}))\otimes Mes(G({\mathbb R}))^*$ est invariant par cet automorphisme.
 
 Puisque l'automorphisme respecte clairement la stabilit\'e, (6) r\'esulte de (5). Il est clair que l'espace $D_{orb,unip}(G({\mathbb R}))\otimes Mes(G({\mathbb R}))^*$ est invariant par l'automorphisme. Soit ${\bf G}'$ une donn\'ee endoscopique elliptique et relevante de $G$, avec $G'\not=G$. Soit $\boldsymbol{\delta}\in D^{st}_{tr-unip}({\bf G}')\otimes Mes(G'({\mathbb R}))^*$, posons $\boldsymbol{\gamma}=transfert(\boldsymbol{\delta})$. On doit prouver que $\boldsymbol{\gamma}[r]$ appartient \`a $D_{tr-unip}(G({\mathbb R}))\otimes Mes(G({\mathbb R}))^*$. Comme dans la preuve de [III] 6.7, on montre qu'il existe une constante $c(r)\not=0$ telle que $\boldsymbol{\gamma}[r]=c(r)transfert(\boldsymbol{\delta}[r])$. En utilisant (6) par r\'ecurrence, on a $\boldsymbol{\delta}[r]\in D^{st}_{tr-unip}({\bf G}')\otimes Mes(G'({\mathbb R}))^*$, d'o\`u la conclusion. Cela prouve (5) et (4).

 Supposons que le triplet $(G_{1},G_{2},j_{*})$ est \'equivalent \`a un triplet quasi-\'el\'ementaire. Notons $(G_{0,1},G_{0,2},j_{0*})$ le triplet \'el\'ementaire sur $F_{0}={\mathbb R}$ ou ${\mathbb C}$ tel que $(G_{1},G_{2},j_{*})$ soit \'equivalent  au triplet d\'eduit de $(G_{0,1},G_{0,2},j_{0*})$ par restriction des scalaires de $F_{0}$ \`a ${\mathbb R}$. 
 Aux Levi $M_{1}$ et $M_{2}$    sont associ\'es des Levi $M_{0,1}$ et $M_{0,2}$ de $G_{0,1}$ et $G_{0,2}$.  Supposons l'une des conditions suivantes v\'erifi\'ee:
 
 (7)  $(G_{0,1},G_{0,2},j_{0*})$ est du type (1) de [III] 6.1;
 
 (8) $(G_{0,1},G_{0,2},j_{0,*})$ est de l'un des types (2) ou (3) de [III] 6.1; pour $i=1,2$, si $G_{0,i}$ est de type $B_{n}$, resp. $C_{n}$, les \'el\'ements de $\Sigma^{M_{0,i}}(T_{0,i})$ sont des racines longues, resp. courtes. 
 
Chacune de ces conditions implique (3).  La condition (8) implique que  le groupe $M_{SC}$ d\'efini ci-dessus est isomorphe \`a un produit de groupes $SL_{k}({\mathbb R})$ si $F_{0}={\mathbb R}$, de groupes $SL_{k}({\mathbb C})$ si $F_{0}={\mathbb C}$. 
 
 Soient $B_{1}$ et $B_{2}$ des fonctions comme en [III] 6.4, v\'erifiant toutes deux les hypoth\`eses de [II] 1.8. On a alors un lemme similaire \`a [III] 6.5.
 
 \ass{Lemme}{On suppose v\'erifi\'ee  l'une des conditions (7) ou (8). Pour $i=1,2$, soient $J_{i}\in {\cal J}_{M_{i}}^{G_{i}}(B_{i})$ et $\boldsymbol{\delta}_{i}\in D_{tr-unip}^{st}(M_{i}({\mathbb R}))\otimes Mes(M_{i}({\mathbb R}))^*$. On suppose que $J_{1}$ et $J_{2}$ se correspondent par la bijection entre les ensembles $\Sigma(A_{M_{i}},B_{i})$ et que $\boldsymbol{\delta}_{1}$ et $\boldsymbol{\delta}_{2}$ se correspondent par l'isomorphisme (4) ci-dessus. Alors, pour tout $X_{1}\in \mathfrak{a}_{M_{1}}({\mathbb R})$ en position g\'en\'erale et proche de $0$, les \'el\'ements $\sigma_{J_{1}}^{G_{1}}(\boldsymbol{\delta}_{1},exp(X_{1}))$ et $c_{M_{1},M_{2}}^{G_{1},G_{2}}\sigma_{J_{2}}^{G_{2}}(\boldsymbol{\delta}_{2},exp(j_{*}(X_{1})))$ se correspondent par l'isomorphisme (4). }
 
 Ce lemme sera d\'emontr\'e sous hypoth\`eses en 5.7.  Il faut l'inclure dans notre sch\'ema de r\'ecurrence. On a d\'efini l'entier $N(G_{1},G_{2},j_{*})$ en [III] 6.1. Rappelons que $(G_{2},G_{1},j_{*}^{-1})$ est aussi un triplet quasi-\'el\'ementaire. On pose
 $$N^{max}(G_{1},G_{2},j_{*})=sup(N(G_{1},G_{2},j_{*}),N(G_{2},G_{1},j_{*}^{-1})).$$
Plus simplement, si (7) est v\'erifi\'ee, $N^{max}(G_{1},G_{2},j_{*})=0$. Si (8) est v\'erifi\'ee et que les syst\`emes de racines de $G_{0,1}$ et $G_{0,2}$ sont de type $B_{n}$ ou $C_{n}$, on a $
  N^{max}(G_{1},G_{2},j_{*})=[F_{0}:{\mathbb R}](4n^2-1)$. Soit $N\in {\mathbb N}$. Pour d\'emontrer le lemme relativement \`a un triplet tel que $N^{max}(G_{1},G_{2},j_{*})=N$, on le suppose v\'erifi\'e pour tout triplet  $(G'_{1},G'_{2},j'_{*})$ v\'erifiant des conditions similaires et tel que $N^{max}(G'_{1},G'_{2},j'_{*})<N$. On suppose aussi v\'erifi\'es tous nos r\'esultats concernant des triplets $(KG,K\tilde{G},{\bf a})$ tels que $dim(G_{SC})\leq N$. Par ailleurs, pour d\'emontrer un r\'esultat concernant un tel triplet $(KG,K\tilde{G},{\bf a})$ tel que $dim(G_{SC})=N$, on suppose le lemme ci-dessus v\'erifi\'e pour tout triplet $(G_{1},G_{2},j_{*})$ comme ci-dessus tel que $N^{max}(G_{1},G_{2},j_{*})< N$.
  
  On a
  
  (9) supposons le lemme v\'erifi\'e dans le cas o\`u  les trois conditions suivantes sont satisfaites:  $(G_{1},G_{2},j_{*})$ est quasi-\'el\'ementaire,  (8) est v\'erifi\'ee et $B_{1}$ est constante de valeur $1$;
   alors le lemme est v\'erifi\'e.
  
  Preuve. Comme en [III] 6.7, on montre que l'assertion du lemme est insensible \`a la multiplication de  $j_{*}$, $B_{1}$ ou $B_{2}$ par des constantes. Cela nous ram\`ene au cas o\`u $(G_{1},G_{2},j_{*})$ est quasi-\'el\'ementaire. Si (7) est v\'erifi\'e, $j_{*}$ provient d'un isomorphisme de $G_{1}$ sur $G_{2}$ et l'assertion est tautologique. Si (8) est v\'erifi\'ee, on voit que l'une des fonctions $B_{1}$ ou $B_{2}$ est constante. Il est clair que le lemme pour $(G_{1},G_{2},j_{*})$ est \'equivalent \`a celui pour $(G_{2},G_{1},j_{*}^{-1})$. Quitte \`a \'echanger ces deux triplets et \`a multiplier nos fonctions par des constantes, on peut donc supposer $B_{1}$ constante de valeur $1$. Cela prouve (9).

 \bigskip
 
 \subsection{R\'ealisation conditionnelle du programme de 2.4}
 On consid\`ere un $K$-triplet $(KG,K\tilde{G},{\bf a})$, un $K$-espace de Levi $K\tilde{M}\in {\cal L}(K\tilde{M}_{0})$ et une classe de conjugaison stable semi-simple ${\cal O}$ dans $K\tilde{M}({\mathbb R})$. On va r\'ealiser le programme fix\'e en 2.4 sous l'hypoth\`ese (Hyp) de 2.5, que l'on rappelle:

(Hyp) pour tout $\boldsymbol{\gamma}\in D_{ orb}(K\tilde{M}({\mathbb R}),\omega)\otimes Mes(M({\mathbb R}))^*$ dont le support est form\'e d'\'el\'ements fortement $\tilde{G}$-r\'eguliers et pour tout ${\bf f}\in I(K\tilde{G}({\mathbb R}),\omega)\otimes Mes(G({\mathbb R}))$, on a l'\'egalit\'e
 $$I_{K\tilde{M}}^{K\tilde{G},{\cal E}}(\boldsymbol{\gamma},{\bf f})= I_{K\tilde{M}}^{K\tilde{G}}(\boldsymbol{\gamma},{\bf f}).$$

  Ecrivons $K\tilde{G}=(\tilde{G}_{p})_{p\in \Pi}$, $K\tilde{M}=(\tilde{M}_{p})_{p\in \Pi^M}$. Soit $J\in {\cal J}_{K\tilde{M}}^{K\tilde{G}}$. Pour $p\in \Pi^M$, on peut identifier $J$ \`a un \'el\'ement de ${\cal J}_{\tilde{M}_{p}}^{\tilde{G}_{p}}$. En [II] 3.1, on a associ\'e \`a $J$ un sous-groupe $G_{p,J}\subset G_{p}$ et un sous-espace tordu $\tilde{G}_{p,J}\subset \tilde{G}_{p}$, qui contiennent respectivement $M_{p}$ et $\tilde{M}_{p}$. On voit que la collection $(\tilde{G}_{p,J})_{p\in \Pi^M}$ s'\'etend en un $K$-espace tordu $K\tilde{G}_{J}=(\tilde{G}_{p,J})_{p\in \Pi_{J}}$, dont $K\tilde{M}$ est un $K$-espace de Levi. Signalons que, parce qu'on se s'int\'eressera qu'\`a des distributions induites \`a partir de $K\tilde{M}$, les composantes $\tilde{G}_{p,J}$ pour $p\in \Pi_{J}-\Pi^M$ ne joueront aucun r\^ole. 
 
 Soient $\boldsymbol{\gamma}\in D_{orb}({\cal O},\omega)\otimes Mes(M({\mathbb R}))^*$ et $a\in A_{K\tilde{M}}({\mathbb R})$ en position g\'en\'erale et proche de $1$. On  d\'efinit comme en [II] 3.2 le terme
 $$\rho_{J}^{K\tilde{G}}(\boldsymbol{\gamma},a)\in (D_{g\acute{e}om}({\cal O},\omega)\otimes Mes(M({\mathbb R}))^*)/Ann_{{\cal O}}^{K\tilde{G}}.$$
 Comme dans ce paragraphe, c'est l'image de $\rho_{J}^{K\tilde{G}_{J}}(\boldsymbol{\gamma},a)$ par l'application naturelle
  $$(1) \qquad  (D_{g\acute{e}om}({\cal O},\omega)\otimes Mes(M({\mathbb R}))^*)/Ann_{{\cal O}}^{K\tilde{G}_{J}}\to  (D_{g\acute{e}om}({\cal O},\omega)\otimes Mes(M({\mathbb R}))^*)/Ann_{{\cal O}}^{K\tilde{G}}.$$
  
  Soit ${\bf M}'=(M',{\cal M}',\tilde{\zeta})$ une donn\'ee endoscopique de $(M,\tilde{M},{\bf a})$, elliptique et relevante,   soit ${\cal O}'$ une classe de conjugaison stable semi-simple dans $\tilde{M}'({\mathbb R})$ corespondant \`a ${\cal O}$ et soit $\boldsymbol{\delta}\in D^{st}_{tr-orb}({\bf M}',{\cal O}')\otimes Mes(M'({\mathbb R}))^*$. Pour $a$ comme ci-dessus, on d\'efinit $\rho_{J}^{K\tilde{G},{\cal E}}({\bf M}',\boldsymbol{\delta},a)$ comme en 2.4.
  
  \ass{Proposition}{Le terme $\rho_{J}^{K\tilde{G},{\cal E}}({\bf M}',\boldsymbol{\delta},a)$ est l'image de $\rho_{J}^{K\tilde{G}_{J},{\cal E}}({\bf M}',\boldsymbol{\delta},a)$ par l'application (1).}
  
  Cela sera prouv\'e dans les paragraphes 5.3 et 5.5.
  
  Admettons ce r\'esultat. On va v\'erifier  la condition 2.5(4). Soit donc $J\in {\cal J}_{K\tilde{M}}^{K\tilde{G}}$ un \'el\'ement non-maximal. Soit $\boldsymbol{\gamma}\in D_{tr-orb}({\cal O},\omega)\otimes Mes(M({\mathbb R}))^*$. On l'\'ecrit
  $$(2) \qquad \boldsymbol{\gamma}=\boldsymbol{\gamma}_{orb}+\sum_{i=1,...,n}transfert(\boldsymbol{\delta}_{i})$$
  comme en 2.5(1). On pose
  $$\rho_{J}^{K\tilde{G}}(\boldsymbol{\gamma},a)=\rho_{J}^{K\tilde{G}}(\boldsymbol{\gamma}_{orb},a)+\sum_{i=1,...,n}\rho_{J}^{K\tilde{G},{\cal E}}({\bf M}'_{i},\boldsymbol{\delta}_{i},a).$$
  On doit prouver que cette expression ne d\'epend pas de la d\'ecomposition (2). D'apr\`es ce que l'on a dit ci-dessus pour $\boldsymbol{\gamma}_{orb}$ et d'apr\`es la proposition pour les autres termes, $\rho_{J}^{K\tilde{G}}(\boldsymbol{\gamma},a)$ est l'image de
  $$\rho_{J}^{K\tilde{G}_{J}}(\boldsymbol{\gamma}_{orb},a)+\sum_{i=1,...,n}\rho_{J}^{K\tilde{G}_{J},{\cal E}}({\bf M}'_{i},\boldsymbol{\delta}_{i},a)$$
  par l'application (1). Or $dim(G_{J,SC})<dim(G_{SC})$ puisque $J$ n'est pas maximal. Par r\'ecurrence, l'expression ci-dessus ne d\'epend pas de la d\'ecomposition (1): elle vaut $\rho_{J}^{K\tilde{G}_{J}}(\boldsymbol{\gamma},a)$. Il en est donc de m\^eme de $\rho_{J}^{K\tilde{G}}(\boldsymbol{\gamma},a)$. Cela prouve 2.5(4). Comme on l'a vu en 2.5, cette relation suffit, sous l'hypoth\`ese (Hyp), pour valider le programme de 2.4.
  
  \bigskip
  
  \subsection{Preuve de la proposition 5.2, premier cas}
  On fixe ${\bf M}'$, ${\cal O}'$, $\boldsymbol{\delta}$ comme dans cette proposition.  Pour simplifier, on fixe des mesures de Haar sur tous les groupes qui apparaissent, afin de se d\'ebarrasser des espaces de mesures. On fixe  un \'el\'ement $\epsilon\in {\cal O}'$ tel que $M'_{\epsilon}$ soit quasi-d\'eploy\'e et un diagramme $(\epsilon,B^{M'},T',B^M,T,\eta)$ joignant $\epsilon$ \`a un \'el\'ement $\eta\in {\cal O}$.   L'\'el\'ement $\eta$ appartient \`a une composante $\tilde{M}_{p}({\mathbb R})$ de $K\tilde{M}({\mathbb R})$.  Fixons des donn\'ees auxiliaires $M'_{1}$,...,$\Delta_{1}$ pour ${\bf M}'$ et un \'el\'ement $\epsilon_{1}\in \tilde{M}'_{1}({\mathbb R})$ au-dessus de $\epsilon$. On note ${\cal O}'_{1}$ la classe de conjugaison stable de $\epsilon_{1}$ dans $\tilde{M}'_{1}({\mathbb R})$. On peut identifier $\boldsymbol{\delta}$ \`a un \'el\'ement $\boldsymbol{\delta}_{1}\in D^{st}_{tr-orb,\lambda_{1}}({\cal O}'_{1}) $.
Supposons d'abord
  
  (1) $A_{M'_{\epsilon}}\not=A_{M'}$.
  
  Notons $\tilde{R}'$ le commutant de $A_{M'_{\epsilon}}$ dans $\tilde{M}'$ et ${\cal O}_{\tilde{R}'}$ la classe de conjugaison stable de $\epsilon$ dans $\tilde{R}'({\mathbb R})$. On note $\tilde{R}'_{1}$ l'image r\'eciproque de $\tilde{R}'$ dans $\tilde{M}'_{1}$ et ${\cal O}'_{\tilde{R}'_{1}}$ la classe de conjugaison stable de $\epsilon_{1}$ dans $\tilde{R}'_{1}({\mathbb R})$. 
  D'apr\`es le lemme  4.4, on peut fixer $\boldsymbol{\delta}_{\tilde{R}'_{1}}\in D_{tr-orb,\lambda_{1}}^{st}({\cal O}'_{\tilde{R}'_{1}})$ de sorte que $\boldsymbol{\delta}_{1}$ soit l'image de $\boldsymbol{\delta}_{\tilde{R}'_{1}}$ par induction de $\tilde{R}'_{1}$ \`a $\tilde{M}'_{1}$.  Du diagramme fix\'e se d\'eduit un homomorphisme
  $\xi:T^{\theta,0}\to T'$.
  On note par anticipation $A_{\tilde{R}_{p}}$ la composante neutre de l'image r\'eciproque de $A_{M'_{\epsilon}}$ par cet homomorphisme.  Dans la composante $\tilde{M}_{p}$, on a fix\'e un espace de Levi minimal $\tilde{M}_{p,0}$. Quitte \`a conjuguer $\eta$, on peut supposer $A_{\tilde{R}_{p}}\subset A_{\tilde{M}_{p,0}}$. 
   On note $\tilde{R}_{p}$ le commutant de $A_{\tilde{R}_{p}}$ dans $\tilde{M}_{p}$. C'est un espace de Levi semi-standard, qui se compl\`ete en un $K$-espace de Levi $K\tilde{R}$. On note ${\cal O}_{\tilde{R}}$ la classe de conjugaison stable de $\eta$ dans $K\tilde{R}({\mathbb R})$. L'espace $\tilde{R}'$ appara\^{\i}t comme l'espace endoscopique d'une donn\'ee endoscopique ${\bf R}'$ de $(KR,K\tilde{R},{\bf a})$. Des donn\'ees auxiliaires fix\'ees pour ${\bf M}'$ se d\'eduisent des donn\'ees auxiliaires pour ${\bf R}'$. On voit que $\boldsymbol{\delta}_{\tilde{R}'_{1}}$ s'identifie \`a un \'el\'ement de $D^{st}_{tr-orb}({\bf R}',{\cal O}'_{\tilde{R}'})$ que l'on note $\boldsymbol{\delta}_{{\bf R}'}$. Pour $J\in {\cal J}_{K\tilde{M}}^{K\tilde{G}}$ et $a\in A_{K\tilde{M}}({\mathbb R})$ en position g\'en\'erale et proche de $1$, on a alors la formule d'induction 2.5(8):
   $$(2) \qquad \rho_{J}^{K\tilde{G},{\cal E}}({\bf M}',\boldsymbol{\delta},a)=\sum_{K\tilde{L}\in {\cal L}^{K\tilde{G}}(K\tilde{R}), J\in {\cal J}_{K\tilde{R}}^{K\tilde{L}}}d_{\tilde{R}}^{\tilde{G}}(\tilde{M},\tilde{L})\rho_{J}^{K\tilde{L},{\cal E}}({\bf R}',\boldsymbol{\delta}_{{\bf R}'},a)^{K\tilde{M}}.$$
   On a la formule parall\`ele
  $$(3) \qquad \rho_{J}^{K\tilde{G}_{J},{\cal E}}({\bf M}',\boldsymbol{\delta},a)=\sum_{K\tilde{L}'\in {\cal L}^{K\tilde{G}_{J}}(K\tilde{R}), J\in {\cal J}_{K\tilde{R}}^{K\tilde{L}'}}d_{\tilde{R}}^{\tilde{G}_{J}}(\tilde{M},\tilde{L}')\rho_{J}^{K\tilde{L}',{\cal E}}({\bf R}',\boldsymbol{\delta}_{{\bf R}'},a)^{K\tilde{M}}.$$ 
  On voit ais\'ement que l'application $K\tilde{L}\mapsto K\tilde{L}_{J}$ est une bijection de l'ensemble des  $K\tilde{L}\in {\cal L}^{K\tilde{G}}(K\tilde{R})$ tels que $ J\in {\cal J}_{K\tilde{R}}^{K\tilde{L}}$ sur l'ensemble des $K\tilde{L}'\in {\cal L}^{K\tilde{G}_{J}}(K\tilde{R})$ tels que $J\in {\cal J}_{K\tilde{R}}^{K\tilde{L}'}$. Pour $K\tilde{L}$ dans l'ensemble de d\'epart, on a
  $$d_{\tilde{R}}^{\tilde{G}}(\tilde{M},\tilde{L})=d_{\tilde{R}}^{\tilde{G}_{J}}(\tilde{M},\tilde{L}_{J}).$$
  En effet, cela r\'esulte des \'egalit\'es ${\cal A}_{\tilde{G}_{J}}={\cal A}_{\tilde{G}}$, ${\cal A}_{\tilde{L}_{J}}={\cal A}_{\tilde{L}}$. L'hypoth\`ese (1) permet d'appliquer par r\'ecurrence la proposition 5.2 \`a tous les termes apparaissant dans les formules ci-dessus. C'est-\`a-dire que, pour tout $K\tilde{L}$ intervenant dans (2), $\rho_{J}^{K\tilde{L},{\cal E}}({\bf R}',\boldsymbol{\delta}_{{\bf R}'},a)$ est l'image de $\rho_{J}^{K\tilde{L}_{J},{\cal E}}({\bf R}',\boldsymbol{\delta}_{{\bf R}'},a)$ par l'application
  $$D_{g\acute{e}om}({\cal O}_{\tilde{R}},\omega)/Ann_{{\cal O}_{\tilde{R}}}^{K\tilde{L}_{J}}\to D_{g\acute{e}om}({\cal O}_{\tilde{R}},\omega)/Ann_{{\cal O}_{\tilde{R}}}^{K\tilde{L}}.$$
  Il en r\'esulte que $\rho_{J}^{K\tilde{L},{\cal E}}({\bf R}',\boldsymbol{\delta}_{{\bf R}'},a)^{K\tilde{M}}$ est l'image de $\rho_{J}^{K\tilde{L}_{J},{\cal E}}({\bf R}',\boldsymbol{\delta}_{{\bf R}'},a)^{K\tilde{M}}$ par l'application 
  $$D_{g\acute{e}om}({\cal O},\omega)/Ann_{{\cal O}}^{K\tilde{G}_{J}}\to D_{g\acute{e}om}({\cal O},\omega)/Ann_{{\cal O}}^{K\tilde{G}}.$$
  Les formules (2) et (3) entra\^{\i}nent que $\rho_{J}^{K\tilde{G},{\cal E}}({\bf M}',\boldsymbol{\delta},a)$ est l'image de $\rho_{J}^{K\tilde{G}_{J},{\cal E}}({\bf M}',\boldsymbol{\delta},a)$ par la m\^eme application. Cela prouve la proposition 5.2 sous l'hypoth\`ese (1). 
  
  \bigskip
  
  \subsection{Comparaison des espaces $K\tilde{G}$ et $K\tilde{G}_{J}$ }
 On conserve les m\^emes notations qu'au d\'ebut du  paragraphe pr\'ec\'edent (mais on n'impose pas l'hypoth\`ese (1) de ce paragraphe).  On a not\'e $p$ l'\'el\'ement de $\Pi$ tel que $\eta\in \tilde{M}_{p}({\mathbb R})$. Pour simplifier, posons simplement $\tilde{M}=\tilde{M}_{p}$, $\tilde{G}=\tilde{G}_{p}$. Fixons une paire de Borel \'epingl\'ee ${\cal E}=(B,T,(E_{\alpha})_{\alpha\in \Delta})$ de $G$ de sorte que $M$ soit standard pour $(B,T)$. On peut supposer que $(B\cap M,T)$ est la paire de Borel pour $M$ figurant dans le diagramme fix\'e en 5.3. On note $\Sigma^G(T)$ l'ensemble des racines de $T$ dans $\mathfrak{g}$ et $\Sigma^G(A_{\tilde{M}})$ celui des racines de $A_{\tilde{M}}$ dans $\mathfrak{g}$. On peut consid\'erer que $\Sigma^G(A_{\tilde{M}})\subset \mathfrak{a}_{\tilde{M}}^*$.  On a une application de restriction
 $$\begin{array}{ccc}\Sigma^G(T)&\to&\Sigma^G(A_{\tilde{M}})\cup\{0\}\\ \alpha&\mapsto&\alpha_{A_{\tilde{M}}}\\ \end{array}$$
 Soit $J\in {\cal J}_{\tilde{M}}^{\tilde{G}}$. Cet \'el\'ement d\'etermine un ${\mathbb Z}$-module $R_{J}$ de rang $a_{\tilde{M}}-a_{\tilde{G}}$ dans $\mathfrak{a}_{\tilde{M}}^*$.  Le groupe $G_{J}$ est engendr\'e par $T$ et les sous-espaces radiciels associ\'es aux $\alpha\in \Sigma^G(T)$ tels que $\alpha_{A_{\tilde{M}}}\in R_{J}$. On note $\Sigma^{G_{J}}(T)$ l'ensemble des racines de $T$ dans $\mathfrak{g}_{J}$. On a donc $\Sigma^{G_{J}}(T)\subset \Sigma^G(T)$. La paire ${\cal E}$ est munie d'un automorphisme $\theta$.    L'ensemble $\Sigma^{G_{J}}(T)$ est invariant par $\theta$.
 Posons $B_{J}=B\cap G_{J}$. Notons $\Delta_{J}$ l'ensemble  des \'el\'ements de $\Sigma^{G_{J}}(T)$ qui sont simples pour l'ordre associ\'e \`a $B_{J}$. Les ensembles $\Delta$ et $\Delta_{J}$ contiennent tous deux l'ensemble $\Delta^M$ des racines simples de $T$ dans $\mathfrak{m}$ pour l'ordre associ\'e \`a $B^M$. Compl\'etons $(B_{J},T)$ en une paire de Borel \'epingl\'ee ${\cal E}_{J}=(B_{J},T,(E_{J,\alpha})_{\alpha\in \Delta_{J}})$ de sorte que $E_{J,\alpha}=E_{\alpha}$ pour $\alpha\in \Delta^M$. 
 
  Soit $\alpha\in \Sigma^{G_{J}}(T)$. On voit que, si $\alpha$ est de type $1$, resp. $2$, en tant qu'\'el\'ement de $\Sigma^G(T)$, $\alpha$ est encore de type $1$, resp. $2$, en tant qu'\'el\'ement de $\Sigma^{G_{J}}(T)$. Par contre, si $\alpha$ est de type $3$ en tant qu'\'el\'ement de $\Sigma^G(T)$, $\alpha$ peut \^etre de type $1$ ou $3$ en tant qu'\'el\'ement de $\Sigma^{G_{J}}(T)$. En effet, il existe $\beta\in \Sigma^G(T)$  qui est de type $2$, de sorte que $\alpha=\beta+\theta^{n_{\beta}/2}(\beta)$. Si $\beta\in \Sigma^{G_{J}}(T)$, $\alpha$ reste de type $3$ dans $\Sigma^{G_{J}}(T)$. Si $\beta\not\in \Sigma^{G_{J}}(T)$, $\alpha$ devient de type $1$ dans $\Sigma^{G_{J}}(T)$. Notons $\Sigma_{irr}^{G_{J}}(T)$ l'ensemble des $\alpha\in \Sigma^{G_{J}}(T)$ qui sont de type $3$ dans $\Sigma^G(T)$ et de type $1$ dans $\Sigma^{G_{J}}(T)$ (l'indice $irr$ \'evoquant peut-\^etre  "irr\'egulier"). Soient $e\in Z(\tilde{G},{\cal E})$ et $e_{J}\in Z(\tilde{G}_{J},{\cal E}_{J})$. On a $e_{J}=t_{J}e$, avec un $t_{J}\in T$.  L'\'egalit\'e $E_{J,\alpha}=E_{\alpha}$ pour $\alpha\in \Delta^M$ entra\^{\i}ne que $t_{J}\in Z(M)$. Montrons que
  
  (1) pour $\alpha\in \Sigma^{G_{J}}(T)$, on a $(N\alpha)(t_{J})=1$ si $\alpha\not\in \Sigma_{irr}^{G_{J}}(T)$ et $(N\alpha)(t_{J})=-1$ si $\alpha\in \Sigma_{irr}^{G_{J}}(T)$. 
  
  Preuve. Notons $\theta=ad_{e}$ et $\theta_{J}=ad_{e_{J}}$. Ce sont des automorphismes respectivement de $G$ et $G_{J}$.  Il s'en d\'eduit des automorphismes de $\Sigma^G(T)$ et $\Sigma^{G_{J}}(T)$. L'automorphisme $\theta_{J}$ de $\Sigma^{G_{J}}(T)$ est la restriction de l'automorphisme $\theta$ de $\Sigma^G(T)$. Soit $\alpha\in \Sigma^{G_{J}}(T)$, notons comme toujours $n_{\alpha}$ le plus petit entier $n\geq1$ tel que $\theta^n(\alpha)=\alpha$ (ou $\theta_{J}^n(\alpha)=\alpha$, c'est pareil). Soit $\mathfrak{u}_{\alpha}$ l'espace radiciel associ\'e \`a $\alpha$. On sait que $\theta^{n_{\alpha}}$ agit sur $\mathfrak{u}_{\alpha}$ par $1$ si $\alpha$, vu comme \'el\'ement de $\Sigma^G(T)$,  est de type $1$ ou $2$, et qu'il agit par $-1$ si $\alpha$ est de type $3$ (cf. [KS] 1.3). Une propri\'et\'e analogue vaut pour $\theta_{J}^{n_{\alpha}}$. On en d\'eduit que $\theta_{J}^{n_{\alpha}}$ co\"{\i}ncide avec $\theta^{n_{\alpha}}$  sur $\mathfrak{u}_{\alpha}$ si $\alpha\not\in \Sigma_{irr}^{G_{J}}(T)$, tandis que $\theta_{J}^{n_{\alpha}}$ co\"{\i}ncide avec $-\theta^{n_{\alpha}}$  sur $\mathfrak{u}_{\alpha}$ si $\alpha\in \Sigma_{irr}^{G_{J}}(T)$. Or, par  construction, $\theta_{J}^{n_{\alpha}}$ co\"{\i}ncide avec $(N\alpha)(t_{J})\theta^{n_{\alpha}}$ sur $\mathfrak{u}_{\alpha}$. L'assertion (1) s'ensuit. 
  
  Passons maintenant aux  groupes duaux. On fixe une paire de Borel \'epingl\'ee $\hat{{\cal E}}=(\hat{B},\hat{T},(\hat{E}_{\hat{\alpha}})_{\alpha\in \Delta})$ de $\hat{G}$. Rappelons qu'il y a une bijection $\alpha\mapsto \hat{\alpha}$ entre $\Sigma^G(T)$ et $\Sigma^{\hat{G}}(\hat{T})$. On identifie $\hat{M}$ au Levi standard associ\'e au sous-ensemble $\Delta^M$ de $\Delta$. On peut identifier $\hat{G}_{J}$ au sous-groupe de $\hat{G}$ engendr\'e par $\hat{T}$ et les sous-espaces radiciels associ\'es aux $\hat{\alpha}$ pour $\alpha\in \Sigma^{G_{J}}(T)$. On pose $\hat{B}_{J}=\hat{G}_{J}\cap \hat{B}$ et on compl\`ete $(\hat{B}_{J},\hat{T})$ en une paire de Borel \'epingl\'ee $\hat{{\cal E}}_{J}=(\hat{B}_{J},\hat{T},(\hat{E}_{J,\hat{\alpha}})_{\alpha\in \Delta_{J}})$ de sorte que $\hat{E}_{J,\hat{\alpha}}=\hat{E}_{\hat{\alpha}}$ si $\alpha\in \Delta^M$. On note $\hat{\theta}$ l'automorphisme de $\hat{G}$ associ\'e \`a $\theta$ qui conserve $\hat{{\cal E}}$ et $\sigma\mapsto \sigma_{G}$ l'action galoisienne qui conserve $\hat{{\cal E}}$. On d\'efinit de fa\c{c}on similaire l'automorphisme $\hat{\theta}_{J}$ et l'action galoisienne $\sigma\mapsto \sigma_{G_{J}}$. L'automorphisme $\hat{\theta}_{J}$ n'est pas en g\'en\'eral la restriction de $\hat{\theta}$ \`a $\hat{G}_{J}$ car cette restriction n'a pas de raison de respecter l'\'epinglage de ce groupe. Mais il existe $s_{J}\in \hat{T}$ tel que $\hat{\theta}_{J}$ soit la restriction de $ad_{s_{J}}\circ \hat{\theta}$.  L'\'egalit\'e $\hat{E}_{J,\hat{\alpha}}=\hat{E}_{\hat{\alpha}}$ si $\alpha\in \Delta^M$ entra\^{\i}ne que $s_{J}\in Z(\hat{M})$. Comme ci-dessus
  
    (2) pour $\alpha\in \Sigma^{G_{J}}(T)$, on a $(N\hat{\alpha})(s_{J})=1$ si $\alpha\not\in \Sigma_{irr}^{G_{J}}(T)$ et $(N\hat{\alpha})(s_{J})=-1$ si $\alpha\in \Sigma_{irr}^{G_{J}}(T)$. 
    
    De m\^eme l'action galoisienne $\sigma\mapsto \sigma_{G_{J}}$ n'est pas en g\'en\'eral la restriction de $\sigma\mapsto \sigma_{G}$.  Prolongeons ces actions galoisiennes en des actions de $W_{{\mathbb R}}$. La preuve de [II] 1.10(10) montre qu'il existe un cocycle $\chi:W_{{\mathbb R}}\to Z(\hat{M})$ de sorte que $w_{G_{J}}$ soit la restriction de $ad_{\chi(w)}\circ w_{G}$ pour tout $w\in W_{{\mathbb R}}$. 

  Introduisons les espaces tordus $(\hat{G}\rtimes W_{{\mathbb R}})\hat{\theta}$, $(\hat{G}_{J}\rtimes W_{{\mathbb R}})\hat{\theta}_{J}$ et $(\hat{M}\rtimes W_{{\mathbb R}})\hat{\theta}^M$, o\`u $\hat{\theta}^M$ est la restriction commune de $\hat{\theta}$ et $\hat{\theta}_{J}$ \`a $\hat{M}$. Il n'y a pas en g\'en\'eral de plongement d'espaces tordus 
  $$(\hat{G}_{J}\rtimes W_{{\mathbb R}})\hat{\theta}_{J} \to (\hat{G}\rtimes W_{{\mathbb R}})\hat{\theta}.$$
  On a toutefois deux plongements
  $$\begin{array}{ccc}&&(\hat{G}\rtimes W_{{\mathbb R}})\hat{\theta}\\ &\nearrow&\\ (\hat{M}\rtimes W_{{\mathbb R}})\hat{\theta}^M&&\\ &\searrow&\\&&(\hat{G}_{J}\rtimes W_{{\mathbb R}})\hat{\theta}_{J}\\ \end{array}$$
  Les restrictions \`a $\hat{T}$ de $\hat{\theta}$, $\hat{\theta}_{J}$ et $\hat{\theta}^M$ co\"{\i}ncident. Pour simplifier, nous noterons $\hat{\theta}$ cette restriction.

  \bigskip
  
  \subsection{Preuve de la proposition 5.2, deuxi\`eme cas}
   On conserve les notations  de 5.3, on suppose maintenant
  
  (1) $A_{M'_{\epsilon}}=A_{M'}$.
  
  Rappelons que ${\bf M}'=(M',{\cal M}',\tilde{\zeta})$. On peut supposer $\tilde{\zeta}=\zeta\hat{\theta}^M$, avec $\zeta\in \hat{T}$. Soient $J\in {\cal J}_{K\tilde{M}}^{K\tilde{G}}$ et $a\in A_{K\tilde{M}}({\mathbb R})$ en position g\'en\'erale et proche de $1$. Par d\'efinition
  $$(2) \qquad \rho_{J}^{K\tilde{G},{\cal E}}({\bf M}',\boldsymbol{\delta},a)=\sum_{s\in \zeta Z(\hat{M})^{\Gamma_{{\mathbb R}},\hat{\theta}}/Z(\hat{G})^{\Gamma_{{\mathbb R}},\hat{\theta}}}i_{\tilde{M}'}(\tilde{G},\tilde{G}'(s\hat{\theta}))$$
  $$\sum_{J'\in {\cal J}_{\tilde{M}'}^{\tilde{G}'(s\hat{\theta})}(B ); J'\mapsto J} transfert(\sigma_{J'}^{{\bf G}'(s\hat{\theta})}(\boldsymbol{\delta},\xi(a))),$$
   $$(3) \qquad \rho_{J}^{K\tilde{G}_{J},{\cal E}}({\bf M}',\boldsymbol{\delta},a)=\sum_{t\in \zeta Z(\hat{M})^{\Gamma_{{\mathbb R}},\hat{\theta}}/Z(\hat{G}_{J})^{\Gamma_{{\mathbb R}},\hat{\theta}}}i_{\tilde{M}'}(\tilde{G}_{J},\tilde{G}'_{J}(t\hat{\theta}_{J}))$$
   $$\sum_{J'\in {\cal J}_{\tilde{M}'}^{\tilde{G}_{J}'(t\hat{\theta}_{J})}(B_{J} ); J'\mapsto J} transfert(\sigma_{J'}^{{\bf G}'_{J}(t\hat{\theta}_{J})}(\boldsymbol{\delta},a)) .$$
   
   Pour simplifier, on a not\'e uniform\'ement $B$ les fonctions $B^{\tilde{G}}_{{\cal O}'}$ qui interviennent dans la premi\`ere somme et $B_{J}$ les analogues de la deuxi\`eme somme. 
   Rappelons que les sommes en $J'$ sont vides ou r\'eduites \`a un seul \'el\'ement.
   
     Soit $s\in \zeta Z(\hat{M})^{\Gamma_{{\mathbb R}},\hat{\theta}}/Z(\hat{G})^{\Gamma_{{\mathbb R}},\hat{\theta}}$.  Introduisons le syst\`eme   de racines $\Sigma^{G'(s\hat{\theta})_{\epsilon}}(T')$ de  $G'(s\hat{\theta})_{\epsilon}$ relativement au tore $T'$.  Comme en [II] 1.8, on note  $\Sigma^{G'(s\hat{\theta})_{\epsilon}}(T',B)$ l'ensemble des $B(\beta)^{-1}\beta$ pour $\beta\in \Sigma^{G'(s\hat{\theta})_{\epsilon}}(T')$. Un \'el\'ement  $\beta$ de $ \Sigma^{G'(s\hat{\theta})_{\epsilon}}(T')$ ou de $\Sigma^{G'(s\hat{\theta})_{\epsilon}}(T',B)$ peut s'identifier \`a un \'el\'ement de $\mathfrak{t}^{\theta,*}$ via l'isomorphisme $\xi:\mathfrak{t}^{\theta}\to \mathfrak{t}'$.  Notons $\Sigma^{G'(s\hat{\theta})_{\epsilon},J}(T',B)$ le sous-ensemble des $\beta\in \Sigma^{G'(s\hat{\theta})_{\epsilon}}(T',B)$ tels que $ \beta_{A_{\tilde{M}}}\in R_{J}$. Notons $\Sigma^{G'(s\hat{\theta})_{\epsilon},J}(T')$ le sous-ensemble des $\beta\in \Sigma^{G'(s\hat{\theta})_{\epsilon}}(T')$ tels que $B(\beta)^{-1}\beta\in \Sigma^{G'(s\hat{\theta})_{\epsilon},J}(T',B)$. Par d\'efinition, il existe $J'\in {\cal J}_{\tilde{M}'}^{\tilde{G}'(s\hat{\theta})}(B )$ tel que $J'\mapsto J$ si et seulement s'il existe $\beta_{1},...,\beta_{n}\in \Sigma^{G'(s\hat{\theta})_{\epsilon},J}(T',B)$, de sorte $n=a_{\tilde{M}}-a_{\tilde{G}}$ et que la famille $ (\beta_{i,A_{\tilde{M}}})_{i=1,...,n}$
      engendre   $R_{J}$. Cette description entra\^{\i}ne

     (4) supposons qu'il existe  $J'\in {\cal J}_{\tilde{M}'}^{\tilde{G}'(s\hat{\theta})}(B )$ tel que $J'\mapsto J$; alors ${\bf G}'(s\hat{\theta})$ est elliptique.
   
   En effet, pour une famille   $(\beta_{i})_{i=1,...,n}$ comme ci-dessus, la famille $ (\beta_{i,A_{\tilde{M}}})_{i=1,...,n}$ est lin\'eairement ind\'ependante, donc
     les restrictions \`a $\mathfrak{a}_{M'_{\epsilon}}$ des $\beta_{i}$ sont lin\'eairement ind\'ependantes. A fortiori, $a_{\tilde{M}}-a_{\tilde{G}}=n\leq  a_{M'_{\epsilon}}-  a_{G'(s\hat{\theta})} $. Par l'hypoth\`ese (1) et l'ellipticit\'e de ${\bf M}'$, on a $a_{M'_{\epsilon}}=a_{\tilde{M}}$, d'o\`u    $a_{G'(s\hat{\theta})}\leq a_{\tilde{G}}$. Cette in\'egalit\'e est forc\'ement une \'egalit\'e, donc ${\bf G}'(s\hat{\theta})$ est elliptique.

   Supposons qu'il existe  $J'\in {\cal J}_{\tilde{M}'}^{\tilde{G}'(s\hat{\theta})}(B )$ tel que $J'\mapsto J$. Soit $J_{s}$ l'unique \'el\'ement qui s'envoie sur $J$, introduisons le groupe $G'(s\hat{\theta})_{\epsilon,J_{s}}$ comme en [II] 3.3. On a alors $\Sigma^{G'(s\hat{\theta})_{\epsilon},J}(T')=\Sigma^{G'(s\hat{\theta})_{\epsilon,J_{s}}}(T')$.
     
     Notons $t$ l'image de $s$ par l'application naturelle
     $$\zeta Z(\hat{M})^{\Gamma_{{\mathbb R}},\hat{\theta}}/Z(\hat{G})^{\Gamma_{{\mathbb R}},\hat{\theta}}\to \zeta Z(\hat{M})^{\Gamma_{{\mathbb R}},\hat{\theta}}/Z(\hat{G}_{J})^{\Gamma_{{\mathbb R}},\hat{\theta}}.$$
   Par analogie avec ce qui pr\'ec\`ede, on d\'efinit les ensembles $\Sigma^{G'_{J}(t\hat{\theta}_{J})_{\epsilon}}(T') $ et $\Sigma^{G'_{J}(t\hat{\theta}_{J})_{\epsilon}}(T',B_{J})$.
   Montrons que
   
   (5)  on a l'\'egalit\'e $  \Sigma^{G'(s\hat{\theta})_{\epsilon},J}(T',B)=\Sigma^{G'_{J}(t\hat{\theta}_{J})_{\epsilon}}(T',B_{J})$.

  On doit une fois de plus rappeler la description du syst\`eme de racines $\Sigma^{G'(s\hat{\theta})_{\epsilon}}(T')$  ainsi que la description de la fonction $B $. On \'ecrit $\eta=\nu e$, avec $\nu\in T$ et $e\in {\cal Z}(\tilde{G},{\cal E})$. D'apr\`es [W2] 3.3 et [II] 1.11, le syst\`eme de racines est la r\'eunion des ensembles suivants
   
   (a) les $N\alpha$ pour $\alpha\in \Sigma^G(T)$ de type $1$ tels que $N\alpha(\nu)=1$ et $N\hat{\alpha}(s)=1$; on a $B(N\alpha)=n_{\alpha}$;
   
    (b) les $2N\alpha$ pour $\alpha\in \Sigma^G(T)$ de type $2$ tels que $N\alpha(\nu)=1$ et $N\hat{\alpha}(s)=1$; on a $B(2N\alpha)=2n_{\alpha}$;
   
  (c) les $2N\alpha$ pour $\alpha\in \Sigma^G(T)$ de type $2$ tels que $N\alpha(\nu)=-1$ et $N\hat{\alpha}(s)=1$;  on a $B(2N\alpha)=n_{\alpha}$;
   
    (d) les $N\alpha$ pour $\alpha\in \Sigma^G(T)$ de type $3$ tels que $N\alpha(\nu)=1$ et $N\hat{\alpha}(s)=-1$; on a $B(N\alpha)=2n_{\alpha}$.
    
  Evidemment, les $\alpha$ parcourent ici les orbites dans $\Sigma^G(T)$ pour l'action du groupe d'automorphismes engendr\'e par $\theta$.  Rappelons qu'\`a tout $\alpha_{2}$ de type $2$ est associ\'e une racine $\alpha_{3}=\alpha_{2}+\theta^{n_{\alpha_{2}}/2}(\alpha_{2})$ de type $3$. On a $N\alpha_{2}=N\alpha_{3}$, $N\hat{\alpha}_{2}=N\hat{\alpha}_{3}$ mais $n_{\alpha_{2}}=2n_{\alpha_{3}}$. Cette correspondance se quotiente en une bijection entre orbites de type $2$ et orbites de type (3). Ainsi, les cas (c) et (d) peuvent \^etre remplac\'es par
  
  (c') les $2N\alpha$ pour $\alpha\in \Sigma^G(T)$ de type $3$ tels que $N\alpha(\nu)=-1$ et $N\hat{\alpha}(s)=1$;  on a $B(2N\alpha)=2n_{\alpha}$;
   
 (d') les $N\alpha$ pour $\alpha\in \Sigma^G(T)$ de type $2$ tels que $N\alpha(\nu)=1$ et $N\hat{\alpha}(s)=-1$; on a $B(N\alpha)=n_{\alpha}$.  
 
 Pour tout $\alpha\in \Sigma^G(T)$, notons $\alpha_{res}$ sa restriction \`a $\mathfrak{t}^{\theta}$. Pour $\beta=N\alpha$, l'\'el\'ement $n_{\alpha}^{-1}\beta$ co\"{\i}ncide avec $\alpha_{res}$ comme  forme lin\'eaire sur $\mathfrak{t}^{\theta}$. De m\^eme, pour $\beta=2N\alpha$, l'\'el\'ement $(2n_{\alpha})^{-1}\beta$ co\"{\i}ncide avec $\alpha_{res}$. Dire que la restriction \`a $\mathfrak{a}_{\tilde{M}}$ appartient \`a $R_{J}$ revient \`a dire que $\alpha_{A_{\tilde{M}}}\in R_{J}$. On voit alors que le syst\`eme $  \Sigma^{G'(s\hat{\theta})_{\epsilon},J}(T',B)$ est  form\'e des $\alpha_{res}$ pour $\alpha\in \Sigma^G(T)$ tel que $\alpha_{A_{\tilde{M}}}\in R_{J}$ et tel que l'une des conditions suivantes soit v\'erifi\'ee

   (e)  $\alpha$ de type $1$, $N\alpha(\nu)=1$ et $N\hat{\alpha}(s)=1$;  
     
    (f)  $\alpha$   de type $2$,  $N\alpha(\nu)=1$ et $N\hat{\alpha}(s)=\pm 1$; 
       
  (g)  $\alpha$ de type $3$, $N\alpha(\nu)=-1$ et $N\hat{\alpha}(s)=1$.

      L'ensemble  $\Sigma^{G'_{J}(t\hat{\theta}_{J})_{\epsilon}}(T')$ se d\'ecrit par des conditions analogues \`a (a), (b), (c'), (d'). Le syst\`eme $\Sigma^G(T)$ doit \^etre remplac\'e par $\Sigma^{G_{J}}(T)$, c'est-\`a-dire le sous-ensemble des $\alpha\in \Sigma^G(T)$ tels que $\alpha_{A_{\tilde{M}}}\in R_{J}$. Les types $1$, $2$ et $3$ sont relatifs \`a ce syst\`eme et, pour  les distinguer des pr\'ec\'edents, on note ces types $1_{J}$, $2_{J}$ et $3_{J}$. Le terme $s$ ne change pas. Par contre, l'\'ecriture $\eta=\nu e$ est remplac\'ee par $\eta=\nu_{J}e_{J}$, avec $\nu_{J}\in T$ et $e_{J}\in {\cal Z}(\tilde{G}_{J},{\cal E}_{J})$. On a $\nu_{J}=\nu t_{J}^{-1}$, o\`u $t_{J}$ v\'erifie    5.4(1). Le calcul se poursuit et on obtient que $\Sigma^{G'_{J}(t\hat{\theta}_{J})_{\epsilon}}(T',B_{J})$ est  form\'e des $\alpha_{res}$ pour $\alpha\in \Sigma^G(T)$ tel que $\alpha_{A_{\tilde{M}}}\in R_{J}$ et tel que l'une des conditions suivantes soit v\'erifi\'ee

     $(e)_{J}$  $\alpha$ de type $1_{J}$, $N\alpha(t_{J}^{-1}\nu)=1$ et $N\hat{\alpha}(s)=1$;  
     
    $(f)_{J}$  $\alpha$   de type $2_{J}$,  $N\alpha(t_{J}^{-1}\nu)=1$ et $N\hat{\alpha}(s)=\pm 1$; 
       
  $(g)_{J}$  $\alpha$ de type $3_{J}$, $N\alpha(t_{J}^{-1}\nu)=-1$ et $N\hat{\alpha}(s)=1$.

     Soit $\alpha \in \Sigma^G(T)$ tel que $\alpha_{A_{\tilde{M}}}\in R_{J}$. Supposons d'abord que son type soit le m\^eme dans $\Sigma^G(T)$ et dans $\Sigma^{G_{J}}(T)$. D'apr\`es 5.4(1), on a $N\alpha(t_{J})=1$. Les conditions (e), (f) , (g) sont alors respectivement \'equivalentes \`a $(e)_{J}$, $(f)_{J}$, $(g)_{J}$. Supposons maintenant que le type de $\alpha$ change. On a vu qu'alors $\alpha$ est de type (3) et de type $(1)_{J}$. D'apr\`es 5.4(1), on a $N\alpha(t_{J})=-1$. Mais alors, les conditions (g) et $(e)_{J}$ (qui sont les seules pouvant concerner $\alpha$) sont \'equivalentes. Cela prouve l'\'egalit\'e (5).

    Notons ${\cal Z}$ l'ensemble des $s\in \zeta Z(\hat{M})^{\Gamma_{{\mathbb R}},\hat{\theta}}/Z(\hat{G})^{\Gamma_{{\mathbb R}},\hat{\theta}}$ tels que $i_{\tilde{M}'}(\tilde{G},\tilde{G}'(s\hat{\theta}))\not=0$ (c'est-\`a-dire que ${\bf G}'(s\hat{\theta})$ est elliptique) et qu'il existe  $J'\in {\cal J}_{\tilde{M}'}^{\tilde{G}'(s\hat{\theta})}(B )$ tel que $J'\mapsto J$. Notons  ${\cal Z}_{J}$ l'ensemble des $t\in \zeta Z(\hat{M})^{\Gamma_{{\mathbb R}},\hat{\theta}}/Z(\hat{G}_{J})^{\Gamma_{{\mathbb R}},\hat{\theta}}$ tels que $i_{\tilde{M}'}(\tilde{G}_{J},\tilde{G}_{J}'(t\hat{\theta}_{J}))\not=0$   et qu'il existe  $J'\in {\cal J}_{\tilde{M}'}^{\tilde{G}_{J}'(s\hat{\theta}_{J})}(B _{J})$ tel que $J'\mapsto J$. Montrons que
   
   (6) ${\cal Z}$ est l'image r\'eciproque de ${\cal Z}_{J}$ par la projection naturelle
   $$\zeta Z(\hat{M})^{\Gamma_{{\mathbb R}},\hat{\theta}}/Z(\hat{G})^{\Gamma_{{\mathbb R}},\hat{\theta}}\to \zeta Z(\hat{M})^{\Gamma_{{\mathbb R}},\hat{\theta}}/Z(\hat{G}_{J})^{\Gamma_{{\mathbb R}},\hat{\theta}}.$$
   
   L'assertion (4) montre que, dans la d\'efinition de ${\cal Z}$, la condition que ${\bf G}'(s\hat{\theta})$ est elliptique est superflue: elle est entra\^{\i}n\'ee par l'existence de $J'$.  La condition $s\in {\cal Z}$ \'equivaut donc \`a ce qu'il  existe $\beta_{1},...,\beta_{n}\in \Sigma^{G'(s\hat{\theta})_{\epsilon},J}(T',B)$, de sorte $n=a_{\tilde{M}}-a_{\tilde{G}}$ et que la famille $ (\beta_{i,A_{\tilde{M}}})_{i=1,...,n}$
      engendre   $R_{J}$.
   De m\^eme, la condition $t\in {\cal Z}_{J}$ \'equivaut \`a ce qu'il  existe $\beta_{1},...,\beta_{n}\in \Sigma^{G_{J}'(t\hat{\theta}_{J})_{\epsilon}}(T',B_{J})$, de sorte $n=a_{\tilde{M}}-a_{\tilde{G}}$ et que la famille $ (\beta_{i,A_{\tilde{M}}})_{i=1,...,n}$
      engendre   $R_{J}$. L'assertion (5) montre que la condition pour $s$ \'equivaut \`a la condition pour $t$, o\`u $t$ est l'image de $s$ par la projection naturelle. Cela prouve (6).
   
 Posons $\hat{Z}=Z(\hat{G})\cap \hat{T}^{\Gamma_{{\mathbb R}},\hat{\theta},0}$.  Soit $s\in {\cal Z}$, notons $t$ son image dans ${\cal Z}_{J}$. Le groupe $\hat{Z}$ est d'indice fini dans $ Z(\hat{G}'(s\hat{\theta}))^{\Gamma_{{\mathbb R}},\hat{\theta}}$ comme dans $Z(\hat{G}'_{J}(t\hat{\theta}_{J}))^{\Gamma_{{\mathbb R}},\hat{\theta}}$. On pose
 $$[Z(\hat{G}'_{J}(t\hat{\theta}_{J}))^{\Gamma_{{\mathbb R}},\hat{\theta}}:Z(\hat{G}'(s\hat{\theta}))^{\Gamma_{{\mathbb R}},\hat{\theta}}]=[Z(\hat{G}'_{J}(t\hat{\theta}_{J}))^{\Gamma_{{\mathbb R}},\hat{\theta}} :\hat{Z}]
   [ Z(\hat{G}'(s\hat{\theta}))^{\Gamma_{{\mathbb R}},\hat{\theta}}:\hat{Z}]^{-1}.$$
 Montrons que  
    
   (7)  on a l'\'egalit\'e
   $$i_{\tilde{M}'}(\tilde{G}_{J},\tilde{G}_{J}'(t\hat{\theta}_{J}))[Z(\hat{G}'_{J}(t\hat{\theta}_{J}))^{\Gamma_{{\mathbb R}},\hat{\theta}}  
   :Z(\hat{G}'(s\hat{\theta}))^{\Gamma_{{\mathbb R}},\hat{\theta}}]=[Z(\hat{G}_{J})^{\Gamma_{{\mathbb R}},\hat{\theta}}:Z(\hat{G})^{\Gamma_{{\mathbb R}},\hat{\theta}}]i_{\tilde{M}'}(\tilde{G},\tilde{G}'(s\hat{\theta})).$$
   
   On a un diagramme commutatif
   $$\begin{array}{ccc}Z(\hat{M})^{\Gamma_{{\mathbb R}},\hat{\theta}}/Z(\hat{G})^{\Gamma_{{\mathbb R}},\hat{\theta}}&\stackrel{A}{\to}&Z(\hat{M})^{\Gamma_{{\mathbb R}},\hat{\theta}}/Z(\hat{G}_{J})^{\Gamma_{{\mathbb R}},\hat{\theta}}\\ B\downarrow\,&&C \downarrow\,\\ Z(\hat{M}')^{\Gamma_{{\mathbb R}}}/\hat{Z}&\stackrel{D}{\to}&Z(\hat{M}')^{\Gamma_{{\mathbb R}}}/Z(\hat{G}_{J}'(t\hat{\theta}_{J}))^{\Gamma_{{\mathbb R}}}\\ E\downarrow\,&&\\ Z(\hat{M}')^{\Gamma_{{\mathbb R}}}/Z(\hat{G}'(s\hat{\theta}))^{\Gamma_{{\mathbb R}}}&&\\ \end{array}$$
   Toutes les fl\`eches sont surjectives et de noyaux finis. On en d\'eduit les \'egalit\'es
   $$\vert ker(A)\vert \vert ker(C)\vert =\vert ker(B)\vert \vert ker(D)\vert =\vert ker(EB)\vert \vert ker(E)\vert ^{-1}\vert ker(D)\vert .$$
   Par d\'efinition, on a
   $$i_{\tilde{M}'}(\tilde{G}_{J},\tilde{G}_{J}'(t\hat{\theta}_{J}))=\vert ker(C)\vert ^{-1},$$
   $$[Z(\hat{G}'_{J}(t\hat{\theta}_{J}))^{\Gamma_{{\mathbb R}},\hat{\theta}}  
   :Z(\hat{G}'(s\hat{\theta}))^{\Gamma_{{\mathbb R}},\hat{\theta}}]=\vert ker(D)\vert \vert ker(E)\vert ^{-1},$$
   $$[Z(\hat{G}_{J})^{\Gamma_{{\mathbb R}},\hat{\theta}}:Z(\hat{G})^{\Gamma_{{\mathbb R}},\hat{\theta}}]=\vert ker(A)\vert ,$$
   $$i_{\tilde{M}'}(\tilde{G},\tilde{G}'(s\hat{\theta}))=\vert ker(EB)\vert ^{-1}.$$
   L'assertion (7) s'en d\'eduit.

 Pour $s\in {\cal Z}$, resp. $t \in {\cal Z}_{J}$, notons $J_{s}$, resp. $J_{t}$, l'unique \'el\'ement de ${\cal J}_{\tilde{M}'}^{\tilde{G}'(s\hat{\theta})}(B)$, resp. de ${\cal J}_{\tilde{M}'}^{\tilde{G}'_{J}(t\hat{\theta}_{J})}(B_{J})$, qui s'envoie sur $J$. 
 
 \ass{Lemme}{Soit $s\in {\cal Z}$, notons $t$ son image dans ${\cal Z}_{J}$. On a l'\'egalit\'e  
   $$[Z(\hat{G}'_{J}(t\hat{\theta}_{J}))^{\Gamma_{{\mathbb R}}}:Z(\hat{G}'(s\hat{\theta}))^{\Gamma_{{\mathbb R}}}]transfert(\sigma_{J_{s}}^{{\bf G}'(s\hat{\theta})}(\boldsymbol{\delta},\xi(a)))=transfert(\sigma_{J_{t}}^{{\bf G}'_{J}(t\hat{\theta}_{J})}(\boldsymbol{\delta},\xi(a))).$$}
   
   Nous prouverons ce lemme au paragraphe suivant. Admettons-le et achevons la preuve de la proposition  5.2. La formule (2) se r\'ecrit
   $$\rho_{J}^{K\tilde{G},{\cal E}}({\bf M}',\boldsymbol{\delta},a)=\sum_{s\in {\cal Z}}i_{\tilde{M}'}(\tilde{G},\tilde{G}'(s\hat{\theta}))transfert(\sigma_{J_{s}}^{{\bf G}'(s\hat{\theta})}(\boldsymbol{\delta},\xi(a))).$$
   En utilisant (7) et (8), on obtient
   $$\rho_{J}^{K\tilde{G},{\cal E}}({\bf M}',\boldsymbol{\delta},a)=[Z(\hat{G}_{J})^{\Gamma_{{\mathbb R}},\hat{\theta}}:Z(\hat{G})^{\Gamma_{{\mathbb R}},\hat{\theta}}]^{-1}\sum_{s\in {\cal Z}}i_{\tilde{M}'}(\tilde{G}_{J},\tilde{G}_{J}'(t\hat{\theta}_{J}))transfert(\sigma_{J_{t}}^{{\bf G}'_{J}(t\hat{\theta}_{J})}(\boldsymbol{\delta},\xi(a))),$$
   o\`u $t$ est la projection de $s$ dans ${\cal Z}_{J}$. 
   Le nombre d'\'el\'ements du noyau de cette projection ${\cal Z}\to {\cal Z}_{J}$ est \'egal \`a  
   $[Z(\hat{G}_{J})^{\Gamma_{{\mathbb R}},\hat{\theta}}:Z(\hat{G})^{\Gamma_{{\mathbb R}},\hat{\theta}}]$. La formule ci-dessus devient
   $$\rho_{J}^{K\tilde{G},{\cal E}}({\bf M}',\boldsymbol{\delta},a)=\sum_{t\in {\cal Z}_{J}}i_{\tilde{M}'}(\tilde{G}_{J},\tilde{G}_{J}'(t\hat{\theta}_{J}))transfert(\sigma_{J_{t}}^{{\bf G}'_{J}(t\hat{\theta}_{J})}(\boldsymbol{\delta},\xi(a))).$$
   Mais le membre de droite est \'egal \`a celui de la formule (3), donc est \'egal \`a $\rho_{J}^{K\tilde{G}_{J},{\cal E}}({\bf M}',\boldsymbol{\delta},a)$. Cela prouve la proposition 5.2.
   
   \bigskip
   
   \subsection{Preuve du lemme 5.5}
   
   On fixe $s\in {\cal Z}$, on note $t$  son image dans ${\cal Z}_{J}$.  Pour simplifier, on pose ${\bf G}'={\bf G}' (s\hat{\theta})$ et ${\bf G}'_{J}={\bf G}'_{J}(t\hat{\theta}_{J})$. On fixe des donn\'ees auxiliaires   $G'_{1}$,...,$\Delta_{1}$  pour la donn\'ee ${\bf G}'$ et $G'_{J,1}$,...,$\Delta_{J,1}$ pour la donn\'ee ${\bf G}'_{J}$. On note $\tilde{M}'_{1}$ et $\tilde{M}'_{J,1}$, resp. $\tilde{T}'_{1}$ et $\tilde{T}'_{J,1}$,  les images r\'eciproques de $\tilde{M}'$, resp. $\tilde{T}'$, dans $\tilde{G}'_{1}$ et $\tilde{G}'_{J,1}$. On fixe des \'el\'ements $\epsilon_{1}\in \tilde{T}'_{1}({\mathbb R})$ et $\epsilon_{J,1}\in \tilde{T}'_{J,1}({\mathbb R})$ se projetant sur $\epsilon$. On note ${\cal O}'_{1}$ la classe de conjugaison stable de $\epsilon_{1}$ dans $\tilde{M}'_{1}({\mathbb R})$ et ${\cal O}'_{J,1}$ celle de $\epsilon_{J,1}$ dans $\tilde{M}'_{J,1}({\mathbb R})$. Notons $M'_{\epsilon,sc}$ l'image r\'eciproque de $M'_{\epsilon}$ dans $G'_{\epsilon,SC}$. Puisque $G' _{\epsilon,SC}=G'_{1,\epsilon_{1},SC}$, c'est aussi l'image r\'eciproque $M'_{1,\epsilon_{1},sc}$ de $M'_{1,\epsilon_{1}}$ dans $G'_{1,\epsilon_{1},SC}$. 
  On a une suite d'homomorphismes
  $$(1) \qquad D_{unip}^{st}(M'_{\epsilon,SC}({\mathbb R}))\stackrel{\iota^*_{M'_{\epsilon,SC},M'_{\epsilon,sc}}}{\to}D_{unip}^{st}(M'_{\epsilon,sc}({\mathbb R}))
   \simeq D_{unip}^{st}(M'_{1,\epsilon_{1},sc}({\mathbb R}))$$
   $$\stackrel{\iota^*_{M'_{1,\epsilon_{1},sc},M'_{1,\epsilon_{1}}}}{\to}D_{ unip}^{st}(M'_{1,\epsilon_{1}}({\mathbb R})) \stackrel{desc_{\epsilon_{1}}^{st,\tilde{M}'_{1},*}}{\to}D_{g\acute{e}om}^{st}({\cal O}'_{1})\to D_{g\acute{e}om,\lambda_{1}}^{st}(\tilde{M}'_{1}({\mathbb R}),{\cal O}'_{1})\simeq D_{g\acute{e}om}^{st}({\bf M}',{\cal O}').$$
  Introduisons le groupe $\Xi_{\epsilon}=Z_{M'}(\epsilon)/M'_{\epsilon}$. Son sous-groupe d'invariants $\Xi_{\epsilon}^{\Gamma_{{\mathbb R}}}$ agit sur chacun des espaces ci-dessus, de la fa\c{c}on suivante. Il agit trivialement sur $D_{g\acute{e}om}^{st}({\bf M}',{\cal O}')$ et $D_{g\acute{e}om,\lambda_{1}}^{st}(\tilde{M}'_{1}({\mathbb R}),{\cal O}'_{1})$. Posons $ C_{1,\sharp}=\{c\in C_{1}({\mathbb R}); c\epsilon_{1}\in {\cal O}'_{1}\}$. Soit $x\in Z_{M'}(\epsilon)$ dont l'image dans $\Xi_{\epsilon}$ soit fixe par $\Gamma_{{\mathbb R}}$. Alors $ad_{x}(\epsilon_{1})=c_{1}(x)\epsilon_{1}$ o\`u $c_{1}(x)\in C_{1,\sharp}$. L'appication $x\mapsto c_{1}(x)$ se quotiente en un homomorphisme $\Xi_{\epsilon}^{\Gamma_{{\mathbb R}}}\to C_{1,\sharp}$. Le groupe $C_{1,\sharp}$ agit par multiplication sur ${\cal O}'_{1}$, donc aussi sur $D_{g\acute{e}om}^{st}({\cal O}'_{1})$. On tord cette action par la restriction \`a $C_{1,\sharp}$ du caract\`ere $\lambda_{1}$. Via l'homomorphisme pr\'ec\'edent, on obtient une action de $\Xi_{\epsilon}^{\Gamma_{{\mathbb R}}}$ sur $D_{g\acute{e}om}^{st}({\cal O}'_{1})$. Pour $x$ comme ci-dessus, l'application $ad_{x}$ pr\'eserve $M'_{1,\epsilon_{1}}$. Ce groupe est quasi-d\'eploy\'e. Fixons-en une paire de Borel \'epingl\'ee d\'efinie sur ${\mathbb R}$. Quitte \`a multiplier $x$ par un \'el\'ement de $M'_{\epsilon}$, on peut supposer que $ad_{x}$ respecte cette paire de Borel \'epingl\'ee. Alors $ad_{x}$ est un automorphisme de $M'_{1,\epsilon_{1}}$ qui est d\'efini sur ${\mathbb R}$. On obtient ainsi une action de $\Xi_{\epsilon}^{\Gamma_{{\mathbb R}}}$ sur $M'_{1,\epsilon_{1}}$ par automorphismes d\'efinis sur ${\mathbb R}$. D'o\`u aussi une action sur $D_{ unip}^{st}(M'_{1,\epsilon_{1}}({\mathbb R}))$. On tord cette action par le caract\`ere $x\mapsto \lambda_{1}(c_{1}(x))$ et on obtient l'action cherch\'ee de $\Xi_{\epsilon}^{\Gamma_{{\mathbb R}}}$ sur $D_{ unip}^{st}(M'_{1,\epsilon_{1}}({\mathbb R}))$. La m\^eme connstruction d\'efinit des actions de ce groupe sur   $ D_{unip}^{st}(M'_{1,\epsilon_{1},sc}({\mathbb R}))$, $D_{unip}^{st}(M'_{\epsilon,sc}({\mathbb R}))$ et $D_{unip}^{st}(M'_{\epsilon,SC}({\mathbb R}))$. Les homomorphismes de la suite ci-dessus sont \'equivariants pour les actions de $\Xi_{\epsilon}^{\Gamma_{{\mathbb R}}}$. En notant par un exposant  les sous-espaces d'invariants, on obtient une suite d'homomorphismes
    $$(2) \quad D_{unip}^{st}(M'_{\epsilon,SC}({\mathbb R}))^{\Xi_{\epsilon}^{\Gamma_{{\mathbb R}}}}\stackrel{\iota^*_{M'_{\epsilon,SC},M'_{\epsilon,sc}}}{\to}D_{unip}^{st}(M'_{\epsilon,sc}({\mathbb R}))^{\Xi_{\epsilon}^{\Gamma_{{\mathbb R}}}}
   \simeq D_{unip}^{st}(M'_{1,\epsilon_{1},sc}({\mathbb R}))^{\Xi_{\epsilon}^{\Gamma_{{\mathbb R}}}}$$
   $$\stackrel{\iota^*_{M'_{1,\epsilon_{1},sc},M'_{1,\epsilon_{1}}}}{\to}D_{ unip}^{st}(M'_{1,\epsilon_{1}}({\mathbb R}))^{\Xi_{\epsilon}^{\Gamma_{{\mathbb R}}}} \stackrel{desc_{\epsilon_{1}}^{st,\tilde{M}'_{1},*}}{\to}D_{g\acute{e}om}^{st}({\cal O}'_{1})^{\Xi_{\epsilon}^{\Gamma_{{\mathbb R}}}}$$
   $$\to D_{g\acute{e}om,\lambda_{1}}^{st}(\tilde{M}'_{1}({\mathbb R}),{\cal O}'_{1})\simeq D_{g\acute{e}om}^{st}({\bf M}',{\cal O}').$$
    Maintenant, les homomorphismes de cette suite sont injectifs. En effet, les homomorphismes  $\iota^*_{M'_{\epsilon,SC},M'_{\epsilon,sc}}$ et $\iota^*_{M'_{1,\epsilon_{1},sc},M'_{1,\epsilon_{1}}}$ sont injectifs d'apr\`es 3.3 (4), m\^eme sans passer aux invariants. L''application $desc_{\epsilon_{1}}^{st,\tilde{M}'_{1},*}$ ci-dessus est injective d'apr\`es [I] 4.8. Enfin, l'injectivit\'e de l'homomorphisme $D_{g\acute{e}om}^{st}({\cal O}'_{1})^{\Xi_{\epsilon}^{\Gamma_{{\mathbb R}}}}\to D_{g\acute{e}om,\lambda_{1}}^{st}(\tilde{M}'_{1}({\mathbb R}),{\cal O}'_{1})$ se v\'erifie ais\'ement sur les d\'efinitions. 
    
      D'apr\`es la d\'efinition de 2.1, notre distribution $\boldsymbol{\delta}\in D_{ tr-orb}^{st}({\bf M}',{\cal O}')$ se remonte en un \'el\'ement $\boldsymbol{\delta}_{1}\in D_{ tr-orb}^{st}({\cal O}'_{1})$.  On voit que l'action que l'on a d\'efinie de $\Xi_{\epsilon}^{\Gamma_{{\mathbb R}}}$ sur $D_{g\acute{e}om}^{st}({\cal O}'_{1})$ respecte le sous-espace $D_{ tr-orb}^{st}({\cal O}'_{1})$. Quitte \`a moyenner par cette action, on peut supposer $\boldsymbol{\delta}_{1}\in D_{tr-orb}^{st}({\cal O}'_{1})^{\Xi_{\epsilon}^{\Gamma_{{\mathbb R}}}}$. D'apr\`es le lemme 4.3(ii), $\boldsymbol{\delta}_{1}$ se rel\`eve en un \'el\'ement $\boldsymbol{\delta}_{\epsilon_{1}}\in D_{tr-unip}(M'_{1,\epsilon_{1}}({\mathbb R}))$. Pour la m\^eme raison, on peut supposer $\boldsymbol{\delta}_{\epsilon_{1}}\in D_{tr-unip}(M'_{1,\epsilon_{1}}({\mathbb R}))^{\Xi_{\epsilon}^{\Gamma_{{\mathbb R}}}}$. D'apr\`es le lemme 3.3, $\boldsymbol{\delta}_{\epsilon_{1}}$ se rel\`eve en un \'el\'ement $\boldsymbol{\delta}_{sc}\in D^{st}_{tr-unip}(M'_{\epsilon,sc}({\mathbb R}))$ dont on peut encore supposer qu'il est invariant par $\Xi_{\epsilon}^{\Gamma_{{\mathbb R}}}$. Enfin, toujours d'apr\`es le lemme 3.3, $\boldsymbol{\delta}_{sc}$ se rel\`eve en un \'el\'ement 
     $\boldsymbol{\delta}_{SC}\in D^{st}_{tr-unip}(M'_{\epsilon,SC}({\mathbb R}))$, dont on peut encore supposer qu'il est invariant par $\Xi_{\epsilon}^{\Gamma_{{\mathbb R}}}$. L'injectivit\'e des homomorphismes de la suite (10) entra\^{\i}ne que toutes ces distributions sont uniquement d\'etermin\'ees par $\boldsymbol{\delta}$. 
      
      L'\'el\'ement $J_{s}\in {\cal J}_{\tilde{M}'}^{\tilde{G}'}(B)$ s'identifie \`a un \'el\'ement  de ${\cal J}_{\tilde{M}'_{1}}^{\tilde{G}'_{1}}(B)$, puis \`a un \'el\'ement de ${\cal J}_{M'_{1,\epsilon_{1}}}^{G'_{1,\epsilon_{1}}}(B)$, puis \`a un \'el\'ement de ${\cal J}_{M'_{\epsilon,sc}}^{G'_{\epsilon,SC}}(B)$.        Les formalit\'es que l'on a pass\'ees entra\^{\i}nent que $\sigma_{J_{s}}^{{\bf G}'}(\boldsymbol{\delta},\xi(a))$ est l'image naturelle dans $D_{g\acute{e}om}^{st}({\bf M}',{\cal O}')/Ann_{{\cal O}'}^{st,{\bf G}'}$ de l'\'el\'ement $\sigma_{J_{s}}^{\tilde{G}'_{1}}(\boldsymbol{\delta}_{1},a_{1})$, o\`u $a_{1}$ est un \'el\'ement quelconque de $A_{M'_{1}}({\mathbb R})$ se projetant sur $\xi(a)\in A_{M'}({\mathbb R})$. L'hypoth\`ese (1) de 5.5 permet d'appliquer la relation  4.6 (9): $\sigma_{J_{s}}^{\tilde{G}'_{1}}(\boldsymbol{\delta}_{1},a_{1})$ est l'image par $desc_{\epsilon_{1}}^{st,\tilde{M}'_{1},*}$ de $e_{\tilde{M}_{1}}^{\tilde{G}_{1}'}(\epsilon_{1})\sigma_{J_{s}}^{G'_{1,\epsilon_{1}}}(\boldsymbol{\delta}_{\epsilon_{1}},a_{1})$. Un calcul facile, cf. [III] 7.1(4), montre que $e_{\tilde{M}_{1}}^{\tilde{G}_{1}'}(\epsilon_{1})=e_{\tilde{M}'}^{\tilde{G}'}(\epsilon)$. 
      D'apr\`es le lemme 3.4, $\sigma_{J_{s}}^{G'_{1,\epsilon_{1}}}(\boldsymbol{\delta}_{\epsilon_{1}},a_{1})$ est l'image par $\iota^*_{M'_{1,\epsilon_{1},sc},M'_{1,\epsilon_{1}}}$ de $\sigma_{J_{s}}^{G'_{\epsilon,SC}}(\boldsymbol{\delta}_{sc},a')$, o\`u $a'$ est un \'el\'ement quelconque de $A_{M'_{\epsilon,sc}}({\mathbb R})$ ayant m\^eme projection que $a_{1}$ dans $A_{M'_{1}}({\mathbb R})/A_{G'_{1}}({\mathbb R})$, ou encore ayant m\^eme projection que $\xi(a)$ dans $A_{M'}({\mathbb R})/A_{G'}({\mathbb R})$.        En r\'esum\'e, $\sigma_{J_{s}}^{{\bf G}'}(\boldsymbol{\delta},\xi(a))$ provient par la suite (1) de l'\'el\'ement $e_{\tilde{M}'}^{\tilde{G}'}(\epsilon) \sigma_{J_{s}}^{G'_{\epsilon,SC}}(\boldsymbol{\delta}_{sc},a')\in D_{unip}^{st}(M'_{\epsilon,sc}({\mathbb R}))/Ann_{unip}^{st,G'_{\epsilon,SC}}$.

  Des constructions analogues valent pour la donn\'ee ${\bf G}'_{J}$. On note cette fois $M'_{\epsilon,J,sc}$ l'image r\'eciproque de $M'_{\epsilon}$ dans $G'_{J,\epsilon,SC}$. 
De fa\c{c}on analogue \`a (1), on a  une suite
 $$(3) \qquad D_{unip}^{st}(M'_{\epsilon,SC}({\mathbb R}))\stackrel{\iota^*_{M'_{\epsilon,SC},M'_{\epsilon,J,sc}}}{\to}D^{st}_{unip}(M'_{\epsilon,J,sc}({\mathbb R}))
 \simeq D_{unip}^{st}(M'_{J,1,\epsilon_{J,1},sc}({\mathbb R}))$$
 $$\stackrel{\iota^*_{M'_{J,1,\epsilon_{J,1},sc},M'_{J,1,\epsilon_{J,1}}}}{\to}D_{ unip}^{st}(M'_{J,1,\epsilon_{J,1}}({\mathbb R}))\stackrel{desc_{\epsilon_{J,1}}^{st,\tilde{M}'_{J,1},*}}{\to}D_{g\acute{e}om}^{st}({\cal O}'_{J,1})$$
  $$\to D_{g\acute{e}om,\lambda_{J,1}}^{st}(\tilde{M}'_{J,1}({\mathbb R}),{\cal O}'_{J,1})\simeq D_{g\acute{e}om}^{st}({\bf M}',{\cal O}').$$
  On a aussi une suite analogue \`a (2)
  $$(4) \qquad D_{unip}^{st}(M'_{\epsilon,SC}({\mathbb R}))^{\Xi_{\epsilon}^{\Gamma_{{\mathbb R}}}}\stackrel{\iota^*_{M'_{\epsilon,SC},M'_{\epsilon,J,sc}}}{\to}D^{st}_{unip}(M'_{\epsilon,J,sc}({\mathbb R}))^{\Xi_{\epsilon}^{\Gamma_{{\mathbb R}}}}
 \simeq D_{unip}^{st}(M'_{J,1,\epsilon_{J,1},sc}({\mathbb R}))^{\Xi_{\epsilon}^{\Gamma_{{\mathbb R}}}}$$
 $$\stackrel{\iota^*_{M'_{J,1,\epsilon_{J,1},sc},M'_{J,1,\epsilon_{J,1}}}}{\to}D_{ unip}^{st}(M'_{J,1,\epsilon_{J,1}}({\mathbb R}))^{\Xi_{\epsilon}^{\Gamma_{{\mathbb R}}}}\stackrel{desc_{\epsilon_{J,1}}^{st,\tilde{M}'_{J,1},*}}{\to}D_{g\acute{e}om}^{st}({\cal O}'_{J,1})^{\Xi_{\epsilon}^{\Gamma_{{\mathbb R}}}}$$
  $$\to D_{g\acute{e}om,\lambda_{J,1}}^{st}(\tilde{M}'_{J,1}({\mathbb R}),{\cal O}'_{J,1})\simeq D_{g\acute{e}om}^{st}({\bf M}',{\cal O}').$$
  L'\'el\'ement $\boldsymbol{\delta}$ se remonte successivement en \'el\'ements $\boldsymbol{\delta}_{J,1}$, $\boldsymbol{\delta}_{\epsilon_{J,1}}$, $\boldsymbol{\delta}_{J,sc}$ et $\boldsymbol{\delta}_{J,SC}$. Le m\^eme calcul que plus haut montre que $\sigma_{J_{t}}^{{\bf G}_{J}'}(\boldsymbol{\delta},\xi(a))$ provient par la suite (3) de l'\'el\'ement $e_{\tilde{M}'}^{\tilde{G}_{J}'}(\epsilon)i_{J_{t}}^{G'_{J,\epsilon,SC}}\sigma_{J_{t}}^{G'_{J,\epsilon,SC}}(\boldsymbol{\delta}_{J,sc},a'_{J})\in D_{unip}^{st}(M'_{\epsilon,J,sc}({\mathbb R}))/Ann_{unip}^{st,G'_{J,\epsilon,SC}}$, o\`u $a'_{J}$ est un \'el\'ement quelconque de $A_{M'_{\epsilon,J, sc}}({\mathbb R})$ ayant m\^eme projection   que $\xi(a)$ dans $A_{M'}({\mathbb R})/A_{G'_{J}}({\mathbb R})$.

  Les compos\'es des  suites (1) et (3) ne sont pas \'egaux. L'isomorphisme compos\'e
  $$ D_{g\acute{e}om,\lambda_{1}}^{st}(\tilde{M}'_{1}({\mathbb R}),{\cal O}'_{1})\simeq D_{g\acute{e}om}^{st}({\bf M}',{\cal O}')\simeq D_{g\acute{e}om,\lambda_{J,1}}^{st}(\tilde{M}'_{J,1}({\mathbb R}),{\cal O}'_{J,1})$$
  est donn\'e par la fonction de transition $\tilde{\lambda}$ entre les deux s\'eries de donn\'ees auxiliaires pour ${\bf M}'$ d\'eduites des donn\'ees auxiliaires pour ${\bf G}'$ et pour ${\bf G}'_{J}$. En g\'en\'eral, cet isomorphisme est plus difficile \`a calculer que dans le cas non-archim\'edien. Quand des distributions font intervenir des d\'eriv\'ees dans les directions centrales de $M'_{1}$ ou $M'_{J,1}$, la formule fait intervenir des d\'eriv\'ees de cette fonction de transition. Mais ce n'est pas le cas pour les distributions provenant de $M'_{\epsilon,SC}$. Pour celles-ci, le calcul est le m\^eme que dans le cas archim\'edien. On obtient que les compos\'es des suites (1) et (3) sont proportionnels, 
  la constante de proportionnalit\'e \'etant la valeur $\tilde{\lambda}(\epsilon_{1},\epsilon_{J,1})$.  On note simplement $c$ cette valeur. Pr\'ecis\'ement, pour $\boldsymbol{\tau}_{SC}\in  D_{unip}^{st}(M'_{\epsilon,SC}({\mathbb R}))$, l'image de $\boldsymbol{\tau}_{SC}$ par (1)  est \'egale \`a l'image de $c\boldsymbol{\tau}_{SC}$ par  (3). D'autre part, dans les suites (2) et (4) interviennent des actions de $\Xi_{\epsilon}^{\Gamma_{{\mathbb R}}}$ sur $D_{unip}^{st}(M'_{\epsilon,SC}({\mathbb R}))$. Ces actions sont les m\^emes d'apr\`es [I] 4.8(2). Puisque les suites (2) et (4) sont injectives, on en d\'eduit l'\'egalit\'e $\boldsymbol{\delta}_{J,SC}=c\boldsymbol{\delta}_{SC}$. 
  
Posons
$$C= e_{\tilde{M}'}^{\tilde{G}_{J}'}(\epsilon)e_{\tilde{M}'}^{\tilde{G}'}(\epsilon)^{-1} [Z(\hat{G}'_{J})^{\Gamma_{{\mathbb R}}}:
Z(\hat{G}')^{\Gamma_{{\mathbb R}}}]^{-1}.$$

Supposons d\'emontr\'ee l'assertion suivante, o\`u $a'$ et $a'_{J}$ sont comme ci-dessus:
  
  (5) il existe un \'el\'ement $\boldsymbol{\tau}\in D_{unip}^{st}(M'_{\epsilon,SC}({\mathbb R}))$ dont l'image dans $D_{unip}^{st}(M'_{\epsilon,sc}({\mathbb R}))/Ann_{unip}^{st,G'_{\epsilon,SC}}$ soit $\sigma_{J_{s}}^{G'_{\epsilon,SC}}(\boldsymbol{\delta}_{sc},a')$ et dont l'image dans 
$D_{unip}^{st}(M'_{\epsilon,J,sc}({\mathbb R}))/Ann_{unip}^{st,G'_{J,\epsilon,SC}}$ soit $c^{-1}C\sigma_{J_{t}}^{G'_{J,\epsilon,SC}}(\boldsymbol{\delta}_{J,sc},a'_{J})$.

Fixons un tel $\boldsymbol{\tau}$ et notons $\underline{\boldsymbol{\tau}}$ son image dans $D^{st}_{g\acute{e}om}({\bf M}',{\cal O}')$ par la suite (1).  C'est aussi l'image de $c\boldsymbol{\tau}$ par la suite (3). Alors
 $\sigma_{J_{s}}^{{\bf G}'}(\boldsymbol{\delta},\xi(a))$ est l'image de  $e_{\tilde{M}'}^{\tilde{G}'}(\epsilon)\underline{\boldsymbol{\tau}}$ modulo $Ann_{{\cal O}'}^{{\bf G}'}$, tandis que $\sigma_{J_{t}}^{{\bf G}_{J}'}(\boldsymbol{\delta},\xi(a))$  est l'image de  $C^{-1}e_{\tilde{M}'}^{\tilde{G}_{J}'}(\epsilon)\underline{\boldsymbol{\tau}}$ modulo $Ann_{{\cal O}'}^{{\bf G}'_{J}}$. Donc $transfert(\sigma_{J_{s}}^{{\bf G}'}(\boldsymbol{\delta},\xi(a)))$, resp. $transfert(\sigma_{J_{t}}^{{\bf G}_{J}'}(\boldsymbol{\delta},\xi(a)))$, est l'image de $e_{\tilde{M}'}^{\tilde{G}'}(\epsilon)transfert(\underline{\boldsymbol{\tau}})$, resp. $C^{-1}e_{\tilde{M}'}^{\tilde{G}_{J}'}(\epsilon)transfert(\underline{\boldsymbol{\tau}})$, dans $D_{g\acute{e}om}({\cal O})/Ann_{{\cal O}}^{\tilde{G}}$.  
Compte tenu de la d\'efinition de $C$, cela d\'emontre le lemme 5.5.

D\'emontrons (5). Introduisons le groupe $(G'_{\epsilon,SC})_{J_{s}}$ et son rev\^etement simplement connexe que l'on note simplement $G^1$. Notons $M^1$ l'image r\'eciproque de $M'_{\epsilon,sc}$ dans $G^1$ et $T^1$ le sous-tore maximal de $M^1$ qui se projette dans $T'$. Posons $G^2=G'_{J,\epsilon,SC}$, $M^2=M'_{\epsilon,J,sc}$ et notons $T^2$ le sous-tore maximal de $M^2$ qui se projette dans $T'$. L'homomorphisme 
$$D_{unip}^{st}(M'_{\epsilon,SC}({\mathbb R}))\stackrel{\iota^*_{M'_{\epsilon,SC},M'_{\epsilon,sc}}}{\to}D_{unip}^{st}(M'_{\epsilon,sc}({\mathbb R}))$$
est le compos\'e de la suite
$$D_{unip}^{st}(M'_{\epsilon,SC}({\mathbb R}))\stackrel{\iota^*_{M'_{\epsilon,SC},M^1}}{\to} D_{unip}^{st}(M^1({\mathbb R}))\stackrel{\iota^*_{M^1,M'_{\epsilon,sc}}}{\to}D_{unip}^{st}(M'_{\epsilon,sc}({\mathbb R})).$$
En utilisant successivement 3.2(3) et le lemme 3.4, on voit que  $\sigma_{J_{s}}^{G'_{\epsilon,SC}}(\boldsymbol{\delta}_{sc},a')$ est l'image par le deuxi\`eme homomorphisme ci-dessus de $i_{J_{s}}^{G'_{\epsilon,SC}}\sigma_{J_{s}}^{G^1}(\boldsymbol{\delta}^1,a^1)$, o\`u $\boldsymbol{\delta}^1=\iota^*_{M'_{\epsilon,SC},M^1}(\boldsymbol{\delta}_{SC})$ et $a^1$ est est un \'el\'ement quelconque de $A_{M^1}({\mathbb R})$ ayant m\^eme projection que  $\xi(a)$ dans $A_{M'}({\mathbb R})/A_{G'}({\mathbb R})$.
Posons
$$C'=C(i_{J_{s}}^{G'_{\epsilon,SC}})^{-1}$$
et  $\boldsymbol{\delta}^2=c^{-1}\boldsymbol{\delta}_{J,sc}=\iota^*_{M'_{\epsilon,SC},M^2}(\boldsymbol{\delta}_{SC})$. Pour $i=1,2$,  on note $a^{i}$ un \'el\'ement de $A_{M^{i}}({\mathbb R})$ qui a m\^eme projection que $\xi(a)$ dans $A_{M'}({\mathbb R})/A_{G'}({\mathbb R}) =A_{M'}({\mathbb R})/A_{G'_{J}}({\mathbb R})$.
  L'assertion (5) r\'esulte de

(6)  il existe un \'el\'ement $\boldsymbol{\tau}\in D_{unip}^{st}(M'_{\epsilon,SC}({\mathbb R}))$ dont l'image dans $D_{unip}^{st}(M^1({\mathbb R}))/Ann_{unip}^{G^1}$ soit $\sigma_{J_{s}}^{G^1}(\boldsymbol{\delta}^1,a^1)$ et dont l'image dans $D_{unip}^{st}(M^2({\mathbb R}))/Ann_{unip}^{G^2}$ soit $C'\sigma_{J_{t}}^{G^2}(\boldsymbol{\delta}^2,a^2)$.

Dans la preuve de 5.5(5), on a calcul\'e les ensembles de racines $\Sigma^{G^{i}}(T^{i})$ pour $i=1,2$. On peut de m\^eme calculer les ensembles de coracines associ\'es $\check{\Sigma}^{G^{i}}(T^{i})$, en utilisant les formules de [W2] 3.3. Donnons le r\'esultat. Pour $\alpha\in \Sigma^G(T)$, notons $\check{\alpha}^{res}$ l'image de $\check{\alpha}\in \mathfrak{t}$ dans $\mathfrak{t}'$ par l'application $\xi$. On obtient que les \'el\'ements de  $\check{\Sigma}^{G^1}(T^1)$  sont de la forme $c^1(\alpha)\check{\alpha}^{res}$, o\`u $\alpha\in \Sigma^G(T)$ est  tel que $\alpha_{A_{\tilde{M}}}\in R_{J}$ et $c^1(\alpha)\in {\mathbb Q}^{\times}$. Les $\alpha$ sont soumis \`a l'une des conditions suivantes et le terme $c^1(\alpha)$ est d\'ecrit dans chaque cas:

(a) $\alpha$ de type $1$, $N\alpha(\nu)=1$, $N\hat{\alpha}(s)=1$; $c^1(\alpha)=1$;

(b)  $\alpha$ de type $2$, $N\alpha(\nu)=1$, $N\hat{\alpha}(s)=1$; $c^1(\alpha)=1$;

(c) $\alpha$ de type $3$, $N\alpha(\nu)=-1$, $N\hat{\alpha}(s)=1$; $c^1(\alpha)=1/2$;

(d) $\alpha$ de type $2$, $N\alpha(\nu)=1$, $N\hat{\alpha}(s)=-1$; $c^1(\alpha)=2$.

Une description analogue vaut pour $\check{\Sigma}^{G^{2}}(T^2)$. Les types $1$, $2$, $3$ sont remplac\'es par $1_{J}$, $2_{J}$, $3_{J}$, l'\'el\'ement $\nu$ est remplac\'e par $t_{J}^{-1}\nu$ et la fonction $c^2$ est donn\'ee par les m\^emes formules que ci-dessus. En utilisant la relation 5.4(1), on voit que la seule diff\'erence entre les ensembles $\check{\Sigma}^{G^1}(T^1)$ et $\check{\Sigma}^{G^2}(T^2)$ provient des racines v\'erifiant (c) et telles que $\alpha$ soit de type $1_{J}$. Dans ce cas, $\check{\Sigma}^{G_{1}}(T^1)$ contient $\check{\alpha}^{res}/2$ tandis que $\check{\Sigma}^{G_{2}}(T^2)$ contient $\check{\alpha}^{res}$. Puisque $\mathfrak{t}^1$ et $\mathfrak{t}^2$, vus comme sous-espaces de $\mathfrak{t}'$, sont engendr\'es par les ensembles de coracines, on en d\'eduit d\'ej\`a que $\mathfrak{t}^1=\mathfrak{t}^2$. On note $j_{*}:\mathfrak{t}^1\to \mathfrak{t}^2$ et $j^*:\mathfrak{t}^{2*}\to \mathfrak{t}^{1*}$ les identit\'es. En reprenant les calculs de 5.5, on voit que, de m\^eme, la seule diff\'erence entre les ensembles de racines  $\Sigma^{G^{1}}(T^1)$ et $\Sigma^{G^{2}}(T^2)$ provient des racines v\'erifiant (c) et telles que $\alpha$ soit de type $1_{J}$. Dans ce cas, $\Sigma^{G^1}(T^1)$ contient $2N\alpha$ tandis que $\Sigma^{G^2}(T^2)$ contient $N\alpha$. On d\'efinit une fonction $b:\Sigma^{G^2}(T^2)\to {\mathbb Q}_{>0}$ qui vaut $1$ sauf sur les racines $N\alpha$ pr\'ec\'edentes, pour lesquelles $b(N\alpha)=1/2$. On obtient que $\Sigma^{G^1}(T^1)$, resp. $\check{\Sigma}^{G^1}(T^1)$, est form\'e des $b(\beta_{2})^{-1}j^*(\beta_{2})$, resp.   $b(\beta_{2})j_{*}^{-1}(\check{\beta}_{2})$, pour $\beta_{2}\in \Sigma^{G^2}(T^2)$. Cela montre que $(G^1,G^2,j_{*})$ est un triplet endoscopique non standard. On a une fonction $B$ sur $\Sigma^{G^1}(T^1)$ et une fonction $B_{J}$ sur $\Sigma^{G^2}(T^2)$. En reprenant les formules de 5.5, on voit que ces fonctions sont reli\'ees comme en [III] 6.4. C'est-\`a-dire, soit $\beta_{2}\in \Sigma^{G^{2}}(T^2)$, notons $\beta_{1}=b(\beta_{2})^{-1}j^*(\beta_{2})$ l'\'el\'ement de $\Sigma^{G^1}(T^1)$ qui lui correspond; on a alors $B(\beta_{1})=\frac{B_{J}(\beta_{2})}{b(\beta_{2})}$. 

Les sous-ensembles $\Sigma^{M^1}(T^1)$ et  $\Sigma^{M^2}(T^2)$ se d\'ecrivent comme ci-dessus, en rempla\c{c}ant la condition $\alpha\in \Sigma^G(T)$ par $\alpha\in \Sigma^M(T)$. Puisque $\Sigma^M(T)\subset \Sigma^{G_{J}}(T)$, une racine dans $\Sigma^M(T)$ est de m\^eme type dans $G$ et $G_{J}$. Cela entra\^{\i}ne que la fonction $b$ vaut $1$ sur $\Sigma^{M^2}(T^2)$. Autrement dit, la correspondance endoscopique non standard se restreint en la correspondance naturelle entre les ensembles de racines de $M^1$ et $M^2$ (celle qui provient de l'identification de ces ensembles de racines \`a celui de $M_{\epsilon,SC}$). 

On peut d\'ecomposer notre triplet endoscopique non standard en produit de triplets $(G^1_{i},G^2_{i},j_{*,i})$ pour $i=1,...,m$, chacun d'eux \'etant  \'equivalent  \`a un triplet quasi-\'el\'ementaire. Les Levi $M^1$ et $M^2$ et les tores $T^1$ et $T^2$  se d\'ecomposent conform\'ement en produits $\prod_{i=1,...,m}M^1_{i}$  etc.... On note $b_{i}$ la restriction de la fonction $b$ \`a $\Sigma^{G^2_{i}}(T^2_{i})$. On va prouver que

(7) pour tout $i=1,...,m$, les donn\'ees $(G^1_{i},G^2_{i},j_{*,i})$, $M^1_{i}$, $M^2_{i}$ et $b_{i}$ v\'erifient les conditions du lemme 5.1 et on a l'in\'egalit\'e $N^{max}(G^1_{i},G^2_{i},j_{*,i})< dim(G_{SC})$. 

Un raisonnement analogue \`a celui de la preuve du lemme 6.1 de [III] nous ram\`ene au cas o\`u le groupe $G_{AD}$ est simple. Il y a un cas particulier: celui o\`u $G_{AD}$ est de type $A_{2n}$ et o\`u l'action de $\theta$ sur ce syst\`eme est l'automorphisme  non trivial. Hors de ce cas, il n'y a pas de racines de type $3$. La d\'efinition de $b$ entra\^{\i}ne alors que cette fonction est constante de valeur $1$.  Donc $j_{*}$ provient d'un isomorphisme de $G^1$ sur $G^2$ et l'assertion est claire. Consid\'erons le cas particulier ci-dessus: $G_{AD}$ est de type $A_{2n}$ et l'action de $\theta$ est non triviale. Le groupe $G_{J}$ peut se r\'ealiser comme intersection de commutants dans $G$ d'\'el\'ements de $A_{\tilde{M}}$. Dans un groupe de type $A_{2n}$, un tel commutant est de type $A_{n_{1}}\times... A_{n_{h}}$. De plus, le syst\`eme de racines de $G_{J}$ est stable par $\theta$. Il en r\'esulte que ce syst\`eme de racines, muni de son automorphisme $\theta$, est produit de sous-syst\`emes de l'un des types suivants

(e) $A_{2n'}$ muni de l'automorphisme non trivial;

(f) $A_{2n'-1}$ muni de l'automorphisme non trivial;

(g) $A_{n'}\times A_{n'}$ muni de la permutation des deux facteurs. 

Quant \`a l'action galoisienne, l'\'el\'ement non trivial de $\Gamma_{{\mathbb R}}$ ne peut agir sur le syst\`eme $A_{2n}$ que par l'identit\'e ou par $\theta$. Il en r\'esulte que chacun des sous-syst\`emes ci-dessus est stable par cette action. Rappelons que le groupe $G^2$ est d\'eduit de $\tilde{G}_{J}$ par la suite d'op\'erations suivantes: on passe de $\tilde{G}_{J}$ \`a un groupe $G'_{J}(t\hat{\theta}_{J})$, on passe de celui-ci au commutant $G'_{J}(t\hat{\theta})_{\epsilon}$, on passe ensuite au rev\^etement simplement connexe. Il est clair que ces op\'erations se d\'ecomposent selon la d\'ecomposition ci-dessus du syst\`eme de racines de $G_{J}$. C'est-\`a-dire que, si on fixe $i\in \{1,...,m\}$, il existe une composante de l'un des types (e), (f), (g) ci-dessus de sorte que $G^2_{i}$ soit une composante irr\'eductible (sur ${\mathbb R}$) d'un groupe issu par le m\^eme proc\'ed\'e que $G^2$ \`a partir de cette composante. Comme plus haut, l'assertion \`a prouver pour les donn\'ees index\'ees par $i$ est claire si $b_{i}$ est constante. Par construction de la fonction $b$, la fonction $b_{i}$ est constante sauf si la composante en question contient des racines de type $1$ provenant de racines de type $3$ dans $G$.  Les racines de type $3$ dans $G$ sont celles qui sont fix\'ees par $\theta$. Une telle racine n'intervient pas dans une composante de type (g). Elle peut intervenir dans une composante de type (e), mais alors elle y est encore de type $3$. Il reste les composantes de type (f). Supposons que $G^2_{i}$ soit issu d'une telle composante.    Les racines de type $1$ de cette composante  provenant de racines de type $3$ dans $G$ sont exactement celles qui sont fix\'ees par $\theta$. Par endoscopie  tordue, on cr\'ee des groupes de syst\`emes de racines de type $B_{p}\times D_{q} \times A_{r_{1}}\times...\times A_{r_{k}}$. En passant \`a un commutant, les types $D_{q}$ ou $A_{r}$ ne cr\'eent que des syst\`emes de m\^eme type. Or, d'apr\`es la classification des donn\'ees endoscopiques non standard \'el\'ementaires, cf. [III] 6.1, de tels types ne peuvent intervenir que dans des donn\'ees \'el\'ementaires "triviales". On obtient la conclusion pour notre triplet $(G^1_{i},G^2_{i},j_{*,i})$ sauf si celui-ci provient d'une composante $B_{p}$ ci-dessus. En passant \`a un commutant dans une telle composante, on obtient un groupe de m\^eme type que ci-dessus, c'est-\`a-dire $B_{p'}\times D_{q'} \times A_{r'_{1}}\times...\times A_{r'_{k}}$. Par le m\^eme argument de classification, on obtient la conclusion pour notre triplet $(G^1_{i},G^2_{i},j_{*,i})$ sauf si le syst\`eme de racines de $G^2_{i}$ est la composante $B_{p'}$. Supposons qu'il en soit ainsi.    On v\'erifie ais\'ement que les racines de $A_{2n'+1}$ fix\'ees par $\theta$ cr\'eent des racines courtes dans la composante $B_{p}$. Cela passe \`a la composante $B_{p'}$ ci-dessus. Il en r\'esulte que les racines de $G^2_{i}$ sur lesquelles $b$ ne vaut pas $1$ sont les racines courtes de $G^2_{i}$.  Sur celles-ci, $b$ vaut $1/2$. Puisque $b$ vaut $1$ sur les racines dans $M^2_{i}$, ce groupe ne contient pas de racines courtes. D'autre part, d'apr\`es la classification de [III] 6.1, ou bien le triplet $(G^1_{i},G^2_{i},j_{*,i})$ est trivial, ou bien $G^1_{i}$ est de type $C_{p'}$ (la premi\`ere possibilit\'e est d'ailleurs exclue sauf si $p'=1$ puisque, d'apr\`es la description ci-dessus, $b$ n'est pas constante si $p'>1$). Cela montre que les donn\'ees index\'ees par $i$ v\'erifient les hypoth\`eses de 5.1. Par ailleurs, il r\'esulte des d\'efinitions que $N^{max}(G^1_{i},G^2_{i},j_{*,i})=4(p')^2-1$. On a n\'ecessairement $p'\leq n'$ (o\`u $n'$ est l'entier associ\'e \`a la composante de type (f) fix\'ee). On a aussi $n'\leq n$. Puisque $dim(G_{SC})=(2n+1)^2-1$, on en d\'eduit l'in\'egalit\'e $N^{max}(G^1_{i},G^2_{i},j_{*,i})< dim(G_{SC})$. Cela v\'erifie (7).

 Remarquons que, dans l'\'egalit\'e (6), les faits que les points $a^1$ et $a^2$ soient proches de $1$ et que $j_{*}:\mathfrak{t}^1\to \mathfrak{t}^2$ soit l'identit\'e (ces espaces \'etant vus comme sous-espaces de $\mathfrak{t}'$) entra\^{\i}nent que  l'on peut identifier ces deux points $a^1$ et $a^2$. Ou encore, avec les notations du lemme 5.1, on a $a^1=exp(X)$ et $a^2=exp(j_{*}(X))$, pour un $X\in \mathfrak{a}_{M^1}({\mathbb R})$. 
D'apr\`es (7), on peut appliquer ce lemme 5.1.  
Celui-ci nous dit que  $\sigma_{J_{s}}^{G^1}(\boldsymbol{\delta}^1,a^1)$ s'envoie sur $c_{M^{1},M^{2}}^{G^1,G^2}\sigma_{J_{t}}^{G^2}(\boldsymbol{\delta}^2,a^2)$ par la correspondance
$$D_{unip}^{st}(M^1({\mathbb R}))/Ann_{unip}^{G^1,st}\simeq D_{unip}^{st}(M^2({\mathbb R}))/Ann_{unip}^{G^2,st}.$$
Comme on l'a vu en 5.1, on a un diagramme commutatif
$$\begin{array}{ccccc}&&D_{unip}^{st}(M_{\epsilon,SC}({\mathbb R}))&&\\ &\swarrow&&\searrow&\\ D_{unip}^{st}(M^1({\mathbb R}))&&\to&&D_{unip}^{st}(M^2({\mathbb R}))\\ \end{array}$$
D'apr\`es le corollaire 3.6(ii), l'\'el\'ement $\sigma_{J_{s}}^{G^1}(\boldsymbol{\delta}^1,a^1)$ provient d'un \'el\'ement $\boldsymbol{\tau}\in D_{unip}^{st}(M_{\epsilon,SC}({\mathbb R}))$. Donc $\boldsymbol{\tau}$ s'envoie sur  $c_{M^{1},M^{2}}^{G^1,G^2}\sigma_{J_{t}}^{G^2}(\boldsymbol{\delta}^2,a^2)$. Pour d\'emontrer (6) et le lemme 5.5, il reste \`a prouver l'\'egalit\'e
$$(8) \qquad c_{M^{1},M^{2}}^{G^1,G^2}=C'.$$

Puisque $b$ prend pour valeurs $1$ et $1/2$, $j_{*}^{-1}$ envoie le ${\mathbb Z}$-module engendr\'e par $\check{\Sigma}^{G^2}(T^2)$ dans le ${\mathbb Z}$-module engendr\'e par $\check{\Sigma}^{G^1}(T^1)$. Puisque les groupes $G^1$ et $G^2$ sont simplement connexes, l'homomorphisme $j_{*}^{-1}$ se rel\`eve en un homomorphisme $T^2\to T^1$. Dualement, on a un homomorphisme $\hat{T}^1\to \hat{T}^2$. Il se restreint en un homomorphisme
$$\hat{j}:Z(\hat{M}^1)^{\Gamma_{{\mathbb R}}}\to Z(\hat{M}^2)^{\Gamma_{{\mathbb R}}}.$$
En se reportant \`a la d\'efinition de [III] 6.4, on v\'erifie l'\'egalit\'e $c_{M^{1},M^{2}}^{G^1,G^2}=\vert ker(\hat{j})\vert^{-1} $.  Notons $\hat{M}'_{\epsilon,ad}$ le dual du groupe $M'_{\epsilon,sc}$ introduit plus haut. Notons aussi $\hat{H}$ le dual de $(G'_{\epsilon,SC})_{J_{s}}$. Alors $\hat{M}^1=\hat{M}'_{\epsilon,ad}/Z(\hat{H})$. Consid\'erons le diagramme commutatif
$$\begin{array}{ccc}Z(\hat{M}')^{\Gamma_{{\mathbb R}}}/\hat{Z}&=&Z(\hat{M}')^{\Gamma_{{\mathbb R}}}/\hat{Z}\\ \downarrow&&\downarrow\\ Z(\hat{M}')^{\Gamma_{{\mathbb R}}}/Z(\hat{G}')^{\Gamma_{{\mathbb R}}}&&Z(\hat{M}')^{\Gamma_{{\mathbb R}}}/Z(\hat{G}_{J}')^{\Gamma_{{\mathbb R}}}\\ \downarrow&&\downarrow\\ Z(\hat{M}'_{\epsilon})^{\Gamma_{{\mathbb R}}}/Z(\hat{G}'_{\epsilon})^{\Gamma_{{\mathbb R}}}&&Z(\hat{M}'_{\epsilon})^{\Gamma_{{\mathbb R}}}/Z(\hat{G}'_{J,\epsilon})^{\Gamma_{{\mathbb R}}}\\ \parallel&&\\ Z(\hat{M}'_{\epsilon,ad})^{\Gamma_{{\mathbb R}}}&&\\ \downarrow&&\parallel\\ Z(\hat{M}'_{\epsilon,ad})^{\Gamma_{{\mathbb R}}}/Z(\hat{H})^{\Gamma_{{\mathbb R}}}&&\\ \parallel&&\\ Z(\hat{M}^1)^{\Gamma_{{\mathbb R}}}&\stackrel{\hat{j}}{\to}&Z(\hat{M}^2)^{\Gamma_{{\mathbb R}}}\\ \end{array}$$
Toutes les fl\`eches sont surjectives et \`a noyaux finis. Calculons le nombre d'\'el\'ements du noyau de l'application compos\'ee
$$Z(\hat{M}')^{\Gamma_{{\mathbb R}}}/\hat{Z}\to Z(\hat{M}^2)^{\Gamma_{{\mathbb R}}}.$$
En utilisant le chemin  nord-est du diagramme, et en appliquant la d\'efinition de [III] 4.3, on obtient
$$[Z(\hat{G}'_{J})^{\Gamma_{{\mathbb R}}}:\hat{Z}] e_{\tilde{M}'}^{\tilde{G}'_{J}}(\epsilon)^{-1}.$$ 
En utilisant le chemin sud-ouest, on obtient
$$[Z(\hat{G}')^{\Gamma_{{\mathbb R}}}:\hat{Z}]e_{\tilde{M}'}^{\tilde{G}'}(\epsilon)^{-1}\vert Z(\hat{H})^{\Gamma_{{\mathbb R}}}\vert (c_{M^{1},M^{2}}^{G^1,G^2})^{-1}.$$
En utilisant la d\'efinition de [III] 1.2, on voit que $\vert Z(\hat{H})^{\Gamma_{{\mathbb R}}}\vert =(i_{J_{s}}^{G'_{\epsilon,SC}})^{-1}$. L'assertion (8) r\'esulte alors de la d\'efinition de $C'$. Cela ach\`eve la preuve. 

\bigskip

\subsection{Preuve du lemme 5.1}
On consid\`ere un triplet $(G_{1},G_{2},j_{*})$ muni de diverses donn\'ees comme en 5.1. Gr\^ace \`a 5.1(9), on suppose de plus que notre triplet est quasi-\'el\'ementaire, que la condition 5.1(8) est v\'erifi\'ee et que $B_{1}$ est constante de valeur $1$.  Comme on l'a remarqu\'e en 5.1, la condition 5.1(8)  implique que le rev\^etement simplement connexe commun  $M_{SC}$ des groupes d\'eriv\'es de $M_{1}$ et $M_{2}$ est isomorphe \`a un produit de groupes $SL_{k}({\mathbb R})$ si $F_{0}={\mathbb R}$, de groupes $SL_{k}({\mathbb C})$ si $F_{0}={\mathbb C}$. Soit $H$ un groupe quasi-d\'eploy\'e  tel que $H_{SC}$ soit isomorphe  \`a un produit de groupes $SL_{k}({\mathbb R})$ ou $SL_{k}({\mathbb C})$. On a

(1) $D_{tr-unip}^{st}(H({\mathbb R}))=D_{orb,unip}^{st}(H({\mathbb R}))$.

Preuve. Le lemme 3.3 nous ram\`ene au cas $H=H_{SC}$. Il suffit donc de traiter les cas $H=SL_{k}({\mathbb R})$ ou $H=SL_{k}({\mathbb C})$. Le m\^eme lemme nous ram\`ene aux cas $H=GL_{k}({\mathbb R})$ ou $H=GL_{k}({\mathbb C})$. Mais alors, ces groupes n'ont pas de donn\'ees endoscopiques elliptiques autres que la donn\'ee "maximale" ${\bf H}$. La d\'efinition de 2.1 montre que $D_{tr-orb}(H({\mathbb R}))=D_{orb}(H({\mathbb R}))$ et l'assertion s'ensuit. $\square$

Reprenons la preuve de [III] 7.7. On introduit le triplet $(G,\tilde{G},{\bf a})$ associ\'e \`a $(G_{1},G_{2},j_{*})$ comme en [III] 6.2. On fixe un \'el\'ement $\eta\in \tilde{G}({\mathbb R})$ qui conserve une paire de Borel \'epingl\'ee ${\cal E}=(B,T,(E_{\alpha})_{\alpha\in \Delta})$ de $G$ d\'efinie sur ${\mathbb R}$. De cette paire se d\'eduit une paire de Borel \'epingl\'ee de $G_{\eta}$ d\'efinie sur ${\mathbb R}$. On peut identifier $G_{1}$ \`a $G_{\eta}$ de sorte que le Levi $M_{1}$ de $G_{1}$ devienne un Levi de $G_{\eta}$ standard pour cette paire. On note $\tilde{M}$ le commutant de $A_{M_{1}}$ dans $\tilde{G}$.   On a $\eta\in \tilde{M}({\mathbb R})$, $M_{\eta}=M_{1}$  et $M$ est standard pour ${\cal E}$.  
D\'ecrivons plus concr\`etement ces objets. Il y a quatre cas.

(a) $G_{1}=Sp(2n)$ et $M_{1}$ est un Levi isomorphe \`a $GL(n_{1})\times...\times GL(n_{k})$, avec $n_{1}+...+n_{k}=n$. Dans ce cas, $G=SL(2n)$ et $ad_{\eta}$ est l'automorphisme ext\'erieur habituel de ce groupe. On v\'erifie que $M$ est le Levi standard de $G$ de blocs $n_{1}\times ...\times n_{k}\times n_{k}\times...\times n_{1}$.

(b) $G_{1}=Spin(2n+1)$ et $M_{1}$ est l'image r\'eciproque dans ce groupe d'un Levi de $SO(2n+1)$ isomorphe \`a $GL(n_{1})\times...\times GL(n_{k})$, avec $n_{1}+...+n_{k}=n$.  Pr\'ecis\'ement, $M_{1}$ est isomorphe au groupe des $(x_{1},...,x_{k},t)\in GL(n_{1})\times...\times GL(n_{k})\times GL(1)$ tels que $det(x_{1})...det(x_{k})t^2=1$. 
 Dans ce cas, $G=Spin(2n+2)$. L'action de $O(2n+2)$ sur $SO(2n+2)$ se rel\`eve en une action sur $G$ 
 et on peut r\'ealiser $ad_{\eta}$ comme la conjugaison par une sym\'etrie \'el\'ementaire qui est un \'el\'ement de $O(2n+2)$ de d\'eterminant $-1$. On v\'erifie que $M$ est l'image r\'eciproque dans $G$ d'un Levi de $SO(2n+2)$ isomorphe \`a $GL(n_{1})\times...\times GL(n_{k})\times GL(1)$.   Comme ci-dessus, $M$ est isomorphe au groupe des $(x_{1},...,x_{k},x_{k+1},t)\in GL(n_{1})\times...\times GL(n_{k})\times GL(1)\times GL(1)$ tels que $det(x_{1})...det(x_{k})x_{k+1}t^2=1$. Mais on peut faire dispara\^{\i}tre l'\'el\'ement $x_{k+1}$ et $M$ est simplement isomorphe $GL(n_{1})\times...\times GL(n_{k})\times GL(1)$.

(c) Les objets sont d\'eduits par restriction des scalaires de ${\mathbb C}$  \`a ${\mathbb R}$ des objets du cas (a) d\'efinis sur ${\mathbb C}$. 

(d) Les objets sont d\'eduits par restriction des scalaires de ${\mathbb C}$  \`a ${\mathbb R}$ des objets du cas (b) d\'efinis sur ${\mathbb C}$.

On doit inclure le triplet $(G,\tilde{G},{\bf a})$ dans un $K$-triplet $(KG,K\tilde{G},{\bf a})$. L'espace de Levi $\tilde{M}$ s'\'etend en un $K$-espace de Levi $K\tilde{M}$. On a en fait $K\tilde{M}=\tilde{M}$. En effet, on sait que $K\tilde{M}$ est lui-m\^eme un $K$-espace. D'apr\`es la d\'efinition de [I] 1.11, ce $K$-espace est r\'eduit \`a une seule composante connexe pouvu que $H^1(\Gamma_{{\mathbb R}};M_{SC})$ soit trivial. C'est le cas d'apr\`es la description ci-dessus: le groupe $M_{SC}$ est produit de groupes $SL(m)$ ou de groupes d\'eduits par restriction des scalaires de $SL(m)$ sur ${\mathbb C}$. Comme en [III] 5.1, on note ${\cal Y}^M$ l'ensemble des $y\in M$ tels que $y\sigma(y)^{-1}\in I^M_{\eta}$. Rappelons que, parce que $G$ est simplement connexe, $Z_{G}(\eta)$ est connexe et cette propri\'et\'e perdure lorsqu'on passe \`a un Levi. Il en r\'esulte que $I^M_{\eta}=M_{\eta}$. Consid\'erons  l'ensemble de doubles classes
  $$M_{\eta}\backslash {\cal Y}^M/M({\mathbb R}).$$
  Montrons que
  
  (2) cet ensemble de doubles classes est r\'eduit \`a un \'el\'ement.
  
  Preuve.  On sait que  cet ensemble est en bijection avec le noyau de l'application
  $$ H^1(\Gamma_{{\mathbb R}};M_{\eta})\to H^1(\Gamma_{{\mathbb R}},M).$$
  Il suffit de prouver que $H^1(\Gamma_{{\mathbb R}};M_{\eta})=\{1\}$. Rappelons que $M_{\eta}=M_{1}$. Dans les cas (c) et (d), le groupe $M_{1}$ est complexe et $H^1(\Gamma_{{\mathbb R}};M_{1})=\{1\}$. Dans le cas (a), $M_{1}$ est un produit de groupes $GL(m)$, d'o\`u la m\^eme conclusion. Dans le cas (b), posons $L=GL(n_{1})\times...\times GL(n_{k})$. On a une suite exacte
  $$1\to \{\pm 1\}\to M_{1}\to L\to 1$$
  D'o\`u une suite exacte
  $$M_{1}({\mathbb R})\to L({\mathbb R})\to H^1(\Gamma_{{\mathbb R}};\{\pm 1\})\to  H^1(\Gamma_{{\mathbb R}};M_{1})\to H^1(\Gamma_{{\mathbb R}};L)=\{1\}$$
  Puisque $ H^1(\Gamma_{{\mathbb R}};\{\pm 1\})$ a deux \'el\'ements, il suffit de prouver que la premi\`ere application de la suite ci-dessus n'est pas surjective. Mais il r\'esulte de la description de $M_{1}$ que cette image est form\'ee des $(x_{1},...,x_{k})\in L({\mathbb R})$ tels que $det(x_{1})...det(x_{k})>0$. Cela prouve (2). $\square$
  
  On peut prendre comme ensemble de repr\'esentants de notre ensemble de doubles classes l'ensemble $\dot{{\cal Y}}^M=\{1\}$. 
  
On introduit les donn\'ees endoscopiques maximales   ${\bf G}'=(G',\hat{G}_{\hat{\theta}}\rtimes W_{{\mathbb R}},\hat{\theta})$ de $(KG,K\tilde{G},{\bf a})$ et  ${\bf M}'=(M',\hat{M}_{\hat{\theta}}\rtimes W_{{\mathbb R}},\hat{\theta})$ de $(KM,K\tilde{M},{\bf a})=(M,\tilde{M},{\bf a})$.   Remarquons que ${\bf G}'$ est aussi  la donn\'ee ${\bf G}'(\hat{\theta})$  d\'eduite de ${\bf M}'$ et de l'\'el\'ement $\tilde{s}=\hat{\theta}$. 
Comme en [III] 6.3, l'\'el\'ement $\eta\in \tilde{G}({\mathbb R})$ d\'etermine un \'el\'ement $\epsilon\in {\cal Z}(\tilde{G}')^{\Gamma_{{\mathbb R}}}$. Si l'on remplace les espaces ambiants $\tilde{G}$ et $\tilde{G}'$ par $\tilde{M}$ et $\tilde{M}'$, on obtient \'evidemment le m\^eme \'el\'ement $\epsilon$. 
On  fixe un diagramme $(\epsilon,B^{M'},T',B^M,T,\eta)$, o\`u $B^M=B\cap M$.   Remarquons que   $G_{2}=G'_{\epsilon}=G'$ et $M_{2}=M'_{\epsilon}=M'$. 

 Comme d'habitude, on n\'eglige les espaces de mesures. On dispose d'\'el\'ements $\boldsymbol{\delta}_{1}$ et $\boldsymbol{\delta}_{2}$. On peut identifier $\boldsymbol{\delta}_{2}$ \`a un \'el\'ement de $D_{tr-unip}^{st}(M'_{\epsilon}({\mathbb R}))$.   En fait, d'apr\`es (1), il appartient \`a $D_{orb,unip}^{st}(M'_{\epsilon}({\mathbb R}))$. On pose
$$\boldsymbol{\delta}=desc_{\epsilon}^{st,M',*} (\boldsymbol{\delta}_{2}).$$
C'est un \'el\'ement de $D^{st}_{orb}({\cal O}')$, o\`u ${\cal O}'$ est la classe de conjugaison stable de $\epsilon$ dans $\tilde{M}'({\mathbb R})$. Posons $\boldsymbol{\tau}=transfert(\boldsymbol{\delta})$. C'est un \'el\'ement de $D_{g\acute{e}om}({\cal O})$, o\`u ${\cal O}$ est la classe de conjugaison stable de $\eta$ dans $\tilde{M}({\mathbb R})$. Montrons que

 (3) $\boldsymbol{\tau}$ appartient \`a $D_{orb}({\cal O})$.

Preuve. On a calcul\'e $\boldsymbol{\tau}$ \`a la fin de la preuve de [III] 7.1. Avec les notations de cette r\'ef\'erence, on a
  $$\boldsymbol{\tau}=\sum_{y\in \dot{{\cal Y}}^M}c^M[y] desc_{\eta[y]}^{\tilde{M},*} \circ transfert_{y}(\boldsymbol{\delta}_{1}).$$
  Remarquons que l'application $\iota^*_{M_{\eta[y],sc},M_{\eta[y]}}$ qui figure en [III] 7.1 dispara\^{\i}t puisque les groupes $G_{\eta[y]}$ sont  ici simplement connexes.   L'\'el\'ement $\boldsymbol{\delta}_{1}$ appartient \`a $D_{orb}^{st}({\cal O}')$, c'est-\`a-dire que c'est une combinaison lin\'eaire stable d'int\'egrales orbitales. L'homomorphisme $desc_{\eta[y]}^{\tilde{M},*}$   envoie une combinaison lin\'eaire d'int\'egrales orbitales sur une telle combinaison lin\'eaire.  En g\'en\'eral, les applications $transfert_{y}$ ne v\'erifient pas cette propri\'et\'e.   Mais elles les v\'erifient dans notre cas particulier car $\dot{{\cal Y}}^M$ est r\'eduit \`a $\{1\}$ et $transfert_{1}$ est l'identit\'e. $\square$

  Les \'el\'ements $J_{1}$ et $J_{2}$ du lemme 5.1 d\'eterminent  un \'el\'ement $J\in {\cal J}_{K\tilde{M}}^{K\tilde{G}}$. Soit $a\in A_{K\tilde{M}}({\mathbb R})$ un \'el\'ement en position g\'en\'erale et proche de $1$. D'apr\`es (3), on peut d\'efinir un \'el\'ement $\rho_{J}^{K\tilde{G}}(\boldsymbol{\tau},a)\in D_{g\acute{e}om}({\cal O})/Ann_{{\cal O}}^{K\tilde{G}}$, o\`u ${\cal O}$ est la classe de conjugaison stable de $\theta^*$ dans $K\tilde{M}({\mathbb R})=\tilde{M}({\mathbb R})$. On peut aussi l'\'el\'ement $\rho_{J}^{K\tilde{G},{\cal E}}({\bf M}',\boldsymbol{\delta},a)$ du m\^eme espace. Nos hypoth\`eses de r\'ecurrence pos\'ees en 5.1 nous autorisent \`a appliquer les r\'esultats de 2.4: on a l'\'egalit\'e $\rho_{J}^{K\tilde{G}}(\boldsymbol{\tau},a)=\rho_{J}^{K\tilde{G},{\cal E}}({\bf M}',\boldsymbol{\delta},a)$. A partir de cette \'egalit\'e, la preuve de [III] 7.7 entra\^{\i}ne l'\'egalit\'e du lemme 5.1. Il faut toutefois v\'erifier que cette preuve s'applique. En inspectant cette preuve, on voit que la seule chose \`a v\'erifier et que l'hypoth\`ese (2) de [III] 7.1 est v\'erifi\'ee pour les triplets $(\bar{G}'(\bar{s})_{SC},G'(\tilde{s})_{\epsilon,SC},j_{*})$ qui apparaissent et qui sont diff\'erents de notre triplet $(G_{1},G_{2},j_{*})$ de d\'epart. Pour un tel triplet, le groupe $\bar{G}'(\bar{s})$ est un groupe endoscopique de $G_{\eta}$, d\'eduit de la fa\c{c}on habituelle du groupe endoscopique "maximal" du Levi $M_{\eta}$ de $G_{\eta}$. Autrement dit, il se d\'eduit de $G_{1}$ et du groupe endoscopique "maximal" de $M_{1}$.    Les groupes endoscopiques de $G_{1}$ \'etant bien connus, on voit que les propri\'et\'es suivantes sont v\'erifi\'ees:

 - $(\bar{G}'(\bar{s})_{SC},G'(\tilde{s})_{\epsilon,SC},j_{*})$  est produit de triplets \'equivalents \`a des triplets quasi-\'el\'ementaires  $(G_{1,i},G_{2,i},j_{*i})$, pour $i=1,...,m$, qui sont de l'un des types (1), (2) ou (3) de [III] 6.1; le Levi $\bar{M}'\simeq M_{1}$ de $\bar{G}'(\bar{s})$ d\'etermine pour chaque $i$ un Levi $M_{1,i}$ de $G_{1,i}$;
 
 - si un triplet  $(G_{1,i},G_{2,i},j_{*i})$ est de type (2), resp. (3), les racines dans $M_{1,i}$ sont longues, resp. courtes;
  
  - si $(\bar{G}'(\bar{s})_{SC},G'(\tilde{s})_{\epsilon,SC},j_{*})$ est
diff\'erent de $(G_{1},G_{2},j_{*})$, on a pour tout $i$ l'in\'egalit\'e $N^{max}(G_{1,i},G_{2,i},j_{*i})<N^{max}(G_{1},G_{2},j_{*})$.

Les deux premi\`eres propri\'et\'es impliquent que le triplet  $(\bar{G}'(\bar{s})_{SC},G'(\tilde{s})_{\epsilon,SC},j_{*})$ v\'erifie les conditions de 5.1. On voit qu'alors, l'hypoth\`ese (2) de [III] 7.1 \'equivaut au lemme 5.1. La derni\`ere propri\'et\'e ci-dessus et nos hypoth\`eses de r\'ecurrence assurent que ce lemme est v\'erifi\'e pour les triplets $(\bar{G}'(\bar{s})_{SC},G'(\tilde{s})_{\epsilon,SC},j_{*})$ 
diff\'erents de $(G_{1},G_{2},j_{*})$. La d\'emonstration de [III] 7.7 s'applique donc bien. Cela ach\`eve la preuve. 

\bigskip
       
 \section{Un r\'esultat d'approximation}
 
 \bigskip
 
 \subsection{Un espace de germes de fonctions}
 Dans toute la section, $(G,\tilde{G},{\bf a})$ est quasi-d\'eploy\'e et \`a torsion int\'erieure. 
  Soit $\tilde{M}$ un espace de Levi de $\tilde{G}$. Soit $\tilde{L}\in {\cal L}(\tilde{M})$. Une racine $\alpha\in \Sigma^G(A_{M})$ (consid\'er\'ee comme une forme lin\'eaire sur ${\cal A}_{M}$), se d\'ecompose en $\alpha_{L}+\alpha^L$, o\`u $\alpha_{L}\in {\cal A}_{L}^*$ et $\alpha^L\in {\cal A}^{L,*}$.  Consid\'erons l'ensemble des formes lin\'eaires $\alpha_{L}$,  pour  $\tilde{L}\in {\cal L}(\tilde{M})$ et $\alpha\in \Sigma^G(A_{M})$. On note  $V_{M}^G$ le sous-ensemble des \'el\'ements non nuls. On note $U_{M}^G$ le sous-ensemble des $H\in {\cal A}_{M}$ tels que $\alpha(H)\not=0$ pour tout $\alpha\in V_{M}^G$. C'est le compl\'ementaire dans ${\cal A}_{M}$ d'un ensemble fini d'hyperplans. 
  Notons $\underline{{\cal V}}_{M}^G$ l'espace des fonctions  sur $U_{M}^G$ , qui sont combinaisons lin\'eaires de fonctions
 $$H\mapsto \prod_{i=1,...,n}log(\vert exp(r_{i}\alpha_{i}(H))-exp(-r_{i}\alpha_{i}(H))\vert _{{\mathbb R}})$$
 o\`u les $\alpha_{i}$ appartiennent \`a  $ V_{M}^G$ et les $r_{i}$ sont des r\'eels non nuls. On consid\`ere les \'el\'ements de $\underline{{\cal V}}_{M}^G$ comme des fonctions d\'efinies presque partout sur ${\cal A}_{M}$. Remarquons que ces fonctions sont invariantes par translations par ${\cal A}_{G}$. On peut aussi bien consid\'erer qu'elles sont d\'efinies sur ${\cal A}_{M}^G$. 
 
 Appelons domaine ad\'equat   dans ${\cal A}_{M}$ l'intersection d'un voisinage ouvert de $0$ dans ${\cal A}_{M}$ avec l'ensemble des $H\in {\cal A}_{M}$ qui v\'erifient la condition $\vert \alpha(H)\vert _{{\mathbb R}}>c \vert\vert  H\vert \vert $ pour tout $\alpha\in V_{M}^G $, o\`u $c>0$ est un r\'eel fix\'e (on note ici $\vert \vert .\vert \vert $  la norme euclidienne fix\'ee sur ${\cal A}_{M}$) . Soit $u$ un germe de fonction d\'efini presque partout au voisinage de $0$ dans ${\cal A}_{M}$. On dit qu'il est faiblement \'equivalent \`a $0$ s'il existe $r>0$ et si, pour tout domaine ad\'equat, il existe $C>0$ tel que $\vert u(H)\vert\leq C\vert \vert H\vert \vert ^r$ pour tout $H$ dans le domaine et assez proche de $0$. On dit que deux germes $u$ et $u'$ sont faiblement \'equivalents si et seulement si $u-u'$ est faiblement  \'equivalent \`a $0$.
 
 {\bf Remarque.} Ces d\'efinitions d\'ependent de $G$ mais cela ne nous g\^enera pas.
 
 \bigskip
 
 Une fonction dans $ \underline{{\cal V}}_{M}^G$ peut \^etre faiblement \'equivalente \`a $0$. On note $\underline{{\cal V}}_{0,M}^G$ le sous-espace des \'el\'ements de ${\cal V}_{M}^G$ qui sont faiblement \'equivalents \`a $0$ et on pose ${\cal V}_{M}^G=\underline{{\cal V}}_{M}^G/\underline{{\cal V}}_{0,M}^G$. 
 Remarquons que, si $u$ est une fonction d\'efinie presque partout sur ${\cal A}_{M}$ qui est faiblement \'equivalente \`a $0$ et si $v\in\underline{ {\cal V}}_{\tilde{M}}^{\tilde{G}}$, alors $uv$ est faiblement \'equivalente \`a $0$: le produit d'une fonction \`a croissance logarithmique et d'une  fonction d\'ecroissante en $\vert \vert H\vert \vert ^r$ est d\'ecroissante en $\vert \vert H\vert \vert ^{r-\epsilon}$ pour tout $\epsilon>0$.
  
  \bigskip
 
 \subsection{Approximation des int\'egrales orbitales pond\'er\'ees invariantes}
  On note  $Ind_{\tilde{M}}^{\tilde{G}}$ l'homomorphisme d'induction
 $$\begin{array}{ccc} D_{orb}(\tilde{M}({\mathbb R}))\otimes Mes(M({\mathbb R}))^*&\to & D_{orb}(\tilde{G}({\mathbb R}))\otimes Mes(G({\mathbb R}))^*\\ \boldsymbol{\gamma}&\mapsto& \boldsymbol{\gamma}^{\tilde{G}}\end{array}.$$
  On fixe une r\'eunion finie ${\cal O}$ de classes de conjugaison semi-simples par $M({\mathbb R})$ dans $\tilde{M}({\mathbb R})$.  
  Consid\'erons un \'el\'ement 
 $$\zeta \in {\cal V}_{M}^G\otimes Ind_{\tilde{M}}^{\tilde{G}}(D_{orb}({\cal O})\otimes Mes(M({\mathbb R}))^*).$$ 
 Relevons-le en un \'el\'ement $\underline{\zeta}\in \underline{{\cal V}}_{M}^G\otimes  Ind_{\tilde{M}}^{\tilde{G}}(D_{orb}({\cal O}^{\tilde{G}})\otimes Mes(G({\mathbb R}))^*)$. On peut l'\'evaluer en un point $H\in U_{M}^G$,
on obtient un \'el\'ement $\underline{\zeta}(H)\in Ind_{\tilde{M}}^{\tilde{G}}(D_{orb}({\cal O})\otimes Mes(M({\mathbb R}))^*)$.  La fonction $H\mapsto \underline{\zeta}(H)$ ne d\'epend du choix du rel\`evement $\underline{\zeta}$ qu'\`a faible \'equivalence pr\`es (en un sens similaire \`a celui du paragraphe pr\'ec\'edent). Une telle \'equivalence importera peu, on notera donc simplement $H\mapsto \zeta(H)$ cette fonction. 
 On rappelle que, si $\boldsymbol{\gamma}\in D_{orb}({\cal O})$ et $a\in A_{M}({\mathbb R})$, on peut d\'efinir la distribution $a\boldsymbol{\gamma}$. Sa valeur sur  une fonction $\varphi\in C_{c}^{\infty}(\tilde{M}({\mathbb R}))$ est \'egale \`a la valeur de $\boldsymbol{\gamma}$ sur la fonction $^{a}\varphi$, o\`u $^{a}\varphi(\delta)=\varphi(a\delta)$ pour tout $\delta\in \tilde{M}({\mathbb R})$.   Fixons d\'esormais un syst\`eme de fonctions $B$ sur $\tilde{G}({\mathbb R})$. 
 
 \ass{Proposition}{ Il existe un unique homomorphisme
  $$\xi^{\tilde{G}}(B):D_{orb}({\cal O})\otimes Mes(M({\mathbb R}))^*\to {\cal V}_{M}^G\otimes Ind_{\tilde{M}}^{\tilde{G}}( D_{orb}({\cal O})\otimes Mes(M({\mathbb R}))^*)$$
  tel que, pour tout $\boldsymbol{\gamma}\in D_{orb}({\cal O})\otimes Mes(M({\mathbb R}))^*$ et tout ${\bf f}\in C_{c}^{\infty}(\tilde{G}({\mathbb R})\otimes Mes(G({\mathbb R}))$, le germe de la fonction d\'efinie presque partout  sur ${\cal A}_{M}$ qui, \`a $H\in {\cal A}_{M}$, associe
 $$\sum_{\tilde{L}\in {\cal L}(\tilde{M}}I_{\tilde{L}}^{\tilde{G}}(exp(H_{L})\xi^{\tilde{L}}(\boldsymbol{\gamma},B,H^L),{\bf f})$$
 soit faiblement \'equivalent au germe constant de valeur
 $$I_{\tilde{M}}^{\tilde{G}}(\boldsymbol{\gamma},B,{\bf f}).$$}
 
{\bf Remarques.} (1) Comme souvent, il est implicite que les fonctions $\xi^{\tilde{L}}(B)$ intervenant pour $\tilde{L}\not=\tilde{G}$ ont \'et\'e d\'etermin\'ees par r\'ecurrence par le m\^eme \'enonc\'e appliqu\'e en rempla\c{c}ant $\tilde{G}$ par $\tilde{L}$.

 (2) Pour $\tilde{M}=\tilde{G}$, on a $ {\cal V}_{M}^{M}={\mathbb C}$. L'application $\xi^{\tilde{M}}(B)$ est l'identit\'e. 
 
 \bigskip
 Preuve de l'unicit\'e. On peut supposer par r\'ecurrence que les $\xi^{\tilde{L}}$ sont d\'etermin\'es pour $\tilde{L}\not=\tilde{G}$. Pour un tel $\tilde{L}$, la classe d'\'equivalence faible de la fonction 
$$H\mapsto I_{\tilde{L}}^{\tilde{G}}(exp(H_{L})\xi^{\tilde{L}}(\boldsymbol{\gamma},B,H^L),{\bf f})$$
est bien d\'etermin\'ee. En effet,  $\xi^{\tilde{L}}(\boldsymbol{\gamma},B,H^L)$ est uniquement d\'etermin\'ee modulo une combinaison lin\'eaire de fonctions $H\mapsto u(H^L)\boldsymbol{\gamma}_{\tilde{L}}$, o\`u $u\in \underline{{\cal V}}_{0,M}^L$ et $\boldsymbol{\gamma}_{\tilde{L}}\in D_{orb}({\cal O}^{\tilde{L}})\otimes Mes(L({\mathbb R}))^*$. Or, pour de telles donn\'ees,  la relation 2.4(5) (qui est v\'erifi\'ee dans notre situation quasi-d\'eploy\'ee et \`a torsion int\'erieure) entra\^{\i}ne que la fonction
$$H\mapsto u(H^L)I_{\tilde{L}}^{\tilde{G}}(exp(H_{L})\boldsymbol{\gamma}_{\tilde{L}},{\bf f})$$
est faiblement \'equivalente \`a $0$. 
  Alors la formule de l'\'enonc\'e d\'etermine la classe d'\'equivalence faible de la fonction
 $$H\mapsto I^{\tilde{G}}(\xi^{\tilde{G}}(\boldsymbol{\gamma},B,H),{\bf f}).$$
 Fixons 
  une base $(v_{i})_{i\in I}$ d'un suppl\'ementaire de $\underline{{\cal V}}_{0,M}^G$ dans $\underline{{\cal V}}_{M}^G$. On peut relever $\xi^{\tilde{G}}(\boldsymbol{\gamma},B)$ de fa\c{c}on unique en un terme
 $$\underline{\xi}^{\tilde{G}}(\boldsymbol{\gamma},B)=\sum_{i\in I}v_{i}\otimes \boldsymbol{\gamma}_{i},$$
 o\`u les $\boldsymbol{\gamma}_{i}$ appartiennent \`a $Ind_{\tilde{M}}^{\tilde{G}}(D_{orb}({\cal O})\otimes Mes(M({\mathbb R}))^*)$. Alors $I^{\tilde{G}}(\xi^{\tilde{G}}(\boldsymbol{\gamma},B,H),{\bf f})$ est faiblement \'equivalent \`a
 $$\sum_{i\in I}v_{i}(H)I^{\tilde{G}}(\boldsymbol{\gamma}_{i},{\bf f}).$$
 Donc  la classe d'\'equivalence faible de cette somme est bien d\'etermin\'ee. Par d\'efinition de la base $(v_{i})_{i\in I}$, cela entra\^{\i}ne 
 que les coefficients $I^{\tilde{G}}(\boldsymbol{\gamma}_{i},{\bf f})$ sont bien d\'etermin\'es. Cela \'etant vrai pour tout ${\bf f}$, les distributions $\boldsymbol{\gamma}_{i}$ sont uniquement d\'etermin\'ees.      Donc $\underline{\xi}^{\tilde{G}}(\boldsymbol{\gamma},B)$ est uniquement d\'etermin\'e, ce qui prouve l'unicit\'e. 
 
 Preuve de l'existence. Par lin\'earit\'e, il suffit de traiter le cas o\`u $\boldsymbol{\gamma}$ est l'int\'egrale orbitale associ\'ee \`a un \'el\'ement $\gamma\in \tilde{M}({\mathbb R})$ de partie semi-simple dans ${\cal O}$ et \`a une mesure de Haar sur  $M_{\gamma}({\mathbb R})$. Pour deux espaces de Levi $\tilde{L}$ et $\tilde{L}'$ tels que
 $\tilde{M}\subset\tilde{L}\subset \tilde{L}'$, on sait d\'efinir la fonction $H\mapsto r_{\tilde{L}}^{\tilde{L}'}(\gamma,exp(H),B)$ pour un point $H$ en position g\'en\'erale dans ${\cal A}_{M}$, cf. [II] 1.9. Elle est d'ailleurs invariante par translations par ${\cal A}_{L'}$, donc d\'efinie pour $H$ en position g\'en\'erale dans ${\cal A}_{M}^{L'}$. En fait, on a montr\'e en [II] 1.7(9) qu'elle s'\'etendait par continuit\'e au voisinage de tout point de ${\cal A}_{L}^{L'}$ en position g\'en\'erale (dans cette r\'ef\'erence, il n'y avait pas de syst\`eme de fonctions $B$ mais le r\'esultat s'\'etend \`a ce cas). On peut donc d\'efinir $r_{\tilde{L}}^{\tilde{L}'}(\gamma,exp(H_{L}),B)$ pour tout $H\in {\cal A}_{M}$ en position g\'en\'erale. Plus pr\'ecis\'ement, le r\'esultat de [II] 1.7(9) affirme que l'on a l'\'egalit\'e
 $$(3) \qquad r_{\tilde{L}}^{\tilde{L}'}(\gamma,exp(H_{L}),B)=r_{\tilde{L}}^{\tilde{L}'}(\gamma',exp(H_{L}),B)$$
 o\`u $\gamma'$ est un \'el\'ement quelconque de l'orbite induite de $\gamma$ \`a $\tilde{L}$.
  Cela implique que cette fonction est d\'efinie pour $H\in U_{M}^G$ proche de $0$. D\'efinissons une fonction $H\mapsto \bar{r}_{\tilde{M}}^{\tilde{G}}(\gamma,H,B)$  sur $U_{M}^G$ par la formule de r\'ecurrence
 $$(4)\qquad \bar{r}_{\tilde{M}}^{\tilde{G}}(\gamma,H,B)=-\sum_{\tilde{L}\in{\cal L}(\tilde{M}),\tilde{L}\not=\tilde{G}}\bar{r}_{\tilde{M}}^{\tilde{L}}(\gamma,H^L,B)r_{\tilde{L}}^{\tilde{G}}(\gamma,exp(H_{L}),B).$$
 On v\'erifie facilement que c'est un \'el\'ement de ${\cal V}_{M}^G$.

Cela \'etant, on pose
 $$\xi^{\tilde{G}}(\boldsymbol{\gamma},H,B)=(-1)^{a_{M}-a_{G}}\bar{r}_{\tilde{M}}^{\tilde{G}}(\gamma,H,B)\boldsymbol{\gamma}.$$
 Pour prouver que cette d\'efinition satisfait la condition de l'\'enonc\'e, on doit calculer le germe de l'expression
$$\sum_{\tilde{L}\in {\cal L}(\tilde{M})}I_{\tilde{L}}^{\tilde{G}}(exp(H_{L})\xi^{\tilde{L}}(\boldsymbol{\gamma},B,H^L),{\bf f}).$$ 
Ceci n'est autre que
$$(5) \qquad \sum_{\tilde{L}\in {\cal L}(\tilde{M})}(-1)^{a_{M}-a_{L}}\bar{r}_{\tilde{M}}^{\tilde{L}}(\gamma,H^L,B)I_{\tilde{L}}^{\tilde{G}}(exp(H_{L})\boldsymbol{\gamma}^{\tilde{L}},{\bf f}).$$
Soit $\tilde{L}\in {\cal L}(\tilde{M})$. Fixons une famille $(\gamma_{i})_{i=1,...,m}$ de repr\'esentants des classes de conjugaison par $L({\mathbb R})$ dans l'orbite induite par $\gamma$. On peut d\'ecomposer $I_{\tilde{L}}^{\tilde{G}}(exp(H_{L})\boldsymbol{\gamma}^{\tilde{L}},{\bf f})$ en somme de $I_{\tilde{L}}^{\tilde{G}}(exp(H_{L})\boldsymbol{\gamma}_{i},{\bf f})$, o\`u $\boldsymbol{\gamma}_{i}$ est l'int\'egrale orbitale associ\'ee \`a $\gamma_{i}$ et une certaine mesure sur $L_{\gamma_{i}}({\mathbb R})$.   On a montr\'e en [II] 3.2(1) que le germe de $I_{\tilde{L}}^{\tilde{G}}(exp(H_{L})\boldsymbol{\gamma}_{i},{\bf f})$ \'etait \'equivalent \`a
$$\sum_{\tilde{R}\in {\cal L}(\tilde{L})}(-1)^{a_{L}-a_{R}}r_{\tilde{L}}^{\tilde{R}}(\gamma_{i},exp(H_{L}),B)I^{\tilde{G}}_{\tilde{R}}(\boldsymbol{\gamma}_{i}^{\tilde{R}},B,{\bf f}).$$
L'\'equivalence utilis\'ee dans cette r\'ef\'erence n'\'etait pas la m\^eme qu'ici, mais elle \'etait plus forte d'apr\`es les d\'efinitions. La m\^eme assertion vaut donc pour notre \'equivalence faible.
La relation (3) nous dit que les coefficients $r_{\tilde{L}}^{\tilde{R}}(\gamma_{i},exp(H_{L}),B)$ sont ind\'ependants de $i$ et valent $r_{\tilde{L}}^{\tilde{R}}(\gamma,exp(H_{L}),B)$. On peut regrouper les expressions ci-dessus et on obtient que $I_{\tilde{L}}^{\tilde{G}}(exp(H_{L})\boldsymbol{\gamma}^{\tilde{L}},{\bf f})$ est faiblement \'equivalent \`a
$$\sum_{\tilde{R}\in {\cal L}(\tilde{L})}(-1)^{a_{L}-a_{R}}r_{\tilde{L}}^{\tilde{R}}(\gamma,exp(H_{L}),B)I_{\tilde{R}}^{\tilde{G}}(\boldsymbol{\gamma}^{\tilde{R}},B,{\bf f}).$$
L'\'equivalence faible se conserve par multiplication par une fonction de ${\cal V}_{M}^G$.  Donc l'expression (5) est faiblement \'equivalente \`a
$$\sum_{\tilde{R}\in {\cal L}(\tilde{M})}(-1)^{a_{M}-a_{R}}X^{\tilde{R}}(H)I_{\tilde{R}}^{\tilde{G}}(\boldsymbol{\gamma}^{\tilde{R}},B,{\bf f}),$$
o\`u
$$X^{\tilde{R}}(H)=\sum_{\tilde{L}\in {\cal L}^{\tilde{R}}(\tilde{M})}\bar{r}_{\tilde{M}}^{\tilde{L}}(\gamma,H^L,B)r_{\tilde{L}}^{\tilde{R}}(\gamma,exp(H_{L}),B).$$
Mais la d\'efinition (4) entra\^{\i}ne que $X^{\tilde{M}}(H)=1$ tandis que $X^{\tilde{R}}(H)$ est \'equivalent \`a $0$ si $\tilde{R}\not=\tilde{M}$. Donc (5) est faiblement \'equivalent \`a $I_{\tilde{M}}^{\tilde{G}}(\boldsymbol{\gamma},B,{\bf f})$. $\square$

  \bigskip
 
 \subsection{Approximation des int\'egrales orbitales pond\'er\'ees invariantes stables}
 On fixe une r\'eunion finie ${\cal O}$ de classes de conjugaison stable semi-simples dans $\tilde{M}({\mathbb R})$.   En fixant pour un instant les mesures, on d\'efinit  l'espace  
 $$D^{st}_{tr-orb}(\tilde{G}({\mathbb R}),{\cal O}) =D_{tr-orb}^{st}(\tilde{G}({\mathbb R})) \cap Ind_{\tilde{M}}^{\tilde{G}}(D_{g\acute{e}om}^{st}({\cal O})).$$
 
  \ass{Proposition}{ Il existe un unique homomorphisme
  $$\xi^{\tilde{G},st}(B):D_{tr-orb}^{st}({\cal O})\otimes Mes(M({\mathbb R}))^*\to {\cal V}_{M}^G\otimes D^{st}_{tr-orb}(\tilde{G}({\mathbb R}),{\cal O})\otimes Mes(G({\mathbb R}))^*$$
  tel que, pour tout $\boldsymbol{\delta}\in D_{tr-orb}^{st}({\cal O})\otimes Mes(M({\mathbb R}))^*$ et tout ${\bf f}\in C_{c}^{\infty}(\tilde{G}({\mathbb R})\otimes Mes(G({\mathbb R}))$, le germe de la fonction d\'efinie presque partout  sur ${\cal A}_{M}$ qui, \`a $H\in {\cal A}_{M}$, associe
 $$\sum_{\tilde{L}\in {\cal L}(\tilde{M}}S_{\tilde{L}}^{\tilde{G}}(exp(H_{L})\xi^{\tilde{L},st}(\boldsymbol{\delta},B,H^L),{\bf f})$$
 soit faiblement \'equivalent au germe constant de valeur
 $$S_{\tilde{M}}^{\tilde{G}}(\boldsymbol{\delta},B,{\bf f}).$$}
 
 {\bf Remarques.} (1)  Pour $H$ en position g\'en\'erale, la distribution $exp(H_{L})\xi^{\tilde{L},st}(\boldsymbol{\delta},B,H^L)$ est $\tilde{G}$-\'equisinguli\`ere (c'est-\`a-dire support\'ee par des \'el\'ements $\gamma\in \tilde{L}({\mathbb R})$ tels que $L_{\gamma}=G_{\gamma}$), donc l'int\'egrale $S_{\tilde{L}}^{\tilde{G}}(exp(H_{L})\xi^{\tilde{L},st}(\boldsymbol{\delta},B,H^L),{\bf f})$ est bien d\'efinie.
 
 (2) Si $\tilde{M}=\tilde{G}$, $\xi^{\tilde{M},st}(B)$ est l'inclusion naturelle.
 
 \bigskip
 L'assertion d'unicit\'e se d\'emontre comme pour la proposition pr\'ec\'edente.  L'existence sera d\'emontr\'ee  en 6.5.  

\bigskip 
 \subsection{Approximation des int\'egrales orbitales pond\'er\'ees invariantes associ\'ees aux \'el\'ements de $D_{tr-orb}(\tilde{M}({\mathbb R}))\otimes Mes(M({\mathbb R}))^*$}
 On fixe une r\'eunion finie ${\cal O}$ de classes de conjugaison stable semi-simples dans $\tilde{M}({\mathbb R})$.   En fixant pour un instant les mesures, on d\'efinit l'espace  
  $$D_{tr-orb}(\tilde{G}({\mathbb R}),{\cal O}) =D_{tr-orb}(\tilde{G}({\mathbb R}))\cap Ind_{\tilde{M}}^{\tilde{G}}(D_{g\acute{e}om}({\cal O}) ).$$
  \ass{Proposition}{ Il existe un unique homomorphisme
  $$\xi^{\tilde{G}}(B):D_{tr-orb}({\cal O})\otimes Mes(M({\mathbb R}))^*\to {\cal V}_{M}^G\otimes D_{tr-orb}(\tilde{G}({\mathbb R}),{\cal O})\otimes Mes(G({\mathbb R}))^*  $$
  tel que, pour tout $\boldsymbol{\gamma}\in D_{tr-orb}({\cal O})\otimes Mes(M({\mathbb R}))^*$ et tout ${\bf f}\in C_{c}^{\infty}(\tilde{G}({\mathbb R})\otimes Mes(G({\mathbb R}))$, le germe de la fonction d\'efinie presque partout  sur ${\cal A}_{M}$ qui, \`a $H\in {\cal A}_{M}$, associe
 $$\sum_{\tilde{L}\in {\cal L}(\tilde{M})}I_{\tilde{L}}^{\tilde{G}}(exp(H_{L})\xi^{\tilde{L}}(\boldsymbol{\gamma},B,H^L)^{\tilde{L}},{\bf f})$$
 soit faiblement \'equivalent au germe constant de valeur
 $$I_{\tilde{M}}^{\tilde{G}}(\boldsymbol{\gamma},B,{\bf f}).$$}
 
 {\bf Remarque.} Si $\tilde{M}=\tilde{G}$, $\xi^{\tilde{M}}(B)$ est l'inclusion naturelle. 
 \bigskip
 
 Preuve. L'unicit\'e se d\'emontre comme en 6.2. D\'emontrons l'existence. Par lin\'earit\'e, on peut supposer que $\boldsymbol{\gamma}$ est une int\'egrale orbitale appartenant \`a $D_{orb}({\cal O})\otimes Mes(M({\mathbb R}))^*$ ou qu'il existe une donn\'ee endoscopique elliptique et relevante ${\bf M}'$ de $(M,\tilde{M})$, avec $M'\not=M$, qu'il existe une classe de conjugaison stable ${\cal O}'$ dans $\tilde{M}'({\mathbb R})$ correspondant \`a une classe dans ${\cal O}$  et qu'il existe $\boldsymbol{\delta}\in D_{tr-orb}^{st}({\bf M}',{\cal O}')\otimes Mes(M'({\mathbb R}))^*$, de sorte que $\boldsymbol{\gamma}=transfert(\boldsymbol{\delta})$. Le premier cas est trait\'e par la proposition 6.2. Traitons le second. On \'ecrit ${\bf M}'=(M',{\cal M}',\zeta)$. D'apr\`es 2.4(2) (qui est valide dans notre situation), on a l'\'egalit\'e
 $$(1) \qquad I_{\tilde{M}}^{\tilde{G}}(\boldsymbol{\gamma},B,{\bf f})=I_{\tilde{M}}^{\tilde{G},{\cal E}}({\bf M}',\boldsymbol{\delta},B,{\bf f})=\sum_{s\in \zeta Z(\hat{M})^{\Gamma_{{\mathbb R}}}/Z(\hat{G})^{\Gamma_{{\mathbb R}}}}i_{\tilde{M}'}(\tilde{G},\tilde{G}'(s))S_{{\bf M}'}^{{\bf G}'(s)}(\boldsymbol{\delta},B,{\bf f}^{{\bf G}'(s)}).$$
 On peut appliquer   la proposition 6.3 par r\'ecurrence \`a chacun des termes intervenant ici. Comme toujours, on doit d\'evelopper quelques formalit\'es pour adapter cette proposition au cas de donn\'ees endoscopiques. Le seul point \`a expliquer est que la multiplication $(H,\boldsymbol{\delta})\mapsto exp(H)\boldsymbol{\delta}$ a un sens dans $D_{g\acute{e}om}^{st}({\bf M}')\otimes Mes(M'({\mathbb R}))^*$, pour $H\in {\cal A}_{M}\simeq {\cal A}_{M'}$. En effet, introduisons des donn\'ees auxiliaires $M'_{1}$,...,$\Delta_{1}$. 
 Identifions $\boldsymbol{\delta}$ \`a un \'el\'ement $\boldsymbol{\delta}_{1}$ de $D_{g\acute{e}om,\lambda_{1}}^{st}(\tilde{M}'_{1}({\mathbb R}))\otimes Mes(M'({\mathbb R}))^*$. L'\'el\'ement $H\in {\cal A}_{M}\simeq {\cal A}_{M'}$ se rel\`eve en un \'el\'ement $H_{1}\in {\cal A}_{M'_{1}}$ et la distribution $exp(H_{1})\boldsymbol{\delta}_{1}$ est bien d\'efinie. On a introduit en [IV] 2.1 un caract\`ere $\lambda_{\mathfrak{A}_{M'_{1}}}$ de $\mathfrak{A}_{M'_{1}}=exp( {\cal A}_{M'_{1}})$. La distribution $\lambda_{\mathfrak{A}_{M'_{1}}}(exp(H_{1}))exp(H_{1})\boldsymbol{\delta}_{1}$ ne d\'epend pas du rel\`evement $H_{1}$ de $H$. 
  En effet,  on ne peut modifier $H_{1}$ que par un \'el\'ement de ${\cal A}_{C_{1}}$, donc $exp(H_{1})$ par un \'el\'ement $c\in \mathfrak{A}_{C_{1}}$. Or $c\boldsymbol{\delta}_{1}=\lambda_{1}(c)^{-1}\boldsymbol{\delta}_{1}$ et $\lambda_{1}$ co\"{\i}ncide sur $\mathfrak{A}_{C_{1}}$ avec $\lambda_{\mathfrak{A}_{M'_{1}}}$. D'o\`u l'invariance cherch\'ee. On v\'erifie que, si l'on change de donn\'ees auxiliaires, les distributions $\lambda_{\mathfrak{A}_{M'_{1}}}(exp(H_{1}))exp(H_{1})\boldsymbol{\delta}_{1}$ se recollent en un unique \'el\'ement de 
   $D_{g\acute{e}om}^{st}({\bf M}')\otimes Mes(M'({\mathbb R}))^*$, que l'on note $exp(H)\boldsymbol{\delta}$.
 On applique la proposition. Chaque terme $S_{{\bf M}'}^{{\bf G}'(s)}(\boldsymbol{\delta},B,{\bf f}^{{\bf G}'(s)})$ se d\'eveloppe en une somme index\'ee par des espaces de Levi $\tilde{L}'_{s}\in {\cal L}^{\tilde{G}'(s)}(\tilde{M}')$. Comme on l'a d\'ej\`a vu plusieurs fois, un tel espace de Levi donne naissance \`a un espace de Levi $\tilde{L}\in {\cal L}(\tilde{M})$. Le terme $s$ d\'efinit une donn\'ee endoscopique ${\bf L}'(s)$ de $(L,\tilde{L})$ de sorte que $\tilde{L}'_{s}\simeq \tilde{L}'(s)$. Pour simplifier les notations, on anticipe cela dans la formule suivante.
On obtient que le membre de droite de (1) est faiblement \'equivalent \`a la fonction sur ${\cal A}_{M}$ qui \`a $H$ associe
 $$\sum_{s\in \zeta Z(\hat{M})^{\Gamma_{{\mathbb R}}}/Z(\hat{G})^{\Gamma_{{\mathbb R}}}}i_{\tilde{M}'}(\tilde{G},\tilde{G}'(s))\sum_{\tilde{L}'_{s}\in {\cal L}^{\tilde{G}'(s)}(\tilde{M}')}S_{{\bf L}'(s)}^{{\bf G}'(s)}(exp(H_{L})\xi^{{\bf L}'(s),st}(\boldsymbol{\delta},B,H^L),{\bf f}^{{\bf G}'(s)}).$$
 En regroupant les couples $(s,\tilde{L}'_{s})$ selon l'espace de Levi $\tilde{L}$ associ\'e, on obtient comme en [II] 3.7 que l'expression ci-dessus est \'egale \`a
 $$\sum_{\tilde{L}\in {\cal L}(\tilde{M})}\sum_{s\in \zeta Z(\hat{M})^{\Gamma_{{\mathbb R}}}/Z(\hat{L})^{\Gamma_{{\mathbb R}}}}i_{\tilde{M}'}(\tilde{L},\tilde{L}'(s))\sum_{t\in s Z(\hat{L})^{\Gamma_{{\mathbb R}}}/Z(\hat{G})^{\Gamma_{{\mathbb R}}}}i_{\tilde{L}'(s)}(\tilde{G},\tilde{G}'(t))$$
 $$S_{{\bf L}'(s)}^{{\bf G}'(t)}(exp(H_{L})\xi^{{\bf L}'(s),st}(\boldsymbol{\delta},B,H^L),{\bf f}^{{\bf G}'(t)})$$
 ou encore
 $$\sum_{\tilde{L}\in {\cal L}(\tilde{M})}\sum_{s\in \zeta Z(\hat{M})^{\Gamma_{{\mathbb R}}}/Z(\hat{L})^{\Gamma_{{\mathbb R}}}}i_{\tilde{M}'}(\tilde{L},\tilde{L}'(s))I_{\tilde{L}}^{\tilde{G},{\cal E}}({\bf L}'(s),exp(H_{L})\xi^{{\bf L}'(s),st}(\boldsymbol{\delta},B,H^L),{\bf f}).$$
 Les distributions $exp(H_{L})\xi^{{\bf L}'(s),st}(\boldsymbol{\delta},H^L,B)$ sont $\tilde{G}$-\'equisinguli\`eres pour $H$ en position g\'en\'erale et on peut appliquer la proposition 1.13. Cela nous permet de remplacer ci-dessus 
 $$I_{\tilde{L}}^{\tilde{G},{\cal E}}({\bf L}'(s),exp(H_{L})\xi^{{\bf L}'(s),st}(\boldsymbol{\delta},B,H^L),{\bf f})$$
  par 
  $$I_{\tilde{L}}^{\tilde{G}}(transfert(exp(H_{L})\xi^{{\bf L}'(s),st}(\boldsymbol{\delta},B,H^L)),{\bf f}).$$
   On a
 $$transfert((exp(H_{L})\xi^{{\bf L}'(s),st}(\boldsymbol{\delta},B,H^L))=exp(H_{L})(transfert(\xi^{{\bf L}'(s),st}(\boldsymbol{\delta},B,H^L))).$$
La somme se r\'ecrit
 $$(2)\qquad \sum_{\tilde{L}\in {\cal L}(\tilde{M})}I_{\tilde{L}}^{\tilde{G}}(exp(H_{L})\xi^{\tilde{L}}(\boldsymbol{\delta},B,H^L),{\bf f}),$$
 o\`u
 $$(3) \qquad  \xi^{\tilde{L}}(\boldsymbol{\delta},B,H^L)=\sum_{s\in \zeta Z(\hat{M})^{\Gamma_{{\mathbb R}}}/Z(\hat{L})^{\Gamma_{{\mathbb R}}}}i_{\tilde{M}'}(\tilde{L},\tilde{L}'(s))transfert(\xi^{{\bf L}'(s),st}(\boldsymbol{\delta},B,H^L)).$$
 Pour tout $s$ apparaissant ci-dessus, le terme $\xi^{{\bf L}'(s),st}(\boldsymbol{\delta},B)$ appartient \`a ${\cal V}_{M'}^{L'(s)}\otimes D^{st}_{tr-orb}({\bf L}'(s),{\cal O}')\otimes Mes(L'(s;{\mathbb R}))^*$. On a l'inclusion naturelle ${\cal V}_{M'}^{L'(s)}\subset {\cal V}_{M}^L$. En reprenant les d\'efinitions, on voit que le transfert envoie $D^{st}_{tr-orb}({\bf L}'(s),{\cal O}')\otimes Mes(L'(s;{\mathbb R}))^*$ dans $D_{tr-orb}(\tilde{L}({\mathbb R}),{\cal O})\otimes Mes(L({\mathbb R}))^*$. Donc $\xi^{\tilde{L}}(\boldsymbol{\delta},B)$ est bien un \'el\'ement de ${\cal V}_{M}^L\otimes D_{tr-orb}(\tilde{L}({\mathbb R}),{\cal O})\otimes Mes(L({\mathbb R}))^*$ comme on le voulait. On a montr\'e que le membre de droite de (1) \'etait faiblement \'equivalent \`a (2). C'est la relation de l'\'enonc\'e. $\square$

\bigskip
\subsection{Preuve de la proposition 6.3}
Comme on l'a d\'ej\`a vu plusieurs fois, la preuve est essentiellement la m\^eme
que la pr\'ec\'edente.  On part de la formule
$$S_{\tilde{M}}^{\tilde{G}}(\boldsymbol{\delta},B,f)=I_{\tilde{M}}^{\tilde{G}}(\boldsymbol{\delta},B,f)-\sum_{s\in Z(\hat{M})^{\Gamma_{{\mathbb R}}}/Z(\hat{G})^{\Gamma_{{\mathbb R}}}, s\not=1}i_{\tilde{M}}(\tilde{G},\tilde{G}'(s))S_{{\bf M}}^{{\bf G}'(s)}(\boldsymbol{\delta},B,{\bf f}^{{\bf G}'(s)}).$$
On applique la proposition 6.4 au premier terme et la proposition 6.3 par r\'ecurrence aux autres.  On obtient que $S_{\tilde{M}}^{\tilde{G}}(\boldsymbol{\delta},B,f)$ est faiblement \'equivalente \`a la fonction sur ${\cal A}_{M}$ qui \`a $H$ associe
$$(1) \qquad \sum_{\tilde{L}\in {\cal L}(\tilde{M})}\left(I_{\tilde{L}}^{\tilde{G}}(exp(H_{L})\xi^{\tilde{L}}(\boldsymbol{\delta},B,H^L),{\bf f})-X_{\tilde{L}}(H)\right),$$
o\`u
$$X_{\tilde{L}}(H)=\sum_{s\in Z(\hat{M})^{\Gamma_{{\mathbb R}}}/Z(\hat{L})^{\Gamma_{{\mathbb R}}}}i_{\tilde{M}}(\tilde{L},\tilde{L}'(s))\sum_{t\in s Z(\hat{L})^{\Gamma_{{\mathbb R}}}/Z(\hat{G})^{\Gamma_{{\mathbb R}}}, t\not=1}i_{\tilde{L}'(s)}(\tilde{G},\tilde{G}'(t))$$
$$S_{{\bf L}'(s)}^{{\bf G}'(t)}(exp(H_{L})\xi^{{\bf L}'(s),st}(\boldsymbol{\delta},B,H^L),{\bf f}^{{\bf G}'(t)}).$$
Supposons d'abord $\tilde{L}\not=\tilde{G}$.  D'apr\`es notre hypoth\`ese de r\'ecurrence, le terme $\xi^{\tilde{L},st}(\boldsymbol{\delta},B,H^L)$ est d\'ej\`a d\'efini. On voit qu'ajouter \`a $X_{\tilde{L}}(H)$ le terme $S_{\tilde{L}}^{\tilde{G}}(exp(H_{L})\xi^{\tilde{L},st}(\boldsymbol{\delta},B,H^L),{\bf f})$ revient \`a supprimer la restriction $t\not=1$ dans la formule ci-dessus. La somme en $t$ devient alors 
$$I_{\tilde{L}}^{\tilde{G},{\cal E}}({\bf L}'(s),exp(H_{L})\xi^{{\bf L}'(s),st}(\boldsymbol{\delta},B,H^L),{\bf f}).$$
En appliquant la proposition 1.13, c'est aussi 
$$I_{\tilde{L}}^{\tilde{G}}(exp(H_{L})transfert(\xi^{{\bf L}'(s),st}(\boldsymbol{\delta},B,H^L)),{\bf f}).$$
Alors
$$S_{\tilde{L}}^{\tilde{G}}(exp(H_{L})\xi^{\tilde{L},st}(\boldsymbol{\delta},B,H^L),{\bf f})+X_{\tilde{L}}(H)=I_{\tilde{L}}^{\tilde{G}}(exp(H_{L})\underline{\xi}^{\tilde{L}}(\boldsymbol{\delta},B,H^L),{\bf f}),$$
o\`u
$$ \underline{\xi}^{\tilde{L}}(\boldsymbol{\delta},B,H^L)=\sum_{s\in Z(\hat{M})^{\Gamma_{{\mathbb R}}}/Z(\hat{L})^{\Gamma_{{\mathbb R}}}}i_{\tilde{M}}(\tilde{L},\tilde{L}'(s))transfert(\xi^{{\bf L}'(s),st}(\boldsymbol{\delta},B,H^L)).$$
On applique la d\'efinition 6.4(3), en se rappelant qu'il s'agissait d'une \'egalit\'e dans ${\cal V}_{M}^L$, c'est-\`a-dire \`a faible \'equivalence pr\`es. Elle  implique  que $ \xi^{\tilde{L}}(\boldsymbol{\delta},B,H^L)$ et $ \underline{\xi}^{\tilde{L}}(\boldsymbol{\delta},B,H^L)$ sont faiblement \'equivalents.  Alors le terme index\'e par $\tilde{L}$ de la formule (1) est faiblement \'equivalent \`a  $S_{\tilde{L}}^{\tilde{G}}(exp(H_{L})\xi^{\tilde{L},st}(\boldsymbol{\delta},B,H^L),{\bf f})$.

Supposons maintenant $\tilde{L}=\tilde{G}$. La d\'efinition se simplifie:
$$X_{\tilde{G}}(H)=\sum_{s\in Z(\hat{M})^{\Gamma_{{\mathbb R}}}/Z(\hat{G})^{\Gamma_{{\mathbb R}}},s\not=1}i_{\tilde{M}}(\tilde{G},\tilde{G}'(s))S^{{\bf G}'(s)}(\xi^{{\bf G}'(s),st}(\boldsymbol{\delta},B,H),{\bf f}^{{\bf G}'(s)})$$
$$=\sum_{s\in Z(\hat{M})^{\Gamma_{{\mathbb R}}}/Z(\hat{G})^{\Gamma_{{\mathbb R}}},s\not=1}i_{\tilde{M}}(\tilde{G},\tilde{G}'(s))I^{\tilde{G}}(transfert(\xi^{{\bf G}'(s),st}(\boldsymbol{\delta},B,H)),{\bf f})$$
$$=I^{\tilde{G}}(\xi^{\tilde{G}}(\boldsymbol{\delta},B,H),{\bf f})-I^{\tilde{G}}(\underline{\xi}^{\tilde{G},st}(\boldsymbol{\delta},B,H),{\bf f}),$$
o\`u on a pos\'e
$$(2) \qquad \underline{\xi}^{\tilde{G},st}(\boldsymbol{\delta},B,H)=\xi^{\tilde{G}}(\boldsymbol{\delta},B,H)$$
$$-\sum_{s\in Z(\hat{M})^{\Gamma_{{\mathbb R}}}/Z(\hat{G})^{\Gamma_{{\mathbb R}}},s\not=1}i_{\tilde{M}}(\tilde{G},\tilde{G}'(s))transfert(\xi^{{\bf G}'(s),st}(\boldsymbol{\delta},B,H)).$$
Alors le terme index\'e par $\tilde{G}$ dans la formule (1) est \'egal \`a $I^{\tilde{G}}(\underline{\xi}^{\tilde{G},st}(\boldsymbol{\delta},B,H),{\bf f})$. Cette formule devient
$$(3) \qquad I^{\tilde{G}}(\underline{\xi}^{\tilde{G},st}(\boldsymbol{\delta},B,H),{\bf f})+\sum_{\tilde{L}\in {\cal L}(\tilde{M}),\tilde{L}\not=\tilde{G}}S_{\tilde{L}}^{\tilde{G}}(exp(H_{L})\xi^{\tilde{L},st}(\boldsymbol{\delta},B,H^L),{\bf f}).$$
La distribution $\xi^{\tilde{G}}(\boldsymbol{\delta},B,H)$ appartient \`a $D_{tr-orb}(\tilde{G}({\mathbb R}))\otimes Mes(G({\mathbb R}))^*$. Les distributions $\xi^{{\bf G}'(s),st}(\boldsymbol{\delta},B,H)$ appartiennent \`a $D_{tr-orb}^{st}({\bf G}'(s))\otimes Mes(G'(s;{\mathbb R}))^*$. Leur transfert appartient  \`a $D_{tr-orb}(\tilde{G}({\mathbb R}))\otimes Mes(G({\mathbb R}))^*$. La formule (2) montre que $\underline{\xi}^{\tilde{G},st}(\boldsymbol{\delta},B,H)$ appartient aussi \`a cet espace. Pour la m\^eme raison, $\underline{\xi}^{\tilde{G},st}(\boldsymbol{\delta},B,H)$ appartient \`a $Ind_{\tilde{M}}^{\tilde{G}}(D_{g\acute{e}om}({\cal O})\otimes Mes(M({\mathbb R}))^*)$. Donc  $\underline{\xi}^{\tilde{G},st}(\boldsymbol{\delta},B,H)$ appartient \`a $D_{tr-orb}(\tilde{G}({\mathbb R}),{\cal O})\otimes Mes(G({\mathbb R}))^*$.
Fixons un suppl\'ementaire $\underline{{\cal W}}$ de $\underline{{\cal V}}_{0,M}^G$ dans $\underline{{\cal V}}_{M}^G$.  Il existe un unique \'el\'ement $\xi^{\tilde{G},st}(\boldsymbol{\delta},B)\in \underline{{\cal W}}\otimes D_{tr-orb}(\tilde{G}({\mathbb R}),{\cal O})\otimes Mes(G({\mathbb R}))^*$ qui soit faiblement \'equivalent \`a $\underline{\xi}^{\tilde{G},st}(\boldsymbol{\delta},B)$. Dans la formule (3), on peut remplacer $\underline{\xi}^{\tilde{G},st}(\boldsymbol{\delta},B,H)$ par $\xi^{\tilde{G},st}(\boldsymbol{\delta},B,H)$.
Supposons que ${\bf f}$ soit instable, c'est-\`a-dire que son image  dans $SI(\tilde{G}({\mathbb R}))\otimes Mes(G({\mathbb R}))$ soit nulle. En vertu du th\'eor\`eme 1.4 et des propri\'et\'es \'enonc\'ees en 2.4,  les termes $S_{\tilde{M}}^{\tilde{G}}(\boldsymbol{\delta},B,{\bf f})$ et $S_{\tilde{L}}^{\tilde{G}}(exp(H_{L})\xi^{\tilde{L},st}(\boldsymbol{\delta},B,H^L),{\bf f})$ sont nuls. On obtient que $I^{\tilde{G}}(\xi^{\tilde{G},st}(\boldsymbol{\delta},B,H),{\bf f})$ est faiblement \'equivalent \`a $0$. On peut \'ecrire
$$\xi^{\tilde{G},st}(\boldsymbol{\delta},H,B)=\sum_{i=1,...,n}u_{i}(H)\boldsymbol{\gamma}_{i},$$
o\`u $(u_{i})_{i=1,...,n}$ est une famille d'\'el\'ements lin\'eairement ind\'ependants de $\underline{{\cal W}}$ et les $\boldsymbol{\gamma}_{i}$ appartiennent \`a $D_{tr-orb}(\tilde{G}({\mathbb R}),{\cal O})\otimes Mes(G({\mathbb R}))^*$. Alors
$$I^{\tilde{G}}(\xi^{\tilde{G},st}(\boldsymbol{\delta},B,H),{\bf f})=\sum_{i=1,...,n}u_{i}(H)I^{\tilde{G}}(\boldsymbol{\gamma}_{i},{\bf f}) .$$
Par d\'efinition de $\underline{{\cal W}}$, ceci ne peut \^etre faiblement \'equivalent \`a $0$ que si tous les coefficients sont nuls. Cela \'etant vrai pour tout ${\bf f}$ instable, les distributions $\boldsymbol{\gamma}_{i}$ sont stables. Gr\^ace au lemme [I] 5.13, elles sont induites d'\'el\'ements de $D_{g\acute{e}om}^{st}({\cal O})\otimes Mes(M({\mathbb R}))^*$ et appartiennent donc \`a $D_{tr-orb}^{st}(\tilde{G}({\mathbb R}),{\cal O})\otimes Mes(G({\mathbb R}))^*$. Donc $\xi^{\tilde{G},st}(\boldsymbol{\delta},B)$ prend ses valeurs dans l'espace voulu. Puisque ces valeurs sont stables, on peut remplacer le premier terme de (3) par $S^{\tilde{G}}(\xi^{\tilde{G},st}(\boldsymbol{\delta},B,H),{\bf f})$. Alors (3) devient la formule de l'\'enonc\'e de la proposition 6.3. Cela ach\`eve la d\'emonstration. $\square$

\bigskip

\section{Le cas des groupes complexes}
Consid\'erons tr\`es bri\`evement le cas o\`u le corps de base n'est plus ${\mathbb R}$ mais ${\mathbb C}$.  On consid\`ere donc un triplet $(G,\tilde{G},{\bf a})$ d\'efini sur ${\mathbb C}$. Tous les r\'esultats de l'article restent valables. Pour le voir, on peut reprendre les d\'emonstrations et constater qu'elles valent aussi bien dans le cas complexe. Certaines se simplifient: par exemple, les $K$-espaces disparaissent. On peut aussi appliquer les r\'esultats du cas r\'eel au triplet r\'eel $(G_{{\mathbb R}},\tilde{G}_{{\mathbb R}},{\bf a}_{{\mathbb R}})$ d\'eduit du triplet initial par restriction des scalaires. Cette deuxi\`eme m\'ethode perturbe les hypoth\`eses de r\'ecurrence. Pour $(G,\tilde{G},{\bf a})$, celles-ci concernent des triplets complexes $(G',\tilde{G}',{\bf a}')$ v\'erifiant en tout cas $dim(G'_{SC})\leq dim(G_{SC})$. Pour  $(G_{{\mathbb R}},\tilde{G}_{{\mathbb R}},{\bf a}_{{\mathbb R}})$, elles concernent des triplets r\'eels $(G',\tilde{G}',{\bf a}')$ v\'erifiant $dim(G'_{SC})\leq dim(G_{{\mathbb R},SC})=2dim(G_{SC})$. En fait, on n'applique ces hypoth\`eses de r\'ecurrence qu'\`a des triplets d\'eduits de $(G_{{\mathbb R}},\tilde{G}_{{\mathbb R}},{\bf a}_{{\mathbb R}})$ par des op\'erations naturelles: passage \`a un Levi, \`a une composante neutre d'un commutant, \`a l'espace d'une donn\'ee endoscopique etc... On constate que toutes ces op\'erations construisent des objets de m\^eme type que $(G_{{\mathbb R}},\tilde{G}_{{\mathbb R}},{\bf a}_{{\mathbb R}})$, c'est-\`a-dire d\'eduits par restriction des scalaires d'objets d\'efinis sur ${\mathbb C}$. On voit ainsi que les hypoth\`eses de r\'ecurrence pos\'ees pour  $(G,\tilde{G},{\bf a})$ suffisent \`a assurer la validit\'e des raisonnements pour le triplet  $(G_{{\mathbb R}},\tilde{G}_{{\mathbb R}},{\bf a}_{{\mathbb R}})$.

 \bigskip 
{\bf Bibliographie}

[A1]  J. Arthur: {\it Germ expansions for real groups}, pr\'epublication (2004)

[A2] ------------: {\it  The local behaviour of weighted orbital integrals}, Duke Math. J. 56 (1988), p. 223-293

[A3] ------------: {\it The trace formula in invariant form}, Annals of Math. 114 (1981), p. 1-74

[A4]-------------: {\it Parabolic transfer for real groups}, J. AMS 21 (2008), p. 171-234

[A5]-------------: {\it  On the transfer of distributions: weighted orbital integrals}, Duke Math. J. 99 (1999), p. 209-283

[W1] J.-L. Waldspurger: {\it La formule des traces locale tordue}, pr\'epublication 2012

[W2] -----------------------: {\it L'endoscopie tordue n'est pas si tordue}, Memoirs AMS 908 (2008)

[I] -------------------------: {\it Stabilisation de la formule des traces tordue I: endoscopie tordue sur un corps local}, pr\'epublication 2014

[II] -------------------------: {\it  Stabilisation de la formule des traces tordue II: int\'egrales orbitales et endoscopie sur un corps local non-archim\'edien; d\'efinitions et \'enonc\'es des r\'esultats}, pr\'epublication 2014
 
 [III] -------------------------: {\it Stabilisation de la formule des traces tordue III: int\'egrales orbitales et endoscopie sur un corps local non-archim\'edien; r\'eductions et preuves}, pr\'epublication 2014
 
 [IV] -------------------------: {\it Stabilisation de la formule des traces tordue IV: transfert spectral archim\'edien}, pr\'epublication 2014

 \bigskip
 
 Institut de math\'ematiques de Jussieu-CNRS
 
 2 place Jussieu 75005 Paris
 
 e-mail: waldspur@math.jussieu.fr

\end{document}